\documentclass[12pt,twoside]{amsbook}
\usepackage{cite}
\usepackage{latexsym}
\usepackage{layout}
\usepackage{amsfonts}
\usepackage{amssymb}
\usepackage{epsfig}
\usepackage{color}
\usepackage{tikz}

\usetikzlibrary{patterns}
\usetikzlibrary{intersections}
\usetikzlibrary{decorations.markings}

\newtheorem{thm}{Theorem}[chapter]

\newtheorem{lem}[thm]{Lemma}
\newtheorem{cor}[thm]{Corollary}

\newtheorem{prp}[thm]{Proposition}

\makeindex

\newsavebox\dotbox
\sbox{\dotbox}{\(\displaystyle\bigodot\)}
\DeclareMathOperator*{\bigcdot}{\raisebox{0pt}[\ht\dotbox][\dp\dotbox]{\(\boldsymbol{\cdot}\)}}
\usepackage{afterpage}

\usepackage{afterpage}

\usepackage{fancyhdr}
\begin{document}

\title{Quadratic Residues and Non-Residues: Selected Topics}
\author{Steve Wright 

Department of Mathematics and Statistics

Oakland University

Rochester, Michigan 48309

U.S.A.

e-mail: wright@oakland.edu} 
\maketitle
\markboth{}{}

\frontmatter
\vspace*{3.5cm}
\pagestyle{plain} 
\begin{center}
\textbf{\emph{For Linda} }   
\end{center}
\newpage
\afterpage{\null\newpage}
\thispagestyle{empty}

\newpage
{\pagestyle{plain}
\tableofcontents}

\afterpage{\null\newpage}
\thispagestyle{empty}

\chapter*{Preface}
{\pagestyle{plain}

Although number theory as a coherent mathematical subject started with the work of Fermat in the 1630's, modern number theory, i.e., the systematic and mathematically rigorous development of the subject from fundamental properties of the integers, began in 1801 with the appearance of Gauss' landmark treatise \emph{Disquisitiones Arithmeticae}\index{\emph{Disqusitiones Arithmeticae}} [19]. A major part of the \emph{Disquisitiones} deals with quadratic residues and non-residues: if $p$ is an odd prime, an integer $z$ is a quadratic residue  (respectively, a quadratic non-residue) of $p$ if there is (respectively, is not) an integer $x$ such that $x^2 \equiv z \mod p$. As we shall see, quadratic residues arise naturally as soon as one wants to solve the general quadratic congruence\index{quadratic congruence} $ax^2+bx+c \equiv 0\ \textrm{mod}\ m,\ a \not \equiv 0\ \textrm{mod}\ m$, and this, in fact, motivated some of the interest which Gauss himself had in them. Beginning with Gauss' fundamental contributions, the study of quadratic residues and non-residues has subsequently led directly to many of the key ideas and techniques that are used everywhere in number theory today, and the primary goal of these lecture notes is to use this study as a window through which to view the development of some of those ideas and techniques. In pursuit of that goal, we will employ methods from elementary, analytic, and combinatorial number theory, as well as methods from the theory of algebraic numbers.

In order to follow these lectures most profitably, the reader should have some familiarity with the basic results of elementary number theory. An excellent source for this material (and much more) is the text [30] of Kenneth Ireland and Michael Rosen, \emph{A Classical Introduction to Modern Number Theory}. A feature of this text that is of particular relevance to what we discuss is Ireland and Rosen's treatment of quadratic and higher-power residues, which is noteworthy for its elegance and completeness, as well as for its historical perspicacity. We will in fact make use of some of their work in Chapters 3 and 7.

Although not absolutely necessary, some knowledge of algebraic number theory will also be helpful for reading these notes. We will provide complete proofs of some facts about algebraic numbers and we will quote other facts without proof. Our reference for proof of the latter results is the classical treatise of Erich Hecke [27], \emph{Vorlesungen $\ddot{u}$ber die Theorie der Algebraischen Zahlen}, in the very readable English translation by G. Brauer and J. Goldman. About Hecke's text Andr$\acute{\textrm{e}}$ Weil ([58], foreword) had this to say: ``To improve upon Hecke, in a treatment along classical lines of the theory of algebraic numbers, would be a futile and impossible task." We concur enthusiastically with Weil's assessment and highly recommend Hecke's book to all those who are interested in number theory.

We next offer a brief overview of what is to follow. The notes are arranged in a series of ten chapters. Chapter 1, an introduction to the subsequent chapters, provides some motivation for the study of quadratic residues and non-residues by consideration of what needs to be done when one wishes to solve the general quadratic congruence mentioned above. We briefly discuss the contents of the  \emph{Disquisitiones Arithmeticae}, present some biographical information about Gauss, and also record some basic results from elementary number theory that will be used frequently in the sequel. Chapter 2 provides some useful facts about quadratic residues and non-residues upon which the rest of the chapters are based. Here we also describe a procedure which provides a strategy for solving what we call the \emph{Basic Problem}: if $d$ is an integer, find all primes $p$ such that $d$ is a quadratic residue of $p$. The Law of Quadratic Reciprocity is the subject of Chapter 3. We present seven proofs of this fundamentally important result (five in Chapter 3, one in Chapter 7, and one in Chapter 8), which focus primarily (but not exclusively) on the ideas used in the proofs of quadratic reciprocity which Gauss discovered. Chapter 4 discusses some interesting and important applications of quadratic reciprocity, having to do with the solution of the Basic Problem from Chapter 2 and with the structure of  the finite subsets $S$ of the positive integers possessing at least one of the following two properties: for infinitely many primes $p$, $S$ is a set of quadratic residues of $p$, or for infinitely many primes $p$, $S$ is a set of quadratic non-residues of $p$. Here the fundamental contributions of Dirichlet to the theory of quadratic residues enters our story and begins a major theme that will play throughout the rest of our work. Chapter 4 concludes with an interesting application of quadratic residues in modern cryptology, to so-called zero-knowledge or minimal-disclosure proofs. The use of transcendental methods in the theory of quadratic residues, begun in Chapter 4, continues in Chapter 5 with the study of the zeta function of an algebraic number field and its application to the solution of some of the problems taken up in Chapter 4. Chapter 6 gives elementary proofs of some of the results in Chapter 5 which obviate the use made there of the zeta function. The question of how quadratic residues and non-residues of a prime $p$ are distributed among the integers $1, 2,\dots,p-1$ is considered in Chapter 7, and there we highlight additional results and methods due to Dirichlet which employ the basic theory of $L$-functions attached to Dirichlet characters determined by certain moduli. Because of the importance that  positivity of the values at $s=1$ of Dirichlet $L$-functions plays in the proof of the results of Chapter 7, we present in Chapter 8 a discussion and proof of Dirichlet's class-number formula as a way to defenitively explain why the  values at $s=1$ of $L$-functions are positive. In Chapter 9 the occurrence of quadratic residues and non-residues as arbitrarily long arithmetic progressions is studied by means of some ideas of Harold Davenport [5] and some techniques in combinatorial number theory developed in recent work of the author [62], [63]. A key issue that arises in our approach to this problem is the estimation of certain character sums over the field of $p$ elements, $p$ a prime, and we address this issue by using some results of A. Weil [57] and G. I. Perel'muter [44].  Our discussion concludes with Chapter 10, where the Central Limit Theorem from the theory of probability and a theorem of Davenport and Paul Erd$\ddot{\textrm{o}}$s [6] are used to provide evidence for the contention that as the prime $p$ tends to infinity, quadratic residues of $p$ are distributed randomly throughout certain subintervals of the set $\{1, 2,\dots,p-1\}$.

These notes are an elaboration of the contents of a special-topics-in-mathematics course that was offered during the Summer semesters of 2014 and 2015 at Oakland University. I am very grateful to my colleague Meir Shillor for suggesting that I give such a course, and for thereby providing me with the impetus to think about what such a course would entail. I am also very grateful to my colleagues Eddie Cheng and Serge Kruk, the former for giving me very generous and  valuable assistance with numerous LaTeX issues which arose during the preparation of the manuscript, and the latter for formatting all of the figures in the text. I thank my students Saad Al Najjar, Amelia McIlvenna and Julian Venegas for reading an early version of the notes and offering several insightful comments which were very helpful to me.  My sincere and heartfelt appreciation is also tendered to the anonymous referees for many comments and suggestions which resulted in a very substantial improvement in both the content and exposition of these notes. Finally, and above all others, I am grateful beyond words to my dear wife Linda for her unstinting love, support, and encouragement; this humble missive is dedicated to her.

\newpage
\afterpage{\null\newpage}
\thispagestyle{empty}

\pagestyle{fancy}
\mainmatter
\chapter{Introduction: Solving the General Quadratic Congruence Modulo a Prime}
The purpose of this chapter is to define quadratic residues and non-residues and to use the solution of the general quadratic congruence modulo a prime to indicate one reason why the study of quadratic residues and non-residues is interesting and important. This is done in section 1. The primary source for essential information about quadratic residues is the  \emph{Disquisitiones Arithmeticae} of Carl Friedrich Gauss, one of the most important books about number theory ever written. Because of its singular prominence for number theory and also for what we will do in these lecture notes, the contents of the \emph{Disquisitiones} are discussed briefly in section 2, and some biographical facts about Gauss are also presented. Notation and terminology that will be employed throughout the sequel are recorded in section 3, as well as a few basic facts from elementary number theory that will be used frequently in subsequent work.

\section{Linear and Quadratic Congruences}
One of the central problems of number theory, both ancient and modern, is finding solutions (in the integers) of polynomial equations with integer coefficients in one or more variables. In order to motivate our study, consider the equation
\[
ax \equiv b\ \textrm{mod}\ m,
\]
a linear equation in the unknown integer $x$. Elementary number theory provides an algorithm for determining exactly when this equation has a solution, and for finding all such solutions, which essentially involves nothing more sophisticated than the Euclidean algorithm (see Proposition 1.4 below and the comments after it).

When we consider what happens for the general quadratic congruence
\begin{equation*}
ax^2+bx+c \equiv 0\ \textrm{mod}\ m,\ a \not \equiv 0\ \textrm{mod}\ m \tag{1},
\]
things get more complicated. In order to see what  the issues are, note first that

\begin{eqnarray*}
 (2ax+b)^2&\equiv&b^2-4ac \ \textrm{mod}\ 4am\\
&\textrm{iff}&4a^2x^2+4abx+4ac \equiv 0  \ \textrm{mod}\ 4am\\
&\textrm{iff}&4a(ax^2+bx+c) \equiv 0\ \textrm{mod}\ 4am\\
&\textrm{iff}&ax^2+bx+c \equiv 0\ \textrm{mod}\ m . 
\end{eqnarray*}
Hence (1) has a solution if and only if
\begin{equation*}
2ax \equiv s-b\ \textrm{mod}\ 4am, \tag{2} \]
where $s$ is a solution of
\begin{equation*}
s^2 \equiv b^2-4ac\ \textrm{mod}\ 4am \tag{3}.
\]
Now (2) has a solution if and only if $s-b$ is divisible by $2a$, the greatest common divisor of $2a$ and $4am$, and so it follows that (1) has a solution if and only if (3) has a solution $s$ such that $s-b$ is divisible by $2a$. We have hence reduced the solution of (1) to finding solutions $s$ of (3) which satisfy an appropriate divisibility condition.

Our attention is therefore focused on the following problem: if $n$ and $z$ are integers with $n\geq 2$, find all solutions $x$ of the congruence 
\begin{equation*}
x^2 \equiv z\ \textrm{mod}\ n. \tag{4} \]
Let\[
n=\prod_{i=1}^k p_i^{\alpha_i} \]
be the prime factorization of $n$, and let $\Sigma_i$ denote the set of all solutions of the congruence
\[
x^2 \equiv z\ \textrm{mod}\ p_i^{\alpha_i},\ i=1,\dots,k.\]
Let $s=(s_1,\dots,s_k) \in \Sigma_1 \times \cdots \times \Sigma_k$, and let $\sigma(s)$ denote the simultaneous solution, unique mod $n$, of the system of congruences
\[
x \equiv s_i\ \textrm{mod}\  p_i^{\alpha_i},\  i=1,\dots,k,\]
obtained via the Chinese remainder theorem (Theorem 1.3 below). It is then not difficult to show that the set of all solutions of (4) is given precisely by the set
\[
\{\sigma(s): s \in  \Sigma_1 \times \cdots \times \Sigma_k \}.\]
Consequently (4), and hence also (1), can be solved if we can solve the congruence
\begin{equation*}
x^2 \equiv z\ \textrm{mod}\ p^{\alpha}, \tag{5}\]
where $p$ is a fixed prime and $\alpha$ is a fixed positive integer.

In articles 103 and 104 of  \emph{Disquisitiones Arithmeticae} [19], Gauss gave a series of beautiful formulae which completely solve (5) for all primes $p$ and exponents $\alpha$. 
In order to describe them, let $\sigma \in \{0, 1,\dots,p^{\alpha}-1\}$ denote a solution of (5). 

I. Suppose first that $z$ is not divisible by $p$. If $p=2$ and $\alpha=1$ then $\sigma=1$. If $p$ is odd or $p=2=\alpha$ then $\sigma$ has exactly two values $\pm \sigma_0$. If $p=2$ and $\alpha>2$ then $\sigma$ has exactly four values  $\pm \sigma_0$ and $\pm \sigma_0+2^{\alpha-1}$. 

II. Suppose next that $z$ is divisible by $p$ but not by $p^{\alpha}$. If (5) has a solution it can be shown that the multiplicity of $p$ as a factor of $z$ must be even, say $2\mu$, and so let $z=z_1p^{2\mu}$. Then $\sigma$ is given by the formula
\[
\sigma^{\prime}p^{\mu}+ip^{\alpha-\mu},\ i\in \{0, 1,\dots,p^{\mu}-1\},\]
where $\sigma^{\prime}$ varies over all solutions, determined according to I, of the congruence
\[
x^2 \equiv z_1\ \textrm{mod}\ p^{\alpha-2\mu}.\]

III. Finally if $z$ is divisible by $p^{\alpha}$, and if we set $\alpha=2k$ or $\alpha=2k-1$, depending on whether $\alpha$ is even or odd, then $\sigma$ is given by the formula
\[
ip^k,\ i\in \{0,\dots,p^{\alpha-k}-1\}.\]

We will focus on the most important special case of (5), namely when $p$ is odd and $\alpha=1$, i.e., the congruence
\[
x^2 \equiv z\ \textrm{mod}\ p \tag{6}\]
(note that when $p$ is odd, the solutions of (5) in cases I and II are determined by the solutions of (6) for certain values of $z$). The first thing to do here is to observe that the ring determined by the congruence classes of integers mod $p$ is a field, and so (6) has at most two solutions. We have that $x \equiv 0$ mod $p$ is the unique solution of (6) if and only if $z$ is divisible by $p$, and if $s_0 \not \equiv 0$ mod $p$ is a solution of (6) then so is $-s_0$, and $s_0  \not \equiv -s_0$ mod $p$ because $p$ is an odd prime. These facts are motivation for the following definition:

\vspace{0.3cm}
\emph{Definition}. If $p$ is an odd prime and $z$ is an integer not divisible by $p$, then $z$ is a \emph{quadratic residue} ( respectively, \emph{quadratic non-residue}) \emph{of} $p$ if there is (respectively, is not) an integer $x$ such that $x^2 \equiv z\ \textrm{mod}\ p$.  
\vspace{0.3cm}

As a consequence of our previous discussion and Gauss' solution of (5), solutions of (1) will exist only if (among other things) for each (odd) prime factor $p$ of $4am$, the discriminant $b^2-4ac$ of $ax^2+bx+c$ is either divisible by $p$ or is a quadratic residue of $p$. This remark becomes even more emphatic if the modulus $m$ in (1) is a single odd prime $p$. In that case,
\begin{equation*}
(2ax+b)^2 \equiv b^2-4ac\ \textrm{mod}\ p\ \textrm{iff}\ ax^2+bx+c \equiv 0\ \textrm{mod}\ p, \]
from whence the next proposition follows immediately:
\begin{prp}
Let p be an odd prime.The congruence
\begin{equation*}
ax^2+bx+c \equiv 0\ \textnormal{mod}\ p,\ a\not \equiv 0 \ \textnormal{mod}\ p, \tag{7}\]
has a solution if and only if
\begin{equation*}
x^2 \equiv b^2-4ac \ \textnormal{mod}\ p\tag{8}\]
has a solution, i.e., if and only if either $b^2-4ac$ is divisible by $p$ or $b^2-4ac$ is a quadratic residue of $p$. Moreover, if $(2a)^{-1}$ is the multiplicative inverse of $2a\ \textnormal{mod}\ p$ (which exists because $p$ does not divide $2a$; see Proposition $1.2$ below) then the solutions of $(7)$ are given precisely by the formula
\[
x \equiv (\pm s-b)(2a)^{-1}\ \textnormal{mod}\ p,
\]
where $\pm s$ are the solutions of $(8)$.
\end{prp}
\noindent We take it as self-evident that the solution of the general quadratic congruence (1) is one of the most fundamental and most important problems in the theory of Diophantine equations in two variables.  By virtue of Proposition 1.1 and the discussion which precedes it, quadratic residues and non-residues play a pivotal role in the determination of the solutions of (1). We hope that the reader will now agree: the study of quadratic residues and non-residues is important and interesting!
\section{The \emph{Disquisitiones Arithmeticae}}
As we pointed out in the preface, modern number theory, that is, the systematic and mathematically rigorous development of the subject from fundamental properties of the integers, began in 1801 with the publication of Gauss' great treatise, the \emph{Disquisitiones Arithmeticae}. Because of its first appearance here in our story, and especially because it plays a dominant role in that story, we will now briefly discuss some of the most important aspects of  the  \emph{Disquisitiones}. The book consists of seven sections divided into 366 articles. The first three sections are concerned with establishing basic results in number theory such as the Fundamental Theorem of Arithmetic (proved rigorously here for the first time), Fermat's little theorem, primitive roots and indices, the Chinese remainder theorem, and Wilson's theorem. Perhaps the most important innovation in the \emph{Disquisitiones} is Gauss' introduction in Section I, and its systematic use throughout the rest of the book, of the concept of modular congruence. Gauss shows how modular congruence can be used to give the study of divisibility of the integers a comprehensive and systematic \emph{algebraic} formulation, thereby greatly increasing the power and diversity of the techniques in number theory that were in use up to that time. Certainly one of the most striking examples of the power of modular congruence is the use that Gauss made of it in Section IV in his investigation of quadratic residues, which culminates in the first complete and correct proof of the Law of Quadratic Reciprocity, the most important result of that subject. We will have much more to say about quadratic reciprocity in Chapter 3. 

By far the longest part of the \emph{Disquisitiones}, over half of the entire volume, is taken up by Section V, which contains Gauss' deep and penetrating analysis of quadratic forms. If $(a, b, c)\in \mathbb{Z}\times \mathbb{Z}\times \mathbb{Z}$, the function defined on $\mathbb{Z}\times \mathbb{Z}$ by
\[
(x, y)\rightarrow ax^2+bxy+cy^2\]
is called a (binary) \emph{quadratic form}. We frequently repress the functional dependence on $(x, y)$ and hence also refer to the polynomial $ax^2+bxy+cy^2$ as a quadratic form. In Section V, Gauss first develops a way to systematically classify quadratic forms according to their number-theoretic properties and then investigates in great depth the algebraic and number theoretic structure of quadratic forms using the properties of this classification. Next, he defines an operation, which he called the \emph{composition of forms}, on the set of forms $ax^2+bxy+cy^2$ whose discriminate $b^2-4ac$ is a fixed value. The set of all quadratic forms whose discriminants have a fixed value supports a basic equivalence relation which partitions this set of forms into a finite number of equivalence classes (section 12, Chapter 3), and composition of forms turns this set of equivalence classes  into what we would today call a group. Of course Gauss did not call it that, as the group concept was not formulated until much later; however, he did make essential use of the group \emph{structure} which the composition of forms possesses. Gauss then proceeds to use composition of forms together with the methods that he developed previously to establish additional results concerning the algebraic and number-theoretic structure of forms. We will see in section 12 of Chapter 3 how some of the elements of Gauss' theory of quadratic forms arises in our study of quadratic reciprocity, and we will also make use of some of the basic theory of quadratic forms in Chapter 8.

Section VI of the \emph{Disquisitiones} is concerned with applications of primitive roots, quadratic residues, and quadratic forms to the structure of rational numbers, the solution of modular square-root problems, and primality testing and prime factorization of the integers. Finally, in Section VII, Gauss presents his theory of cyclotomy, the study of divisions of the circle into congruent arcs, which culminates in his famous theorem on the determination of the regular $n$-gons which can be constructed using only a straightedge and compass.

It is appropriate at this juncture to say a few words about Gauss himself. Carl Friedrich Gauss was born in 1777 in Brunswick, a city in the north of present-day Germany, lived for most of his life in G$\ddot{\textrm{o}}$ttingen, and died there in 1855.  Gauss' exceptional mathematical talent was clear from a very early age. Because he was such a gifted and promising young student, Gauss was introduced in 1791 to the Duke of Brunswick, who became a prominant patron and supporter of Gauss for many years (see, in particular, the  dedication in the \emph{Disquisitiones} which Gauss addressed to the Duke).  In 1795, Gauss matriculated at the University of G$\ddot{\textrm{o}}$ttingen, left in 1798 without obtaining a degree, and was granted a doctoral degree from the University of Helmstedt in 1799. After his celebrated computation of the orbit of the asteroid Ceres in 1801 Gauss was appointed director of the newly opened observatory at the University of G$\ddot{\textrm{o}}$ttingen in 1807, which position he held for the rest of his life. In addition to his ground-breaking work in number theory, Gauss made contributions of fundamental importance to geometry (differential geometry and non-Euclidean geometry), analysis  (elliptic functions, elliptic integrals, and the theory of infinite series), physics (potential theory and geomagnetism), geodesy and astronomy (celestial mechanics and the computation of the orbits of celestial bodies), and statistics and probability (the method of least squares and the normal distribution).

\section{Notation, Terminology, and Some Useful Elementary Number Theory}
We now fix some notation and terminology that will be used repeatedly throughout the sequel. The letter $p$ will always denote a generic odd prime, the letter $q$, unless otherwise specified, will denote a generic prime (either even or odd), $P$ is the set of all primes, $\mathbb{Z}$ is the set of all integers, $\mathbb{Q}$ is the set of all rationals, and $\mathbb{R}$ is the set of all real numbers. If $m, n \in \mathbb{Z}$ with $m \leq n$ then $[m, n]$ is the set of all integers at least $m$ and no more than $n$, listed in increasing order, $[m, \infty)$ is the set of all integers exceeding $m-1$, also listed in increasing order, and gcd$(m, n)$ is the greatest common divisor of $m$ and $n$. If $n \in [2, \infty)$ then $U(n)$ will denote the set $\{m\in [1, n-1]: \gcd(m, n)=1\}$. If $z$ is an integer then $\pi(z)$ will denote the set of all prime factors of $z$. If $A$ is a set then $|A|$ will denote the cardinality of $A$, $2^A$ is the set of all subsets of $A$, and $\emptyset$ denotes the empty set. Finally, we will refer to a quadratic residue or quadratic non-residue as simply a residue or non-residue; all other residues of a modulus $m \in [2, \infty)$ will always be called \emph{ordinary} residues. In particular, the minimal non-negative ordinary residues modulo $m$ are the elements of the set $[0, m-1]$.

We also recall some facts from elementary number theory that will be useful in what follows. For more information about them consult any standard text on elementary number theory, e.g., Ireland and Rosen [30] or K. Rosen [48].
 
 If $m$ is a positive integer and $a \in \mathbb{Z}$, recall that an \emph{inverse of $a$ modulo $m$} is an integer $\alpha$ such that $a \alpha \equiv 1 \ \textrm{mod}\ m$.
 
 \begin{prp}
 If $m$ is a positive integer and $a \in \mathbb{Z}$ then $a$ has an inverse modulo $m$ if and only if $\gcd(a, m)=1$. Moreover, the inverse is relatively prime to $m$ and is unique modulo $m$.
 \end{prp}
\begin{thm}
$($Chinese remainder theorem$)$. If $m_1,\dots,m_r$ are pairwise relatively prime positive integers and $(a_1,\dots,a_r)$ is an $r$-tuple of integers
then the system of congruences
\[
x \equiv a_i \ \textnormal{mod}\ m_i,\ i=1,\dots,r,
\]
has a simultaneous solution that is unique modulo $\prod_{i=1}^r m_i$. Moreover, if
\[ 
M_k= \prod _{i \not=k}\ m_i,
\]
and if $y_k$ is the inverse of $M_k \mod m_k$ $($which exits because $\gcd(m_k, M_k)=1)$  then the solution is given by
\[
x \equiv \sum _{k=1}^r \ a_kM_ky_k \ \textnormal{mod}\ \prod_{i=1}^r m_i.
\]
\end{thm}
 
Recall that a \emph{linear Diophantine equation} is an equation of the form
\[
ax+by=c,
\]
where $a, b,$ and $c$ are given integers and $x$ and $y$ are integer-valued unknowns.
\begin{prp}
Let a, b, and c be integers and let $\gcd(a,b)=d$. The Diophantine equation $ax+by=c$ has a solution if and only if d divides c. If d divides c then there are infinitely many solutions $(x, y)$, and if $(x_0, y_0)$ is a particular solution then all solutions are given by
\[
x=x_0+(b/ d)n,\ y=y_0-(a/d)n,\ n \in \mathbb{Z}.
\] 
\end{prp}

Given the Diophantine equation $ax+by=c$ with $c$ divisible by $d=\gcd(a, b)$, the Euclidean algorithm can be used to easily find a particular solution $(x_0, y_0)$. Simply let $k=c/d$ and use the Euclidean algorithm to find integers $m$ and $n$ such that $d=am+bn$; then $(x_0, y_0)=(km, kn)$ is a particular solution, and all solutions can then be found by using Proposition 1.4. The simple first-degree congruence $ax \equiv b$ mod $m$ can thus be easily solved upon the observation that this congruence has a solution $x$ if and only if the  Diophantine equation $ax+my=b$ has the solution $(x, y)$ for some $y\in \mathbb{Z}$. 
\newpage
\chapter{Basic Facts}

In this chapter, we lay the foundations for all of the work that will be done in subsequent chapters. Section 1 defines the Legendre symbol and verifies its basic properties, proves Euler's criterion, and deduces some corollaries which will be very useful in many situations in which we will find ourselves. Motivated by the solutions of a quadratic congruence modulo a prime which we discussed in Chapter 1, we formulate what we will call the Basic Problem and the Fundamental Problem for Primes in section 2. In section 3, we state and prove Gauss' Lemma for residues and non-residues and use it to solve the Fundamental Problem for the prime 2. 
\section{The Legendre Symbol, Euler's Criterion, and other Important Things}

In this section, we establish some fundamental facts about residues and non-residues that will be used repeatedly throughout the rest of these notes.
\begin{prp}
In every complete system of ordinary residues modulo $p$, there are exactly $(p-1)/2$ quadratic residues. 
\end{prp}

\emph{Proof}. It suffices to prove that in $[1, p-1]$ there are exactly $(p-1)/2$ quadratic residues. Note first that $1^2, 2^2,\dots, (\frac{p-1}{2})^2$ are all incongruent mod $p$ (if $1 \leq i, j< p/2$ and $i^2 \equiv j^2\ \textrm{mod}\ p$ then $i \equiv j$ hence $i=j$ or $i \equiv -j$, i.e., $i+j \equiv 0$. But $2 \leq i+j<p$, and so $i+j \equiv 0$ is impossible).

Let $\mathcal{S}$ denote the set of minimal non-negative ordinary residues mod $p$ of $1^2, 2^2,\dots, (\frac{p-1}{2})^2$. The elements of $\mathcal{S}$ are quadratic residues of $p$ and $|\mathcal{S}| = (p-1)/2.$ Suppose that $n \in [1, p-1]$ is a quadratic residue of $p$. Then there exists $r \in [1, p-1]$ such that $r^2 \equiv n$. Then $(p-r)^2 \equiv r^2 \equiv n$ and $\{r, p-r\} \cap [1, (p-1)/2] \not= \emptyset$. Hence $n \in \mathcal{S}$, whence $\mathcal{S}=$ the set of quadratic residues of $p$ inside $[1, p-1]$.$\hspace{12cm} \textrm{QED}   $ 

\emph{Remark}. The proof of Proposition 2.1 provides a way to easily find, at least in principle, the residues of any prime $p$. Simply calculate the integers $1^2, 2^2,\dots, (\frac{p-1}{2})^2$ and then reduce mod $p$. The integers that result from this computation are the residues of $p$ inside $[1, p-1]$. This procedure also finds the modular square 
roots $x$ of a residue $r$ of $p$, i.e., the solutions to the congruence $x^2\equiv r$ mod $p$. For example, in just a few minutes on a hand-held calculator, one finds that the residues of 17 are 1, 2, 4, 8, 9, 13, 15, and 16, with corresponding modular square roots $\pm1,\ \pm6,\ \pm2,\ \pm5,\ \pm3,\ \pm8,\ \pm7$, and $\pm$4, and the residues of 37 are 1, 3, 4, 7, 9, 10, 11, 12, 16, 21, 25, 26, 27, 28, 30, 33, 34, and 36, with corresponding modular square roots $\pm1,\ \pm15,\ \pm2,\ \pm9,\ \pm3,\ \pm11,\ \pm14,\ \pm7,\ \pm4,\ \pm13,\ \pm5,\ \pm10,\ \pm8,\ \pm18,\ \pm17,\ \pm12,\ \pm16$, and $\pm6$. Of course, for large $p$, this method quickly becomes impractical for the calculation of residues and modular square roots, but see section 9 of Chapter 4 for a practical and efficient way to perform these calculations for large values of $p$.

N.B. In the next proposition, all residues and non-residues are taken with respect to a fixed prime $p$.

\begin{prp}
$(i)$ The product of two residues is a residue.

$(ii)$ The product of a residue and a non-residue is a non-residue.

$(iii)$ The product of two non-residues is a residue.
\end{prp}
 \emph{Proof}. $(i)$ If $\alpha, \alpha^{\prime}$ are residues then $x^2 \equiv \alpha, y^2 \equiv \alpha^{\prime}$ imply that $(xy)^2 \equiv \alpha \alpha^{\prime} \ \textrm{mod}\ p$. 

 $(ii)$ Let $\alpha$ be a fixed residue. The integers $0, \alpha,\dots, (p-1)\alpha$ are incongruent mod $p$, hence are a complete system of ordinary residues mod $p$. If $R$ denotes the set of all residues in $[1, p-1]$ then by Proposition 2.2$(i)$, $\{\alpha r: r \in R\}$ is a set of residues of cardinality $(p-1)/2$, hence Proposition 2.1 implies that there are no other residues among $\alpha, 2\alpha,\dots, (p-1)\alpha$, i.e., if $\beta \in [1, p-1] \setminus R$ then $\alpha \beta$ is a non-residue. Statement $(ii)$ is an immediate consequence of this.

$(iii)$ 
Suppose that $\beta$ is a non-residue. Then $0, \beta, 2\beta,\dots, (p-1)\beta$ is a complete system of ordinary residues mod $p$, and by Proposition 2.2$(ii)$ and Proposition 2.1, $\{\beta r: r \in R\}$ is a set of non-residues and there are no other non-residues among $\beta, 2\beta,\dots, (p-1)\beta$. Hence $\beta^{\prime} \in [1, p-1] \setminus R$ implies that $\beta\beta^{\prime}$ is a residue. Statement $(iii)$ is an immediate consequence of this.$\hspace{14.8cm} \textrm{QED}$

The following definition introduces the most important piece of mathematical technology that we will use to study residues and non-residues. 

\vspace{0.3cm}
\emph{Definition}. The \emph{Legendre symbol} $\chi_p$ \emph{of} $p$ is the function $\chi_p: \mathbb{Z} \rightarrow [-1, 1]$ defined by
\[
\chi_p(n)=\left\{\begin{array}{rl}0,& \textrm{if $p$ divides $n$,}\\
1,& \textrm{if gcd($p$, $n$) = 1 and $n$ is a residue of $p$,}\\
-1,& \textrm{if gcd($p$, $n$) = 1 and $n$ is a non-residue of $p$.}\\\end{array}\right.
\]
\vspace{0.3cm}

The next proposition asserts that $\chi_p$ is a completely multiplicative arithmetic function of period $p$. This fact will play a crucial role in much of our subsequent work.
\begin{prp}
$(i)\ \chi_p(n)=0$ if and only if $p$ divides $n$, and if $m \equiv n \ \textnormal{mod}\ p$ then $\chi_p(m)= \chi_p(n)$ $(\chi_p$ is of period $p)$.

$(ii)$ For all $m, n \in \mathbb{Z}, \chi_p(mn)=\chi_p(m) \chi_p(n)$ $(\chi_p$ is completely multiplicative$)$.
\end{prp}

\emph{Proof}. $(i)$ If $m \equiv n \mod p$ then $p$ divides $m$ (respectively, $m$ is a residue/non-residue of $p$) if and only if $p$ divides $n$ (respectively, $n$ is a residue/non-residue of $p$). Hence $\chi_p(m)= \chi_p(n)$.

$(ii)$ $\chi_p(mn)=0$ if and only if $p$ divides $mn$ if and only if $p$ divides $m$ or $n$ if and only if $\chi_p(m)=0$ or $\chi_p(n)=0$ if and only if $\chi_p(m) \chi_p(n)=0$.

Because $\chi_p(n^2)=\big(\chi_p(n)\big)^2$, we may assume that $m\not=n$. Then $\chi_p(mn)=1$ (respectively, $\chi_p(mn)=-1$) if and only if gcd$(mn, p)=1$ and $mn$ is a residue (respectively, $mn$ is a non-residue) of $p$ if and only if gcd$(m, p)=1=$ gcd$(n, p)$ and, by Proposition 2.2, $m$ and $n$ are either both residues or both non-residues of $p$ (respectively, $\{m, n\}$ contains a residue and a non-residue of $p$) if and only if $\chi_p(m) \chi_p(n)=1$ (respectively, $\chi_p(m) \chi_p(n)=-1$).$\hspace{14.8cm}\textrm{QED}$

 \emph{Remark on notation}. As a consequence of Proposition 2.3, $\chi_p$ defines a homomorphism of the group of units in the ring $\mathbb{Z}/p\mathbb{Z}$ into the circle group, i.e., $\chi_p$ is a \emph{character} of the group of units. This is the reason why we have chosen the character-theoretic notation $\chi_p(n)$ for the Legendre symbol, instead of the more traditional notation $\displaystyle\left(\frac{n}{p} \right)$. When $p$ is replaced by an arbitrary integer $m \geq 2$, we will have more to say later (see section 4 of Chapter 4) about characters on the group of units in the ring $\mathbb{Z}/m\mathbb{Z}$ and their use in what we will study here. 

The next result determines the quadratic character of $-1$.
\begin{thm}
\[
\chi_p(-1)=\left\{\begin{array}{rl}1,& \textnormal{if $p \equiv 1 \ \textnormal{mod}\ 4$,}\\
-1,& \textnormal{if $p \equiv -1 \ \textnormal{mod}\ 4$ .}\\\end{array}\right.
\]
\end{thm}
 This theorem is due to Euler [16], who proved it in 1760. It is of considerable importance in the history of number theory because in 1795, the young Gauss (at the ripe old age of 18!) rediscovered it.  Gauss was so struck by the beauty and depth of this result that, as he testifies in the preface to  \emph{Disquisitiones Arithmeticae} [19],  ``I concentrated on it all of my efforts in order to understand the principles on which it depended and to obtain a rigorous proof. When I succeeded in this I was so attracted by these questions that I could not let them be." Thus began Gauss' work in number theory that was to revolutionize the subject.
 
 \emph{Proof of Theorem} 2.4. The proof that we give is Euler's own. It is based on 
 \begin{thm}
 $($Euler's criterion$)$ If $a \in \mathbb{Z}$ and $\gcd(a, p)=1$ then
 \[
 \chi_p(a) \equiv a^{(p-1)/2} \ \textnormal{mod}\ p.
 \]
\end{thm}

If we apply Euler's criterion with $a=-1$ then
\[
\chi_p(-1) \equiv (-1)^{(p-1)/2} \ \textnormal{mod}\ p.
 \]
Hence $\chi_p(-1)-(-1)^{(p-1)/2}$ is either 0 or $\pm2$ and is divisible by $p$, hence
\[
\chi_p(-1)=(-1)^{(p-1)/2}, 
\]
and so $\chi_p(-1)=1$ (respectively, $-1$) if and only if $(p-1)/2$ is even (respectively, odd) if and only if $p \equiv 1 \ \textnormal{mod}\ 4$ (respectively, $p \equiv -1 \ \textrm{mod}\ 4$). This verifies Theorem 2.4.

\emph{Proof of Theorem} 2.5. This is an interesting application of Wilson's theorem, which asserts that
\begin{equation*}
\textrm{if}\ q\ \textrm{is a prime then}\ (q-1)! \equiv -1 \ \textrm{mod}\ q, \tag{*}
\]
and was in fact first stated by Abu Ali al-Hasan ibn al-Haytham in 1000 AD, over 750 years before it was attributed to John Wilson, whose name it now bears. We will use Wilson's theorem to first prove Theorem 2.5; after that we then verify Wilson's theorem.

Suppose that $\chi_p(a)=1$, and so $x^2 \equiv a \mod p$ for some $x \in \mathbb{Z}$. Note now that $1= \textrm{gcd}(a, p)$ implies that $1= \textrm{gcd}(x^2, p)$, and so $1= \textrm{gcd}(x, p)$ ($p$ is prime!), hence by Fermat's little theorem,
\[
a^{(p-1)/2} \equiv (x^2)^{(p-1)/2}=x^{p-1} \equiv 1 \ \textrm{mod}\ p.
\]
 
Suppose that $\chi_p(a)=-1$, i.e., $a$ is a non-residue.  For each $i \in [1, p-1]$, there exists $j\in [1, p-1]$ uniquely determined by $i$, such that
\[
ij \equiv a \ \textrm{mod}\ p
\]
($\mathbb{Z}/p\mathbb{Z}$ is a field) and $i \not= j$ because $a$ is  a non-residue. Hence we can group the integers $1,\dots, p-1$ into $(p-1)/2$ pairs, each pair with a product $\equiv a \mod p$. Multiplying all of these pairs together yields
\[
(p-1)! \equiv a^{(p-1)/2} \ \textrm{mod}\ p,
\]
 and so $(*)$ implies that
 \[
-1 \equiv a^{(p-1)/2} \ \textrm{mod}\ p.  \]  

$\hspace{15cm} \textrm{QED}$

\vspace{0.3cm}

\emph{Proof of Wilson's theorem}. The implication $(*)$ is clearly valid when $q=2$, so assume that $q$ is odd. Use Proposition 1.2 to find for each integer $a \in [1, q-1]$ an integer $\bar{a} \in [1, p-1]$ such that $a\bar{a}\equiv 1$ mod $q$. The integers 1 and $q-1$ are the only integers in $[1, q-1]$ that are their own inverses mod $q$, hence we may group the integers from 2 through $q-2$ into $(q-3)/2$ pairs with the product of each pair congruent to 1 mod $q$. Hence
\[
2\cdot3\cdots(q-3)(q-2)\equiv 1\ \textrm{mod}\ q.\]
Multiplication of both sides of this congruence by $q-1$ yields
\[
(q-1)!=1\cdot2\cdots(q-1) \equiv q-1 \equiv -1\ \textrm{mod}\ q.\]
$\hspace{15.6cm} \textrm{QED}$ 

\emph{Remark}. The converse of Wilson's theorem is also valid.

\section{The Basic Problem and the Fundamental Problem for a Prime}
From our discussion in Chapter 1, if $d$ is the discriminant of $ax^2+bx+c$ and if neither $a$ nor $d$ is divisible by $p$ then
\[
ax^2+bx+c \equiv 0 \ \textrm{mod}\ p
\]
has a solution if and only if $d$ is a residue of $p$. This motivates what we will call the

\vspace{0.3cm}
\textbf{ Basic Problem}. If $d \in \mathbb{Z}$, for what primes $p$ is $d$ a quadratic residue of $p$?
\vspace{0.3cm}

We now present a strategy for solving this problem which employs Proposition 2.3 as the basic tool. Things can be stated precisely and concisely if we use the following

\emph{Notation}: if $z \in \mathbb{Z}$, let
\[
X_{\pm}(z)=\{p: \chi_p(z)= \pm 1\},
\]
\[
\pi_{\textrm{odd}}(z) (\textrm{resp.,}\ \pi_{\textrm{even}}(z))=\{q \in \pi(z): q\ \textrm{has odd (resp., even) multiplicity in}\ z\}.
\]
We point out here how to read the $\pm$ signs. If $\pm$ signs occur simultaneously in different places in an equation, formula, definition, etc., then the $+$ sign is meant to be taken simultaneously in all occurrences of $\pm$ and then the $-$ sign is also to be taken simultaneously in all occurrences of $\pm$. For example, the equation $X_{\pm}(z)=\{p: \chi_p(z)= \pm 1\}$ above asserts that $X_+(z)=\{p: \chi_p(z)=1\}$ and $X_{-}(z)=\{p: \chi_p(z)= -1\}$. We follow this convention in the sequel.

Suppose first that $d > 0$, with gcd$(d, p)=1$. If $\pi_{\textrm{odd}}(d)= \emptyset$ then $d$ is a square, so $d$ is trivially a residue of $p$. Hence assume that $\pi_{\textrm{odd}}(d) \not= \emptyset$. Proposition 2.3 implies that
\[
\chi_p(d)= \prod_{q \in \pi_{\textrm{odd}}(d)}\ \chi_p(q).
\]
Hence
\begin{equation*}
\chi_p(d)=1\ \textrm{iff}\ |\{q \in \pi_{\textrm{odd}}(d): \chi_p(q)=-1\}|\ \textrm{is \emph{even}}. \tag{1}
\]
Let
\[
\mathcal{E}=\{E \subseteq \pi_{\textrm{odd}}(d): |E|\ \textrm{is even} \}.
\]
If $E \in \mathcal{E}$, let $R_{E}$ denote the set of all $p$ such that 
\[
\chi_p(q)=\left\{\begin{array}{rl}-1, & \textrm{if $q \in E$,}\\
1,& \textrm{if $q \in \pi_{\textrm{odd}}(d) \setminus E$.}\\\end{array}\right.
\]
Then $(1)$ implies that
\begin{equation*}
X_+(d)= \Big( \bigcup_{E \in \mathcal{E}}\ R_E \Big) \setminus \pi_{\textrm{even}}(d), \tag{2}
\]
and this union is pairwise disjoint. Moreover
\begin{equation*}
R_E=\Big( \bigcap_{q \in E}\ X_-(q) \Big) \cap \Big( \bigcap_{q \in \pi_{\textrm{odd}}(d) \setminus E}\ X_+(q) \Big). \tag{3}
\]

Suppose next that $d<0$. Then $d=(-1)(-d)$, hence
\begin{equation*}
\chi_p(d)= \prod_{q \in \{-1\} \cup \pi_{\textrm{odd}}(d)}\ \chi_p(q). \tag{4}
\]
If we let
\[
\mathcal{E}_{-1}=\{E \subseteq \{-1\} \cup \pi_{\textrm{odd}}(d): |E|\ \textrm{is even} \} ,
\]
then by applying (4) and an argument similar to the one that yielded (2) and (3) for $X_+(d), d>0$, we also deduce that for $d<0$,
\begin{equation*}
X_+(d)= \Big( \bigcup_{E \in \mathcal{E}_{-1}}\ R_E \Big) \setminus \pi_{\textrm{even}}(d), \tag{5}
\]
where
\begin{equation*}
R_E=\Big( \bigcap_{q \in E}\ X_-(q) \Big) \cap \Big( \bigcap_{q \in(\{-1\} \cup \pi_{\textrm{odd}}(d)) \setminus E}\ X_+(q) \Big), E \in \mathcal{E}_{-1}. \tag{6}
\]

In order to show more concretely how this strategy for the solution of the basic problem is implemented, suppose as an example that we wish to determine $X_+(\pm126)$. First, factor $\pm126$ as $\pm2 \cdot 3^2 \cdot 7$. It follows from this factorization that

\vspace{0.3cm}
$\pi_{\textrm{odd}}(\pm126)= \{2, 7 \},\ \pi_{\textrm{even}}(\pm126)= \{3\},$

\vspace{0.3cm}
\noindent hence

\vspace{0.3cm}
$\mathcal{E}=\{ \emptyset, \{2, 7\} \},\ \mathcal{E}_{-1}= \{\emptyset, \{-1, 2\}, \{-1, 7\}, \{2, 7\} \}.$

\vspace{0.3cm}
\noindent It now follows from (2) and (3) that
\begin{eqnarray*}
X_+(126)&=& (R_{\emptyset} \cup R_{\{2, 7\}}) \setminus \{3\}\\
&=& \Big( \big(X_+(2) \cap X_+(7)\big) \cup \big(X_-(2) \cap X_-(7)\big) \Big) \setminus \{3\},
\end{eqnarray*}
and from (5) and (6) that
\begin{eqnarray*}
X_+(-126)&=& (R_{\emptyset} \cup R_{\{-1, 2\}} \cup R_{\{-1, 7\}} \cup R_{\{2, 7\}}) \setminus \{3\} \\
&=& \Big( \big(X_+(-1) \cap X_+(2) \cap X_+(7)\big) \cup \big(X_-(-1) \cap X_-(2) \cap X_+(7)\big)\\
&\cup&\big(X_-(-1) \cap X_+(2) \cap X_-(7)\big) \cup \big(X_+(-1) \cap X_-(2) \cap X_-(7)\big) \Big) \setminus \{3\}.\\\
\end{eqnarray*}
In order to finish this calculation of $X_+(\pm126)$, we must now calculate $X_+(2) \cap X_+(7),\ X_-(2) \cap X_-(7),\ X_+(-1) \cap X_+(2) \cap X_+(7),\ X_-(-1) \cap X_-(2) \cap X_+(7),\ X_-(-1) \cap X_+(2) \cap X_-(7),$ and $X_+(-1) \cap X_-(2) \cap X_-(7)$, which in turn requires the calculation of $X_{\pm}(-1),\ X_{\pm}(2)$, and $X_{\pm}(7)$. Theorem 2.4 and formulae (2),(3),(5), and (6) hence reduce the solution of the Basic Problem to the solution of the
\vspace{0.3cm}

\textbf{Fundamental Problem for Primes}. If $q$ is \emph{prime}, calculate $X_{\pm}(q)$. 

\vspace{0.3cm}
The Fundamental Problem for odd primes and Basic Problem and will be completely solved in sections 1 and 2 of Chapter 4. The Fundamental Problem for the prime 2 will be completely solved in the next section. 
\section{Gauss' Lemma and the Fundamental Problem for the Prime $2$}

The next theorem, along with Theorems 2.4 and 2.5, will be used many times in our subsequent work.
\begin{thm}
$\chi_p(2)= (-1)^{(p^2-1)/8}.$
\end{thm}
Theorem 2.6 solves the Fundamental Problem for the prime 2. It is easy to see that $(p^2-1)/8$ is even (odd)
if and only if $p \equiv 1$ or 7 mod 8 ($p \equiv $ 3 or 5 mod 8). Hence
\[
X_+(2)= \{p: p \equiv 1\ \textrm{or}\ 7 \  \textrm{mod}\ 8\},
\]
\[
X_-(2)= \{p: p \equiv 3\ \textrm{or}\ 5 \ \textrm{mod}\ 8\}.
\]

 The proof of Theorem 2.6 will use a basic result in the theory of quadratic residues called Gauss' lemma (this lemma was first used by Gauss in his third proof of the Law of Quadratic Reciprocity [20], which proof we will present in Chapter 3). We will first state Gauss' lemma, then use it to prove Theorem 2.6, and then we will prove Gauss' lemma. 

Toward that end, then, let $a \in \mathbb{Z}$, $\gcd(a, p)=1$. Consider the minimal positive ordinary residues mod $p$ of the integers $a,\dots, \frac{1}{2} (p-1)a$. None of these ordinary residues is $p/2$, as $p$ is odd, and they are all distinct as $\gcd(a, p)=1$, hence let
 \[
 u_1,\dots, u_s\ \textrm{be those ordinary residues that are}\ > p/2,
 \]
 \[
  v_1,\dots, v_t\ \textrm{be those ordinary residues that are}\ < p/2.
  \]
N.B. $s+t= \frac{1}{2} (p-1)$. We then have
\begin{thm}
$($Gauss' lemma$)$
\[
\chi_p(a)=(-1)^s.
\]
\end{thm}

\emph{Proof of Theorem} 2.6. Let $\sigma$ be the the number of minimal positive ordinary residues $\mod p$ of the integers in the set
\begin{equation*}
1 \cdot2, 2 \cdot2,\dots, \frac{1}{2} (p-1) \cdot 2 \tag{7} 
\]
that exceed $p/2$. Gauss' lemma implies that
\[
\chi_p(2)=(-1)^{\sigma}.
\]
Because each integer in (7) is less than $p$, $\sigma=$ the number of integers in the set (7) that exceed $p/2$. An integer $2j, j \in [1, (p-1)/2]$ does not exceed $p/2$ if and only if $1 \leq j \leq p/4$, hence the number of integers in (7) that do not exceed $ p/2$ is $[p/4]$, where $[x]$ denotes the greatest integer not exceeding $x$.  Hence
\[
\sigma=\frac{p-1}{2} - \Big[ \frac{p}{4} \Big].
\]
To prove Theorem 2.6, it hence suffices to prove that
\begin{equation*}
 \textrm{ for all odd integers}\ n, \frac{n-1}{2} - \Big[ \frac{n}{4} \Big] \equiv \frac{n^2-1}{8} \ \textrm{mod}\ 2 \tag{8}.
\]
To see this, note first that the congruence in (8) is true for a particular integer $n$ if and only if it is true for $n+8$, because
\[ 
\frac{(n+8)-1}{2}- \Big[ \frac{n+8}{4} \Big]= \frac{n-1}{2}+4- \Big( \Big[ \frac{n}{4} \Big] +2 \Big) \equiv \frac{n-1}{2}-\Big[ \frac{n}{4} \Big] \ \textrm{mod}\ 2,
\]
\[
\frac{(n+8)^2-1}{8}= \frac{n^2-1}{8}+2n+8 \equiv \frac{n^2-1}{8} \ \textrm{mod}\ 2.\]
Thus (8) holds if and only if it holds for $n= \pm1, \pm3$, and it is easy to check that (8) holds for these values of $n$. $\hspace{12.3cm} \textrm{QED}$

\vspace{0.3cm}
\emph{Proof of Theorem }2.7. Let $u_i, v_i$ be as defined before the statement of Gauss' lemma. We claim that
\begin{equation*}
\{p-u_1,\dots, p-u_s, v_1,\dots,v_t\}=[1, \frac12 (p-1)]. \tag{9}
\]
To see this, note first that if $i \not= j$ then $v_i \not= v_j, u_i \not= u_j$ hence $p-u_i \not= p-u_j$.  It is also true that $p-u_i \not= v_j$ for all $i, j$; otherwise $p \equiv a(k+l)$ mod $p$, where $2 \leq k+l \leq \frac{p-1}{2}+ \frac{p-1}{2}=p-1,$ which is impossible because gcd$(a, p)=1$. Hence
\begin{equation*}
|\{p-u_1,\dots, p-u_s, v_1,\dots,v_t\}|=s+t= \frac{p-1}{2}.\tag{10}
\]
But $0<v_i<p/2$ implies that $0<v_i \leq (p-1)/2$ and $p/2<u_i<p$, hence $0<p-u_i \leq (p-1)/2$, and so
\begin{equation*}
\{p-u_1,\dots, p-u_s, v_1,\dots,v_t\} \subseteq [1, \frac{1}{2} (p-1)]. \tag{11}
\]
As $|[1, \frac{1}{2} (p-1)]|= \frac{1}{2} (p-1)$, (9) follows from (10) and (11).

It follows from (9) that
\[
\prod_1^s\ (p-u_i) \prod_1^t\ v_i= \Big(\frac{p-1}{2} \Big)!.
\]
Because
\[
p-u_i \equiv -u_i \ \textrm{mod}\ p
\]
we conclude from the preceding equation that
\begin{equation*}
(-1)^s\prod_1^s\ u_i \prod_1^t\ v_i \equiv \Big(\frac{p-1}{2} \Big)! \ \textrm{mod}\ p. \tag{12}
\]
Because $u_1,\dots,u_s, v_1,\dots,v_t$ are the least positive ordinary residues of $a, \dots, \frac{1}{2} (p-1)a$, it is a consequence of (12) that
\begin{equation*}
(-1)^sa^{(p-1)/2} \Big(\frac{p-1}{2} \Big)! \equiv  \Big(\frac{p-1}{2} \Big)! \ \textrm{mod}\ p. \tag{13}
\]
But $p$ and $(\frac{p-1}{2})!$ are relatively prime, and so (13) implies that
\[
(-1)^sa^{(p-1)/2} \equiv 1 \ \textrm{mod}\ p
\]
i.e.,
\[
a^{(p-1)/2} \equiv (-1)^s \ \textrm{mod}\ p.
\]
By Euler's criterion (Theorem 2.5),
\[
a^{(p-1)/2} \equiv \chi_p(a) \ \textrm{mod}\ p,
\]
hence
\[
\chi_p(a) \equiv (-1)^s \ \textrm{mod}\ p.
\]
It follows that $\chi_p(a)-(-1)^s$ is either 0 or $\pm2$ and is also divisible by $p$ and so
\[
\chi_p(a)=(-1)^s.
\] 
\hspace{15.4cm} \textrm{QED}

We now need to solve the Fundamental Problem for odd primes. This will be done in Chapter 4 by using a result which Gauss called the \emph{theorema aureum}, the ``golden theorem", of number theory. We will discuss that result extensively in the next chapter.

\chapter{Gauss' \emph{Theorema Aureum}: the Law of Quadratic Reciprocity}

Proposition 1.1 of Chapter 1 shows that the solution of the general second-degree congruence $ax^2+bx+c\equiv 0$ mod $p$ for an odd prime $p$ can be reduced to the solution of the congruence $x^2\equiv b^2-4ac $ mod $p$, and we also saw how the solution of  $x^2\equiv n$ mod $m$ for a composite modulus $m$ can be reduced by way of Gauss' algorithm to the solution of $x^2\equiv q$ mod $p$ for prime numbers $p$ and $q$. In this chapter, we will discuss a remarkable theorem known as the \emph{Law of Quadratic Reciprocity}, which provides a very powerful method for determining the solvability of congruences of this last type.  The theorem states that if $p$ and $q$ are distinct odd primes then the congruences $x^2\equiv q$ mod $p$ and  $x^2\equiv p$ mod $q$ are either both solvable or both not solvable, unless $p$ and $q$ are both congruent to 3 mod 4, in which case one is solvable and the other is not. As a simple but nevertheless striking example of the power of this theorem, suppose one wants to know if $x^2\equiv 5$ mod 103 has any solutions. Since 5 is not congruent to 3 mod 4, the quadratic reciprocity law asserts that $x^2\equiv 5$ mod 103 and $x^2\equiv$ 103 mod 5 are both solvable or both not. But solution of the latter congruence reduces to $x^2\equiv 3$ mod 5, which clearly has no solutions. Hence neither does $x^2\equiv 5$ mod 103.

The first rigorous proof of the Law of Quadratic Reciprocity is due to Gauss. He valued this theorem so much that he referred to it as the \emph{theorema aureum}, the golden theorem, of number theory, and in order to acquire a deeper understanding of its content and implications, he searched for various proofs of the theorem, eventually discovering eight different ones. After discussing what type of mathematical principle a reciprocity law might seek to encapsulate in section 1 of this chapter, stating the Law of Quadratic Reciprocity precisely in section 2, and discussing some of the mathematical history which led up to it in section 3, we follow Gauss' example by presenting five different proofs of quadratic reciprocity in the remaining 10 sections. Each of these proofs is chosen to highlight the ideas behind the techniques which Gauss himself employed and to indicate how some of the more modern approaches to quadratic reciprocity are inspired by the work of Gauss. For a more detailed summary of what we do in sections 5-13, consult section 4.

\section{What is a reciprocity law?}

\vspace{0.4cm}
We will motivate why we would want an answer to the question entitling this section by first asking this question: what positive integers $n$ are the sum of two squares? This is an old problem that was solved by Fermat in 1640. We can reduce to the case when $n$ is prime by first observing, as Fermat did, that if a prime number $q$ divides a sum of two squares, neither of which is divisible by $q$, then $q$ is the sum of two squares. Using the identity
\[
(a^2+b^2)(c^2+d^2)=(ad-bc)^2+(ac+bd)^2,\] 
which shows that the property of being the sum of two squares is preserved under multiplication, it can then be easily shown that $n$ is the sum of two squares if and only if $n$ is either a square or each oddprime factor of $n$ of odd multiplicity is the sum of two squares. Because $2=1^2+1^2$, we hence need only consider odd primes $p$.

As we mentioned before, $p$ is the sum of two squares if it divides the sum of two squares and neither of the squares are divisible by $p$, and so we are looking for integers $a$ and $b$ such that
\[
a^2+b^2 \equiv 0\ \textrm{mod}\ p\]
and
\[
a\not \equiv 0 \not \equiv b\ \textrm{mod}\ p.\]
After multiplying the first congruence by the square of the inverse of $b$ mod $p$, it follows that $p$ is the sum of two squares if and only if the congruence
\[
x^2+1\equiv 0 \ \textrm{mod}\ p\]
has a solution, i.e., $-1$ is a residue of $p$. We now invoke Theorem 2.4 of Chapter 2, which asserts that $-1$ is a residue of $p$ if and only if $p\equiv$ 1 mod 4, to conclude that a positive integer $n$ is the sum of two squares if and only if either $n$ is a square or each odd prime factor of $n$ of odd multiplicity is congruent to 1 mod 4.

Another way of saying that the congruence $x^2+1\equiv 0 \ \textrm{mod}\ p$ has a solution is to say that the polynomial $x^2+1$ factors over $\mathbb{Z}/p\mathbb{Z}$ as $(x+c)(x-c)$ for some (nonzero) $c\in \mathbb{Z}$, i.e., $x^2+1$ \emph{splits} over $\mathbb{Z}/p\mathbb{Z}$ (in the remainder of this section, we follow the exposition as set forth in the very nice paper of B. F. Wyman [64]). Our previous discussion hence shows that the problem of deciding when an integer is the sum of two squares comes down to deciding when a certain monic polynomial with integer coefficients splits over $\mathbb{Z}/p\mathbb{Z}$. It is therefore of considerable interest to further study this splitting phenomenon. For that purpose, we will start more generally with a polynomial $f(x)$ with integral coefficients that is irreducible over $\mathbb{Q}$, and for an odd prime $p$, we let $f_p(x)$ denote the polynomial over $\mathbb{Z}/p\mathbb{Z}$ obtained from $f(x)$ by reducing all of its coefficients modulo $p$. We will say that $f(x)$ \emph{splits modulo p} if $f_p(x)$ is the product of distinct linear factors over  $\mathbb{Z}/p\mathbb{Z}$, and if $f(x)$ splits modulo $p$, we will call $p$ a \emph{slitting modulus of} $f(x)$. 

Suppose now that $f(x)=ax^2+bx+c$ is a quadratic polynomial, the case that is of most interest to us here. If $p$ is an odd prime then the congruence
\[
f(x)\equiv 0\ \textrm{mod}\ p\]
has either 0, 1, or 2 solutions, which occur, according to Proposition 1.1 of Chapter 1, if the discriminant $b^2-4ac$ of $f(x)$ is, respectively, a non-residue of $p$, is divisible by $p$, or is a residue of $p$. Because this congruence also has  exactly 2 solutions if and only if $f(x)$ splits modulo  $p$, it follows that $f(x)$ splits modulo $p$ if and only if the discriminant of $f(x)$ is a residue of $p$. We saw before that $x^2+1$ splits modulo $p$ if and only if $p\equiv$ 1 mod 4, and using Theorem 2.6 from Chapter 2 in a similar manner, one can prove that $x^2-2$ splits modulo $p$ if and only if $p\equiv$ 1 mod 8. Another amusing example, which we will let the reader work out, asserts that $x^2+x+1$ splits modulo $p$ if and only if $p\equiv$ 1 mod 3.

In light of these three examples, we will now, for a fixed prime $q$, look for the splitting moduli of $x^2-q$. We wish to determine these moduli by means of congruence conditions that are similar to the conditions which described the splitting moduli of $x^2+1$, $x^2-2$ and $x^2+x+1$. If $p$ is a prime then, over $\mathbb{Z}/p\mathbb{Z}$, $x^2-q$ is the square of a linear polynomial only if $p=q$, and also if $p=2$, hence we may assume that $p$ is an odd prime distinct from $q$. It follows that $x^2-q$ splits modulo $p$ if and only if $q$ is a quadratic residue of $p$, i.e., the Legendre symbol $\chi_p(q)$ is 1. Hence we must find a way to calculate  $\chi_p(q)$ as $p$ varies over the odd primes.

This translation of the splitting modulus problem for $x^2-q$ does not really help much. The Legendre symbol $\chi_p$ is not easy to evaluate directly, and changing the value of $p$ would require a direct calculation to begin again from scratch. Because there are infinitely many primes, this approach to the problem quickly becomes unworkable.

A way to possibly overcome this difficulty is to observe that in this problem $q$ is fixed while the prime $p$ varies, and so if it was possible to somehow use $\chi_q(p)$ in place of $\chi_p(q)$ then only one Legendre symbol would be required. Moreover, the values of  $\chi_q(p)$ are determined only by the ordinary residue class of $p$ modulo $q$, and so we would also have open the possibility of calculating the splitting moduli of $x^2-q$ in terms of ordinary residue classes determined in some way by $q$, as per the descriptions of the splitting moduli in our three examples. This suggests looking for a computationally efficient relationship between $\chi_p(q)$ and $\chi_q(p)$, i.e., is there a useful \emph{reciprocal} relation between the residues (respectively, non-residues) of $p$ and the residues (respectively, non-residues) of $q$? The answer: yes there is, and it is given by the \emph{Law of Quadratic Reciprocity}, one of the fundamental principles of elementary number theory and one of the most powerful tools that we have for analyzing the behavior of residues and non-residues. As we will see (in Chapter 4), it completely solves the problem of determining the splitting moduli of any quadratic polynomial by means of congruence conditions which depend only on the discriminant of the polynomial. 

 We will begin our study of quadratic reciprocity in the next section, but before we do that, it is natural to wonder if there is a  similar principle which can be used to study the splitting moduli of polynomials of degree larger than 2. Using what we have discussed for quadratic polynomials as a guide,  we will say that a polynomial $f(x)$ with integer coefficients \emph{satisfies a reciprocity law} if its splitting moduli are determined solely by congruence conditions which depend only on $f(x)$. This way of formulating these higher-degree reciprocity laws is the main reason that we used the idea of splitting moduli of polynomials in the first place. 
 
 As it turns out, higher reciprocity laws exit for many polynomials. A particularly nice class of examples are provided by the set of cyclotomic polynomials. There is a cyclotomic polynomial corresponding to each integer $n \geq 2$, defined by a primitive $n$-th root of unity, say $\zeta_n=\exp(2\pi i/n)$. The number $\zeta_n$ is algebraic over $\mathbb{Q}$, and is the root of a unique irreducible monic polynomial $\Phi_n(x)$ with integer coefficients of degree $\varphi(n)$, where $\varphi$ denotes Euler's totient function. The polynomial $\Phi_n(x)$ is the \emph{n-th cyclotomic polynomial}. For example, if $n=q$ is prime, one can show that
 \[
 \Phi_q(x)=1+x+\cdots +x^{q-1}\]
 (Chapter 3, section 8 ). The degree of $\Phi_n(x)$ is at least 4 when $n\geq 7$, and $\Phi_n(x)$ satisfies the following very nice reciprocity law (for a proof, consult Wyman [64]):
 \begin{thm}
 $($a cyclotomic reciprocity law$)$. The prime p is a splitting modulus of  $\Phi_n(x)$ if and only if $p\equiv 1\ \textnormal{mod}\ n$.
 \end{thm}
 
It transpires that not every polynomial with integer coefficients satisfies a reciprocity law, but there is an elegant way to characterize the polynomials with rational coefficients which do satisfy one. If $f(x)$ is a polynomial of degree $n$ with coefficients in $\mathbb{Q}$ then $f(x)$ has $n$ complex roots, counted according to multiplicity, and these roots, together with $\mathbb{Q}$, generate a subfield of the complex numbers that we will denote by $K_f$. The set of all field automorphisms of $K_f$ forms a group under the operation of composition of automorphisms. The \emph{Galois group of} $f(x)$ is defined to be the subgroup of all automorphisms $\sigma$ of $K_f$ which fix each rational number, i.e., $\sigma(r)=r$ for all $r\in \mathbb{Q}$. The next theorem neatly characterizes in terms of their Galois groups the polynomials which satisfy a reciprocity law.
\begin{thm}
$($Existence of Reciprocity Laws$)$. If $f(x)$ is a polynomial with integer coefficients and is irreducible over $\mathbb{Q}$ then $f(x)$ satisfies a reciprocity law if and only if the Galois group of $f(x)$ is abelian.
\end{thm}

The polynomials $x^4+4x^2+2$, $x^4-10x^2+4$, $x^4-2$, and $x^5-4x+2$ are all irreducible over $\mathbb{Q}$, and one can show that their Galois groups are, respectively, the cyclic group of order 4, the Klein 4-group, the dihedral group of order 8, and the symmetric group on 5 symbols (Hungerford [29], section V.4).  We hence conclude from Theorem 3.2 that $x^4+4x^2+2$ and $x^4-10x^2+4$ satisfy a reciprocity law, but $x^4-2$ and $x^5-4x+2$ do not.

Two natural questions now arise: how do you prove Theorem 3.2, and if you have an irreducible polynomial with integer coefficients and an abelian Galois group, how do you find the congruence conditions which determine its splitting moduli? The answers to these questions are far beyond the scope of what we will do in these lecture notes, because they make use of essentially all of the machinery of class field theory over the rationals. We will not even attempt an explanation of what class field theory is, except to say that it originated in a program to find reciprocity laws which are similar in spirit to the reciprocity laws for polynomials that we have discussed here, but which are valid in much greater generality. This program, which began with the work of Gauss on quadratic reciprocity, was eventually completed in the 1920's and 30's by Tagaki, E. Artin, Furtw$\ddot {\textrm{a}}$ngler, Hasse, and Chevalley. We now turn to the theorem which inspired all of that work.
 
\vspace{0.4cm}
\section{The Law of Quadratic Reciprocity}

\vspace{0.4cm}
\begin{thm}
$($Law of Quadratic Reciprocity $($LQR$))$ If $p$ and $q$ are distinct odd primes then
\[
\chi_p(q) \chi_q(p)=(-1)^{\frac{1}{2}(p-1) \frac{1}{2}(q-1)}.
\]
\end{thm}

We will begin our study of the the LQR by unpacking the information that is encoded in the elegant and efficient way by which Theorem 3.3 states it. Note first that if $n \in \mathbb{Z}$ is odd then $\frac{1}{2}(n-1)$ is even (respectively, odd) if and only if $n \equiv 1$ mod 4 (respectively, $n \equiv 3$ mod 4). Hence
\[
\chi_p(q) \chi_q(p)=1\ \textrm{iff}\ p\ \textrm{or}\ q \equiv 1 \ \textrm{mod}\ 4,
\]
\[
 \chi_p(q) \chi_q(p)=-1\ \textrm{iff}\ p \equiv q \equiv 3 \ \textrm{mod}\ 4,
 \]
 i.e., 
\[
\chi_p(q)= \chi_q(p)\ \textrm{iff}\ p\ \textrm{or}\ q \equiv 1 \ \textrm{mod}\ 4,
\]
\[
 \chi_p(q)=- \chi_q(p)\ \textrm{iff}\ p \equiv q \equiv 3 \ \textrm{mod}\ 4.
 \]
Thus the LQR states that
\[
\textrm{if $p$ or $q$} \equiv 1 \ \textrm{mod}\ 4\ \textrm{then $p$ is a residue of $q$ if and only if $q$ is a residue of $p$},
\]
and
\[
\textrm{if}\ p \equiv q \equiv 3 \ \textrm{mod}\ 4\ \textrm{then $p$ is a residue of $q$ if and only if $q$ is a non-residue of $p$.}
\]
This is why the theorem is called the law of quadratic \emph{reciprocity}. The classical quotient notation for the Legendre symbol makes the reciprocity typographically explicit: in that notation,  the conclusion of Theorem 3.3 reads as $\displaystyle\left( \frac{p}{q}\right)\left(\frac{q}{p}\right)=(-1)^{\frac{1}{2}(p-1) \frac{1}{2}(q-1)}$.

We next illustrate the usefulness of the LQR in determining whether or not a specific integer is or is not the residue of a specific prime. We can do no better than taking the example which Dirichet used himself in his landmark text \emph{Vorlesungen $\ddot{\textrm{u}}$ber Zahlentheorie} [12]. We wish to know whether 365 is a residue of the prime 1847. The first step is to factor 365 $=5\cdot73$, so that
\[
\chi_{1847}(365)=\chi_{1847}(5)\ \chi_{1847}(73).\]
Because $5\equiv$ 1 mod 4, the LQR implies that
\[
\chi_{1857}(5)=\chi_5(1857)\]
and as $1857\equiv$ 2 mod 5, it follows that
\[
\chi_{1857}(5)=\chi_5(2)=-1.\]
Since $73\equiv$ 1 mod 4, it follows in the same manner from the LQR and the fact that $1847\equiv$ 22 mod 73 that
\[
\chi_{1847}(73)=\chi_{73}(1847)=\chi_{73}(22)=\chi_{73}(2)\ \chi_{73}(11).\]
But now $73\equiv$ 1 mod 8, hence it follows from Theorem 2.6 that
\[
\chi_{73}(2)=1\]
hence
\[
\chi_{1847}(73)=\chi_{73}(11).\]
Using the LQR once more, we have that
\[
\chi_{73}(11)=\chi_{11}(73)=\chi_{11}(7),\]
and because 7 and 11 are each congruent to 3 mod 4, it follows from the LQR that
\[
\chi_{11}(7)=-\chi_7(11)=-\chi_7(4)=-\chi_7(2)^2=-1.\]
Consequently,
\[
\chi_{1847}(73)=\chi_{73}(11)=\chi_{11}(7)=-1,\]
and so finally,
\[
\chi_{1847}(365)=\chi_{1847}(5)\ \chi_{1847}(73)=(-1)(-1)=1.\]
Thus 365 is a residue of 1847; in fact
\[
(\pm496)^2=246016=365+133\cdot1847.\]

Quadratic reciprocity can also be used to calculate the splitting moduli of polynomials of the form $x^2-q$, $q$ a prime, as we alluded to in section 1 above. For example, let $q=5$. Then the residues of 5 are 1 and 4 and so 
\[
\chi_5(1)=\chi_5(4)=1\] 
and 
\[
\chi_5(2)=\chi_5(3)=-1.\]
Hence
\[
\chi_5(p)=1\ \textrm{iff}\ p\equiv 1\ \textrm{or}\ 4\ \textrm{mod}\ 5.\]
Because $5\equiv 1$ mod 4, it follows from the LQR that
\[
\chi_p(5)=\chi_5(p),\]
hence
\[
5\ \textrm{is a residue of}\ p\ \textrm{iff} \ p\equiv 1\ \textrm{or}\ 4\ \textrm{mod}\ 5.\]
Consequently, $x^2-5$ splits modulo $p$ if and only if $p$ is congruent to either 1 or 4 mod 5. 

For a different example, take $q=11$. Then calculation of the residues of 11 shows that 
\[
\chi_{11}(p)=1\ \textrm{iff}\ p\equiv \ \textrm{1, 3, 4, 5, or 9}\ \textrm{mod}\ 11.\]  
We have that $11\equiv$ 3 mod 4, hence by the LQR,
\[
\chi_p(11)=\pm\chi_{11}(p),\]
with the sign determined by the equivalence class of $p$ mod 4. For example, if $p=23$ then $23\equiv 1$ mod 11 and $23 \equiv 11\equiv 3$ mod 4, hence the LQR implies that 
\[
\chi_{23}(11)=-\chi_{11}(23)=-\chi_{11}(1)=-1,\]
while if $p=89$ then $89\equiv$ 1 mod 11 and $89\equiv$ 1 mod 4, and so the LQR implies in this case that
\[
\chi_{98}(11)=\chi_{11}(89)=\chi_{11}(1)=1.\]
Use of the Chinese remainder theorem shows that the value of $\chi_p(11)$ depends on the equivalence class of $p$ modulo $4\cdot 11=44$, and after a few more calculations we see that 
\[
\chi_p(11)=1\ \textrm{iff}\ p\equiv \ \textrm{1, 5, 7, 9, 19, 25, 35, 37, 39, or 43 }\ \textrm{mod}\ 44.\]  
Thus $x^2-11$ splits modulo $p$ if and only if $p \equiv $ 1, 5, 7, 9, 19, 25, 35, 37, 39, or 43 mod 44. We will have much more to say about the utility of quadratic reciprocity in Chapter 4, but these examples already give a good indication of how the LQR makes computation of residues and non-residues much easier.

\vspace{0.4cm}
\section{Some History}

\vspace{0.4cm}

At the end of section 1, we indicated very briefly that many important and far-reaching developments in number theory can trace their genesis to the Law of Quadratic Reciprocity. Thus it is instructive to discuss the history of some of the ideas in number theory which led up to it. In order to do that, we will follow the account of that story as presented in Lemmermeyer [38]. Lemmermeyer's book contains a wealth of information about the development of reciprocity laws in many of their various manifestations, with a penetrating and comprehensive analysis, both mathematical and historical, of the circle of ideas, techniques, and approaches which have been brought to bear on that subject.

The first foreshadowing of quadratic reciprocity appears in the work of Fermat.  Fermat's results on the representation of integers as the sum of two squares, as we saw in section 1, leads directly to the problem of determining the quadratic character of $-1$, which was solved in Chapter 2 by Theorem 2.4. Fermat also studied the representation of primes by quadratic forms of the form $x^2+ny^2$, for $n=\pm2, \pm3$, and $-5$. One can show that when $n=2$ or $3$ then a prime $p$ which divides $x^2+ny^2$ but divides neither $x$ nor $y$ is of the form $a^2+nb^2$ for a pair of integers $a$ and $b$. It now follows from this fact, using the same reasoning that we employed in section 1, that $p$ can be represented by either  $x^2+2y^2$ or $x^2+3y^2$ if and only if $-2$ or, respectively, $-3$, is a residue of $p$. We hence see that the quadratic character of $\pm2$ and $\pm3$ is also implicit in Fermat's work on quadratic forms. 

Euler apparently began to seriously study Fermat's work when he started his mathematical correspondence with Christian Goldbach in 1729. As a result of that study, Euler became interested in the divisors of integers which are represented by quadratic forms $nx^2+my^2$, which eventually led him after several years to the Law of Quadratic Reciprocity. The LQR was first conjectured by Euler [15] in an equivalent form in 1744, based on extensive numerical evidence, but he could not prove it. Research done by Lagrange during the years from 1773 to 1775, in particular his work on a general theory of binary quadratic forms, inspired Euler to return to the study of quadratic residues, and in a paper [17] published in 1783 after his death, Euler gave, still without proof, a formulation of the LQR that is very close to that which is used most commonly today. 

Euler's original formulation of quadratic reciprocity in 1744  can be stated (using more modern notation) as follows:
\begin{thm}
Suppose that $p$ is an odd prime and a is a positive integer not divisible by p. If q is a prime such that $p\equiv \pm q\ \textnormal{mod}\ 4a$ then $\chi_p(a)=\chi_q(a)$.
\end{thm}

\noindent This says that the value of the Legendre symbol $\chi_p(a)$ depends only on the ordinary residue class of $p$ modulo $4a$, and that the value of  $\chi_p(a)$ is the same for all primes $p$ with a fixed remained $r$ or $4a-r$ when divided by $4a$.

In section 6 below, we will deduce the LQR by proving Theorem 3.4 directly by means of Gauss' Lemma (Theorem 2.7), and so it 
makes sense to verify that the LQR is equivalent to Theorem 3.4, which is what we will do now.
\begin{prp}
The Law of Quadratic Reciprocity is equivalent to Theorem $3.4$.
\end{prp}

\emph{Proof}. Suppose that Theorem 3.4 is true. Assume that $p>q$ are odd primes. We need to verify that
\begin{equation*}
\chi_p(q) \chi_q(p)=(-1)^{\frac{1}{2}(p-1)\frac{1}{2}(q-1)}\tag{1}.\]

Suppose first that $p\equiv q$ mod 4, and let $a=(p-q)/4>0$.  Then we have that
\[
\chi_q(p)=\chi_q(p-q)=\chi_p(a)\]
and
\[
\chi_p(q)=\chi_p(-1)\chi_p(p-q)=\chi_p(-1)\chi_p(a).\]
Because $p$ does not divide $a$, it follows from Theorem 3.4 that
\[ 
\chi_p(a)=\chi_q(a),\]
and so
\[
\chi_q(p)\chi_p(q)=\chi_p(-1).\]
Equation (1) is now a consequence of this equation, Theorem 2.4, and the assumption that $p\equiv q$ mod 4.

On the other hand, if $p\not \equiv q$ mod 4 then $p\equiv -q$ mod 4, hence we set $a=(p+q)/4>0$ and so deduce that
\[
\chi_q(p)=\chi_q(p+q)=\chi_q(a)\]
and
\[
\chi_p(q)=\chi_p(p+q)=\chi_p(a).\]
As $p$ does not divide $a$, we conclude by way of Theorem 3.4 that
\[
\chi_q(p)\chi_p(q)=1,\]
and (1) follows from this equation because  $p\equiv -q$ mod 4.

In order to verify the converse, we assume that the LQR is valid and then for an odd prime $p$, a positive integer $a$ not divisible by $p$, and a prime $q$ for which $p\equiv \pm q\ \textnormal{mod}\ 4a$, we seek to verify that
\begin{equation*}
\chi_p(a)=\chi_q(a)\tag{2}.\]

By virtue of the multiplicativity of the Legendre symbol (Proposition 2.3$(ii)$ of Chapter 2), we immediately reduce to the case  that $a$ is prime. Suppose first that $a=2$. Then $p\equiv \pm q$ mod 8, hence (2) is true by Theorem 2.6. 

Now let $a$ be a fixed odd prime, and assume first that $p\equiv q$ mod $4a$. Then $a$ is neither $p$ nor $q$, $p\equiv q$ mod $a$, and $p\equiv q$ mod $4$. Hence 
\[
\chi_a(q)\chi_a(p)=1\]
and
\[
\frac{p-1}{2}+\frac{q-1}{2}\equiv 0\ \textrm{mod}\ 2.\]
It hence follows from the LQR that
\begin{eqnarray*}
\chi_p(a)\chi_q(a)&=&(-1)^{\frac{1}{2}(a-1)\cdot \frac{1}{2}(p-1)}(-1)^{\frac{1}{2}(a-1)\cdot \frac{1}{2}(q-1)} \chi_a(p)\chi_a(q)\\
&=&(-1)^{\frac{1}{2}(a-1)\big({\frac{p-1}{2}}+\frac{q-1}{2}\big)}\\
&=&1,
\end{eqnarray*}
i.e., (2) is valid.

Assume next that  $p\equiv -q$ mod $4a$. Then $p\equiv -q$ mod 4 and $p\equiv -q$ mod $a$, hence 
\[
\frac{p-1}{2}+\frac{q-1}{2}\equiv 1\ \textrm{mod}\ 2,\]
and, because of Theorem 2.4,
\[
\chi_a(p)\chi_a(q)=(-1)^{\frac{1}{2}(a-1)}.\]
As $a$ is again neither $p$ nor $q$, from the LQR it therefore follows that
\begin{eqnarray*}
\chi_p(a)\chi_q(a)&=&(-1)^{\frac{1}{2}(a-1)\cdot \frac{1}{2}(p-1)}(-1)^{\frac{1}{2}(a-1)\cdot \frac{1}{2}(q-1)} \chi_a(p)\chi_a(q)\\
&=&(-1)^{\frac{1}{2}(a-1)\big(1+{\frac{p-1}{2}}+\frac{q-1}{2}\big)}\\
&=&1
\end{eqnarray*}
in this case as well. $\hspace{11.7cm}\ \textrm{QED}$

The formulation of the LQR which appears in modern texts, including the one before the reader, was introduced by Legendre [36] in 1785. In this paper and in his influential book [37], Legendre discussed the LQR at length and in depth. In particular, he gave a proof of the LQR that depended on the assumption of primes which satisfy certain auxiliary conditions, but as Legendre was unable to rigorously verify   these conditions, as he himself admitted ([36], p. 520), his argument is not complete. Interestingly enough, one of the auxiliary conditions posited the existence of infinitely many primes in any arithmetic progression whose initial term and common difference are relatively prime, a very important result which Dirichlet would prove in 1837, and which we will have much more to say about later (see section 4   of Chapter 4). Legendre introduced his symbol in [36] as a particularly elegant way to state the LQR.

Lagrange's role in the history of quadratic reciprocity also needs to be mentioned, and comes from his work [33], [34] in refining and generalizing Euler's work on the representation of integers by quadratic forms.  For example, Lagrange proved that if $p$ is a prime which is congruent to either 1 or 9 mod 20 then $p=x^2+5y^2$, a result that was conjectured without proof by Euler, and he also showed that if $p$ and $q$ are primes each of which are congruent to 3 or 7 mod 20 then $pq=x^2+5y^2$. From these two results the quadratic character of $-5$ can be deduced. Lagrange also used quadratic residues to study the existence of nontrivial solutions to certain Diophantine equations of the form $ax^2+by^2+cz^2=0$.

The important role which the theory of quadratic residues played in the work of Euler and Lagrange on quadratic forms and the extensive discussion of the LQR by Legendre all contributed to making the proof of the LQR  one of the major unsolved problems of number theory in the eighteenth century. The first rigorous and correct proof was discovered by Gauss in 1796. He considered this result one of his greatest contributions to mathematics, returning to it again and again throughout his career. Gauss eventually found eight different proofs of the LQR and published six of them (although the two unpublished proofs, usually referred to as proofs 7 and 8, are not complete, according to D. Gr$\ddot{\textrm{o}}$ger [26]; I am grateful to an anonymous referee  for this reference). A major goal of Gauss' later work in number theory was to generalize quadratic reciprocity to higher powers, in particular to cubic and biquadratic (fourth-power) residues, and he sought out different ways of verifying quadratic reciprocity in order to gain as much insight as he could into ways to achieve that goal. We will discuss with a bit more detail in the next section the six proofs of the LQR which Gauss published. 

The establishment of generalizations of quadratic reciprocity that covered arbitrary power residues, the so-called \emph{higher reciprocity laws}, was a major theme of number theory in the nineteenth century and led to many of the most important advances in the subject during that time. Further generalizations to number systems extending beyond the integers, in particular and most importantly, to rings of algebraic integers in algebraic number fields (see section 8 of this chapter and section 1 of Chapter 5 for the relevant definitions), was a major theme of twentieth-century number theory and led to many of the most important advances during \emph{that} time. For an especially apt example of this latter development which closely follows the theme of this chapter, we direct the reader's attention to Hecke's penetrating analysis of quadratic reciprocity in an arbitrary algebraic number field ([27], Chapter VIII). 
\section{Proofs of the Law of Quadratic Reciprocity}
The Law of Quadratic Reciprocity has inspired more proofs, by far, than any other theorem of number theory. Lemmermeyer [38] documents 196 different proofs which have appeared up to the year 2000 and the web site at http://www.rzuser.uni-heidelberg.de/$\sim$hb3/

\noindent fchrono.html summarizes 266 proofs, with still more being found (I am grateful to an anonymous referee  for this web-site reference). Many of these arguments have been used as a basis for developing new methods and opening up new areas of study in number theory. In this book we will give seven different proofs of the LQR, with our primary emphasis being on the six proofs which Gauss himself published. As an introduction to what we will do in that regard, we will briefly discuss next these six proofs of Gauss.

Gauss' first proof of the LQR, which involved an extremely long and complicated induction argument, was published in \emph{Disquisitiones Arithmeticae} ([19], articles 135-145) (for a very readable account of a simplified version of Gauss' argument, see Dirichlet [12], Chapter 3, sections 48-51). It is similar to a proof attempted by Legendre, in that it also requires an ``auxiliary" prime. The complexity of Gauss' argument stems from the necessity of rigorously establishing the existence of this prime, and the formidable technical calculations which Gauss had to use there caused his argument to garner little attention for many years thereafter. However, those calculations were found to be useful in the development of algebraic $K$-theory in the 1970's; in fact a proof of quadratic reciprocity can be deduced from certain results in the $K$-theory of the rational numbers (Rosenberg [49], Theorem 4.4.9 and Corollary 4.4.10)! 

Gauss' second proof also appeared in the \emph{Disquisitiones} (article 262) and uses his genus theory of quadratic forms, a classification of forms that is closely related to Lagrange's classification of quadratic forms by means of unimodular substitutions (see section 12 of this chapter and also section 3 of chapter 8). The main point of the argument here is the verification of an inequality for the number of genera. This proof has a very nice modern formulation using the concept of narrow equivalence of ideals in a quadratic number field, and we will present this argument in sections 10, 11, and 12 below. 

Gauss' Lemma (Theorem 2.7) is the main idea underlying the third and fifth proofs ([20] and [22], respectively) which Gauss gave of the LQR. In the next section, we will present a simplification of Gauss' third proof due to Eisenstein; it is one of the most elegant and elementary proofs that we have of quadratic reciprocity and it has become the standard argument. In section 6 , Theorem 3.4 will be deduced directly from Gauss' Lemma, thereby proving the LQR after an appeal to Proposition 3.5; this approach is in a spirit similar to the idea use by Gauss in his fifth proof.

We make our first contact in these notes with the theory of algebraic numbers in our account in sections 7, 8 and 9 of Gauss' sixth proof [23] of the LQR. The main idea in this proof and Gauss' fourth proof [21] is the employment of quadratic Gauss sums (without and, respectively, with the determination of the signs). Gauss did not himself make use of algebraic number theory in any modern sense but we will use it in order to more clearly arrange the rather ingenious calculations that are a signal feature of these  techniques. We end in section 13 with a proof of quadratic reciprocity by means of the Galois theory of cyclotomic fields, an approach which foreshadows the way to more general formulations of reciprocity laws that are the subject of research in number theory today.

\section{A Proof of Quadratic Reciprocity via Gauss' Lemma}

The first proof that we will give of Theorem 3.3 is a simplification, due to Eisenstein, of Gauss' third proof [20]. It is by now the standard argument and uses an ingenious application of Theorem 2.7 (Gauss' lemma). Theorem 2.7 enters the reasoning by way of

 \begin{lem}
If $a \in \mathbb{Z}$ is odd and $\gcd(a, p)=1$ then
\[
\chi_p(a)=(-1)^{S(a, p)},
\]
where
\[
S(a, p)= \sum_{k=1}^{\frac{1}{2}(p-1)}\ \Big[\frac{ka}{p} \Big].
\]
\end{lem}

 Assume Lemma 3.6 for the time being, with its proof to come shortly.  

We begin our first proof of Theorem 3.3 by outlining the strategy of the argument. Let $p$ and $q$ be distinct odd primes and consider the set $L$ of points $(x, y)$ in the plane, where $x, y \in [1, \infty), 1 \leq x \leq \frac{1}{2}(p-1),$ and $1 \leq y \leq \frac{1}{2}(q-1)$, i.e., the set of \emph{lattice points} inside the rectangle with corners at $(0, 0), (0,  \frac{1}{2}(q-1)), (\frac{1}{2}(p-1), 0), (\frac{1}{2}(p-1),  \frac{1}{2}(q-1))$.

Let $l$ be the line with equation $qx=py$. To prove Theorem 3.3, one shows first that 
\begin{equation*}
\textrm{no point of $L$ lies on $l$}. \tag{3}
\]
Hence
\begin{eqnarray*}
L&=&\ \textrm{set of all points of $L$ which lie below $l$} \cup\ \textrm{set of all points of $L$ which lie above $l$}\\
&=& L_1 \cup L_2,
\end{eqnarray*}
consequently
\begin{equation*}
\frac{1}{2}(p-1) \frac{1}{2}(q-1)=|L|= |L_1| + |L_2|. \tag{4}
\]
This geometry is illustrated in Figure 1.
 \begin{figure}[h]
  \centering
  \begin{tikzpicture}
    \draw[name path=xaxis, thick, ->] (0,0) -- (5,0);
    \draw[name path=yaxis, thick, ->] (0,0) -- (0,6);
    \draw (0,5) node[anchor=east] {$\frac12 (q-1)$};
    \draw (3,0) node[anchor=north] {$\frac12 (p-1)$};
    \foreach \x in {1,2,3} {
      \draw (\x, -2pt) -- (\x, 2pt);
      \foreach \y in {1,2,3,4,5} {
        \draw (-2pt, \y) -- (2pt, \y);
        \fill (\x, \y) circle (1pt);
      }
    }
    \draw[thick, ->] (0,0) -- (3.8, 11*3.8 / 7) node[anchor=west] {$py = qx$};
  \end{tikzpicture}
  \caption{The lattice-point count}
  \label{fig:3.1}
\end{figure}
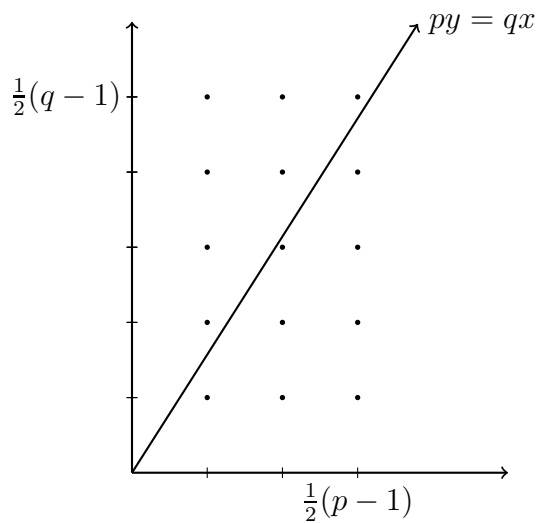

 The next step is to 
 \begin{equation*}
 \textrm{count the number of points in $L_1$ and $L_2$}. \tag{5} 
 \]
 The result is
 \[
 |L_1|=S(q, p),\ |L_2|=S(p, q),
 \]
 hence from (4),
 \[
 \frac{1}{2}(p-1) \frac{1}{2}(q-1)=S(q, p)+S(p, q).
 \]
 It then follows from Lemma 3.3 that
 \[
 (-1)^{\frac{1}{2}(p-1) \frac{1}{2}(q-1)}=(-1)^{S(q, p)}(-1)^{S(p, q)}=\chi_p(q) \chi_q(p),
 \] 
which is the conclusion of Theorem 3.3. Thus, we need to verify (3), implement (5), and prove Lemma 3.6.

\emph{Verification of} (3). Suppose that $(x, y) \in L$ satisfies $qx=py$. Then $q$, being prime, must divide either $p$ or $y$. Because $p$ is prime and $q \not= p$, $q$ must divide $y$, which is not possible because $1 \leq y \leq \frac{1}{2}(q-1)<q$.

\emph{Implementation of} (5). 
\begin{eqnarray*}
L_1&=&\{(x, y) \in L: qx>py\}\\
&=&\{(x, y) \in L: 1 \leq x \leq \frac{1}{2}(p-1), 1 \leq y < \frac{qx}{p} \}\\
&=&\bigcup_{1 \leq x \leq \frac{1}{2}(p-1)}\ \{(x, y): 1 \leq y \leq \Big[\frac{qx}{p} \Big] \},
\end{eqnarray*}
and this union is pairwise disjoint. Hence
\[
|L_1|=\sum_{x=1}^{\frac{1}{2}(p-1)}\ \Big[\frac{qx}{p} \Big]=S(q, p).\]
In Figure 1, $L_1$ is the set of lattice points which lie below the line $py=qx$.

On the other hand,
\begin{eqnarray*}
L_2&=&\{(x, y) \in L: qx<py\}\\
&=&\{(x, y) \in L: 1 \leq y \leq \frac{1}{2}(q-1), 1 \leq x< \frac{py}{q} \}\\
&=&\bigcup_{1 \leq y \leq \frac{1}{2}(q-1)}\ \{(x, y): 1 \leq x \leq \Big[\frac{py}{q} \Big] \},
\end{eqnarray*}
hence
\[
|L_2|=\sum_{y=1}^{\frac{1}{2}(q-1)}\ \Big[\frac{py}{q} \Big]=S(p, q).
\]
Here, $L_2$ is the set of lattice points in Figure 1 which lie above the line $py=qx$.
 
 Note that this part of the proof contains no number theory but is instead a purely geometric lattice-point count. All of the number theory is concentrated in the proof of Lemma 3.6, which is still to come. Indeed, that is the main idea in Gauss' third proof: divide the argument into two parts, a \emph{number-theoretic part} (Lemma 3.6) and a \emph{ geometric part} (the lattice-point count). Coupling geometry to number theory is a very powerful method for proving things, which Gauss pioneered in much of his work.

\emph{Proof of Lemma} 3.6. We set up shop in order to apply Gauss' lemma: take the minimal positive ordinary residues mod $p$ of the integers $a,\dots, \frac{1}{2} (p-1)a$, observe as before that none of these ordinary residues is $p/2$, as $p$ is odd, and they are all distinct as $\gcd(a, p)=1$, hence let
 \[
 u_1,\dots, u_s\ \textrm{be those ordinary residues that are}\ > p/2,
 \]
 \[
  v_1,\dots, v_t\ \textrm{be those ordinary residues that are}\ < p/2.
  \]
By the division algorithm, for each $j \in [1, \frac{1}{2}(p-1)]$,
\[
ja=p \Big[ \frac{ja}{p} \Big]+r_j,
\]
\[
r_j=\textrm{a}\ u_k\ \textrm{or a}\ v_l.
\]
Adding these equations together, we obtain
\begin{equation*}
a\sum_{j=1}^{\frac{1}{2}(p-1)}\ j=p\sum_{j=1}^{\frac{1}{2}(p-1)}\  \Big[\frac{ja}{p} \Big]+ \sum_{j=1}^s\ u_j+ \sum_{j=1}^t\ v_j. \tag{6}
\]
Next, recall from (9) of Chapter 2 that
\[
\{p-u_1,\dots, p-u_s, v_1,\dots,v_t\}=[1, \frac{1}{2} (p-1)].
\]
Hence
\begin{equation*}
\sum_{j=1}^{\frac{1}{2}(p-1)}\ j=sp-\sum_{j=1}^s\ u_j+ \sum_{j=1}^t\ v_j. \tag{7}
\]
Subtracting  (7) from (6) yields
\[
(a-1)\sum_{j=1}^{\frac{1}{2}(p-1)}\ j=pS(a, p)-sp+2\sum_{j=1}^s\ u_j.
\]
Hence
\[
p(S(a, p)-s)\ \textrm{is even ($a$ is odd!)}, 
\]
and so 
\[
S(a, p)-s \ \textrm{is even ($p$ is odd!)},
\]
whence
\[
(-1)^{S(a, p)}=(-1)^s.
\]
Gauss' lemma now implies that
\[
\chi_p(a)=(-1)^s,
\]
and so
\[
\chi_p(a)=(-1)^{S(a, p)}.\]
$\hspace{15.6cm} \textrm{QED}$

\section{Another Proof of Quadratic Reciprocity via Gauss' Lemma}

The proof of Theorem 3.3 that we will do here is similar to the first proof that we presented in section 5, in that it also uses Gauss' Lemma as a crucial tool. However, the strategy for the argument in this section  is different from the one that was used in section 5; we will deduce Theorem 3.4, Euler's version of Theorem 3.3, directly from Gauss' Lemma, and then use the equivalence of Theorem 3.4 and Theorem 3.3 that was established in Proposition 3.5 to conclude that Theorem 3.3 is valid (the details of the following argument, as well as those of the proof of Proposition 3.5, are taken from some lecture notes of F. Lemmermeyer posted at www.fen.bilkent.edu.tr/$\sim$franz/nt/ch6.pdf).

We proceed to implement this strategy. Let $p>q$ be odd primes, let $a$ be a positive integer not divisible by $p$, and suppose that $p\equiv \pm q$ mod $4a$. We wish to show that
\[
\chi_q(a)=\chi_p(a).\]

Suppose that $p\equiv q$ mod $4a$(when  $p\equiv  -q$ mod $4a$, a straightforward modification of the following argument, which we will leave to the interested reader, can be used to also verify that $\chi_q(a)$ and $\chi_p(a)$ are the same). In order to apply Gauss' Lemma to calculate $\chi_p(a)$, we need to find (the parity of) the number $s$ of all integers $a, 2a,\dots, \big((p-1)/2\big)a$ which have positive minimal ordinary residue modulo $p$ between $p/2$ and $p$. If
\[
(x)=x-[x]\]
denotes the fractional part of a real number $x$, then we observe that $s$ is the cardinality of the set
\[
P(a)=\left\{r\in [1, (p-1)/2]: \left(arp^{-1}\right)>\frac{1}{2}\right\}.\]
This set is the pairwise disjoint union of sets $P_m(a), m\in[0, a-1]$, where $P_m(a)$ consists of all integers $r\in [1, (p-1)/2]$ such that 
\[
\frac{m}{2}<\frac{ar}{p}<\frac{m+1}{2}p\ \textrm{and}\ \left(arp^{-1}\right)>\frac{1}{2}.\]
From this definition of $P_m(a)$, we conclude that $P_m(a)$ is empty when $m$ is even, and we also observe that the condition $(arp^{-1})> \frac{1}{2}$ is automatically satisfied when $m$ is odd.
Hence
\[
P_m(a)=\left\{r\in \mathbb{Z}: \frac{m}{2a}p<r<\frac{m+1}{2a}p\right\},\]
when $m$ is odd, and
\[
s=\sum_{0\leq m<a,\ m\ \textrm{odd}} |P_m(a)|.\]

Similarly, it follows that
\[
\chi_q(a)=(-1)^t,\]
where
\[
t=\sum_{0\leq m<a,\ m\ \textrm{odd}} |Q_m(a)|,\]
\[
Q_m(a)=\left\{r\in \mathbb{Z}: \frac{m}{2a}q<r<\frac{m+1}{2a}q\right\}.\]

We have that $p-q=4an$, for some positive integer $n$, and so when $m$ is odd,
\begin{eqnarray*}
P_m(a)&=&\left\{r\in \mathbb{Z}: \frac{m}{2a}p<r<\frac{m+1}{2a}p\right\}\\
&=&\left\{r\in \mathbb{Z}: \frac{m}{2a}(q+4an)<r<\frac{m+1}{2a}(q+4an)\right\}\\
&=&\left\{r\in \mathbb{Z}: \frac{m}{2a}q+2mn<r<\frac{m+1}{2a}q+2mn+2n\right\}\\
&=&\left\{r^{\prime}\in \mathbb{Z}: \frac{m}{2a}<r^{\prime}<\frac{m+1}{2a}q+2n\right\},
\end{eqnarray*}
where $r^{\prime}=r-2mn$. Hence
\[
|P_m(a)|=|Q_m(a)|+2n.\]
We conclude that $s\equiv t$ mod 2, whence
\[
\chi_p(a)=(-1)^s=(-1)^t=\chi_q(a).\]
$\hspace{15.5cm}\ \textrm{QED}$ 
\section{A Proof of Quadratic Reciprocity via Gauss Sums: Introduction}

The third proof of Theorem 3.3 that we will give is a version of Gauss' sixth proof.  It uses ingenious calculations based on some basic facts from algebraic number theory, and anticipates some important techniques that we will use later to study various properties of residues and non-residues in greater depth.

Gauss' sixth proof of quadratic reciprocity [23] appeared in 1818. He mentions in the introduction to this paper that for years he had searched for a method that would generalize to the cubic and bi-quadratic cases and that finally his untiring efforts were crowned with success. The purpose of publishing this sixth proof, he states, was to bring to a close this part of the higher arithmetic dealing with quadratic residues and to say, in a sense, farewell.  Our third proof of Theorem 3.3 is a reworking of Gauss' argument from [23] using some basic facts from the theory of algebraic numbers. We start first with a rather detailed discussion of the algebraic number theory that will be required; this is the content of the next section. This information is then used in section 9 to prove Theorem 3.3, following the development given in Ireland and Rosen [30], sections 6.2 and 6.3.

\section{Algebraic Number Theory}
In the introduction to his great memoirs [24] and [25] on biquadratic reciprocity, Gauss asserts that the theory of quadratic residues had reached such a state of perfection that no more improvement was necessary. However, he said ``The theory of biquadratic residues is by far more difficult". After struggling with this problem for a long time, he had been able to derive satisfactory results in only a few special cases, with the proofs being so difficult that he realized ``$\dots$ the previously accepted principles of arithmetic are in no way sufficient for the foundations of a general theory, that rather such a theory necessarily demands that to a certain extent the domain of the higher arithmetic needs to be endlessly enlarged$\dots$". In modern language, Gauss is calling for a theory of algebraic numbers. He began to answer that call himself in the second paper [25], where he used the subring of the complex numbers $\mathbb{Z}+i\mathbb{Z}$, what is now called the Gaussian integers, to formulate a precise statement of the Law of Biquadratic Reciprocity. Subsequently, Gauss' call has been so resoundingly answered by the work of Dirichlet, Dedekind, Kummer, Kronecker, Hilbert, and many others, that today the theory of algebraic numbers is indispensable in virtually all areas of number theory. In this section we establish some basic facts from the theory of algebraic numbers which will be required in our third proof of Theorem 3.3. This information will also play an important role in further developments in subsequent chapters of these notes.

Let $\textbf{C}$ denote the complex numbers.
\vspace{0.3cm}

\emph{Definition}. A \emph{field of complex numbers} is a nonzero subfield of $\textbf{C}$.

\vspace{0.3cm}
N. B. Every field of complex numbers contains the field $\mathbb{Q}$ of rational numbers.

\vspace{0.3cm}
\emph{Notation}: if $A$ is a commutative ring then $A[x]$ will denote the ring of all polynomials in $x$ with coefficients in $A$.
\vspace{0.3cm}

\emph{Definitions}. Let $F$ be a field of complex numbers. A complex number $\theta$ is \emph{algebraic over $F$} if there exists $f \in F[x]$ such that $f \not \equiv 0$ and $f(\theta)=0$. If $\theta$ is algebraic over $F$, let 
\[
M(\theta)=\{f \in F[x]: f\ \textrm{ is monic and}\ f(\theta)=0\}
\]
(N.B. $M(\theta) \not = \emptyset)$. An element of $M(\theta)$ of smallest degree is a \emph{minimal polynomial of $\theta$ over $F$}.
\begin{prp}
The minimal polynomial of a complex number algebraic over a field of complex numbers $F$ is unique and irreducible over $F$.
\end{prp}

\emph{Proof}. Let $r$ and $s$ be minimal polynomials of the number $\theta$ algebraic over $F$. Use the division algorithm in $F[x]$ to find $d, f \in F[x]$ such that
\[
r=ds+f,\ f\equiv 0\ \textrm{or degree of $f<$ degree of s}.\]
Hence
\[
f(\theta)=r(\theta)-d(\theta)s(\theta)=0.\]
If $f \not \equiv 0$ then, upon dividing $f$ by its leading coefficient, we get a monic polynomial over $F$ of lower degree than $s$ and not identically 0 which has $\theta$ as a root, which is not possible because $s$ is a minimal polynomial of $\theta$ over $F$. Hence $f \equiv 0$ and so $s$ divides $r$ over $F$, Similarly, $r$ divides $s$ over $F$. Hence $r=\alpha s$ for some $\alpha \in F$, and as $r$ and $s$ are both monic, $\alpha=1$, and so $r=s$.  This proves that the minimal polynomial is unique.

To show that the minimal polynomial $m$ is irreducible over $F$, suppose that $m=rs$, where $r$ and $s$ are non-constant elements of $F[x]$. Then the degrees of $r$ and $s$ are both less than the degree of $m$, and $\theta$ is a root of either $r$ or $s$. Hence a constant multiple of either $r$ or $s$ is a monic polynomial in $F[x]$ having $\theta$ as a root and is of degree less than the degree of $m$, contradicting the minimality of the degree of $m$. $\hspace{6cm} \ \textrm{QED}$

\vspace{0.3cm}

\emph{Definition}. Let $\theta$ be algebraic over $F$. The \emph{degree of $\theta$ over $F$} is the degree of the minimal polynomial of $\theta$ over $F$.
\vspace{0.3cm}
\begin{lem}
If $\theta \in \textbf{C}$, $F$ is a field of complex numbers, and $f \in F[x]$ is monic, irreducible over $F$, and $f(\theta)=0$ then f is the minimal polynomial of $\theta$ over F.
\end{lem}

\emph{Proof}. Let $m$ be the minimal polynomial of $\theta$ over $F$. The division algorithm in $F[x]$ implies that there exits $q$, $r \in F[x]$ such that
\[
f=qm+r,\ r \equiv 0 \ \textrm{or degree of $r$ $< $ degree of $m$}.
\]
But\[r(\theta)=f(\theta)-q(\theta)m(\theta)=0,
\]
and so if $r\not \equiv 0$ then we divide $r$ by its leading coefficient to get a monic polynomial over $F$ that is not identically 0, has $\theta$ as a root, and is of degree less than the degree of $m$, which is impossible by the minimality of the degree of $m$. Hence $r \equiv 0$ and so $f=qm$. But $f$ is irreducible over $F$, and so either $q$ or $m$ is constant. If $m$ is constant then $m \equiv 1$ ($m$ is monic), not possible because $m(\theta)=0$. Hence $q$ is constant, and because $f, m$ are both monic, $q \equiv 1$. Hence $f=m. 
\hspace{13cm} \textrm{QED}$      
\vspace{0cm}

We will now discuss two examples that will be of major importance in subsequent work.

(1) Let $m \in \mathbb{Z}\setminus \{1\}$ be square-free, i.e., $m$ does not have a square $\not= 1$ as a factor. Then $\sqrt m$ is irrational, hence $x^2-m$ is irreducible over $\mathbb{Q}$. Lemma 3.8 imllies that  $x^2-m$ is the minimal polynomial of $\sqrt m$ over $\mathbb{Q}$ and so $\sqrt m$ is algebraic over $\mathbb{Q}$ of degree 2.

(2) Let $q$ be a prime and let
\[
\zeta_q= \exp \Big( \frac{2 \pi i}{q}\Big).
\]
Then $\zeta_q^q=1, \zeta_q \not =1$, hence we deduce from the factorization
\[
x^q-1=(x-1) \Big(\sum_{k=0}^{q-1}\ x^k \Big)
\]
that $\zeta_q$ is a root of $\sum_{k=0}^{q-1}\ x^k$.

We claim that $\sum_{k=1}^{q-1}\ x^k$ is irreducible over $\mathbb{Q}$. To see this, note first that a polynomial $f(x)$ is irreducible if and only if $f(x+1)$ is irreducible, because $f(x+1)=g(x)h(x)$ if and only if $f(x)=g(x-1)h(x-1)$. Hence
\[
\sum_{k=0}^{q-1}\ x^k=\frac{x^q-1}{x-1}\ \textrm{is irreducible if and only if $\frac{(x+1)^q-1}{x}$ is irreducible}.
\]
It follows from the binomial theorem that
\[
\frac{(x+1)^q-1}{x}= \sum_{k=1}^q\ \left( \begin{array}{c} q\\ k \end{array} \right) x^{k-1}.
\]
We now recall the following fact about binomial coefficients: if $q$ is a prime then $q$ divides the binomial coefficient $\left( \begin{array}{c} q\\ k \end{array} \right), k=1,\dots,q-1$. Hence
\[
\frac{(x+1)^q-1}{x}=x^{q-1}+q(x^{q-2}+\dots )+q,
\]
and this polynomial is irreducible over $\mathbb{Q}$ by way of
\begin{lem}
$($Eisenstein's criterion$)$ If q is prime and $f(x)= \sum_{k=0}^n\ a_kx^k$ is a polynomial in $\mathbb{Z}[x]$ whose coefficients satisfy: $q$ does not divide $a_n, q^2$ does not divide $a_0$, and $q$ divides $a_k, k=0,1,\dots,n-1,$ then $f(x)$ is irreducible over $\mathbb{Q}$.
\end{lem}
\noindent Thus $\zeta_q$ has minimal polynomial $\sum_{k=0}^{q-1}\ x^k$ and hence is algebraic over $\mathbb{Q}$ of degree $q-1$.

\emph{Proof of Lemma} 3.9. We assert first that if a polynomial $h \in \mathbb{Z}[x]$ does not factor into a product of polynomials in $\mathbb{Z}[x]$ of degree lower than the degree of $h$ then it is irreducible over $\mathbb{Q}$. In order to see this, suppose that $h$ is not constant (otherwise the assertion is trivial) and that $h=rs$, where $r$ and $s$ are polynomials in $\mathbb{Q}[x]$, both not constant and of lower degree than $h$. By clearing denominators and factoring out the greatest common divisors of appropriate integer coefficients, we find integers $a, b, c,$ and polynomials $g, u, v$ in $\mathbb{Z}[x]$ such that $h=ag$, degree of $r$ = degree of $u$, degree of $s$ = degree of $v$, 
\[
abg=cuv,\]
and all of the coefficients of $g$ (respectively, $u$, $v$) are relatively prime, i.e., the greatest common divisor of all of the coefficients of $g$ (respectively, $u$, $v$) is 1.

We claim that the coefficients of the product $uv$ are also relatively prime. Assume this for now. Then $|ab|=$ the greatest common divisor of the coefficients of $abg=$ the greatest common divisor of the coefficients of $cuv=|c|$, hence $ab=\pm c$. But then $h=\pm auv$ and this is a factorization of $h$ as a product of polynomials in $\mathbb{Z}[x]$ of lower degree..

We must now verify our claim.  Suppose that the coefficients of $uv$ have a common prime factor $r$. Let $\mathbb{F}_r$ denote the field of ordinary residue classes mod $r$. If $s \in \mathbb{Z}[x]$ and if we let $\bar{s}$ denote the polynomial in $\mathbb{F}_r[x]$ obtained from $s$ by reducing the coefficients of $s$ mod $r$, then $s \rightarrow \bar{s}$ defines a homomorphism of $Z[x]$ onto $\mathbb{F}_r[x]$. Because $r$ divides all of the coefficients of $uv$, it hence follows that
\[
0=\overline{uv}=\bar{u}\bar{v}\ \textrm{in}\  \mathbb{F}_r[x].\]
Because $\mathbb{F}_r$ is a field, $\mathbb{F}_r[x]$ is an integral domain, hence we conclude from this equation that either $\bar{u}$ or $\bar{v}$ is 0 in $\mathbb{F}_r[x]$, i.e., either all of the coefficients of $u$ or of $v$ are divisible by $r$. This contradicts the fact that the coefficients of $u$ (respectively, $v$) are relatively prime. The assertion that the product of two polynomials in $\mathbb{Z}[x]$ has all of its coefficients relatively prime whenever the coefficients of each polynomial are relatively prime is often referred to as Gauss' lemma, not to be confused, of course, with the statement in Theorem 2.7.

Next suppose that $f(x)=\sum_{k=0}^n a_kx^k \in \mathbb{Z}[x]$ satisfies the hypotheses of Lemma 3.9. By virtue of what we just showed, we need only prove that $f$ does not factor into polynomials of lower degree in $\mathbb{Z}[x]$. Suppose, on the contrary, that
\[
f(x)=\Big(\sum_{k=0}^s b_kx^k\Big)\Big(\sum_{k=0}^t c_kx^k\Big) \]
is a factorization of $f$ in $\mathbb{Z}[x]$ with $b_s \not= 0 \not= c_t$ and $s$ and $t$ both less than $n$. Because $a_0 \equiv 0$ mod $q$, $a_0 \not \equiv 0$ mod $q^2$ and $a_0=b_0c_0$, one element of the set $\{b_0, c_0\}$ is $\not \equiv 0$ mod $q$ and the other is $ \equiv 0$ mod $q$. Assume that $b_0$ is the former element and $c_0$ is the latter. As $a_n \not \equiv 0$ mod $q$ and $a_n=b_sc_t$, it follows that $b_s \not \equiv 0 \not \equiv c_t$ mod $q$. Let $m$ be the smallest value of $k$ such that $c_k \not \equiv 0$ mod $q$. Then $m>0$, hence
\[
a_m=\sum_{j=0}^{m-i} b_j c_{m-j} \]
for some $i \in[0, m-1]$. Because $b_0 \not \equiv 0 \not \equiv c_m$ mod $q$ and $c_{m-1},\dots,c_i$ are all $\equiv 0$ mod $q$, it follows that $a_m \not \equiv 0$ mod $q$, and so $m=n$. Hence $t=n$, contradicting the assumption on $t$ and $n$. $\hspace{15.1cm} \textrm{QED}$

The crucial fact about algebraic numbers that we will need in order to prove the LQR is that the set of all algebraic integers (see the definition after the proof of Theorem 3.10)  forms a subring of the field of complex numbers. The verification of that fact is the goal of the next two results.

For use in the proof of the next theorem, we recall that if $n$ is a positive integer, then the \emph{elementary symmetric polynomials in n variables} are the polynomials in the variables $x_1,\dots,x_n$ defined by
\[
 \sigma_1=\sum_{i=1}^n x_i,\]

\hspace{7cm}  \vdots 

\hspace{6.4cm} $\sigma_i=\textrm{sum of all products of $i$ different $x_j$},$

\hspace{7cm}  \vdots 
\[
\sigma_n=\prod_{i=1}^n x_i.\]
The elementary symmetric polynomials have the property that if $\pi$ is a permutation of the set $[1, n]$ then $\sigma_i(x_{\pi(1)},\dots,x_{\pi(n)})=\sigma_i(x_1,\dots,x_n)$, i.e., $\sigma_i$ is unchanged by any permutation of its variables.

\begin{thm}
If F is a field of complex numbers then the set of all complex numbers algebraic over F is a field of complex numbers which contains F.
\end{thm}

\emph{Proof}. Let $\alpha$ and $\beta$ be algebraic over $F$. We want to show that $\alpha \pm \beta, \alpha \beta,$ and $\alpha/\beta,$ provided that $\beta \not=0$, are all algebraic over $F$. 
We will do this by the explicit construction of  polynomials over $F$ that have these numbers as roots.

Start with $\alpha+\beta$. Let $f$ and $g$ denote the minimal polynomials of, respectively, $\alpha$ and $\beta$, of degree $m$ and $n$, respectively. Let $\alpha_1,\dots,\alpha_m$ and $\beta_1,\dots,\beta_n$ denote the roots of $f$ and $g$ in $\textbf{C}$, with $\alpha_1=\alpha$ and $\beta_1=\beta$. Now consider the polynomial
\begin{equation*}
\prod_{i=1}^m \prod_{j=1}^n (x-\alpha_i-\beta_j)=x^{mn}+\sum_{i=1}^{mn} c_i(\alpha_1,\dots,\alpha_m, \beta_1,\dots,\beta_n)x^{mn-i},\tag{8}
\end{equation*}
where each coefficient $c_i$ is a polynomial in the $\alpha_i$'s and $\beta_j$'s over $F$ (in fact, over $\mathbb{Z}$). We claim that 
\begin {equation*}
c_i(\alpha_1,\dots,\alpha_m, \beta_1,\dots,\beta_n) \in F,\ i=1,\dots,mn. \tag{9} \] 
If this is true then the polynomial $(8)$ is in $F[x]$ and has $\alpha_1+\beta_1=\alpha+\beta$ as a root, whence $\alpha+\beta$ is algebraic over $F$.

In order to verify $(9)$, we will make use of the following result from the classical theory of equations (see Weisner [59], Theorem 49.10). Let $\tau_1,\dots,\tau_m, \sigma_1,\dots,\sigma_n$ denote, respectively, the elementary symmetric polynomials in $m$ and $n$ variables. Suppose that the polynomial $h$ over $F$ in the variables $x_1,\dots,x_m, y_1,\dots,y_n$ has the property that if $\pi$ (respectively, $\nu$) is a permutation of $[1, m]$ (respectively, $[1, n]$) then
\[
h(x_1,\dots,x_m, y_1,\dots,y_n)=h(x_{\pi(1)},\dots,x_{\pi(m)}, y_{\nu(1)},\dots,y_{\nu(n)}),\]
i.e., $h$ remains unchanged when its variables $x_i$ and $y_j$ are permuted amongst themselves. Then there exist a polynomial $l$ over $F$ in the variables $x_1,\dots,x_m, y_1,\dots,y_n$ such that
\begin{eqnarray*}
h(x_1,\dots,&x_m,& y_1,\dots,y_n)\\
&=&l(\tau_1(x_1,\dots,x_m),\dots,\tau_m(x_1,\dots,x_m),\sigma_1( y_1,\dots,y_n),\dots\sigma_n( y_1,\dots,y_n)).
\end{eqnarray*}
Observe next that the left-hand side of $(8)$ remains unchanged when the $\alpha_i$'s and the $\beta_j$'s are permuted amongst themselves (this simply rearranges the order of the factors in the product), and so the same thing is true for each coefficient $c_i$. It thus follows from our result from the theory of equations that there exists a polynomial $l_i$ over $F$ in the variables $x_1,\dots,x_m, y_1,\dots,y_n$ such that
\begin{eqnarray*}
c_i(\alpha_1,\dots,&\alpha_m,& \beta_1,\dots,\beta_n)\\
&=&l_i(\tau_1(\alpha_1,\dots,\alpha_m),\dots,\tau_m(\alpha_1,\dots,\alpha_m),\sigma_1( \beta_1,\dots,\beta_n),\dots,\sigma_n( \beta_1,\dots,\beta_n)).
\end{eqnarray*}
If we can prove that each of the numbers at which $l_i$ is evaluated in this equation is in $F$ then $(9)$ will be verified. Hence it suffices to prove that if $\theta$ is a number algebraic over $F$ of degree $n$, $\theta_1,\dots,\theta_n$ are the roots of the minimal polynomial $m$ of $\theta$ over $F$, and $\sigma$ is an elementary symmetric polynomial in $n$ variables, then $\sigma(\theta_1,\dots,\theta_n) \in F$. But this last statement follows from the fact that all the coefficients of $m$ are in $F$ and
\[
m(x)=\prod_{i=1}^n (x-\theta_i)=x^n+\sum_{i=1}^n (-1)^i \sigma_i(\theta_1,\dots,\theta_n)x^{n-i},\]
where $\sigma_1,\dots,\sigma_n$ are the elementary symmetric polynomials in $n$ variables.

A similar argument shows that $\alpha-\beta$ and $\alpha\beta$ are algebraic over $F$.

Suppose next that $\beta\not=0$ is algebraic over $F$ and let
\[
x^n+\sum_{i=0}^{n-1} a_ix^i\]
be the minimal polynomial of $\beta$ over $F$. Then $1/\beta$ is a root of
\[
1+\sum_{i=0}^{n-1} a_ix^{n-i} \in F[x], \]
and so $1/\beta$ is algebraic over $F$. Then $\alpha/\beta=\alpha \cdot (1/\beta)$ is algebraic over $F$. $\hspace{2.1cm}\ \textrm{QED}$

\vspace{0.3cm}

 \emph{Notation}: $\mathcal{A}(F)$ denotes the field of all complex numbers algebraic over $F$.
\vspace{0.3cm}

\emph{Definition}. An element of $\mathcal{A}(\mathbb{Q})$ is an \emph{algebraic integer} if its minimal polynomial over $\mathbb{Q}$ has all of its coefficients in $\mathbb{Z}$.

\vspace{0.3cm}
 By virtue of examples (1) and (2), all square-free integers and $\exp(2\pi i/q), q$ a prime, are algebraic integers. 

\begin{thm}
The set of all algebraic integers is a subring of $\mathcal{A}(\mathbb{Q})$ containing $\mathbb{Z}$.
\end{thm}

\emph{Proof}. Let $\alpha$ and $\beta$ be algebraic integers. We need to prove that $\alpha \pm \beta$ and $\alpha \beta$ are algebraic integers. This can be done by first observing that the result from the theory of equations that we used in the proof of Theorem 3.10 holds \emph{mutatis mutandis} if the field $F$ there is replaced by the ring $\mathbb{Z}$ of integers (Weisner [59], Theorem 49.9). If we then let
$\alpha_1,\dots,\alpha_m$ and $\beta_1,\dots,\beta_n$ denote the roots of the minimal polynomial over $\mathbb{Q}$ of $\alpha$ and $\beta$, respectively, then the proof of Theorem 3.10, with $F$ replaced in that proof by $\mathbb{Z}$, verifies that $\alpha \pm \beta$ and $\alpha \beta$ are roots of monic polynomials in $\mathbb{Z}[x]$. We now invoke the following fact: if a complex number $\theta$ is the root of a monic polynomial in $\mathbb{Z}[x]$ then it is an algebraic integer.

In order to prove the last statement about $\theta$, let $f \in \mathbb{Z}[x]$ be  monic with $f(\theta)=0$. If $m$ is the minimal polynomial of $\theta$ over $\mathbb{Q}$ then we must show that $m \in \mathbb{Z}[x]$. It follows from the proof of Proposition 3.7 that there is a $q \in \mathbb{Q}[x]$ such that $f=qm$ and so we can find a rational number $a/b$ and $u, v \in \mathbb{Z}[x]$ such that $f= (a/b)uv, u$ (respectively, $v$) is a constant multiple of $m$ (respectively, $q$), and $u$ (respectively, $v$) has all of its coefficients relatively prime.

 We have that 
\[
bf=auv.\]
As $f$ is monic and $u, v \in \mathbb{Z}[x]$, we conclude that $a$ divides $b$ in $\mathbb{Z}$, say $b=ak$ for some $k \in \mathbb{Z}$. Hence 
\[
kf=uv.\]
Because $f \in \mathbb{Z}[x]$, it follows that $k$ is a common factor of all of the coefficients of $uv$. Because of the claim that we verified in the proof of Lemma 3.9, the coefficients of $uv$ are relatively prime, hence $k=\pm 1$, and so
\[
f=\pm uv.\]
As $f$ is monic, the leading coefficient of $u$ is $\pm 1$. But $u$ is a constant multiple of $m$ and $m$ is monic, hence $m=\pm u \in \mathbb{Z}[x]$, and so $\theta$ is an algebraic integer. $\hspace{4.2cm} \textrm{QED}$

\vspace{0.3cm}

\emph{Notation}: $\mathcal{R}$ will denote the ring of algebraic integers.

\vspace{0.3cm}

The proof of Theorem 3.3 that we will give in the next section requires the following simple lemma: 
\begin{lem}
$\mathcal{R} \cap \mathbb{Q}=\mathbb{Z}$.
\end{lem}

\emph{Proof}. If $q \in \mathcal{R} \cap \mathbb{Q}$ then $x-q$ is the minimal polynomial of $q$ over $\mathbb{Q}$, hence $x-q \in \mathbb{Z}[x]$, hence $q \in \mathbb{Z}.\hspace{13.4cm} \textrm{QED}$

\section{Proof of Quadratic Reciprocity via Gauss Sums: Conclusion}

As a warm-up for our third proof of Theorem 3.3, and also as a first illustration of how useful algebraic number theory is to what we will do subsequently, we will reprove Theorem 2.6, which assets that 
\[
\chi_p(2)=(-1)^{\varepsilon},\ \textrm{where $\varepsilon \equiv \frac{p^2-1}{8} \mod 2$},
\]
by using algebraic number theory.
 Let
 \[
 \zeta= e^{\pi/4}= \frac{1}{\sqrt 2}+\frac{1}{\sqrt 2} i.
 \]
 Then 
 \[
 \zeta^{-1}= e^{-\pi/4}= \frac{1}{\sqrt 2}-\frac{1}{\sqrt 2} i,
 \]
hence
\[
\tau=\zeta+\zeta^{-1}= \sqrt 2 \in \mathcal{R},
\]
and so we can work in the ring $\mathcal{R}$ of algebraic integers.

If $p$ is an odd prime then we let 
\[
(p)=\
\textrm{the ideal in $\mathcal{R}$ generated by $p=p\mathcal{R}=\{p\alpha: \alpha \in \mathcal{R}\}$}.
\]
If $\alpha, \beta \in \mathcal{R}$, then we will write
\[
\alpha \equiv \beta \ \textrm{mod}\ p 
\]
if $\alpha-\beta \in (p)$. Euler's criterion (Theorem 2.5) implies that
\[
\tau^{p-1}=(\tau^2)^{(p-1)/2}=2^{(p-1)/2} \equiv \chi_p(2) \ \textrm{mod}\ p,
\]
hence 
\begin{equation*}
\tau^p \equiv \chi_p(2) \tau \ \textrm{mod}\ p. \tag{10}
\]
We now make use of the following lemma, which follows from the binomial theorem and the fact that $p$ divides the binomial coefficient $\left( \begin{array}{c} p\\ k \end{array} \right), k=1,\dots,p-1$.
\begin{lem}
If $\alpha, \beta \in \mathcal{R}$ then
\[
(\alpha+\beta)^p \equiv \alpha^p+\beta^p \ \textnormal{mod}\ p.
\]
\end{lem}
Hence
\begin{equation*}
\tau^p=(\zeta+\zeta^{-1})^p \equiv \zeta^p+\zeta^{-p} \ \textrm{mod}\ p. \tag{11}
\]
 The next step is to calculate $\zeta^p+\zeta^{-p}$. Begin by noting that
 \[
 \zeta^8=1,
 \]
hence if $p \equiv \pm 1 \ \textrm{mod}\ 8$, then
\[
\zeta^p+\zeta^{-p}=\zeta+\zeta^{-1}=\tau,
\]
and if $p \equiv \pm3 \ \textrm{mod}\ 8$, then
\begin{eqnarray*}
\zeta^p+\zeta^{-p}&=&\zeta^3+\zeta^{-3}\\
&=&-(\zeta^{-1}+\zeta)\\
&=&-\tau,
\end{eqnarray*}
where the second line follows from the first because $\zeta^4=-1 implies that \zeta^3=-\zeta^{-1}, and so \zeta^{-3}=-\zeta.$ Hence
\begin{equation*}
\zeta^p+\zeta^{-p}=(-1)^{\varepsilon}\tau,\  \varepsilon \equiv \frac{p^2-1}{8} \ \textrm{mod}\ 2. \tag{12}
\]
 As a consequence of the congruences (10), (11), and (12), it follows that
 \[
 \chi_p(2) \tau \equiv \zeta^p+\zeta^{-p} \equiv (-1)^{\varepsilon}\tau \ \textrm{mod}\ p.
 \]
Multiply this congruence by $\tau$ and use $\tau^2=2$ to derive
\begin{equation*}
2\chi_p(2) \equiv 2(-1)^{\varepsilon} \ \textrm{mod}\ p. \tag{13}
\]
Now this congruence is in $\mathcal{R}$, so there exits $\alpha \in \mathcal{R}$ such that
\[
2\chi_p(2)=2(-1)^{\varepsilon}+\alpha p,
\]
hence
\[
\alpha=\frac{2(\chi_p(2)-(-1)^{\varepsilon})}{p} \in \mathcal{R} \cap \mathbb{Q}=\mathbb{Z}\ (\textrm{by Lemma 3.12}).
\]
Hence (13) is in fact a congruence in $\mathbb{Z}$, and so
\[
\chi_p(2) \equiv (-1)^{\varepsilon} \ \textrm{mod}\ p\ \textrm{in}\ \mathbb{Z},
\]
whence, as before,
\[
\chi_p(2)=(-1)^{\varepsilon}.
\]

This proof of Theorem 2.6 depends on the equation $\tau^2=2$. Can one get a similar equation with an odd prime $p$ replacing the 2 on the right-hand side of this equation? Yes one can, and a proof of the LQR will follow in a similar way from the ring structure of $\mathcal{R}$.

In order to see how that goes, let $\zeta=e^{2 \pi i/p}$ and set
\begin{eqnarray*}
g&=& \sum_{n=0}^{p-1}\ \chi_p(n)\zeta^n,\\
p^{*}&=&(-1)^{(p-1)/2}p.
\end{eqnarray*}
The sum $g$ is called a \emph{Gauss sum}; these sums were first used by Gauss in his famous study of cyclotomy which concluded \emph{Disquisitiones Arithmeticae} ([19], section VII). The analogue of the equation $\tau^2=2$ is given by
\begin{thm}
$g^2=p^{*}.$
\end{thm}

Assume this for now; we deduce LQR from it like so: let $q$ be an odd prime, $q \not= p$. Then
\[
g^{q-1}=(g^2)^{(q-1)/2}=(p^{*})^{(q-1)/2} \equiv \chi_q(p^{*}) \ \textrm{mod}\ q,
\]
where the last equivalence follows from Euler's criterion. Hence
\begin{equation*}
g^q \equiv \chi_q(p^{*})g \ \textrm{mod}\ q, \tag{14}
\]
where this congruence is now in $\mathcal{R}$, because $g \in \mathcal{R}$. If $n \in \mathbb{Z}$ then $\chi_p(n)^q=\chi_p(n)$ because $\chi_p(n) \in [-1, 1]$ and $q$ is odd; consequently Lemma 3.13 implies that
\begin{eqnarray*}
(15)\   g^q&=&\Big(\sum_n\ \chi_p(n) \zeta^n \Big)^q\\
&\equiv& \sum_n\ \chi_p(n)^q \zeta^{qn} \ \textrm{mod}\ q\\
&\equiv& \sum_n\ \chi_p(n) \zeta^{qn} \ \textrm{mod}\ q,
\end{eqnarray*}

We now need 
\begin{lem}
If $a \in \mathbb{Z}$ then
\[
\sum_n\ \chi_p(n) \zeta^{an}= \chi_p(a)g.
\]
\end{lem}

The sum on the left-hand side of this  equation is another Gauss sum. Lemma 3.15 records a very important relation satisfied by Gauss sums; in addition to the use that we make of it here, it will also play an important role in some calculations that are performed in Chapter 7, where we study certain distributions of residues and non-residues. 

Assume Lemma 3.15 for now; this lemma and (15) imply that
\begin{equation*}
g^q \equiv \chi_p(q)g \ \textrm{mod}\ q. \tag{16}
\]
A consequence of (14) and (16) is that
\[
\chi_q(p^{*})g \equiv \chi_p(q)g \ \textrm{mod}\ q.
\]
Multiply by $g$ and use $g^2=p^{*}$ to derive
\[
\chi_q(p^{*})p^{*} \equiv \chi_p(q)p^{*} \ \textrm{mod}\ q,
\]
and then apply Lemma 3.12 and the fact that $\chi_q(p^{*}), \chi_p(q)$ are both $\pm1$ as before to get
\begin{equation*}
\chi_q(p^{*})=\chi_p(q). \tag{17}\]
Theorem 2.4 implies that
\[
\chi_q(-1)=(-1)^{(q-1)/2},
\]
hence, because of (17),
\begin{eqnarray*}
\chi_p(q)&=&\chi_q(-1)^{\frac{1}{2}(p-1)}\chi_q(p)\\
&=&(-1)^{\frac{1}{2}(q-1) \frac{1}{2}(p-1)}\chi_q(p),
\end{eqnarray*}
which is the LQR. 
 
 We must now prove Theorem 3.14 and Lemma 3.15. Since Lemma 3.15 is used in the proof of Theorem 3.14, we verify Lemma 3.15 first.
 
 \emph{Proof of Lemma} 3.15. Suppose first that $p$ divides $a$. Then $\zeta^{an}=1$, for all $n$ and $\chi_p(0)=0$ so
 \[
 \sum_n\ \chi_p(n)\zeta^{an}=\sum_{n=1}^{p-1}\ \chi_p(n).
 \]
Half of the terms of the sum on the right-hand side are 1 and the other half are $-1$ (Proposition 2.1), and so this sum is 0. Because $\chi_p(a)=0$ ($p$ divides $a$), the conclusion of Lemma 3.15 is valid.

Suppose that $p$ does not divide $a$. Then
\begin{equation*}
\chi_p(a) \sum_n\ \chi_p(n) \zeta^{an}= \sum_n\ \chi_p(an) \zeta^{an}. \tag{18}
\]
Observe next that whenever $n$ runs through a complete system of ordinary residues mod $p$, so does $an$, and also that $\chi_p(an)$ and $\zeta^{an}$ depend only on the residue class mod $p$ of $an$. Hence the sum on the right-hand side of (18) is
\[
\sum_{n=0}^{p-1}\ \chi_p(n) \zeta^n=g.
\]
Hence
\[
\chi_p(a) \sum_n\ \chi_p(n) \zeta^{an}=g.
\]
Now multiply through by $\chi_p(a)$ and use the fact that $\chi_p(a)^2=1$, since $p$ does not divide $a$. 

$\hspace{15cm} \textrm{QED}$

\emph{Proof of Theorem }3.14. We must prove that $g^2=p^{*}$. 

Suppose that $\gcd(a, p)=1$ and let
\[
g(a)= \sum_{n=0}^{p-1}\ \chi_p(n) \zeta^{an}.
\]
The idea of this argument is to calculate
\[
\sum_{a=0}^{p-1}\ g(a)g(-a)
\]
in two different ways, equate the expressions resulting from that, and see what happens.

For the first way, use Lemma 3.15 to obtain
\begin{eqnarray*}
g(a)g(-a)&=&\chi_p(a)\chi_p(-a)g^2\\
&=&\chi_p(-a^2)g^2\\
&=&\chi_p(-1)g^2,\ a=1,\dots,p-1,
\end{eqnarray*}
Hence this and the fact that $g(0)=\sum_0^{p-1}\chi_p(n)=0$ imply that
\begin{equation*}
\sum_{a=0}^{p-1}\ g(a)g(-a)=(p-1)\chi_p(-1)g^2. \tag{19}
\]

Now for the second way. We have that
\[
g(a)g(-a)= \sum_{1 \leq x, y \leq p-1}\ \chi_p(x) \chi_p(y) \zeta ^{a(x-y)}.
\]
Hence 
\begin{equation*}
\sum_{a=0}^{p-1}\ g(a)g(-a)=\sum_{1 \leq x, y \leq p-1}\ \chi_p(x) \chi_p(y) \sum_a\ \zeta^{a(x-y)}. \tag{20}
\]
The next step is to calculate 
\[
\sum_a\ \zeta^{a(x-y)}
\]
for fixed $x$ and $y$. If $x \not= y$ then
\[
1 \leq |x-y| \leq p-1
\]
and so $p$ does not divide $x-y$, hence $\zeta^{x-y} \not= 1$, hence 
\[
\sum_a\ \zeta^{a(x-y)}=\frac{\zeta^{(x-y)p}-1}{\zeta^{x-y}-1}=0,\ (\zeta^p=1\ !).
\]
Hence
\begin{equation*}
\sum_a\ \zeta^{a(x-y)}=\left\{\begin{array}{ll}p,\ \textrm{if $x=y$,}\\
0,\ \textrm{if $x \not=y$.}\\\end{array}\right. \tag{21}
\]
Because of (20) and (21), it follows that
\begin{equation*}
\sum_{a=0}^{p-1}\ g(a)g(-a)= (p-1)p. \tag{22}
\]
Equations (19) and (22) imply that
\[
(p-1) \chi_p(-1)g^2=(p-1)p,\]
hence from Theorem 2.4,
\[
g^2=\chi_p(-1)p=(-1)^{(p-1)/2}p.
\]
$\hspace{15.6cm} \textrm{QED}$

As Ireland and Rosen point out, the main idea of Gauss' sixth proof of quadratic reciprocity is to consider the polynomial
\[
f_k(x)=\sum _{t=0}^{p-1} \chi_p(t)x^{kt},\]
and then to show, as Gauss did without using any roots of unity, that\[
1+x+\cdots+x^{p-1}\]
divides
\[
f_1(x)^2-(-1)^{\frac{1}{2}(p-1)}p\ \textrm{and}\ f_q(x)-\chi_p(q)f_1(x).\]
Quadratic reciprocity then follows by noting that 
\[
f_q(x)\equiv f_1(x)^q\ \textrm{mod}\ q.\]
The proof which we have presented in this section amounts to setting $x=\exp(2\pi i/p)$ in this argument and then performing the necessary calculations by means of congruences in the ring of algebraic integers. This observation was also made by Eisenstein and Jacobi, who used it as a stepping stone to the proof of  higher-degree reciprocity laws via Gauss sums.

\section{A Proof of Quadratic Reciprocity via Ideal Theory: Introduction}

The second proof which Gauss gave for quadratic reciprocity rests on his genus theory for quadratic forms, developed (along with his theory of composition of forms) in articles 231-261 of the \emph{Disquisitiones}. A primary goal of genus theory is the determination of the number of genera of quadratic forms with a given discriminant, and Gauss deduced the LQR from a formula which he established for the number of genera. The argument which Gauss used here to deduce the LQR is quite elegant, but the proof of the formula for the number of genera on which it is based is much more involved, as it uses virtually all of the formidable mathematical technology involved in Gauss' development of his composition of forms. Fortunately, the ideas and techniques which Gauss employed can be given a formulation that is much easier to follow using the modern theory of ideals in the ring of algebraic integers in a quadratic number field. This is the approach that we will take in our fourth proof of Theorem 3.3; the necessary information about this ideal theory will be given in the next section, and we will then deduce Theorem 3.3 from that ideal theory in section 12.  
\section{The Structure of Ideals in a Quadratic Number Field}
Let $F$ be a field of complex numbers. With respect to its addition and multiplication, $F$ is a vector space over $\mathbb{Q}$, and we say that $F$ has \emph{degree n $($over $ \mathbb{Q })$} if $n$ is the dimension of $F$ over $\mathbb{Q}$. The following definition singles out the subfields of the complex numbers which will be most important for the number theory which we will be most concerned with in all of what follows.

\vspace{0.3cm}
 \emph{Definition}. $F$ is an \emph{algebraic number field} if the degree of $F$ is finite.

\vspace{0.3cm}
We let $F$ denote an algebraic number field of degree $n$ that will remain fixed in the discussion until indicated otherwise. Because the non-negative integral powers of a nonzero element of $F$ cannot form a set that is linearly independent over $\mathbb{Q}$, every element of $F$ is algebraic over $\mathbb{Q}$; this fact is the primary reason why fields of complex numbers which are finite dimensional over $\mathbb{Q}$ are called \emph{algebraic} number fields. 

The structure of the ideals in the ring $R=\mathcal{R} \cap F$ of all algebraic integers contained in $F$ will play a very important role in many situations in which we will be interested. We begin our discussion of those matters by recalling that if $A$ is a commutative ring with identity then an \emph{ideal of A} is a subring $I$ of $A$ such that $ab \in I$ whenever $a \in A$ and $b \in I$. If $s\in A$ then the set $\{as: a\in A\}$ is an ideal of $A$, called the \emph{principal ideal generated by s}. We will denote this ideal by either $sA$ or $(s)$.  In a slight generalization of the principal-ideal notation which will be useful, for two elements $s$ and $t$ of $A$, we will let $(s, t)$ denote the ideal of $A$ generated by $s$ and $t$, i.e., the set $\{as+bt: (a, b)\in A\times A\}$. This notation should not be confused with the ordered pair  $(s, t)$; the meaning of the notation will always be clear from the context in which it is employed. An ideal $I$ of $A$ is \emph{prime} if $\{0\} \not= I \not= A$ and if $a, b$ are elements of $A$ such that $ab \in I$ then $a \in I$ or $b \in I$. An ideal $M$ of $A$ is \emph{maximal} if $\{0\} \not= M \not= A$ and whenever $I$ is an ideal of $A$ such that $M \subseteq I$ then $M=I$ or $I=A$.  

A basic fact in the theory of commutative rings with identity asserts that all maximal ideals in such rings are prime ideals (Hungerford [29], Theorem III.2.19), with the converse false in general. However, in the ring $R$ of algebraic integers in $F$ this converse is true, i.e.,

$(i)$ an ideal of $R$ is prime if and only if it is maximal.

\noindent This fact will be proved in section 1 of Chapter 5. We will also show in section 1 of Chapter 5 that in $R$,

$(ii)$ if $I$ is a nonzero ideal of $R$ then the quotient ring $R/I$ has finite cardinality,

\noindent and also that

$(iii)$ if $P$ is a prime ideal of $R$ then $P$ contains a unique rational prime $q$ and the cardinality of $R/P$ is $q^n$, for a positive integer $n$ uniquely determined by $P$. The integer $n$ is called the \emph{degree of P}.

By far the most important feature of the structure of proper, nonzero ideals of $R$ is the fact that they can be factored in a unique way as the product of prime ideals. We now explain precisely what this means.
\vspace{0.3cm}

\emph{Definition}. Let $A$ be a commutative ring with identity, $I$, $J$ (not necessarily distinct) ideals of $A$. The \emph{$($ideal$)$ product $IJ$ of $I$ and $J$} is the ideal of $A$ generated by the set of products
\[
\{xy: (x, y) \in I \times J\},\]
i.e., $IJ$ is the smallest ideal of $A$, relative to subset inclusion, which contains this set of products.

\vspace{0.3cm}
One can easily show that $IJ$ consists precisely of all sums of the form $\sum_i x_iy_i$, where $x_i \in I$ and $y_i \in J$, for all $i$. It is also easy to show that the ideal product is commutative and associative. We then have
\begin{thm}
$($Fundamental Theorem of Ideal Theory$)$ Every nonzero, proper ideal I of R is a product of prime ideals and this factorization is unique up to the order of the factors. Moreover, the set of prime ideal factors of I is precisely the set of prime ideals of R which contain I, i.e., the set of prime ideals of R containing I is nonempty and finite, and if $\{P_1,\dots,P_k\}$ is this set then there exist a $k$-tuple $(m_1,\dots,m_k)$ of positive integers, uniquely determined by I, such that $I=P_1^{m_1} \cdots P_k^{m_k}$.
\end{thm} 

Theorem 3.16, one of the most important theorems in algebraic number theory, was proved by R. Dedekind in 1871, and appeared as Supplement X in his famous series of addenda to Dirichlet's landmark text \emph{Vorlesungen $\ddot{u}$ber Zahlentheorie} [12]. Because the proof of Theorem 3.16 requires a rather substantial effort, and also because we will not make use of it at its full strength until section 2 of Chapter 5, we will defer its proof until section 5 of Chapter 5. 

Our fourth proof of Theorem 3.3 will require detailed information about how ideals factor according to Theorem 3.16 for a special class of algebraic number fields. Let $m\not=1$ be a square-free integer. Then $\sqrt m$ is an algebraic integer with minimal polynomial $x^2-m$ over $\mathbb{Q}$. It is not difficult to show that the field of complex numbers generated by $\sqrt m$ over $\mathbb{Q}$, i.e., the smallest subfield of the complex numbers containing $\sqrt m$ and $\mathbb{Q}$, the so-called \emph{quadratic number field determined by m}, is  
\[
\mathbb{Q}(\sqrt m)=\{u+v \sqrt m: (u, v) \in \mathbb{Q} \times \mathbb{Q} \}.\]
It is an immediate consequence of this equation that the set $\{1, \sqrt m\}$ is a vector-space basis of  $\mathbb{Q}(\sqrt m)$ over $\mathbb{Q}$, hence $\mathbb{Q}(\sqrt m)$ has degree 2. With a bit more effort, one can also show that 
\[
\mathcal{R} \cap \mathbb{Q}(\sqrt m)=\{k+n\omega: (k, n) \in \mathbb{Z} \times \mathbb{Z}\},\]
where
\[
\omega=\left\{\begin{array}{rl} \sqrt m\ ,& \textrm{if $m \equiv 2\ \textrm{or}\ 3$ or  mod 4,}\\
\displaystyle \frac{1+\sqrt m}{2}\ ,& \textrm{if $m\equiv 1\ \textnormal{mod}\ 4$ } 
\\\end{array}\right. 
\]
(Hecke [27], pp. 95, 96).

Let $F=\mathbb{Q}(\sqrt m), R=\mathcal{R} \cap F$. The proof of Theorem 3.3 which we will present in the next section requires the determination of the prime-ideal factorization of each ideal $qR$ of $R$, $q$ a rational prime, and the calculation of the degree of each factor, in accordance with the conclusions of Theorem 3.16. This is done in 
\begin{prp}
$($Decomposition law in $\mathbb{Q}(\sqrt m))$ Let p be an odd prime.

$(i)$ If $\chi_p(m)=1$ then $pR$ factors into the product of two distinct prime ideals, each of degree $1$. Moreover, for any $a\in \mathbb{Z}$ such that $a^2\equiv m\ \textnormal{mod}\ p$, we can take the prime-ideal factors of $pR$ to be
\[
(p,\ a+\sqrt m)\  and\ (p,\ a-\sqrt m).\]

$(ii)$ If $\chi_p(m)=0$ then pR is the square of a prime ideal I, and the degree of I is $1$. Moreover, we can take I to be the ideal $(p, \sqrt m)$.

$(iii)$ If $\chi_p(m)=-1$ then pR is prime in R, of degree $2$.

If $m \equiv 1\ \textnormal{mod}\ 8$ then

$(iv)$ $2$R factors into the product of two distinct prime ideals, each of degree $1$.  Moreover, we can take the prime-ideal factors of $2$R to be
\[
 \left(2,\ \frac{1+ \sqrt m}{2}\right)\ and\ \left(2,\ \frac{1-\sqrt m}{2}\right),\]
 
If $m \equiv 5\ \textnormal{mod}\ 8$ then 

$(v)$ $2$R is prime in R of degree $2$. 

If $m \equiv 2\ \textnormal{mod}\ 4 $ then 

$(vi)$ $2$R is the square of a prime ideal I, and the degree of I is $1$. Moreover, we may take $I$ to be the ideal $(2, \sqrt m)$. 

 If $m\equiv 3\ \textnormal{mod}\ 4$ then 
 
 $(vii)$ $2$R is the square of a prime ideal I, and the degree of I is $1$. Moreover, we may take $I$ to be the ideal $(2, 1+\sqrt m)$. 
\end{prp}

\emph{Proof}. Hecke [27], section 29, Theorem 90. $\hspace{7.2cm}\ \textrm{QED}$

Let us return now to the general situation of an algebraic number field with its subring $R$ of algebraic integers. The next bit of mathematical technology that we have need of is the ideal class group of $R$. In order to define this group, we first declare  that the ideals $I$ and $J$ of $R$ are \emph{equivalent}, and write $I\sim J$, if there exist nonzero elements $\alpha$ and $\beta$ of $R$ such that $\alpha I=\beta J$. This defines an equivalence relation on the set of all ideals of $R$, and we refer to the corresponding equivalence classes as the \emph{ideal classes} of $R$. If we let $[I]$ denote the ideal class which contains the ideal $I$ then we can define a multiplication on the set of ideal classes by declaring that the product of $[I]$ and $[J]$ is $[IJ]$. It can be shown that when endowed with this product (which is well-defined), the ideal classes of $R$ form an abelian group, called the \emph{ideal-class group of R} (Hecke [27], section 33). It is easy to see that the set of all principal ideals of $R$, i.e., the set of all ideals of the form $\alpha R, \alpha \in R,$ is an ideal class of $R$, called the \emph{principal class}, and one can prove that the principal class is the identity element of the ideal-class group. It is one of the fundamental theorems of algebraic number theory that the ideal-class group is always \emph{finite} (see Hecke [27], section 33, Theorem 96), and the order of the ideal-class group of $R$ is called the \emph{class number of R}. 

Proposition 3.17, when combined with some additional mathematical technology, can be used effectively to compute ideal-class groups and class numbers for quadratic number fields. We illustrate how things go with three examples. But first, the additional technology that is required.

Let $F= \mathbb{Q}(\sqrt m)$ be a fixed  quadratic number field, with $R=\mathcal{A}\cap F$. A subset $\{\omega_1, \omega_2\}$ of $R$ is an \emph{integral basis of R} if every element $\alpha$ in $R$ can be expressed \emph{uniquely} in the form
\[
\alpha=a_1\omega_1+a_2\omega_2,\ \textrm{where}\ (a_1, a_2) \in \mathbb{Z}\times \mathbb{Z}.\]
The ring $R$ always has an integral basis (Hecke [27], section 22, Theorem 64), so we select one, say $\{\omega_1, \omega_2\}$, let $\alpha \in R$, write $\alpha=a_1\omega_1+a_2\omega_2$ for some $(a_1, a_2) \in \mathbb{Z}\times \mathbb{Z}$, and then set
\[
N(\alpha)=(a_1\omega_1+a_2\omega_2)(a_1\omega_1^{\prime}+a_2\omega_2^{\prime}),\]
where the superscripted primes denote the algebraic conjugate taken over $\mathbb{Q}$. 
The number $N(\alpha)$ defined by this formula is called the \emph{norm of $\alpha$}, it does not depend upon the integral basis of $R$ used to define it, it maps $R$ into $\mathbb{Z}$ and it is multiplicative in the sense that
\[
N(\alpha\beta)=N(\alpha)N(\beta)\ \textrm{for all}\ (\alpha, \beta)\in R\times R.\]
One can show that either the subset $\{1, \frac{1}{2}(1+\sqrt m)\}$ or the subset $\{1, \sqrt m\}$ of $R$ is an integral basis of $R$, if either $m$ is, or, respectively, is not, congruent to 1 mod 4. It follows that if $(a_1, a_2) \in \mathbb{Z}\times \mathbb{Z}$ then either
\[
N\left(a_1+a_2\frac{1+\sqrt m}{2}\right)=a_1^2+a_1a_2+\frac{1-m}{4}a_2^2\]
or
\[
N(a_1+a_2\sqrt m)=a_1^2-ma_2^2,\]
if either $m$ is, or, respectively, is not, congruent to 1 mod 4. We also observe that when $m\equiv 1$ mod 4, then
\[
N\left(a_1+a_2\frac{1+\sqrt m}{2}\right)=N\left(\frac{x+y\sqrt m}{2}\right)=\frac{x^2-my^2}{4},\]
where
\[
x=2a_1+a_2\ \textrm{and}\ y=a_2.\]

 If $I$ is an ideal of $R$ then we extend the definition of the norm of numbers in $R$ to ideals by defining the \emph{norm $N(I)$ of I} to be the cardinality of $R/I$. N.B. It follows from observation $(ii)$ above that $N(I)<+\infty$ for all nonzero ideals $I$ of $R$. We also extend the norm to all elements of $F$ by first noting that any integral basis  $\{\omega_1, \omega_2\}$ of $R$ is a vector-space basis of $F$ over $\mathbb{Q}$, and so if $w\in F$, we chose the uniquely determined element $(r_1, r_2)\in \mathbb{Q}\times \mathbb{Q}$ such that $w=r_1\omega_1+r_2\omega_2$ and then set $N(w)$ equal to $(r_1\omega_1+r_2\omega_2)(r_1\omega_1^{\prime}+r_2\omega_2^{\prime})$ as before. This definition of the norm on $F$ also is independent from the integral basis of $R$ used to define it, it maps $F$ into $\mathbb{Q}$, and it is multiplicative on $F$. 

The following lemma is what we require for the calculation of ideal-class groups and class numbers of quadratic fields. It asserts that the norm on ideals of $R$ is multiplicative with respect to the product of ideals, it explains exactly how the norm on ideals extends the norm on elements of $R$, and it concludes with a very useful inequality that will be used  to limit the ideals in $R$ which can determine elements of the ideal-class group. The interested reader can consult Hecke [27], Theorem 29 for a proof of statement $(i)$ of Lemma 3.18, Hecke [27], pp. 87-88 for a proof of statement $(ii)$, and Marcus [40], Corollary 2, p. 136 for a proof of $(iii)$. It will also be convenient here and subsequently to define the \emph{discriminant of $\mathbb{Q}(\sqrt m)$} as either $m$ or $4m$, if $m$  is or, respectively, is not congruent to $1$  mod 4.

\begin{lem}
$(i)$ If I and J are ideals of R, then
\[
N(IJ)=N(I)N(J).\]

$(ii)$ If $0\not= \alpha \in R$ then the norm of the principal ideal generated by $\alpha$ is $|N(\alpha)|$.

$(iii)$ In each ideal class of R there is an ideal I such that
\[
N(I)\leq \lambda=\frac{1}{2}\left(\frac{4}{\pi}\right)^s\sqrt{|d|},\]
where $d$ is the discriminant of $F$ and $s$ is either $0$, if $m$ is positive, or $1$, if $m$ is negative.

\end{lem}
\vspace{0.4cm}
The constant $\displaystyle{\frac{1}{2}\left(\frac{4}{\pi}\right)^s\sqrt{|d|}}$ is called \emph{Minkowski's constant}, and arises in the study of the geometry of numbers.

In order to gain some insight into the techniques that we will employ to prove Theorem 3.3 in the next section, we will now calculate the ideal-class group and the class number of of three quadratic number fields using Proposition 3.17 and Lemma 3.18.

\vspace{0.4cm}
\emph{Example} 1.

Let $F=\mathbb{Q}(\sqrt 2\ ),\ R=\mathcal{A} \cap F$. Then $s=0$ and $d=8$, hence the value of Minkowski's constant $\lambda$ in Lemma 3.18$(iii)$ is
\[
\frac{1}{2}\sqrt 8<2,\]
and so by Lemma 3.18$(iii)$, every ideal class of $R$ contains an ideal $I$ with $N(I)\leq 1$ hence $|R/I|=N(I)=1$, hence $I=(1)$. Conclusion: $R$ has only one ideal class, the principal class, and so $R$ has class number 1.

\vspace{0.4cm}
\emph{Example} 2

Let $F=\mathbb{Q}(\sqrt {-5}\ ),\ R=\mathcal{A} \cap F=\mathbb{Z}+\sqrt {-5}\ \mathbb{Z}.$ Then $s=1$ and $d=-20$, so $\lambda=\displaystyle{\frac{4\sqrt 5}{\pi}}<3$, hence every ideal class of $R$ contains an ideal $I$ such that $N(I)$ is either 1 or 2.

If $N(I)=1$ then $I=(1)$. Suppose that $N(I)=2$. Then the additive group of $R/I$ has order 2, and so
\[
2(\alpha+I)=I,\ \textrm{for all $\alpha \in R$,}\]
and taking $\alpha=1$, we obtain $2 \in I$. Hence all of the prime factors of $I$ must contain 2, so we factor the ideal $(2)$ by way of Proposition 3.17$(vii)$ as
\[
(2)=(2, 1+\sqrt{-5})^2.\]
It follows that $I$ must be a power $J^k$ of $J=(2, 1+\sqrt{-5})$. The degree of $J$ is 1, hence $N(J)=2$. But then Lemma 3.18$(i)$ implies that $2=N(I)=N(J)^k=2^k$, hence $k=1$ and so $I=J$. Conclusion: there are at most \emph{two} ideal classes of $R$, namely $[(1)]$ and $[J]$.

We claim that  $J$ is not principal. If this claim is true then the ideal-class group of $R$ is $\{[(1)], [J]\}$ and $R$ has class number 2.

In order to verify our claim, suppose there exits $\alpha\in R$ such that $J=(\alpha)$. Lemma 3.18$(ii)$ implies that
\[
|N(\alpha)|=N(J)=2,\]
hence $N(\alpha)=2$ (all nonzero elements of $R$ have positive norm). But there exist $a, b\in \mathbb{Z}$ such that $\alpha=a+b\sqrt{-5}$, hence
\[
a^2+5b^2=N(\alpha)=2,\]
and this is clearly impossible.

\vspace{0.4cm}
\emph{Example} 3

Let $F=\mathbb{Q}(\sqrt {-23}\ ),\ R=\mathcal{A} \cap F=\mathbb{Z}+\left(\frac{1}{2}(1+\sqrt {-23})\right)\mathbb{Z}.$ Then $ s=1$ and  $d=-23$ hence $\lambda=\displaystyle{\frac{2\sqrt {23}}{\pi}}<4$, and so every ideal class contains an ideal with norm 1, 2, or 3. As in example 2, every ideal of norm 2 (respectively, 3) must have all of its prime factors containing 2 (respectively, 3), and so factoring via Proposition 3.17 $(i)$ and $(iv)$, we obtain
\[
(2)=\left(2,\ \frac{1+\sqrt {-23}}{2}\right) \left(2,\ \frac{1-\sqrt {-23}}{2}\right)=I_1I_2,\]
\[
(3)=(3,\ 1+\sqrt {-23}) (3,\ 1-\sqrt {-23})=I_3I_4,\]  
hence the ideals of norm 2 are $I_1,\ I_2$ and the ideals of norm 3 are $I_3,\ I_4$. 

It is easily verified that the elements of $R$ are either of the form $a+b\sqrt{-23}$, where $a, b\in \mathbb{Z}$ or $\frac{1}{2}(a+b\sqrt{-23})$, where $a$ and $b$ are \emph{odd} elements of $\mathbb{Z}$. Hence the norm of an element of $R$ is either $a^2+23b^2$ or $\frac{1}{4}(a^2+23b^2)$ for $a, b\in \mathbb{Z}$, neither of which can be 2 or 3. Hence $I_1,\ I_2,\ I_3,$ and $I_4$ are all not principal. It follows that in order to calculate the ideal-class group of $R$, we must determine the inequivalent ideals among $I_1,\ I_2,\ I_3,$ and $I_4$.

We first look at $I_1$ and $I_4$. $I_1\sim I_4$ if and only if $[I_1][I_4]^{-1}=[(1)]$. But $I_3I_4=(3)\sim (1)$, and so $[I_4]^{-1}=[I_3]$, hence we need to see if $I_1I_3$ is principal. 

Lemma 3.18$(i)$ implies that
\begin{equation*}
N(I_1I_3)=N(I_1)N(I_3)=2\cdot3=6.\tag{23}\]

\emph{Claim}: an ideal $I\not= \{0\}$ of $R$ is principal if and only if there exits $\alpha \in R$ such that $N(\alpha)=N(I)$ and there is a generating set $S$ of $I$ such that $s/\alpha\in R$, for all $s\in S$.

The necessity of this is clear. For the sufficiency, let $\alpha\in R$ satisfy the stated conditions. Then $J=(1/\alpha)I$ is an ideal of $R$ and Lemma 3.18 $(i), (ii)$ imply that
\[
N(I)=N\big((\alpha)J\big)=N(I)N(J).\]
hence $N(J)=1$ and so $J=(1)$, whence $I=(\alpha)$.

So in light of (23), we must look for elements of $R$ of norm 6. If $a, b\in \mathbb{Z}$ then $a^2+23b^2\not= 6$; on the other hand, 
\[
6=\frac{a^2+23b^2}{4}\]
if and only if $a^2=1=b^2$. Hence there are exactly two principal ideals of norm 6: $\left(\frac{1}{2}(1\pm \sqrt{-23}\ )\right)$. Let $\alpha=\displaystyle{\frac{1+\sqrt {-23}}{2}}$. We have that
\[
I_1I_3=\left(6,\ 2+2\sqrt{-23},\ \frac{3}{2}\big(1+\sqrt {-23}\ \big),\ \frac{1}{2}\big(1+\sqrt {-23}\ \big)^2\right).\]
Now divide each of these generators by $\alpha$: you always get an element of $R$. Hence by the claim, $I_1I_3=(\alpha)$, and so $I_1\sim I_4$. Following the same line of reasoning also shows that $I_2\sim I_3$. 

We now claim that $I_1$ is \emph{not} equivalent to $I_2$. Otherwise, $[(1)]=[(2)]=[I_1I_2]=[I_1^2]$, hence there exits $\alpha\in R$ such that $I_1^2=(\alpha)$, and so $N(\alpha)=N(I_1^2)=4$, whence $\alpha=\pm 2$. But then $(2)I_1=I_1^2I_2=(2)I_2$, hence $I_1=I_2$, which contradicts the fact that these ideals are distinct.

It follows that the ideal-class group of $R$ is $\{[(1)], [I_1], [I_2]\}$, and $R$ has class number 3. Since the ideal-class group is of prime order, it is cyclic, and since the order is 3, both ideal classes $[I_1]$ and $[I_2]$ are generators of the group.
\section{Proof of Quadratic Reciprocity via Ideal Theory: Conclusion}

Let $F=\mathbb{Q}(\sqrt m)$ be a quadratic number field. The proof of Theorem 3.3 which we give in this section depends on an equivalence relation defined on the ideals of $R=\mathcal{A}\cap F$ which is similar to, but possibly different from, the equivalence relation $\sim$ which determines the ideal classes of $R$.  If $I$ and $J$ are ideals of $R$ then we declare that \emph{I is equivalent to J in the narrow sense}, and write $I\approx J$, if there exists an element $s\in F$ such that $N(s)>0$ and $I=sJ$. The relation $\approx$ is clearly an equivalence relation, and we will call the corresponding set of equivalence classes \emph{narrow ideal classes of $R$}. Because each element of $F$ is the quotient of two elements of $R$, it follows that $I\approx J$ implies that $I\sim J$, hence each narrow ideal class is a subset of some ideal class. If the norm function $N$ is always positive on $F$, which occurs when $m<0$, then there is no difference between ordinary equivalence of ideals and equivalence in the narrow  sense.

On the other hand, when $m>0$ the difference between these two equivalence relations is mediated by the units in $R$, i.e., the elements of $R$ which have a multiplicative inverse in $R$. It is easy to see that an element $u$ in $R$ is a unit if and only if $N(u)=\pm1$. If $R$ has a unit of norm $-1$ then there is also no difference between the two equivalence relations because if $I=kJ$ for some $k \in F$, then upon multiplication of the element $k$ by a suitable unit of norm $-1$, one can insure that the element $k$ in this equation has positive norm. On the other hand, if there is no unit with a negative norm then it follows easily from this assumption and the fact that $N(\sqrt m)=-m<0$ that every nonzero ideal $I$ of $R$ is not narrowly equivalent to $\sqrt m\ I$. Hence each ideal class $[I]$ in the ordinary sense in the union of exactly two ideal classes in the narrow sense, namely the narrow ideal class containing $I$ and the narrow ideal class containing $\sqrt m\ I$. If we let $h$ denote the class number of $R$ and $h_0$ denote the number of narrow ideal classes of $R$, i.e., $h_0$ is the \emph{narrow class number of R}, it follows that either $h_0=h$ or $h_0=2h$; in particular the number of narrow ideal classes is finite.

The set of all narrow ideal classes can be given an abelian group structure using multiplication of ideals in the same way as we did for the ideal-class  group, with the set of all nonzero principal ideals $(\mu)$ with $N(\mu)>0$ acting as the identity element. We call this group the \emph{narrow ideal-class group of R}. The algebraic relationship between the ideal-class group and the narrow ideal-class group can be described by considering the set $\mathcal{I}$ of all nonzero ideals of $R$ as an abelian semigroup under the ideal product and containing the semigroup $H$ of all nonzero principal ideals of $R$. The quotient semigroup $\mathcal{I}/H$ is then in fact a group which is isomorphic to the ideal-class group. If we let $H_0$ denote the semigroup of all nonzero principal ideals $(\mu)$ with $N(\mu)>0$ then the quotient semigroup $\mathcal{I}/H_0$ is a group which is isomorphic to the narrow ideal-class group.

The next theorem is the basis for the proof of Theorem 3.3 that will be presented here. The proof of the theorem, rather long and technical, employs several results from algebraic number theory that would take us too far afield to clearly explain, so we will be content to cite Hecke [27], proof of Theorem 132, for the details. We will use the theorem to deduce two corollaries from which Theorem 3.3 will follow by an elegant argument. 

The statement of the theorem requires a bit of terminology from the theory of finite abelian groups. If $G$ is such a group whose order exceeds 1 then $G$ is isomorphic to a uniquely determined finite direct sum of cyclic groups of prime-power order. If $q$ is a prime number then the \emph{basis number of q belonging to G} is the number $b(q)$ of cyclic summands of $G$ whose orders are all divisible by $q$. It follows that if $b(q)=0$ then the order of $G$ is not divisible by $q$.

\begin{thm}
If t denotes the number of distinct prime ideals in $R$ which contain the discriminant of F then the basis number of $2$ belonging to the narrow ideal-class group of R is $t-1$.
\end{thm}

\begin{cor}
If the discriminant of $F$ is divisible by a single prime then the order of the narrow ideal-class group is odd.
\end{cor}

\emph{Proof}. If the discriminant of $F=\mathbb{Q}(\sqrt m)$ is divisible by a single prime then either $m=\pm2$ or $m\equiv $ 1 mod 4 and either $m$ or $-m$ is prime. If $m=\pm2$ then the discriminant is $\pm8$ and so from Proposition 3.17$(vi)$, it follows that $(\pm8)=(2)^3$ has only a single prime-ideal factor, i.e., $t=1$. If $m\equiv $ 1 mod 4 and either $m$ or $-m$ is prime  then it follows from Proposition 3.17$(ii)$ that $(m)$ has only a single prime-ideal factor, hence $t=1$ here as well. Theorem 3.19 then implies that the basis number of 2 belonging to the narrow ideal-class group is 0. Hence the order of the narrow ideal-class group is not divisible by 2.$\hspace{15.2cm}\ \textrm{QED}$
\begin{cor}
If the discriminant of F is the product of two positive primes p and q, each congruent to $3\ \textnormal{mod}\ 4$, then either p or q is the norm of an element of R.
\end{cor}
\emph{Proof}. We prove first that the norm of each unit of $R$ is 1. If not, i.e., there is an element $\alpha$ of $R$ with norm $-1$, then there exist rational integers $x$ and $y$ such that
\[
x^2-pqy^2=-4
\]
($pq \equiv 1$ mod 4), hence
\[
-4\equiv x^2\ \textrm{mod}\ pq,\]
and so $-1$ is a quadratic residue of $p$. However, Theorem 2.4 implies that $\chi_p(-1)=-1$, which is not possible.

Now use Proposition 3.17$(ii)$ to factor $(p)$ and $(q)$ as
\begin{equation*}
(p)=I^2,\ (q)=J^2,\tag{24}
\]
where $I$ and $J$ are prime ideals such that
\begin{equation*}
N(I)=p,\ N(J)=q.\tag{25}\]
The discriminant of $\mathbb{Q}(\sqrt{pq})$ is $pq$, hence any prime ideal which is a factor of $(pq)$ must contain either $p$ or $q$, and hence must be equal to either $I$ or $J$. If $(pq)$ has only a single prime factor $Q$, say, then $Q\cap\mathbb{Z}$ is a prime ideal of $\mathbb{Z}$ which contains both $p$ and $q$, which is not possible. It hence follows that $t=2$ in Theorem 3.19, and consequently, the proof of Theorem 3.19 (Hecke [27], pp. 160-162) implies that 
\begin{equation*}
I^{\varepsilon_1}J^{\varepsilon_2}\approx (1),\tag{26}\]
where $\varepsilon_i\in \{0, 1\}, i=1, 2$, and $\varepsilon_1\not=0\not= \varepsilon_2$.

If $\varepsilon_1=\varepsilon_2=1$, then (24) and (26) imply that
\[
(\sqrt{pq})=IJ\approx (1),\]
and so the definition of equivalence in the narrow sense implies that there exists $\alpha \in R$ with $N(\alpha)>0$ such that
\[
(\sqrt{pq})=(\alpha).\]
We conclude the existence of a unit $u$ of $R$ such that $u\sqrt{pq}=\alpha$, hence
\[
N(u)=\frac{N(\alpha)}{N(\sqrt{pq})}=-\frac{N(\alpha)}{pq}<0,\]
which cannot happen because every unit of $R$ has norm 1. Hence either $\varepsilon_1$ or $\varepsilon_2$ is 1 and the other is 0, and consequently from (26) it follows that either $I$ or $J$ is principal and is generated by an element $\gamma$ of $R$ of positive norm. Hence by (25) and Lemma 3.18$(ii)$, either $p$ or $q$ is the norm of $\gamma$. $\hspace{13.5cm}\ \textrm{QED}$

With Corollaries 3.20 and 3.21 in hand, we can now prove Theorem 3.3. Given distinct odd primes $p$ and $q$, we wish to prove that
\begin{equation*}
\chi_p(q) \chi_q(p)=(-1)^{\frac{1}{2}(p-1) \frac{1}{2}(q-1)}\tag{27}.
\]

It will be most convenient to divide the reasoning into the three cases which determine the sign on the right-hand side of (27).

\emph{Case} 1.   

Suppose that $p\equiv q\equiv 1$ mod 4. We will show that $\chi_p(q)$ and $\chi_q(p)$ are simultaneously 1, hence also simultaneously $-1$, hence both sides of (27) are 1.

Assume that $\chi_p(q)=1$. Then according to Proposition 3.17$(i),$ the ideal $(p)$ in $R=\mathcal{R} \cap \mathbb{Q}(\sqrt q)$ factors as $IJ$ with each of the factors having degree 1. If $h_0$ is the narrow class number of $R$ then $I^{h_0}$  is in the narrow principal class of $R$, hence there exists an element $\alpha$ in $R$ of positive norm such that
\begin{equation*}
I^{h_0}=(\alpha).\tag{28}\]
Because there are rational integers $x$ and $y$ such that $\alpha=\frac{1}{2}(x+y\sqrt q)$, we take the norm of both sides of (28) and use the facts that $p\in I$, the degree of $I$ is 1, and $\alpha$ has positive norm to conclude that
\[
p^{h_0}=\frac{x^2-qy^2}{4},\]
from which it follows that
\[
4p^{h_0}\equiv x^2\ \textrm{mod}\ q.\]
This means that $4p^{h_0}$ is a residue of $q$, hence
\[
\chi_q(p)^{h_0}=\chi_q(4p^{h_0})=1.\]
Because $\mathbb{Q}(\sqrt q)$ has discriminant $q$, it follows from Corollary 3.20 that $h_0$ is odd, and so this equation implies that $\chi_q(p)=1$. An interchange of the roles of $p$ and $q$ in this argument also shows that $\chi_q(p)=1$ implies that  $\chi_p(q)=1$.

\emph{Case} 2.

Suppose that $q\equiv 1$ mod 4 and $p\equiv 3$ mod 4. The argument in Case 1 shows that if $\chi_p(q)=1$ then $\chi_q(p)=1$. Hence by Theorem 2.4,
\[
\chi_q(-p)=\chi_q(-1)\chi_q(p)=1\]
Conversely, if $\chi_q(-p)=1$, then we can apply the argument in Case 1, using the field $\mathbb{Q}(\sqrt {-p})$, to obtain $\chi_p(q)=1$. It follows that
\[
\chi_p(q)=\chi_q(-p)=\chi_q(p),\]
and both sides of (27) are again equal to 1.

\emph{Case} 3.

Suppose that $p\equiv q\equiv 3$ mod 4. Applying the same reasoning as we did in Case 1 or 2, it follows that $\chi_q(-p)=1$ implies that $\chi_p(-q)=-1$, but verification of the converse cannot be proved in that way. Instead, we work in the field $\mathbb{Q}(\sqrt{pq})$, in which, according to Corollary 3.21, $p$ or $q$ is the norm of an algebraic integer
$\frac{1}{2}(x+y\sqrt{pq})$. If $p$ is that norm then
\[
4p=x^2-pqy^2.\]
This equation implies that $x$ is divisible by $p$, say $x=ap$, and so it follows that $4=pa^2-qy^2$. Upon observing that $a$ is not divisible by $q$ and $y$ is not divisible by $p$, it hence follows from this equation that
\[
\chi_q(p)=\chi_q(pa^2)=\chi_q(pa^2-qy^2)=\chi_q(4)=1,\]
and similarly,
\[
\chi_p(-q)=1.\]
Hence by Theorem 2.4 again, it follows that
\[
\chi_p(q)=\chi_p(-1)\chi_p(-q)=-1.\]
If $q$ is the prime that is the norm of an element of $R$, the same argument applies to show that $\chi_p(q)$ and $\chi_q(p)$ still have opposite signs. Thus (27) is verified for this final case.$\hspace{1cm}\ \textrm{QED}$

Before we move on to the last proof of Theorem 3.3 that will be presented in this chapter, we will explain how the proof that we just gave is related to Gauss' second proof of quadratic reciprocity. In order to do that we need to describe how to get quadratic forms from ideals of quadratic number fields. In the discussion which follows, we will eschew the proof of the results mentioned; for those please consult Landau [35], Part Four, Chapters I-IV or Hecke [27], section 53.

Recall from Chapter 1 that a (binary) quadratic form is a polynomial in the variables $x$ and $y$ of the form 
\[
ax^2+bxy+cy^2\]
where $(a, b, c)\in \mathbb{Z}\times \mathbb{Z}\times \mathbb{Z}$, and we will denote this form by $[a, b, c]$. The \emph{discriminant of} $[a, b, c]$ is the familiar algebraic invariant $d=b^2-4ac$, and in the classical theory, one presupposes the form is \emph{irreducible}, i.e., it is not a product of linear factors with integer coefficients, so that the discriminant is not a perfect square and is congruent to either 0 or 1 mod 4. A discriminant is \emph{fundamental} if $\gcd(a, b, c)=1$. It can be shown that the set of fundamental discriminants consists of precisely the integers, positive and negative, which are either square-free and congruent to 1 mod 4 or are of the form $4n$, where $n$ is square-free and congruent to 2 or 3 mod 4. Thus the set of fundamental discriminants of quadratic forms coincides with the set of discriminants of quadratic number fields.

For each fundamental discriminant $d$, let $\mathcal{Q}(d)$ denote the set of all irreducible quadratic forms with discriminant $d$, so that $\mathcal{Q}(d)$ consists of all irreducible forms $[a, b, c]$ such that $\gcd(a, b, c)=1$ and $b^2-4ac=d$. It transpires that there is a way to manufacture certain elements of $\mathcal{Q}(d)$ from the ideals in the ring of algebraic integers in $\mathbb{Q}(\sqrt {m(d)} )$, where $m(d)=d$ if $d\equiv 1$ mod 4, and $m(d)=d/4$ if $d\equiv 0$ mod 4. In order to describe this procedure, we first single out the forms in $\mathcal{Q}(d)$ that will arise from it.

The set of quadratic forms that we need is determined by the manner in which quadratic forms represent the integers. If $n$ is an integer and $q(x, y)$ is a quadratic form, we will say that $n$ is \emph{represented by} $q(x, y)$ if there exist integers $x$ and $y$ such that $n=q(x, y)$. The sign of the integers which a given quadratic form $q=[a, b, c]$ in $\mathcal{Q}(d)$ represents depends on the sign of the discriminant $d$. If $d>0$ then $q$ represents both positive and negative integers. If $d<0$ and $a>0$ then $q$ represents no negative integers, and $q(x, y)$ represents 0 only if $x=y=0$. If $d<0$ and $a<0$ then $q$ represents no positive integers, and $q(x, y)$ represents 0 only if $x=y=0$. Hence forms with positive discriminant are called \emph{indefinite} and forms with negative discriminant are called \emph{positive} or \emph{negative definite} if $a$ is, respectively, positive or negative (Hecke [27], section 53, Theorem 153).

Now let $R=\mathcal{R}\cap \mathbb{Q}(\sqrt{m(d)} )$ and let $\mathcal{I}(d)$ denote the set of all nonzero ideals of $R$. For each $I\in \mathcal{I}(d)$, we  choose an integral basis  $\{\alpha, \beta\}$ of $I$ such that $\alpha \beta'-\alpha' \beta=N(I)\sqrt d$ is positive or pure imaginary with positive imaginary part (an integral basis with this property always exits, according to  Hecke [27], p.190). If $(x, y)\in \mathbb{Z}\times \mathbb{Z}$ then we let
\[
Q_I(x, y)=\frac{(x\alpha+y\beta)(x\alpha^{\prime}+y\beta^{\prime})}{N(I)}.\]
One can show that if $I\in \mathcal{I}(d)$ then $Q_I\in \mathcal{Q}(d)$. Moreover, if $d>0$, and $f \in\mathcal{Q}(d)$ then there is an $I\in \mathcal{I}(d)$ such that $Q_I=f$, and if  $d<0$  then for each positive definite form $f \in\mathcal{Q}(d)$, there is an $I\in \mathcal{I}(d)$ such that $Q_I=f$ (Hecke [27], pp.190-192).

The relation of narrow equivalence of ideals in $\mathcal{I}(d)$ has an important connection to an equivalence relation on the set $\mathcal{Q}(d)$. We declare that forms $q(x, y)=ax^2+bxy+cy^2$ and  $q_1(X, Y)=a_1X^2+b_1XY+c_1Y^2$ in $\mathcal{Q}(d)$ are \emph{equivalent} if there is a linear transformation defined by
\[
x=\alpha X+\beta Y,\ y=\gamma X+\delta Y,\]
where $\alpha, \beta, \gamma,$ and $\delta$ are integers satisfying $\alpha \delta-\beta \gamma=1$, such that 
\[
q(\alpha X+\beta Y,\ \gamma X+\delta Y)=q_1(X, Y).\]
These transformations are called \emph{modular substitutions}, and each modular substitution maps $\mathcal{Q}(d)$ bijectively onto $\mathcal{Q}(d)$. It is a classical result of Lagrange that each equivalence class of forms determined by this equivalence relation contains a form $[a, b, c]$ whose coefficients satisfy
\[
|b|\leq |a|\leq |c|,\]
and it follows from this fact that the number of equivalence classes is finite (Landau [35], Theorem 197). There is always at least one form in $\mathcal{Q}(d)$, called the \emph{principal form}, defined by
\[
x^2-\frac{1}{4}dy^2,\ \textrm{if $d\equiv 0$ mod 4,}\]
or
\[
x^2+xy-\frac{1}{4}(d-1)y^2,\ \textrm{if $d\equiv 1$ mod 4,}\]
hence the number of equivalence classes is a positive integer. In what follows, when we speak of a class of quadratic forms, we will mean one of these equivalence classes.

One can now prove the following very important theorem:
\begin{thm}
For each $I\in \mathcal{I}(d)$ the class of the form $Q_I$ does not depend on the integral basis of I used to define it, and ideals $I$ and $J$ in $\mathcal{I}(d)$ are in the same narrow ideal class if and only if $Q_I$ and $Q_J$ are in the same class of forms in $\mathcal{Q}(d)$.
\end{thm}

\emph{Proof}. Hecke [27], section 53, Theorem 154. $\hspace{7.1cm}\ \textrm{QED}$

As we mentioned in section 10, the second proof which Gauss gave for the LQR uses his genus theory of quadratic forms. The genus which a quadratic form belongs to is an equivalence class determined by yet another equivalence relation that Gauss defined on $\mathcal{Q}(d)$. The definition of this equivalence relation is based on a very subtle and detailed analysis of the manner by which a quadratic form represents odd integers, even integers, and the residues and non-residues of primes which divide the discriminant of the form. It is done in such a way that each genus is the union of certain classes of forms determined by modular substitutions. Because of the rather daunting complexity of Gauss' definition, we are instead going to define a notion of genus on the set $\mathcal{I}(d)$, which can be done rather more transparently, and then lift it to $\mathcal{Q}(d)$ via Theorem 3.22. What we will end up with is Gauss' definition of genera of quadratic forms.

Toward that end, we declare that nonzero ideals $I$ and $J$ of $R$ such that $I+(d)=R=J+(d)$ \emph{have the same genus} if there exists a number $\gamma$ in $\mathbb{Q}(\sqrt {m(d)} )$ such that
\[
N(I)\equiv |N(\gamma)|N(J)\ \textrm{mod}\ (d).\]
On the set of nonzero ideals which are relatively prime to $(d)$, i.e., the nonzero ideals whose sum with $(d)$ is $R$, this congruence defines an equivalence relation which partitions the nonzero ideals relatively prime to $(d)$ into \emph{genera}. The genera form an abelian group in the usual way with the identity in this group given by the  the genus containing $(1)$, thus called the \emph{principal genus}, and hence containing the set of all principal ideals which are generated by the elements of $R$ of positive norm. It is not difficult to see that ideals that are equivalent in the narrow sense belong to the same genus if they are all relatively prime to $(d)$; consequently each genus is the union of certain narrow ideal classes. The narrow classes belonging to the principal genus form a subgroup of the narrow ideal-class group, and so if $f$ is the order of this subgroup, $g$ is the number of genera, and $h_0$ is the narrow class number of $R$, then each genus is the union of exactly $f$ narrow ideal classes and $h_0=fg$.

Returning to $\mathcal{Q}(d)$, we take forms $q_1$ and $q_2$ in $\mathcal{Q}(d)$ (both positive definite if $d<0$), choose ideals $I_1$ and $I_2$ in $\mathcal{I}(d)$ such that $q_j=Q_{I_j}$ for $j=1, 2$, and we define the class containing $q_1$ and the class containing $q_2$ to be \emph{in the same genus} if $I_1$ and $I_2$ have the same genus. By virtue of Theorem 3.22, this relation is well-defined, each genus of forms in $\mathcal{Q}(d)$ is the union of $f$ classes of forms, the number of genera is $g$, and the total number of classes of forms (classes of positive definite forms if $d<0$) is $fg$.

The fundamental problem of genus theory is the determination of the number of genera. The solution to this problem is given in the following theorem:
\begin{thm}
If $t$ is the number of distinct prime ideals of $R$ which contain $d$ then the number of genera is $2^{t-1}$. Moreover, a narrow ideal class is the square of a narrow ideal class if and only if it is contained in the principle genus.
\end{thm}

\emph{Proof}. Hecke [27], section 48, Theorem 145.$\hspace{7cm}\ \textrm{QED}$

This theorem, due to Gauss in an equivalent form  ([19], articles 261, 286, 287), is the basis of his second proof of quadratic reciprocity. Actually Gauss only used the fact that the number of genera does not exceed $2^{t-1}$ in his argument, because if $d$ has at most two prime factors, then this inequality implies that there are only at most two genera  in $\mathcal{Q}(d)$, one of which is always the principal genus. Gauss then used his definition of genera to deduce quadratic reciprocity from this fact by means of an argument whose line of reasoning, when converted into the language of ideals, is very similar to the one that is used in the proof of quadratic reciprocity given in this section. Interestingly enough, Gauss then used quadratic reciprocity and his theory of composition of forms to derive the reverse inequality that the number of genera can be no less that $2^{t-1}$.

\section{A Proof of Quadratic Reciprocity via Galois Theory}

The proof of Theorem 3.3 which we gave in the last section shows that quadratic reciprocity results from certain factorization properties of the ideals in a quadratic number field. In this final section of Chapter 3, we derive quadratic reciprocity from the structure of the Galois group of certain cyclotomic number fields.

We preface the argument by recalling some relevant facts from Galois theory. Let $K\subseteq F$ be an inclusion of fields; we say that $F$ is an \emph{extension of K}. An automorphism $\sigma$ of $F$ is \emph{Galois over K} if each element of $K$ is fixed by $\sigma$, i.e., $\sigma(k)=k$ for all $k \in K$. The set of all automorphisms of $F$ which are Galois over $K$ forms a group under composition of automorphisms, called the \emph{Galois group of F over K}, and denoted by $G_K(F)$. Now let $E$ be an \emph{intermediate field}, i.e., a subfield of $F$ which contains $K$, and also let $H$ be a subgroup of $G_K(F)$. We let
\[
E'=\{\sigma \in G_K(F): \sigma(\xi)=\xi,\ \textrm{for all}\ \xi\in E\},\]
\[
H'=\{\xi \in F: \sigma(\xi)=\xi,\ \textrm{for all}\ \sigma \in H\}.\]
It is clear that $E'=G_E(F)$ is a subgroup of $G_K(F)$ and $H'$, the fixed set of $H$, is an intermediate subfield of $F$. The field $F$ is \emph{Galois over $K$} if $G_K(F)'=K$, i.e., the fixed field of $G_K(F)$ is $K$. 

The field $F$ is naturally a vector space over $K$ with respect to the addition and multiplication in $F$; the dimension of $F$ as a vector space over $K$ is called the \emph{degree of F over K} and is denoted by $[F:K]$. $F$ is a \emph{finite extension of} $K$ if the degree of $F$ over $K$ is finite. The following theorem, often called the Fundamental Theorem of Galois Theory, plays, as the name connotes, a central role in the theory of fields.
\begin{thm}
If F is a finite Galois extension of K then the mapping $E\rightarrow E'$ is a bijection of the set of all intermediate fields E onto the set of all subgroups of $G_K(F)$ whose inverse mapping is $H\rightarrow H'$ and which has the following properties:

$(i)$ F is a Galois extension of $E'$and the order of $G_{E'}(F)$ is $[F:E']$. In particular, the order of $G_K(F)$ is $[F:K]$.

$(ii)$ E is a Galois extension of K if and only if $E'$ is a normal subgroup of $G_K(F)$, and if E is Galois over K then $G_K(E)$ is isomorphic to the quotient group $G_K(F)\big/ E'$.
\end{thm}

\emph{Proof}. Hungerford [29], Theorem V.2.5. $\hspace{7.5cm}\ \textrm{QED}$

In particular, if $F$ is a finite Galois extension of $K$ and $G_K(F)$ is \emph{abelian}, i.e., $F$ is an \emph{abelian extension of K}, then all intermediate subfields $E$ are Galois extensions of $K$,
\[
\textrm{the order of}\ G_K(E)\ \textrm{is}\ [E:K]\ \textrm{and}\]
\[
G_K(E)\ \textrm{is isomorphic to}\ G_K(F)\big/G_E(F).\]

We now specialize to the case when $F$ is an algebraic number field, the situation that is of most interest to us. A natural question that arises after one contemplates Theorem 3.24 asks: when is an algebraic number field a Galois extension of $\mathbb{Q}$? The answer involves the concept of a splitting field of a polynomial, an idea that we have already encountered in section 1 of this chapter. 

If $K\subseteq F$ are fields of complex numbers, then $F$ is a \emph{splitting field over K} if $F$ is generated over $K$ by the roots of a polynomial $f(x)\in K[x]$, in other words, $F$ is the smallest subfield of the complex numbers which contains $K$ and all the roots of $f(x)$. In particular, if this holds we also say that $F$ is the \emph{splitting field of $f(x)$ over K}. The next theorem is true in much greater generality, but it will be more than sufficient to meet our needs.
\begin{thm}
An algebraic number field is Galois over $\mathbb{Q}$ if and only if it is the splitting field of a polynomial in $\mathbb{Q}[x]$.
\end{thm}

\emph{Proof}. Hungerford [29], Theorem V.3.11.$\hspace{7.6cm}\ \textrm{QED}$

In order to verify Theorem 3.3, we are going to work in the algebraic number field $\mathbb{Q}(\zeta)$ generated over $\mathbb{Q}$ by the $p$-th root of unity $\zeta=\exp(2\pi i/p)$, for a fixed odd prime $p$, the \emph{cyclotomic number field determined by $p$}. As we saw in the examples from section 8, $\zeta$ is an algebraic integer with minimal polynomial over $\mathbb{Q}$ given by
\[
1+x+\dots+x^{p-1}.\]
It follows from this fact and it is not too difficult to prove that
\begin{equation*}
\mathbb{Q}(\zeta)=\left\{\sum_{i=0}^{p-1}\ r_i\zeta^i: (r_0,\dots,r_{p-1})\in \mathbb{Q}^p\right\}\tag{29}\]
(Hecke [27], section 30, p. 98). It is also true that
\begin{equation*}
R(\zeta)=\mathcal{R}\cap \mathbb{Q}(\zeta)=\left\{\sum_{i=0}^{p-1}\ z_i\zeta^i: (z_0,\dots,z_{p-1})\in \mathbb{Z}^p\right\},\tag{30}\]
although the proof of this requires quite a bit more work, for instance, see Marcus [40], Theorem 10, p. 30.

From the factorization
\[
1+x+\dots+x^{p-1}=\prod_{i=0}^{p-1}\ (x-\zeta^i),\]
 (29), and Theorem 3.25 it follows that $\mathbb{Q}(\zeta)$ is Galois over $\mathbb{Q}$. We proceed to calculate the Galois group $G$ of $\mathbb{Q}(\zeta)$ over $\mathbb{Q}$. Recall that $U(p)$ denotes the group of units in $\mathbb{Z}/p\mathbb{Z}$; this group is cyclic of order $p-1$.
 \begin{prp}
 There is an isomorphism $\theta$ of G onto $U(p)$ such that for $\sigma \in G$,
 \[
 \sigma(\zeta)=\zeta^{\theta(\sigma)}.\]
 \end{prp}

\emph{Proof}. Because $\zeta^p=1$, it follows that $\sigma(\zeta)^p=1$. Hence
\[
\sigma(\zeta)=\zeta^{\theta(\sigma)},\]
where $\theta(\sigma)$ is an integer determined uniquely by $\sigma$ modulo $p$. If $\tau=\sigma^{-1}$ then
\[
\zeta=\tau \sigma(\zeta)=\tau\big(\zeta^{\theta(\sigma)}\big)=\zeta^{\theta(\tau)\theta(\sigma)},\]
hence $\theta(\tau)\theta(\sigma)$ is in the coset in $\mathbb{Z}/p\mathbb{Z}$ containing 1. It follows that $\theta$ maps $G$ into $U(p)$. If $\tau, \sigma \in G$ then
\[
\zeta^{\theta(\tau\sigma)}=(\tau\sigma)(\zeta)=\tau\big(\sigma(\zeta)\big)=\zeta^{\theta(\tau)\theta(\sigma)},\]
hence $\theta(\tau\sigma)\equiv \theta(\tau)\theta(\sigma)$ mod $p$ and so $\theta$ is a homomorphism. If $\theta(\sigma)\equiv 1$ mod $p$ then $\sigma(\zeta)=\zeta$, hence by (29), $\sigma(\alpha)=\alpha$, for all $\alpha\in \mathbb{Q}(\zeta)$, whence $\theta$ is a monomorphism. As $\mathbb{Q}(\zeta)$ is Galois over $\mathbb{Q}$, we have that
\[
\big|U(p)\big|=p-1=[\mathbb{Q}(\zeta):\mathbb{Q}]=|G|,\]
hence $\theta$ is surjective.$\hspace{11.6cm}\ \textrm{QED}$

It is a consequence of Proposition 3.26 that for every integer $a\in \mathbb{Z}$ not divisible by $p$, there exists $\sigma_a\in G$ such that
\[
\sigma_a(\zeta)=\zeta^a,\]
and the map $a\rightarrow \sigma_a$ is the inverse of $\theta$.
\begin{lem}
If q is a prime distinct from p then for all $w\in R(\zeta),\ \sigma_q(w)\equiv w^q$ \textnormal{mod} $qR(\zeta)$.
\end{lem}

\emph{Proof}. From (30), we have that
\[
w=\sum_i\ z_i\zeta^i,\]
where $z_i\in \mathbb{Z}$, for all $i$. Since $\sigma_q(\zeta)=\zeta^q$, it follows from Fermat's little theorem that
\[
\sigma_q(w)\equiv \sum_i\ z_i^{q}\zeta^{qi}\ \textrm{mod}\ qR(\zeta).\]
Because the ring $R(\zeta)/qR(\zeta)$ has characteristic $q$, the $q$-th power map in this ring is additive,
hence
\begin{eqnarray*}
\sigma_q(w)&\equiv&\sum_i\ z_i^q\zeta^{qi}\ \textrm{mod}\ qR(\zeta)\\
&\equiv&\Big(\sum_i\ z_i\zeta^i\Big)^q\ \textrm{mod}\ qR(\zeta)\\
&\equiv& w^q\ \textrm{mod}\ qR(\zeta).
\end{eqnarray*}
$\hspace{15.4cm}\ \textrm{QED}$

The way is now clear to the proof of Theorem 3.3. We begin by looking for a square root of $(-1)^{\frac{1}{2}(p-1)}p=p^{*}$ in $\mathbb{Q}(\zeta)$. This can be found by applying Theorem 3.14, but we want to avoid the use of Gauss sums. Instead, we will find this square root via the equation
\[
p=\prod_{i=1}^{p-1}\ (1-\zeta^i).\]
If the terms corresponding to $i$ and $p-i$ are combined, then
\[
(1-\zeta^i)(1-\zeta^{p-i})=(1-\zeta^i)(1-\zeta^{-i})=-\zeta^{-i}(1-\zeta^i)^2.\]
Hence
\[
p=(-1)^{\frac{1}{2}(p-1)}\zeta^n\prod_{i=1}^{\frac{1}{2}(p-1)}\ (1-\zeta^i)^2,\ \textrm{where}\ n=-\sum_{i=1}^{\frac{1}{2}(p-1)}\ i.\]
Now choose $z\in \mathbb{Z}$ such that $2z\equiv$ 1 mod $p$. Then $\zeta^n=(\zeta^{nz})^2$, and so
\[
p^{*}=\left(\zeta^{nz}\prod_{i=1}^{\frac{1}{2}(p-1)}\ (1-\zeta^i)\right)^2,\]
whence $p^{*}=\tau^2$ for some $\tau\in \mathbb{Q}(\zeta)$.

Let $q$ be an odd prime distinct from $p$. Then
\[
\sigma_q(\tau)^2=\sigma_q(\tau^2)=\sigma_q(p^{*})=p^{*}=\tau^2,\]
hence $\sigma_q(\tau)=\pm \tau$, with the plus sign holding if and only if $\sigma_q$ is in the Galois group of $\mathbb{Q}(\zeta)$ over $\mathbb{Q}(\tau)$. Proposition 3.26 implies that $G$ is cyclic of order $p-1$, hence abelian, and so we conclude from Theorem 3.24 that 
\[
G\big/G_{\mathbb{Q}(\tau)}\big(\mathbb{Q}(\zeta)\big)\ \textrm{is isomorphic to}\ G_{\mathbb{Q}}\big(\mathbb{Q}(\tau)\big)\]
and
\[
\big|G_{\mathbb{Q}}\big(\mathbb{Q}(\tau)\big)\big|=[\mathbb{Q}(\tau):\mathbb{Q}].\]
As the minimal polynomial of $\tau$ over $\mathbb{Q}$ is $x^2-p^{*}$, it follows that $[\mathbb{Q}(\tau):\mathbb{Q}]=2$, and so $G\big/G_{\mathbb{Q}(\tau)}\big(\mathbb{Q}(\zeta)\big)$ is cyclic of order 2. Because $G$ is cyclic of even order, we conclude that $G_{\mathbb{Q}(\tau)}\big(\mathbb{Q}(\zeta)\big)$ consists of all the squares of the elements of $G$. Hence $\sigma_q(\tau)=\tau$ if and only if $\sigma_q$ is a square in $G$. Since the map $a\rightarrow \sigma_a$ is an isomorphism of $U(p)$ onto $G$, it follows that $\sigma_q$ is a square in $G$ if and only if $q$ is a square in $\mathbb{Z}/p\mathbb{Z}$. In other words,
\[
\sigma_q(\tau)=\chi_p(q)\tau.\]

Let $Q$ be a prime ideal of $R(\zeta)$ containing $q$. Lemma 3.27 implies that 
\[
\chi_p(q)\tau=\sigma_q(\tau)\equiv \tau^q\ \textrm{mod}\ Q,\ \textrm{i.e.,}\]
\[
(\chi_p(q)-\tau^{q-1})\tau\in Q.\]
If $\tau\in Q$ then $p=(-1)^{\frac{1}{2}(p-1)}\tau^2\in Q$, and this is not possible because $q$ is the only rational prime which $Q$ contains. We conclude that 
\begin{equation*}
\chi_p(q)\equiv \tau^{q-1}\ \textrm{mod}\ Q.\tag{31}\]
On the other hand, it follows from Euler's criterion that 
\begin{equation*}
\tau^{q-1}=(p^{*})^{\frac{1}{2}(q-1)}\equiv \chi_q(p^{*})\ \textrm{mod}\ q.\tag{32}\]
Because $q\in Q$, it follows from (31) and (32) that
\[
\chi_p(q)\equiv \chi_q(p^{*})\ \textrm{mod}\ Q,\]
hence
\[
\chi_p(q)= \chi_q(p^{*})\]
because $2\notin Q$. The LQR is now an immediate consequence of this equation and Theorem 2.4. $\hspace{14.6cm}\ \textrm{QED}$

The argument which we have given here (taken from [30], section 13.3) shows that quadratic reciprocity results from the fact that the cyclotomic number field $\mathbb{Q}\big(\exp(2\pi i/p)\big)$ is an abelian extension of $\mathbb{Q}$ with a cyclic Galois group. We also know from Theorem 3.2 that an irreducible polynomial $f(x) \in \mathbb{Z}[x]$ satisfies a reciprocity law if and only if the splitting field of $f(x)$ is an abelian extension of $\mathbb{Q}$. This path to higher reciprocity laws by way of the Galois theory of abelian extensions eventually became one of the principal thoroughfares to the creation of class field theory.   

\chapter{Four Interesting Applications of Quadratic Reciprocity}
Gauss called the Law of Quadratic Reciprocity the golden theorem of number theory because, when it is in hand, the study of quadratic residues and non-residues can be pursued to a significantly deeper level. We have already seen some examples of how useful the LQR can be in answering questions about the calculation of specific residues or non-residues. In this chapter, 
we will study four applications of the LQR which illustrate how it can be used to shed further light on interesting properties of residues and non-residues. 

Our first application will use quadratic reciprocity to completely solve the Basic Problem and the Fundamental Problem for Odd Primes that we introduced in Chapter 2. If $z$ is an integer, recall that the Basic Problem is to determine all primes $p$ such that $z$ is a quadratic residue of $p$ and to determine all primes $p$ such that $z$ is a quadratic non-residue of $p$. The Basic Problem must be solved in order to determine when the quadratic congruence $ax^2+bx+c\equiv 0$ mod $p$ has a solution, as we saw in Chapter 1, and it must also be solved in order to determine the splitting moduli of quadratic polynomials, as we explained in section 1 of Chapter 3. Theorems 2.4 and 2.6 solve the Basic Problem for, respectively, $z=-1$ and $z=2$ and in Chapter 2 we also showed how to reduce the solution of the Basic Problem to its solution when $z$ is an odd prime, which we call the Fundamental Problem for Odd Primes. In section 1 of this chapter, the LQR will be used to solve the Fundamental Problem for Odd Primes and this solution will then be used in section 2 to solve the Basic Problem.

The second application, which we will discuss in section 3, employs quadratic reciprocity to investigate when finite, nonempty subsets of the positive integers occur as sets of residues  of infinitely many primes.  In addition to the LQR, the key lemma which we will use to answer that question also employs Dirichlet's theorem on primes in arithmetic progression.  We take the appearance of Dirichlet's theorem here as an opportunity to discuss Dirichlet's proof of that theorem in section 4, because many of the ideas and techniques of his reasoning will be used extensively in much of the work that we will do in subsequent chapters.

If $S$ is a finite, nonempty subset of the positive integers which is a set of residues for infinitely many primes, a natural question that immediately occurs asks: how large is the set of all primes $p$ such that $S$ is a set of residues of $p$? In order to answer that question, one must find a way to accurately measure the size of an infinite set of primes. A good way to make that measurement is provided by the concept of the natural or asymptotic density of a set of primes, which we will discuss in section 5. In section 6, we apply quadratic reciprocity a third time in order to deduce a very nice way to calculate the asymptotic density of the set of all primes $p$ such that $S$ is a set of residues of $p$.

Number theory, and in particular, quadratic residues, has been applied extensively in modern cryptology. As one example of those applications, suppose that you receive an identification number from person $A$ and you want to verify that $A$ validly is in possession of the identification number, i.e., you want to be sure that $A$ really is who he claims to be, without knowing anything else about $A$. Or, for a more mathematical example, $A$ wants to convince you that he knows the prime factors of a very large number, without telling you what the prime factors are. This second example is actually used by smart cards to verify personal identification numbers. In section 7, we will describe methods, known as zero-knowledge or minimum-disclosure proofs, which use quadratic residues to securely verify the identity of someone and to convince someone that you are who you say you are. Jacobi symbols and our fourth application of the LQR are used in sections 8 and 9 to describe and verify an algorithm for fast and efficient computation of Legendre symbols that is required for the calculations in the zero-knowledge proof of section 7.

\section{Solution of the Fundamental Problem for Odd Primes}

We will now use quadratic reciprocity to solve the Fundamental Problem for Odd Primes. Let $q$ be an odd prime, and recall from Chapter 2 that the sets $X_{\pm}(q)$ are defined by
 \[
X_{\pm}(q)=\{p: \chi_p(q)= \pm 1\}.
\]
The Fundamental Problem for Odd Primes requires that the primes $p$ in these sets be found in some explicit and concrete manner.
 
 Let $r_i^+$ (respectively, $r_i^-$), $i=1,\dots, \frac{1}{2}(q-1)$ denote the residues (respectively, non-residues) of $q$ in $[1, q-1]$. Note, as we pointed out in Chapter 2, that the residues and non-residues of $q$ can be found by simply calculating the integers $1^2, 2^2,\dots, (\frac{q-1}{2})^2$ and then reducing mod $q$. The integers that result from this computation are the residues of $q$ inside $[1, q-1]$. We consider the two cases which are determined by whether $q$ is congruent to 1 or 3 mod 4.

 \emph{Case} 1: $q \equiv 1$ mod 4. 
 
 In this case, the LQR implies immediately that
\begin{eqnarray*}
X_{\pm}(q)&=&\{p: \chi_p(q)= \pm 1\}\\
&=&\{p: \chi_q(p)= \pm 1\}\\
&=&\bigcup_{i=1}^{\frac{1}{2}(q-1)}\ \{p: p \equiv r_i^{\pm} \ \textrm{mod}\ q \}.
\end{eqnarray*}
\emph{Example}: $q=17$.

\vspace{0.3cm}
We find that the residues of 17 are 1, 2, 4, 8, 9, 13, 15, and 16 and the non-residues of 17 are 3, 5, 6, 7, 10, 11, 12, and 14. Hence
\[
X_+(17)=\{p: p \equiv 1, 2, 4, 8, 9, 13, 15,\ \textrm{or}\ 16 \ \textrm{mod}\ 17\},
\]
\[
X_-(17)=\{p : p \equiv 3, 5, 6, 7, 10, 11, 12,\ \textrm{or}\ 14\ \textrm{mod}\ 17\}.
\]
(Recall that $p$ always denotes an \emph{odd} prime.)

\emph{Case} 2: $q \equiv 3$ mod 4. 

Note first (from Theorem 2.4) that
\[
X_{\pm}(-1)= \{p: p \equiv \pm1 \ \textrm{mod}\ 4\}.
\]
Hence as a consequence of the LQR,
\begin{equation*}
X_+(q)=\big(X_+(-1) \cap \{p: \chi_q(p)=1\}\big) \cup \big(X_-(-1) \cap \{p: \chi_q(p)=-1\}\big).\tag{1}
\]
Now for $i=1,\dots, \frac{1}{2}(q-1)$, let
\[
x \equiv x_i^{\pm} \mod 4q,\ 1 \leq x_i^{\pm} \leq 4q-1,
\]
be the simultaneous solutions of
\[
x \equiv \pm1 \ \textrm{mod}\ 4,\]
\[
x \equiv r_i^{\pm} \ \textrm{mod}\ q,
\]
obtained from the Chinese remainder theorem (Theorem 1.3). If we set
\[
V(q)=\{x_1^{+},\dots,x_{\frac{1}{2}(q-1)}^{+}, x_1^{-},\dots,x_{\frac{1}{2}(q-1)}^{-}\}
\]
then (1) implies that
\[
X_+(q)= \bigcup_{n \in V(q)}\ \{p:p  \equiv n \ \textrm{mod}\ 4q\}.
\]

In order to calculate $X_-(q)$, recall that $U(4q)$ denotes the set $\{n \in [1, 4q-1]: \gcd(n, 4q)=1\}$ and then observe that 
\[
V(q) \subseteq U(4q),
\]
\[
\{p: p \not= q\}=\bigcup_{n \in U(4q)}\ \{p: p \equiv n \ \textrm{mod}\ 4q\}.
\]
Hence
\begin{eqnarray*}
X_-(q)&=&\{p: p \not= q\} \setminus X_+(q)\\
&=& \bigcup_{n \in U(4q) \setminus V(q)}\ \{p: p \equiv n \ \textrm{mod}\ 4q\}.
\end{eqnarray*}

\emph{Example}: $q=7$.

\vspace{0.3cm}
The residues of 7 are 1, 2, and 4 and the non-residues are 3, 5, and 6. Because of the Chinese remainder theorem, the simultaneous solutions of the congruence pairs
 \[
 p\equiv 1 \mod 4 \ \textrm{and}\ p \equiv 1\ \textrm{mod}\ 7,
\]
 \[
 p\equiv 1 \mod 4 \ \textrm{and}\ p \equiv 2\ \textrm{mod}\ 7,
\]
\[
 p\equiv 1 \mod 4 \ \textrm{and}\ p \equiv 4\ \textrm{mod}\ 7,
\]
\[
 p\equiv -1 \mod 4 \ \textrm{and}\ p \equiv 3\ \textrm{mod}\ 7,
\]
\[
 p\equiv -1 \mod 4 \ \textrm{and}\ p \equiv 5\ \textrm{mod}\ 7,
\]
\[
 p\equiv -1 \mod 4 \ \textrm{and}\ p \equiv 6\ \textrm{mod}\ 7,
\]
are, respectively,
\[
p \equiv 1 \ \textrm{mod}\ 28,
\]
\[
p \equiv 9 \ \textrm{mod}\ 28,
\]
\[
p \equiv 25 \ \textrm{mod}\ 28,
\]
\[
p \equiv 3 \ \textrm{mod}\ 28,
\]
\[
p \equiv 19 \ \textrm{mod}\ 28,
\]
\[
p \equiv 27 \ \textrm{mod}\ 28.
\]
Hence
\[
X_+(7)=\{p: p \equiv 1, 3, 9,19,25,\ \textrm{or}\ 27 \ \textrm{mod}\ 28\}.
\]
We have that
\[
U(28)=\{1,3,5,9,11,13,15,17,19,23,25,27\},
\]
\[
V(7)=\{1,3,9,19,25,27\},
\]
hence,
\[
U(28) \setminus V(7)=\{5,11,13,15,17,23\},
\]
and so
\[
X_-(7)=\{p: p \equiv 5,11,13,15,17,\ \textrm{or}\ 23 \ \textrm{mod}\ 28\}.
\]

\section{Solution of the Basic Problem}

If $d$ is a fixed but arbitrary integer, we recall formulae (2) and (5) for $X_+(d)$ from Chapter 2. Suppose first that $d>0$.
Let
\[
\mathcal{E}=\{E \subseteq \pi_{\textrm{odd}}(d): |E|\ \textrm{is even} \},
\]
where $ \pi_{\textrm{odd}}(d)$ denotes the set of all prime factors of $d$ of odd multiplicity. If $E \in \mathcal{E}$, let $R_{E}$ denote the set of all $p$ such that 
\[
\chi_p(q)=\left\{\begin{array}{ll}-1,\ \textrm{if $q \in E$,}\\
1,\ \textrm{if $q \in \pi_{\textrm{odd}}(d) \setminus E$.}\\\end{array}\right.
\]
Then formula (2) of Chapter 2 is
\[
X_+(d)= \Big( \bigcup_{E \in \mathcal{E}}\ R_E \Big) \setminus \pi_{\textrm{even}}(d), 
\]
where $\pi_{\textrm{even}}(d)$ denotes the set of all prime factors of $d$ of even multiplicity, and this union is pairwise disjoint. Moreover
\[
R_E=\Big( \bigcap_{q \in E}\ X_-(q) \Big) \cap \Big( \bigcap_{q \in \pi_{\textrm{odd}}(d) \setminus E}\ X_+(q) \Big). 
\]

Suppose next that $d<0$, and let
\[
\mathcal{E}_{-1}=\{E \subseteq \{-1\} \cup \pi_{\textrm{odd}}(d): |E|\ \textrm{is even} \}.
\]
Then formula (5) of Chapter 2 is
\[
X_+(d)= \Big( \bigcup_{E \in \mathcal{E}_{-1}}\ R_E \Big) \setminus \pi_{\textrm{even}}(d), 
\]
where
\[
R_E=\Big( \bigcap_{q \in E}\ X_-(q) \Big) \cap \Big( \bigcap_{q \in(\{-1\} \cup \pi_{\textrm{odd}}(d)) \setminus E}\ X_+(q) \Big), E \in \mathcal{E}_{-1}. 
\]

We can now use formula (2) or (5) of Chapter 2 in concert with the solution that we have of the Fundamental Problem for Odd Primes to calculate $X_{\pm}(d)$, thereby solving the Basic Problem. The formulae that we have derived for the calculation of $X_{\pm}(q)$ where $q$ is either $-1$ or a prime show that each of these sets is equal to a union of certain equivalence classes mod 4, 8, an odd prime, or 4 times an odd prime. It follows that when we employ formula (2) or (5) of Chapter 2 to calculate $X_+(d)$, each of the sets $R_E$ occurring in those formulae can hence be calculated by the method of successive substitution, a generalization of the Chinese remainder theorem that can be used to solve simultaneous congruences when the moduli of the congruences are no longer pairwise relatively prime. 

The method of successive substitution works as follows. We have a series of congruences of the form
\begin{equation*}
x \equiv  a_i \ \textrm{mod}\ m_i,\ i=1,\dots,k, \tag{2} 
\end{equation*}
where $(m_1,\dots,m_k)$ is a given $k$-tuple of moduli and $(a_1,\dots,a_k)$ is a given $k$-tuple of integers, which we wish to solve simultaneously. Denoting by lcm$(a, b)$ the least common multiple of the integers $a$ and $b$, one starts with
\begin{prp}
The congruences
\[
x \equiv a_1 \ \textnormal{mod}\ m_1,\  x \equiv a_2 \ \textnormal{mod}\ m_2
\]
have a simultaneous solution if and only if $\gcd(m_1, m_2)$ divides $a_1-a_2$. The solution is unique modulo $\textnormal{lcm}( m_1, m_2)$ and is given by
\[
x \equiv a_1+x_0m_1 \ \textnormal{mod} \ \textnormal{lcm}(m_1, m_2),
\]
where $x_0$ is a solution of
\[
m_1x_0 \equiv a_2-a_1 \ \textnormal{mod}\ m_2.
\]
\end{prp}
\noindent The congruences $(2)$ are then solved by first using Proposition 4.1 to solve the first two congruences in $(2)$, then, if necessary, pairing the solution so obtained with the third congruence in $(2)$ and applying Proposition 4.1 to solve that congruence pair, and continuing in this manner, successively applying Proposition 4.1 to the pair of congruences consisting of the solution obtained from step $i-1$ and the $i$-th congruence in $(2)$. This procedure confirms that $(2)$ has a simultaneous solution if and only if $\gcd(m_i, m_j)$ divides $a_i-a_j$ for all $i$ and $j$, and that the solution is unique modulo the least common multiple of $m_1,\dots,m_k$. Proposition 4.1 is not difficult to verify, and so we will leave that to the interested reader.

Consequently, once the residues and non-residues of each integer in $\pi_{\textrm{odd}}(d)$ are determined, $X_+(d)$ can be calculated by repeated applications of the method of successive substitutions. In particular, one finds a positive integer $m(d)$ and a subset $V(d)$ of $U\big(m(d)\big)$ such that 
\[
X_+(d)= \Big(\bigcup_{n \in V(d)}\ \{p: p \equiv n \ \textrm{mod}\ m(d)\} \Big) \setminus \pi_{\textrm{even}}(d).
\]   
The modulus $m(d)$ is determined like so: if $d>0$ and $\pi_{\textrm{odd}}(d)$ contains neither 2 nor a prime $\equiv 3$ mod 4, then $m(d)$ is the product of all the elements of $\pi_{\textrm{odd}}(d)$; otherwise, $m(d)$ is 4 times this product. 

The formula for $X_-(d)$ can now be obtained from the one for $X_+(d)$ by first observing that as a consequence of the above determination of $m(d)$, 
\[
\pi\big(m(d)\big) \cup \{2\}= \pi_{\textrm{odd}}(d) \cup \{2\},
\]
and so
\[
\pi(d) \cup \{2\}=\pi\big(m(d)\big) \cup \{2\} \cup \pi_{\textrm{even}}(d).\]
Upon  recalling that $P$ denotes the set of all primes, it follows that 
\begin{eqnarray*}
X_-(d)&=&P \setminus \big(X_+(d) \cup \{2\} \cup \pi(d) \big)\\
&=&P \setminus \big(\pi\big(m(d)\big) \cup \{2\} \cup X_+(d) \cup \pi_{\textrm{even}}(d) \big)\\
&=&\big[P\setminus \big(\pi\big(m(d)\big) \cup \{2\} \big) \big] \setminus \big[X_+(d) \cup \pi_{\textrm{even}}(d) \big].\\
\end{eqnarray*}
Because
\[
P \setminus \big(\pi\big(m(d)\big) \cup \{2\} \big)=\bigcup_{n \in U(m(d))}\ \{p: p \equiv n \ \textrm{mod}\ m(d)\}, \]
\[
X_+(d) \cup \pi_{\textrm{even}}(d)= \Big(\bigcup_{n \in V(d)}\ \{p: p \equiv n \ \textrm{mod}\ m(d)\} \Big) \cup \pi_{\textrm{even}}(d),\]
it hence follows that
\[
X_-(d)=\Big( \bigcup_{n \in U(m(d)) \setminus V(d)}\ \{p: p \equiv n \ \textrm{mod}\ m(d)\} \Big) \setminus \pi_{\textrm{even}}(d).
\]

The set $V(d)$ that appears in the formulae which calculate $X_{\pm}(d)$ is obtained from applications of the method of successive substitution to the calculation of each of the sets $R_E$ which appears in equation (2) or (5) of Chapter 2. A natural question which arises asks: are \emph{all} of the integers in $V(d)$ and $U\big(m(d)\big) \setminus V(d)$ which arise from these calculations required for the determination of $X_{\pm}(d)$? The answer is yes, if for each pair of relatively prime positive integers $m$ and $n$, the set $\{z \in \mathbb{Z}: z \equiv n \mod m\}$ contains primes. Remarkably enough, $\{z \in \mathbb{Z}: z \equiv n \mod m\}$ in fact always contains \emph{infinitely many} primes. This is a famous theorem of Dirichlet [10], and the connection of that theorem to the calculation of $X_{\pm}(d)$ was Dirichlet's primary motivation for proving it. Much more is to come (in section 4 below) about Dirichlet's theorem and its use in the study of residues and non-residues.

We next illustrate the procedure which we have described for the solution of the Basic Problem by calculating $X_{\pm}(126)$. From the calculations using this example that we preformed in section 2 of Chapter 2, it follows that
\[
X_+(126)=\Big( \big(X_+(2) \cap X_+(7)\big) \cup \big(X_-(2) \cap X_-(7)\big) \Big) \setminus \{3\}.
\]
hence we must calculate $X_+(2) \cap X_+(7)$ and $X_-(2) \cap X_-(7)$. 

\emph{Calculation of} $X_+(2) \cap X_+(7)$.

Theorem 2.6 implies that
\[
X_+(2)=\{p: p \equiv 1\ \textrm{or }\ 7 \ \textrm{mod}\ 8\},
\]
and we have from the calculation of $X_+(7)$ above that
\[
X_+(7)=\{p: p \equiv 1,3,9,19,25,\ \textrm{or} \ 27 \ \textrm{mod}\ 28\}.
\]
 
In order to calculate $X_+(2) \cap X_+(7)$, we need to solve at most 12 (but in fact exactly six) pairs of simultaneous congruences. We do this by applying Proposition 4.1.
We have that $\gcd(8, 28)=4$, lcm$(8, 28)=56$, and so Proposition 4.1 implies that $X_+(2) \cap X_+(7)$ consists of the union of all odd prime simultaneous solutions of the congruence pairs
\[
x \equiv 1 \mod 8,\ x \equiv 1 \ \textrm{mod}\ 28,
\]
\[
x \equiv 1 \mod 8,\ x \equiv 9 \ \textrm{mod}\ 28,
\]
\[
x \equiv 1 \mod 8,\ x \equiv 25 \ \textrm{mod}\ 28,
\]
\[
x \equiv 7 \mod 8,\ x \equiv 3 \ \textrm{mod}\ 28,
\]
\[
x \equiv 7 \mod 8,\ x \equiv 19\ \textrm{mod}\ 28,
\]
\[
x \equiv 7 \mod 8,\ x \equiv 27 \ \textrm{mod}\ 28,
\]
whose odd prime solutions are, respectively,
\[
p \equiv 1 \ \textrm{mod}\ 56,
\]
\[
p \equiv 9 \ \textrm{mod}\ 56,
\]
\[
p \equiv 25 \ \textrm{mod}\ 56,
\]
\[
p \equiv 31 \ \textrm{mod}\ 56,
\]
\[
p \equiv 47 \ \textrm{mod}\ 56,
\]
\[
p \equiv 55 \ \textrm{mod}\ 56.
\]
 
\emph{Calculation of} $X_-(2) \cap X_-(7)$. 

From Theorem 2.6 and the calculation of $X_-(7)$ above, it follows that
\[
X_-(2)= \{p: p \equiv 3\ \textrm{or}\ 5 \ \textrm{mod}\ 8\},
\]
\[
X_-(7)=\{p: p \equiv 5,11,13,15,17,\ \textrm{or}\ 23 \ \textrm{mod}\ 28\}.
\]
Hence,  again according to Proposition 4.1, $X_-(2) \cap X_-(7)$ consists of the union of all odd prime simultaneous solutions of the congruence pairs
\[
x \equiv 3 \mod 8,\ x \equiv 11 \ \textrm{mod}\ 28,
\]
\[
x \equiv 3 \mod 8,\ x \equiv 15 \ \textrm{mod}\ 28,
\]
\[
x \equiv 3 \mod 8,\ x \equiv 23 \ \textrm{mod}\ 28,
\]
\[
x \equiv 5 \mod 8,\ x \equiv 5 \ \textrm{mod}\ 28,
\]
\[
x \equiv 5 \mod 8,\ x \equiv 13\ \textrm{mod}\ 28,
\]
\[
x \equiv 5 \mod 8,\ x \equiv 17 \ \textrm{mod}\ 28,
\]
whose odd prime solutions are, respectively,
\[
p \equiv 11 \ \textrm{mod}\ 56,
\]
\[
p \equiv 43 \ \textrm{mod}\ 56,
\]
\[
p \equiv 51 \ \textrm{mod}\ 56,
\]
\[
p \equiv 5 \ \textrm{mod}\ 56,
\]
\[
p \equiv 13 \ \textrm{mod}\ 56,
\]
\[
p \equiv 45 \ \textrm{mod}\ 56.
\]

From this calculation of  $X_+(2) \cap X_+(7)$ and  $X_-(2) \cap X_-(7)$, it hence follows that
\[
X_+(126)= \{p   : p \equiv 1,5,9,11,13,25,31,43,45,47,51,\ \textrm{or}\ 55 \ \textrm{mod}\ 56\}.
\]
 
 In order to calculate $X_-(126)$, we simply delete from $U(56)$ the  minimal positive ordinary residues mod 56 that determine $X_+(126)$: the integers resulting from that are 3, 15, 17, 19, 23, 27, 
 29, 33, 37, 39, 41, and 53. Hence 
 \[
 X_-(126)= \{p \not= 3: p \equiv 3,15,17,19,23,27,29,33,37,39,41,\ \textrm{or}\ 53 \ \textrm{mod}\ 56\}.
 \]

\section{Sets of Integers which are Quadratic Residues of Infinitely Many Primes}

In this section we will use the LQR to investigate when a finite non-empty subset of positive integers is the set of residues for infinitely many primes. We start  by looking at singleton sets. Obviously, if $a \in \mathbb{Z}$ is a square then $a$ is a residue of all primes. Is the converse true, i.e., if a positive integer is a residue of all primes, must it be a square?
The answer is yes; in fact a slightly stronger statement is valid:
\begin{thm}
A positive integer is a residue of all but finitely many primes if and only if it is a square.
\end{thm}

This theorem implies that if $S$ is a nonempty finite subset of $[1, \infty)$ then $S$ is a set of residues for all but finitely many primes if and only if every element of $S$ is a square. What if we weaken the requirement that $S$ be a set of residues of all but finitely many primes to the requirement that $S$ be a set of residues for only \emph{infinitely many} primes? Then the somewhat surprising answer is asserted by
\begin{thm}
If S is $\textnormal{any}$ nonempty finite subset of $[1, \infty)$ then S is a set of residues of infinitely many primes.
\end{thm}

 Theorems 4.2 and 4.3 are simple consequences of 
 \begin{lem}
 $($Basic Lemma$)$ If $\Pi=\{p_1,\dots,p_k\}$ is a nonempty finite set of primes and if $\varepsilon: \Pi \rightarrow \{-1, 1\}$ is a fixed function then there exits infinitely many primes $p$ such that
 \[
 \chi_p(p_i)= \varepsilon(p_i),\ i \in [1, k].
 \]
 \end{lem}
 N.B. This lemma asserts that if all of the integers in the set $S$ of Theorem 4.3 are \emph{prime}, then for any pattern of $+1$'s or $-1$'s attached to the elements of $S$, the Legendre symbol $\chi_p$ reproduces that pattern on $S$ for infinitely many primes $p$. Thus the conclusion of Theorem 4.3 can be strengthened considerably when $S$ is a set of primes.
 
  Assume Lemma 4.4 for now. We will use it to first prove Theorems 4.2 and 4.3 and then we will use quadratic reciprocity (and Dirichlet's theorem on primes in arithmetic progression) to prove Lemma 4.4.
  
  \emph{Proof of Theorem} 4.2. 
  
  Suppose that $n \in [1, \infty)$ is not a square. Then $\pi_{\textrm{odd}}(n) \not= \emptyset$ and
  \begin{equation*}
  \chi_p(n)= \prod_{q \in \pi_{\textrm{odd}}(n)}\ \chi_p(q),\ \textrm{for all}\ p \notin \pi(n). \tag{3}
  \]
  Now take any fixed $q_0 \in \pi_{\textrm{odd}}(n)$ and define $\varepsilon: \pi_{\textrm{odd}}(n) \rightarrow \{-1, 1\}$ by
  \[
  \varepsilon(q)=\left\{\begin{array}{rl}-1,& \textrm{if $q=q_0$,}\\
1,& \textrm{if $q \not= q_0$.}\\\end{array}\right. 
\]
Lemma 4.4 implies that there exists infinitely many primes $p$ such that
\[
\chi_p(q)= \varepsilon(q),\ \textrm{for all}\ q \in \pi_{\textrm{odd}}(n),
\]
and so the product in (3), and hence $\chi_p(n)$, is $-1$ for all such $p \notin \pi(n)$.$\hspace{3.4cm} \textrm{QED}$

\emph{Proof of Theorem} 4.3. 

Let $S$ be a fixed nonempty subset of positive integers and let
\[
X=\bigcup_{z \in S}\ \pi_{\textrm{odd}}(z).
\]
We may assume that $X \not= \emptyset$; otherwise all elements of $S$ are squares and Theorem 4.3 is trivially true in that case. Then 
Lemma 4.4 implies that there exists infinitely many primes $p$ such that
\[
\chi_p(q)=1,\ \textrm{for all}\ q \in X,
\]
hence for all such $p$ which are not factors of an element of $S$,
\[
 \chi_p(z)= \prod_{q \in \pi_{\textrm{odd}}(z)}\ \chi_p(q)=1,\ \textrm{for all}\ z \in S.\   
 \] 
$\hspace{15.6cm} \textrm{QED}$ 

 \emph{Proof of Lemma} 4.4. 
 
 It follows from our solution of the Fundamental Problem for all primes (Theorem 2.6 and the calculation of $X_{\pm}(q)$, $q$ an odd prime, in section 1 of this chapter) that Lemma 4.4 is valid when $\Pi$ is a singleton, so assume that $k \geq 2$. We will make use of arithmetic progressions in this argument, and so if $a, b \in [1, \infty)$, let
 \[
 AP(a, b)=\{a+nb: n \in [0, \infty)\}
 \]
 denote the arithmetic progression with initial term $a$ and common difference $b$. We will find the primes that will verify the conclusion of Lemma 4.4 by looking inside certain arithmetic progressions, hence we will need the following theorem, one of the basic results in the theory of prime numbers:
 \begin{thm}
 $($Dirichlet's theorem on primes in arithmetic progression$)$. If $\{a, b\} \subseteq [1, \infty)$ and $\gcd(a, b)=1$ then $AP(a, b)$ contains infinitely many primes.
 \end{thm}

The key ideas in Dirichlet's proof of Theorem 4.5 will be discussed in due course. For now, assume that the elements of the set $\Pi$ in the hypothesis of Lemma 4.4 are ordered as $p_1< \dots <p_k$ and fix $\varepsilon: \Pi \rightarrow \{-1, 1\}$. We need to verify the conclusion of Lemma 4.4 for this $\varepsilon$. Suppose first that $p_1=2$ and $\varepsilon(2)=1$. If $i \in [2, k]$ and $\varepsilon(p_i)=1$, let $k_i=1$, and if $\varepsilon(p_i)=-1$, let $k_i$ be an odd non-residue of $p_i$ such that $\gcd(p_i, k_i)=1$ (if  $\varepsilon(p_i)=-1$ then such a $k_i$ can always be chosen: simply pick any non-residue $x$ of $p_i$ in $[1, p_i-1]$; if $x$ is odd, set $k_i=x$, and if $x$ is even, set $k_i=x+p_i$).

Now, suppose that $i \in [2, k]$ , $p \equiv 1 \mod 8$, and $p \in AP(k_i, 2p_i)$, say $p=k_i+2p_in$, for some $n \in [1, \infty)$. Then LQR implies that
\[
\chi_p(p_i)=\chi_{p_i}(p)=\chi_{p_i}(k_i+2p_in)=\chi_{p_i}(k_i).
\]
It follows from Theorem 2.6 and the choice of $k_i$ that 
\[
\chi_p(2)=1\ \textrm{and}\ \chi_p(p_i)=\varepsilon(p_i).
\]
Hence
\begin{equation*}
\textrm{if}\ p\equiv 1\ \textrm{mod}\  8\ \textrm{and}\ p \in \bigcap_{i=2}^k\ AP(k_i, 2p_i),\ \textrm{then}\ \chi_p(p_i)=\varepsilon(p_i),\ \textrm{for all}\ i \in [1, k]. \tag{4}
\]

We prove next that there are infinitely many primes $\equiv 1\ \textrm{mod}\ 8$ inside $\bigcap_{i=2}^k\ AP(k_i, 2p_i)$. To see this, we first use the fact that each $k_i$ is odd and an inductive construction obtained from solving an appropriate sequence of linear Diophantine equations (Proposition 1.4) to obtain an integer $m$ such that
\begin{equation*}
AP(k_2+2m, 8p_2 \cdots p_k) \subseteq AP(1, 8) \cap \Big( \bigcap_{i=2}^k\ AP(k_i, 2p_i) \Big). \tag{5}
\]
\noindent We then claim that $\gcd(k_2+2m, 8p_2 \cdots p_k)=1$. If this is true then by virtue of Theorem 4.5, we have that $AP(k_2+2m, 8p_2 \cdots p_k)$ contains infinitely many primes $p$, hence for any such $p$, it follows from (4) and (5) that
\begin{equation*}
\chi_p(p_i)=\varepsilon(p_i),\ i \in [1, k], \tag{6}
\]
the conclusion of Lemma 4.4. To verify the claim, assume by way of contradiction that $q$ is a common prime factor of $k_2+2m$ and $8p_2 \cdots p_k$. Then $q \not= 2$ because $k_2$ is odd, hence there is a $j \in [2, k]$ such that $q=p_j$. But (5) implies that there exists $n \in [0, \infty)$ such that
\[
k_2+2m+8p_2 \cdots p_k=k_j+2np_j,
\]
and so $p_j$ divides $k_j$, contrary to the choice of $k_j$.

If $p_1=2$ and $\varepsilon(2)=-1$, a similar argument shows that $\bigcap_{i=2}^k\ AP(k_i, 2p_i)$ contains infinitely many primes $p \equiv 5 \mod 8$, hence (6) is true for all such $p$. If $p_1 \not=2$, simply adjoin 2 to $\Pi$ and repeat this argument.$\hspace{9.8cm} \textrm{QED}$

\section{Intermezzo: Dirichlet's Theorem on Primes in Arithmetic Progression}

In addition to the LQR, Theorem 4.5 also played a key role in the proof of the basic Lemma 4.4, and thus also in the proofs of Theorems 4.2 and 4.3. Because they will play such an important role in our story, we will now discuss the key ingredients of Dirichlet's proof of Theorem 4.5. Dirichlet [10] proved this in 1837 , and it would be hard to overemphasize the importance of this theorem and the methods Dirichlet developed to prove it. As we shall see, he used \emph{analysis}, specifically the theory of infinite series and infinite products of complex-valued functions of a real variable, and in subsequent work [11] also the theory of Fourier series, to discover properties of the primes (for the reader who may benefit from it, we briefly discuss analytic functions, Fourier series, and some of their basic properties in Chapter 7). His use of continuous methods to prove deep results about  discrete sets like the prime numbers was not only a revolutionary insight, but also caused a sensation in the nineteenth century mathematical community.  Dirichlet's results founded the subject of \emph{analytic number theory}, which has become one of the most important areas and a major industry in number theory today. Later (in Chapters 5 and 7) we will also see how Dirichlet used analytic methods to study important properties of residues and non-residues.

Dirichlet is a towering figure in the history of number theory not only because of the many results and methods of fundamental importance which he discovered and developed in that subject but also because of his role as an expositor of that work and the work of Gauss. We have already given an indication of how the work of Gauss, especially the \emph{Disquisitiones Arithmeticae}, brought about a revolutionary transformation in number theory. However, the influence of Gauss' work was rather slow to be realized, due primarily to the difficulty that many of his mathematical contemporaries had in understanding exactly how Gauss had done what he had done in the \emph{Disquisitiones}. Dirichlet is said to have been the first person to completely master the \emph{Disquisitiones}, and legend has it that he was never without a copy of it within easy reach. Many of the results and techniques that Gauss developed in the  \emph{Disquisitiones} were first explained in a more accessible way in Dirichlet's great text [12], the \emph{Vorlesungen $\ddot{\textrm{u}}$ber Zahlentheorie}; John Stillwell, the translator of the \emph{Vorlesungen} into English, called it one of the most important mathematics books of the nineteenth century: the link between Gauss and the number theory of today. If a present-day reader of the \emph{Disquisitiones} finds much of it easier to understand  than a reader in the early days of the nineteenth century did, it is because that modern reader learned number theory the way that Dirichlet first taught  it. 

Now, back to primes in arithmetic progression. In 1737, Euler proved that the series $\sum_ {q \in P} \frac{1}{q}$ diverges and hence deduced Euclid's theorem that there are infinitely many primes. Taking his cue from this result, Dirichlet sought to prove that 
\[
\sum_{p \equiv a\ \textnormal{mod}\ b} \frac{1}{p}
\]
diverges, where $a$ and $b$ are given positive relatively prime integers, thereby showing that the arithmetic progression with constant term $a$ and difference $b$ contains infinitely many primes. To do this, he studied the behavior as $s \rightarrow 1^+$ of the function of $s$ defined by 
\[
\sum_{p \equiv a\ \textnormal{mod}\ b} \frac{1}{p^s}.
\]
This function is difficult to get a handle on; it would be easier if we could replace it by a sum indexed over all of the primes, so consider
\[
\sum_p \delta(p)p^{-s},\ \textrm{where}\ \delta(p)=\left\{\begin{array}{ll}1,\ \textrm{if $p \equiv a\ \textnormal{mod}\ b$,}\\
0,\ \textrm{otherwise.}\\\end{array}\right. 
\]
Dirichlet's profound insight was to replace $\delta(p)$ by certain functions which capture the behavior of $\delta(p)$ closely enough, but which are more amenable to analysis relative to primes in the  ordinary residue classes mod $b$. We now define these functions.

Begin by recalling that if $A$ is a commutative ring with identity 1 then a \emph{unit u of A}  is an element of $A$ that has a multiplicative inverse in $A$, i.e., there exists $v \in A$ such that $uv=1$.
The set of all units of $A$ forms a group under the multiplication of $A$, called the \emph{group of units of A}. Consider now the ring $\mathbb{Z}/b\mathbb{Z}$ of ordinary residue classes of $\mathbb{Z}$ mod $b$. Proposition 1.2 implies that the group of units of $\mathbb{Z}/b\mathbb{Z}$ consists of all ordinary residue classes that are determined by the integers that are relatively prime to $b$. If we hence identify $\mathbb{Z}/b\mathbb{Z}$ in the usual way with the set of ordinary non-negative minimal residues $[0, b-1]$ on which is defined the addition and multiplication induced by addition and multiplication of ordinary residue classes, it follows that 
\[
U(b)=\{n \in [1, b-1]: \gcd(n, b)=1\}
\]
is the group of units of $\mathbb{Z}/b\mathbb{Z}$, and we set
\[
\varphi(b)=|U(b)|;
\]
$\varphi$ is called \emph{Euler's totient function}. 

Let $T$ denote the circle group of all complex numbers of modulus 1, with the group operation defined by ordinary multiplication of complex numbers. A homomorphism of $U(b)$ into $T$ is called a \emph{Dirichlet character modulo b}. We denote by $\chi_0$ the \emph{principal character modulo b}, i.e., the character which sends every element of $U(b)$ to $1 \in T$. If $\chi$ is a Dirichlet character modulo $b$, we extend it to all integers $z$ by setting $\chi(z)= \chi(n)$ if there exists $n \in U(b)$ such that $z \equiv n \ \textrm{mod}\ b$, and setting $\chi(z)=0$, otherwise. It is then easy to verify
\begin{prp}
A  Dirichlet character $\chi$ modulo b is 

$(i)$ of period b, i.e., $\chi(n)=0$ if and only if $\gcd(n, b)>1$ and $\chi(m)=\chi(n)$ whenever $m \equiv n\ \textnormal{mod}\ b$, and is 

$(ii)$ completely multiplicative, i.e., $\chi(mn)=\chi(m) \chi(n)$ for all $m, n \in \mathbb{Z}$.
\end{prp}

We say that a Dirichlet character is \emph{real} if it is real-valued, i.e., its range is  either the set $ \{0, 1\}$ or $[-1, 1]$. In particular the Legendre symbol $\chi_p$ is a real Dirichlet character mod $p$.

For each modulus $b$, the structure theory of finite abelian groups can be used to explicitly construct all Dirichlet characters mod $b$; we will not do this, and instead refer the interested reader to Hecke [27], section 10 or Davenport [6], pp. 27-30. In particular there are exactly $\varphi(b)$ Dirichlet characters mod $b$.

The connection between Dirichlet characters and primes in arithmetic progression can now be made. If $\gcd(a, b)=1$ then Dirichlet showed that
\[
\frac{1}{\varphi(b)} \sum_{\chi}\ \overline{\chi(a)} \chi(p)=\left\{\begin{array}{ll}1,\ \textrm{if $p \equiv a\ \textnormal{mod}\ b$,}\\
0,\ \textrm{otherwise,}\\\end{array}\right. 
\]
where the sum is taken over all Dirichlet characters $\chi$ mod $b$.  These are the so-called \emph{orthogonality relations } for the Dirichlet characters. This equation says that the  characteristic function $\delta(p)$ of the primes in an ordinary equivalence class mod $b$ can be written as a linear combination of Dirichlet characters. Hence
\begin{eqnarray*}
\sum_{p \equiv a\ \textnormal{mod}\ b} \frac{1}{p^s}&=&\sum_p \delta(p)p^{-s}\\
&=&\sum_p\ \Big(\frac{1}{\varphi(b)} \sum_{\chi}\ \overline{\chi(a)} \chi(p) \Big)p^{-s}\\
&=&\frac{1}{\varphi(b)} \sum_p\ p^{-s}+ \frac{1}{\varphi(b)} \sum_{\chi \not= \chi_0}\ \overline{\chi(a)} \Big(\sum_p\ \chi(p)p^{-s} \Big).
\end{eqnarray*}
\noindent After observing that
\[
\lim_{s \rightarrow 1^+}\sum_p p^{-s}= +\infty,\]
Dirichlet deduced immediately from the above equations the following lemma:
\begin{lem}
$\lim_{s \rightarrow 1^+} \sum_{p \equiv a\ \textnormal{mod}\ b} p^{-s}= +\infty$  if for each non-principal Dirichlet character $\chi$ $\textnormal{ mod}$ $b$, $\sum_p \chi(p)p^{-s}$ is bounded as $s \rightarrow 1^+$.
\end{lem}
\noindent Hence Theorem 4.5 will follow if one can prove that
\begin{equation*}
\textrm{for all non-principal Dirichlet characters $\chi\ \textnormal{mod}\ b$}, \sum_p \chi(p)p^{-s}\ \textrm{is bounded as $s \rightarrow 1^+$.} \tag{7}
\]

Let $\chi$ be a given Dirichlet character. In order to verify (7), Dirichlet introduced his next deep insight  into the problem by considering the function
\[
L(s, \chi)=\sum_{n=1}^{\infty} \frac{\chi(n)}{n^s},\ s \in \ \textbf{C},
\]
which has come to be known as the \emph{Dirichlet L-function of $\chi$}. We will prove in Chapter 7 that $L(s, \chi)$ is analytic in the half-plane Re $s>1$, satisfies the infinite-product formula
\[
L(s, \chi)=\prod_{q \in P} \frac{1}{1-\chi(q)q^{-s}},\ \textrm{Re}\ s>1,
\]
the \emph{Euler-Dirichlet product formula}, and is analytic in Re $s>0$ whenever $\chi$ is non-principal. One can take the complex logarithm of both sides of the Euler-Dirichlet product formula to deduce that
\[
\log L(s, \chi)=\sum_{n=2}^{\infty} \frac{\chi(n) \Lambda(n)}{\log n}n^{-s}, \textrm{Re}\ s>1,
\]
where
\[
\Lambda(n)=\left\{\begin{array}{rl} \log q,& \textrm{if $n$  is a power of $q, q \in P$,}\\
0,& \textrm{otherwise.}\\\end{array}\right. 
\]
Using algebraic properties of the character $\chi$ and the function $\Lambda$, Dirichlet proved that (7) is true if
\begin{equation*}
\log L(s, \chi)\ \textrm{is bounded as  $s \rightarrow 1^+$ whenever $\chi$ is non-principal.} \tag{8}
\]
 We should point out that Dirichlet did not use functions of a complex variable in his work, but instead worked only with real values of the variable $s$ (Cauchy's theory of analytic functions of a complex variable, although fully developed by 1825, did not become well-known or commonly employed until the 1840's) . Because $L(s, \chi)$ is continuous on $\textrm{Re}\ s>0$, it follows that
 \[
 \lim_{s \rightarrow 1^+} \log L(s, \chi)= \log L(1, \chi),
 \]
 hence (8) will hold if
 \[
 L(1, \chi) \not= 0 \ \textrm{whenever $\chi$ is non-principal.}
 \]
 We have at last come to the heart of the matter, namely
 \begin{lem}
 If $\chi$ is a non-principal Dirichlet character then $L(1, \chi) \not= 0$.
 \end{lem}
If $\chi$ is not real, Lemma 4.8 is fairly easy to prove, but when $\chi$ is real, this task is much more difficult to do. Dirichlet deduced Lemma 4.8 for real characters by using results from the classical theory of quadratic forms; he established a remarkable formula which calculates $L(1, \chi)$ as the product of a certain parameter and the number of  equivalence classes of quadratic forms (section 12, Chapter 3); because this parameter and the number of equivalence classes are clearly positive, $L(1, \chi)$ must be nonzero. At the conclusion of Chapter 7, we will give an elegant proof of Lemma 4.8 for real characters due to de la Vall$\acute {\textrm{e}}$e Poussin [45], and then in Chapter 8 we will prove Dirichlet's class-number formula for the value of $L(1, \chi)$.

Finally, we note that if $\chi_0$ is the principal character mod $b$ then it is a consequence of the Euler-Dirichlet product formula that
\[
L(s, \chi_0)=\zeta(s)\prod_{q|b}\big(1-q^{-s}\big),
\]
where
\[
\zeta(s)=\sum_{n=1}^{\infty} \frac{1}{n^s}\]
is the Riemann zeta function. 

At this first appearance in our story of $\zeta(s)$, probably the single most important function in analytic number theory, we cannot resist briefly discussing the

\vspace{0.3cm}

\emph{Riemann Hypothesis}: all zeros of $\zeta(s)$ in the strip $0<\textrm{Re}\ s<1$ have real part $\frac{1}{2}$.

\vspace{0.3cm}
\emph{Generalized Riemann Hypothesis $($GRH$)$}: if $\chi$ is a Dirichlet character then all zeros of $L(s, \chi)$ in the strip $0<\textrm{Re}\ s\leq1$ have real part $\frac{1}{2}$.
\vspace{0.3cm}

\noindent Riemann [47] first stated the Riemann Hypothesis (in an equivalent form) in a paper that he published in 1859, in which he derived an explicit formula for the number of primes not exceeding a given real number. By general agreement, verification of the Riemann Hypothesis is the most important unsolved problem in mathematics. One of the most immediate consequences of the truth of the Riemann Hypothesis, and arguably the most significant, is the essentially optimal error estimate for the asymptotic approximation of the cardinality of the set $\{q \in P: q \leq x\}$ given in the Prime Number Theorem (see the statement of this theorem in the next section). This estimate asserts that there is an absolute, positive constant $C$ such that for all $x$ sufficiently large,
 \[
\left|\frac{\big|\{q \in P: q \leq x\} \big|}{\displaystyle \int_2^x \frac{1}{\log t}\ dt}-1\right|\leq  \frac{C}{\sqrt x}\ .\]
The integral $\int_2^x \frac{1}{\log t}\ dt$ appearing in this inequality, the \emph{logarithmic integral of x}, is generally a better asymptotic approximation to the cardinality of $\{q \in P: q \leq x\}$ than the quotient $x/\log x$. Hilbert emphasized the importance of the Riemann Hypothesis in Problem 8 on his famous list of 23 open problems that he presented in 1900 in his address to the second International Congress of Mathematicians. In 2000, the Clay Mathematics Institute (CMI) published a series of seven open problems in mathematics that are considered to be of exceptional importance and have long resisted solution. In order to encourage work on these problems, which have come to be known as the Clay Millennium Prize Problems, for each problem CMI will award to the first person(s) to solve it \$1,000,000 (US). The proof of the Riemann Hypothesis is the second Millennium Prize Problem (as currently listed on the CMI web site).

\section{The Asymptotic Density of Primes}

Theorem 4.3 gives rise to the following natural and interesting question: if $S$ is a nonempty, finite subset of $[1, \infty)$, how large is the necessarily infinite set of primes 
\[
\{p: \chi_p \equiv 1\ \textrm{on}\ S \}\ \textrm{?}
\]
(The meaning of the symbol $\equiv$ used here is as an identity of functions, \emph{not} as a modular congruence; in subsequent uses of this symbol, its meaning will be clear from the context.) To formulate this question precisely, we need a good way to measure the size of an infinite set of primes. This is provided by the concept of the asymptotic density of a set of primes, which we will discuss in this section.

If $\Pi$ is a set of primes and $P$ denotes the set of all primes then the \emph{asymptotic density of $\Pi$ in P} is
\[
\lim_{x\rightarrow +\infty } \frac{\big| \{p \in \Pi:\ p \leq x \}\big|}{\big| \{p \in P:\ p \leq x \}\big|} ,
\]
provided that this limit exists. Roughly speaking, the density of $\Pi$ is the ``proportion" of the set $P$ that is occupied by $\Pi$. Since the asymptotic density of any finite set is clearly 0 and the asymptotic density of any set whose complement in $P$ is finite is clearly 1, only sets of primes which are infinite and have an infinite complement in $P$ are of interest in terms of their asymptotic densities. We can in fact be a bit more precise: recall that if $a(x)$ and $b(x)$ denote positive real-valued functions defined on $(0, +\infty)$, then $a(x)$ \emph{is asymptotic to} $b(x)$ \emph{as} $x \rightarrow +\infty$, denoted by $a(x) \sim b(x)$, if 
\[
\lim_{x\rightarrow +\infty } \frac{a(x)}{b(x)}=1.
\] 
 The Prime Number Theorem (LeVeque, [39], chapter 7; Montgomery and Vaughn, [41], chapter 6) asserts that as $x \rightarrow +\infty$,
 \[
 |\{q \in P: q \leq x\}| \sim \frac{x}{ \log x},
 \]
consequently, if $d$ is the density of $\Pi$ then as $x \rightarrow +\infty$,
\[
|\{q \in \Pi: q \leq x\}| \sim d \frac{x}{\log x}.
\]
Hence the asymptotic density of $\Pi$ provides a way to measure precisely the  ``asymptotic cardinality" of $\Pi$.  

\section{The Density of Primes which have a Given Finite Set of Quadratic Residues}

Theorem 4.3 asserts that if $S$ is a given nonempty finite set of positive integers then the set of primes $\{p: \chi_p \equiv 1\ \textrm{on}\ S\}$ is infinite. In this section, we will prove a theorem which provides a way to calculate the density of the set $\{p: \chi_p \equiv 1\ \textrm{on}\ S\}$. This will be given by a formula which depends on a certain combinatorial parameter that is determined by the prime factors of the elements of $S$. In order to formulate this result, let $F$ denote the Galois field $GF(2)$ of 2 elements, which can be concretely realized as the field  $\mathbb{Z}/2\mathbb{Z}$ of ordinary residue classes mod 2. Let $A \subseteq [1, \infty)$. If $n=|A|$, then we let $F^n$ denote the vector space over $F$ of dimension $n$, arrange the elements $a_1< \dots <a_n$ of $A$ in increasing order, and then define the map $v: 2^A \rightarrow F^n$ like so: if $B \subseteq A$ then
 \[
 \textrm{the $i$-th coordinate of}\ v(B)=\left\{\begin{array}{ll}1,\ \textrm{if $a_i \in B$,}\\
0,\ \textrm{if $a_i \notin B$.}\\\end{array}\right. 
\]
If we recall that $\pi_{\textrm{odd}}(z)$ denotes the set of all prime factors of odd multiplicity of the integer $z$ then we can now state (and eventually prove) the following theorem: 
\begin{thm}
If S is a nonempty, finite subset of $[1, \infty)$,
\[
\mathcal{S}= \{\pi_{\textnormal{odd}}(z): z \in S\},
\]
\[
A=\bigcup_{X \in \mathcal{S}}\ X,
\]
\[
n=|A|,
\]
and
\[
d=\ \textrm{the dimension of the linear span of $v(\mathcal{S})$ in $F^n$},
\]
then the density of $\{p: \chi_p \equiv 1\ \textrm{on}\ S\}$ is $2^{-d}$.
\end{thm}

Theorem 4.9 reduces the calculation of the density of $\{p: \chi_p \equiv 1\ \textrm{on}\ S\}$ to prime factorization of the integers in $S$ and linear algebra over $F$. If we enumerate the nonempty elements of $\mathcal{S}$ as $S_1,\dots,S_m$ (if $\mathcal{S}$ has no such elements then $S$ consists entirely of squares, hence the density is clearly 1) then $d$ is just the rank over $F$ of the $m \times n$ matrix
\[
\left( \begin{array} {ll} v(S_1)(1) \dots v(S_1)(n)\\
 \vdots \hspace{2.1cm}      \vdots\\
v(S_m)(1) \dots v(S_m)(n) \end{array} \right),
\]
where $v(S_i)(j)$ is the $j$-th coordinate of $v(S_i)$. This matrix is often referred to as the \emph{incidence matrix of S}. Because there are only two elementary row (column) operations over $F$, namely row (column) interchange and addition of a row (column) to another row (column), the rank of this matrix is easily calculated by Gauss-Jordan elimination. However, this procedure requires that we first find the prime factors of odd multiplicity of each element of $S$, and that, in general, is not so easy!
 
A few examples will indicate how Theorem 4.9 works in practice. Observe first that if $S$ is a finite set of primes of cardinality $n$, say, then the incidence matrix of $S$ is just the $n\times n$ identity matrix over $F$, hence the dimension of $v(\mathcal{S})$ in $F^n$ is $n$, and so the density of $\{p: \chi_p \equiv 1\ \textrm{on}\ S\}$ is $2^{-n}$. Now chose four primes  $p<q<r<s$, say, and let 
\[
S_1=\{p, pq, qr, rs\}.\]
The incidence matrix of $S_1$ is 
\[
\left( \begin{array} {cccc} 1\ 0\ 0\ 0\\
1\ 1\ 0\ 0\\
0\ 1\ 1\ 0\\
0\ 0\ 1\ 1 \end{array} \right),
\]
which is row equivalent to
\[
\left( \begin{array} {cccc} 1\ 0\ 0\ 0\\
0\ 1\ 0\ 0\\
0\ 0\ 1\ 0\\
0\ 0\ 0\ 1 \end{array} \right).
\]
It follows from Theorem 4.9 that the density of $\{p: \chi_p \equiv 1\ \textrm{on}\ S_1\}$ is $2^{-4}$. If
\[
S_2=\{p, ps, pqr, pqrs\},\]
then the incidence matrix of $S_2$ is
\[
\left( \begin{array} {cccc} 1\ 0\ 0\ 0\\
1\ 0\ 0\ 1\\
1\ 1\ 1\ 0\\
1\ 1\ 1\ 1 \end{array} \right),
\]
which is row equivalent to
\[
\left( \begin{array} {cccc} 1\ 0\ 0\ 0\\
0\ 1\ 1\ 1\\
0\ 0\ 0\ 1\\
0\ 0\ 0\ 0 \end{array} \right),
\]
hence Theorem 4.9 implies that the density of $\{p: \chi_p \equiv 1\ \textrm{on}\ S_2\}$ is $2^{-3}$. Because a 2-dimensional subspace of $F^4$ contains exactly 3 nonzero vectors, it follows that if $S$ consists of 4 nontrivial square-free integers such that $S$ is supported on 4 primes, then the density of  $\{p: \chi_p \equiv 1\ \textrm{on}\ S\}$ cannot be $2^{-2}$. However, for example, if
\[
S_3=\{ps, qr, pqrs\},\]
then the incidence matrix of $S_3$ is
\[
\left( \begin{array} {ccc} 1\ 0\ 0\ 1\\
0\ 1\ 1\ 0\\
1\ 1\ 1\ 1\end{array} \right),
\]
which is row equivalent to
\[
\left( \begin{array} {ccc} 1\ 0\ 0\ 1\\
0\ 1\ 1\ 0\\
0\ 0\ 0\ 0\end{array} \right),
\]
and so the density of $\{p: \chi_p \equiv 1\ \textrm{on}\ S_3\}$ is $2^{-2}$.

We turn now to the
\vspace{0.3cm} 
 
 \emph{Proof of Theorem} 4.9. We first establish a strengthened version of Theorem 4.9 in a special case, and then use it (and another lemma) to prove Theorem 4.9 in general.
 \begin{lem}
 $($Filaseta and Richman $[18]$, Theorem $2$$)$ If $\Pi$ is a nonempty set of primes and $\varepsilon: \Pi \rightarrow \{-1, 1\}$ is a given function then the density of the set $\{p: \chi_p \equiv \varepsilon \ \textrm{on}\ \Pi\}$ is $2^{-|\Pi|}$.
 \end{lem}

\emph{Proof}. Let
\[
X=\{p: \chi_p \equiv \varepsilon \ \textrm{on}\ \Pi\},
\]
\[
K=\ \textrm{product of the elements of $\Pi$}.
\]
If $n \in \mathbb{Z}$ then we let $[n]$ denote the  ordinary residue class mod $4K$ which contains $n$. The proof of Lemma 4.10 can now be outlined in a series of three steps.

\emph{Step} 1. Use the LQR to show that 
\[
X=\bigcup_{n \in U(4K): X \cap [n] \not= \emptyset}\ \{p: p \in [n]\}.
\]

\emph{Step $2$ $($and its implementation$)$} . Here we will make use of the Prime Number Theorem for primes in arithmetic progressions, to wit, if $a \in Z$, $b \in [1, \infty)$, $\gcd(a, b)=1$, and $AP(a, b)$ denotes the arithmetic progression with initial term $a$ and common difference $b$,  then as $x \rightarrow +\infty$,
\[
|\{p \in AP(a, b): p \leq x\}| \sim \frac{1}{\varphi(b)} \frac{x}{ \log x}.\]
For a proof of this important theorem, see either LeVeque [39], section 7.4, or Montgomery and Vaughn, [41], section 11.3. In our situation it asserts that if $n \in U(4K)$ then as $x \rightarrow +\infty$,
\[
|\{p \in [n]: p \leq x\}| \sim \frac{1}{\varphi(4K)} \frac{x}{ \log x}.
\]
From this it follows that
\begin{equation*}
\textrm{the density $d_n$ of $\{p: p \in [n]\}$ is $\frac{1}{\varphi(4K)}$, for all $n \in U(4K)$.} \tag{9}
\end{equation*}
\noindent Because the decomposition of $X$ in Step 1 is pairwise disjoint, (9) implies that
\begin{equation*}
\textrm{density of $X$}=\sum_{n \in U(4K): X \cap [n] \not= \emptyset}\ d_n=\frac{|\{ n \in U(4K): X \cap [n] \not= \emptyset \}|}{\varphi(4K)}. \tag{10}
 \]      
 
 \emph{Step} 3. Use the group structure of $U(4K)$ and the LQR to prove that
 \begin{equation*}
|\{n \in U(4K): X \cap [n] \not= \emptyset\}|= \frac{\varphi(4K)}{2^{|\Pi|}}. \tag{11}
\]
From (10) and (11) it follows that the density of $X$ is $2^{-|\Pi|}$, as desired, hence we need only implement Steps 1 and 3 in order to finish the proof. 

\emph{Implementation of Step} 1. We claim that 
\begin{equation*}
\textrm{if $p, p^{\prime}$ are odd primes and $p \equiv p^{\prime}\ \textrm{mod}\ 4K$ then $\chi_p \equiv \chi_{p^{\prime}}$ on $\Pi$}. \tag{12}
\]
Because $X$ is disjoint from $\{2\} \cup \Pi$ and
\begin{equation*}
P \setminus (\{2\} \cup \Pi)= \bigcup_{n \in U(4K)}\ \{p: p \in [n]\}, \tag{13}
\]
the decomposition of $X$ as asserted in Step 1 follows immediately from (12). 

We verify (12) by using the LQR. Assume that $p \equiv p^{\prime}\ \textrm{mod}\ 4K$ and let $q \in \Pi$. Suppose first that $p$ or $q$ is $\equiv 1\ \textrm{mod}\ 4$. Then $p^{\prime}$ or $q$ is $\equiv 1\ \textrm{mod}\ 4$, and so the LQR implies that
\begin{eqnarray*}
\chi_p(q)&=&\chi_q(p)\\
&=&\chi_q(p^{\prime}+4kK)\ \textrm{for some $k \in \mathbb{Z}$}\\
&=&\chi_q(p^{\prime}),\ \textrm{since $q$ divides $4kK$}\\
&=&\chi_{p^{\prime}}(q).
\end{eqnarray*}
 Suppose next that $p \equiv 3 \equiv q\ \textrm{mod}\ 4$. Then $p^{\prime} \equiv 3\ \textrm{mod}\ 4$ hence it follows from the LQR that
 \[
 \chi_p(q)=-\chi_q(p)=-\chi_q(p^{\prime})=-(-\chi_{p^{\prime}}(q))=\chi_{p^{\prime}}(q).\]
 
 \emph{Implementation of Step} 3. Define the equivalence relation $\sim$ on the set of residue classes $\{[n]: n \in U(4K)\}$ like so: 
 \[
 \textrm{$[n] \sim [n^{\prime}]$ if for all odd primes $p \in [n],\ q \in [n^{\prime}],\ \chi_p \equiv \chi_q$ on $\Pi$}.
 \]
 
 We first count the number of equivalence classes of $\sim$. It is a consequence of (12) that the sets
 \[
 \{q \in \Pi: \chi_p(q)=1\}
 \]
 are the same for all $p \in [n]$, and so we let $I(n)$ denote this subset of $\Pi$. Now if $n \in U(4K)$ and $p \in [n]$ then (13) implies that $ p \notin \Pi$. Hence for all $p \in [n]$, $\chi_p$ takes only the values $\pm1$ on $\Pi$. It follows that
 \[
 \textrm{$[n] \sim [n^{\prime}]$ if and only if $I(n)=I(n^{\prime})$.}
\]
On the other hand, by virtue of Lemma 4.4, if $S \subseteq \Pi$ then there exits infinitely many primes $p$ such that
\[
S= \{q \in \Pi: \chi_p(q)=1\}, 
\]
and so we use (13) to find $n_0 \in U(4K)$ such that $[n_0]$ contains at least one of these primes $p$, hence
\[
S=I(n_0).
\]
We conclude that
\begin{equation*}
\textrm{the number of equivalence classes of}\  \sim\ \textrm{is}\ 2^{|\Pi|}. \tag{14}
\]

 Let $E_n$ denote the equivalence class of $\sim$ which contains $[n]$. We claim that
 \begin{equation*}
 \textrm{multiplication by $n$ maps $E_1$ bijectively onto $E_n$}. \tag{15}
 \]
If this is true then $|E_n|$ is constant as a function of $n \in U(4K)$, hence (14) implies that 
\begin{equation*}
\varphi(4K)=2^{|\Pi|}|E_n|,\ \textrm{for all}\ n \in U(4K). \tag{16}
\]
If we now choose $p \in X$ then there is $n_0 \in U(4K)$ such that $p \in [n_0]$, hence it follows from (12) that
\[
E_{n_0}=\{[n]: X \cap [n] \not= \emptyset\},
\]
and so, in light of (16), 
\[
\varphi(4K)=2^{|\Pi|}|\{n \in U(4K): X \cap [n] \not= \emptyset\}|,\]
which is (11).

It remains only to verify (15). Because $U(4K)$ is a group under the multiplication induced by multiplication of ordinary residue classes mod $4K$, it is clear that multiplication by $n$ on $E_1$ is injective, so we need only prove that $nE_1=E_n$.
 
 We show first that $nE_1 \subseteq E_n$. Let $[n^{\prime}] \in E_1$. We must prove: $[nn^{\prime}] \in E_n$, i.e., $[nn^{\prime}] \sim [n]$, i.e., 
  \begin{equation*}
 \textrm{if $p \in  [nn^{\prime}], q \in [n]$ are odd primes then}\ 
\chi_p \equiv \chi_q\ \textrm{on}\ \Pi. \tag{17}
 \]
 
 In order to verify (17), let $p \in  [nn^{\prime}], q \in [n], p^{\prime} \in [n^{\prime}], q^{\prime} \in [1]$ be odd primes. Because $[n^{\prime}] \sim [1]$,
 \begin{equation*}
 \chi_{p^{\prime}} \equiv \chi_{q^{\prime}}\ \textrm{on}\  \Pi. \tag{18}
 \]
 The choice of $p, q, p^{\prime}, q^{\prime}$ implies that
 \[
 pq^{\prime} \equiv p^{\prime}q\ \textrm{mod}\ 4K.
 \]
This congruence and the LQR when used in an argument similar to the one that was used to prove (12) imply that
\begin{equation*}
\chi_p \chi_{q^{\prime}} \equiv \chi_{p^{\prime}} \chi_q\ \textrm{on}\ \Pi. \tag{19}
\]
Because $\chi_{q^{\prime}}$ and $\chi_{p^{\prime}}$ are both nonzero on $\Pi$, we can use (18) to cancel  $\chi_{q^{\prime}}$ and $\chi_{p^{\prime}}$ from each side of (19) to obtain 
\[
\chi_p \equiv \chi_q\ \textrm{on}\ \Pi.\]

We show next that $E_n \subseteq nE_1$. Let $[n^{\prime}] \in E_n$. The group structure of $U(4K)$ implies that there exits $n_0 \in U(4K)$ such that
\begin{equation*}
[nn_0]=[n^{\prime}], \tag{20}
\]
so we need only show that $[n_0] \in E_1$, i.e., 
\begin{equation*}
\textrm{$\chi_p \equiv \chi_q$ on $\Pi$, for all odd primes $p \in [n_0], q \in [1]$}. \tag{21}
\]

Toward that end, choose odd primes $p^{\prime} \in [n], q^{\prime} \in [n^{\prime}]$. Because $[n] \sim [n^{\prime}]$,
\begin{equation*}
\chi_{p^{\prime}} \equiv \chi_{q^{\prime}}\ \textrm{on}\ \Pi, \tag{22}
\]
and so because of (20), we have that for all $p \in [n_0], q \in [1]$,
\[
pp^{\prime} \equiv qq^{\prime}\ \textrm{mod}\ 4K.
\]
(21) is now a consequence of this congruence, (22), and our previous reasoning. $\hspace{1.5cm} \textrm{QED}$

We will prove Theorem 4.9 by combining Lemma 4.10 with the next lemma, a simple result in enumerative combinatorics.
\begin{lem}
If $A$ is a nonempty finite subset of $[1, \infty), n=|A|, \mathcal{S} \subseteq 2^A, F=$ the Galois field of order $2$, $v: 2^A \rightarrow F^n$ is the map defined at the beginning of this section, and 
\[
\textrm{$d=$ the dimension of the linear span of $v(\mathcal{S})$ in $F^n$},
\]
then the cardinality of the set
\[
\mathcal{N}=\{N \subseteq A: |N \cap S|\ \textrm{is even, for all $S \in \mathcal{S}$}\}
\]
is $2^{n-d}$.
\end{lem}  

\emph{Proof}. Without loss of generality take $A=[1, n]$. Observe first that if $N, T \subseteq A$, then 
\[
|N \cap T|\ \textrm{is even if and only if $\sum_{i=1}\ v(N)(i)v(T)(i)=0$ in $F$}.
\]
Hence there is a bijection of the set of all solutions in $F^n$ of the system of linear equations
\begin{equation*}
\sum_1^n\ v(S)(i)x_i=0, S \in \mathcal{S}, \tag{*}
\]
onto $\mathcal{N}$ given by
\[
(x_1,\dots,x_n) \rightarrow \{i: x_i=1 \}.\]
If $m=|\mathcal{S}|$ and $\sigma: F^n \rightarrow F^m$ is the linear transformation whose representing matrix is the coefficient matrix of the system $(*)$ then
\[
\textrm{the set of all solutions of $(*)$ in $F^n=$ the kernel of $\sigma$}.\]
But $d$ is the rank of $\sigma$ and so the kernel of $\sigma$ has dimension $n-d$. Hence
\[
|\mathcal{N}|=|\textrm{the set of all solutions of $(*)$ in $F^n|=|$kernel of $\sigma|=2^{n-d}$}. \]
$\hspace{15.6cm} \textrm{QED}$

We proceed to prove Theorem 4.9. Let $S, \mathcal{S}, A, n,$ and $d$ be as in the hypothesis of that theorem, let
\[
X=\{p: \chi_p \equiv 1\ \textrm{on}\ S\},\]
\[
\mathcal{N}=\{N \subseteq A: |N \cap S|\ \textrm{is even, for all $S \in \mathcal{S}$}\},
\]
and for each prime $p$, let
\[
N(p)=\{q \in A: \chi_p(q)=-1\}.\]
Then since $X$ is disjoint from $A$,
\begin{eqnarray*}
p \in X &\textrm{iff}& 1= \chi_p(z)= \prod_{q \in \pi_{\textrm{odd}}(z)} \chi_p(q),\ \textrm{for all}\ z \in S,\\
&\textrm{iff}& |N(p) \cap \pi_{\textrm{odd}}(z)|\ \textrm{is even, for all}\ z \in S,\\
&\textrm{iff}& N(p) \in \mathcal{N}.
\end{eqnarray*} 
Hence
\[
X=\bigcup_{N \in \mathcal{N}} \{p: N(p)=N\}\]
and this union is pairwise disjoint. Hence
\[
\textrm{density of}\ X=\sum_{N \in \mathcal{N}}\ \textrm{density of}\ \{p: N(p)=N\}.\]
Lemma 4.10 implies that
\[
\textrm{density of}\ \{p: N(p)=N\}=2^{-n}\ \textrm{for all}\ N \in \mathcal{N},\]
and so
\begin{eqnarray*}
\textrm{density of}\ X&=& 2^{-n} |\mathcal{N}|\\
&=&2^{-n}(2^{n-d}),\ \textrm{by Lemma}\ 4.11\\
&=&2^{-d}. 
\end{eqnarray*}
$\hspace{15.6cm} \textrm{QED}$

The next question which naturally arises asks: what about a version of Theorem 4.9 for quadratic non-residues, i.e., for what finite, nonempty subsets $S$ of $[1, \infty)$ is it true that $S$ is a set of non-residues of infinitely many primes? In contrast to what occurs for residues, this can fail to be true for certain finite subsets $S$ of $[1, \infty)$, and there is a simple obstruction that prevents it from being true. Suppose that there is a subset $T$ of $S$ such that $|T|$ is odd and $\prod_{i \in T} i$ is a square, and suppose that $S$ is a set of non-residues of infinitely many primes. We can then choose $p$ to exceed all of the prime factors of the elements of $T$ and  such that $\chi_p(z)=-1$, for all $z \in T$. Hence 
\[
-1=(-1)^{|T|}=\prod_{i \in T}\ \chi_p(i)=\chi_p \Big(\prod_{i \in T}\ i \Big)=1,\]
a clear contradiction. It follows that the presence of such subsets $T$ of $S$ prevents $S$ from being a set of non-residues of infinitely many primes. The next theorem asserts that those subsets are the only obstructions to $S$ having this property.
\begin{thm}
If $S$ is a finite, nonempty subset of $[1,\infty)$ then $S$ is a set of non-residues of infinitely many primes if and only if for all subsets $T$ of $S$ of odd cardinality, $\prod_{i \in T} i$ is not a square.
\end{thm}
 \noindent This theorem lies somewhat deeper than Theorem 4.3. We will prove it in Chapter 5, where we will once again delve into the theory of algebraic numbers. But before we get to that, we will discuss how to use quadratic residues to design zero-knowledge proofs.

\section{Zero-Knowledge Proofs and Quadratic Residues}

A major issue in modern electronic communication is the secure verification of identification, namely, guaranteeing that the person with whom you are communicating is indeed who you think he is. A typical scenario proceeds as follows: person $P$ sends an electronic message to person $V$ in the form of an identification number. $V$ wants to securely verify that $P$ validly possesses the ID number, without knowing anything more about $P$. Moreover, for security reasons, $P$ does not want $V$ to be able to find out anything about him during the verification procedure, i.e., $V$ is to have \emph{zero knowledge} of $P$. In addition to all of this, $V$ wants to make it virtually impossible for any other person $C$ to use the verification procedure to deceive $V$ into thinking that $C$ is $P$. An identity-verification  algorithm which satisfies all of these requirements is called a \emph{zero-knowledge proof}. 

Zero-knowledge proofs which employ quadratic residues were devised in the 1980's because of the need to maintain security when verifying identification numbers using smart cards, electronic banking and stock transactions, and other similar types of communication. In a zero-knowledge proof there are two parties, the \emph{prover}, a person who wants his identity verified without divulging any other information about himself, and a \emph{verifier}, a person who must be convinced that the prover is who he says he is. The identity of the prover is verified by checking that he has certain secret information that only he possesses. Security is maintained because the procedures used in the zero-knowledge proof guarantee that the probability that someone pretending to be the prover can convince the verifier that she is the prover is extremely small. Moreover, the verifier checks only that the prover is in possession of the secret information, without being able to discover what the secret information is.

We will describe a zero-knowledge proof discovered by Adi Shamir [53] in 1985 (we follow Rosen [48], sections 11.3 and 11.5 for the exposition in this section and in sections 8 and 9 below). The prover $P$ starts by choosing two very large primes $p$ and $q$ such that $p\equiv q \equiv$ 3 mod 4 (to maintain security, these primes should have hundreds of digits), computing $n=pq$, and then sending $n$ to the verifier $V$. Let $I$ be a positive integer that represents particular information, e.g., the personal identification number of $P$. $P$ selects a positive number $c$ such that the integer $w$ obtained by concatenating $I$ with $c$ (the integer obtained by writing the digits of $I$ followed by the digits of $c$) is a quadratic residue modulo $n$, i.e., there is a solution in integers of the congruence $x^2\equiv w$ mod $n$, with $\gcd(x, n)=1$. $P$ sends $w$ to $V$ and then finds a solution $u$ of this congruence. Finding $u$ can easily be done by means of Euler's criterion. In order to see that, note first that $\chi_p(w)=\chi_q(w)=1$ and recall that  $p\equiv q \equiv$ 3 mod 4. Euler's criterion therefore implies that 
\[
w^{\frac{1}{2}(p-1)}\equiv \chi_p(w)=1\ \textrm{mod}\ p,\]
\[
w^{\frac{1}{2}(q-1)}\equiv \chi_q(w)=1\ \textrm{mod}\ q,\]
hence
\[
\big(w^{\frac{1}{4}(p+1)}\big)^2=w^{\frac{1}{2}(p+1)}=w^{\frac{1}{2}(p-1)}\cdot w\equiv w\ \textrm{mod}\ p,\] 
and similarly,
\[
\big(w^{\frac{1}{4}(q+1)}\big)^2\equiv w\ \textrm{mod}\ q.\] 
The prover then finds a solution $u$ of $x^2\equiv w$ mod $n$ by using the Chinese remainder theorem to solve the congruences
\[
u\equiv w^{\frac{1}{4}(p+1)}\ \textrm{mod}\ p,\]
\[
u\equiv w^{\frac{1}{4}(q+1)}\ \textrm{mod}\ q.\]
Of course, in order to find $u$ in this way, one must know the primes $p$ and $q$.

$P$ convinces $V$ that $P$ knows $u$ by using an interactive proof that is composed of iterations of the following four-step cycle:

$(i)$ $P$ choses a random number $r$ and sends $V$ a message containing two integers: $x$, where $x\equiv r^2$ mod $n$, with $\gcd(x, n)=1$ and $0\leq x<n$, and $y$, where $y\equiv w\overline{x}$ mod $n$,  $0\leq x<n$, and $\overline{x}$ denotes the inverse of $x$ modulo $n$.

$(ii)$ $V$ checks that $xy\equiv w$ mod $n$, then choses a random bit $b$ equal to either 0 or 1, and sends $b$ to $P$.

$(iii)$ If $b=0$, $P$ sends $r$ to $V$. If $b=1$ then $P$ calculates $s\equiv u\overline{r}$ mod $n$, $0\leq s<n$, and sends $s$ to $V$.

$(iv)$ $V$ computes the square modulo $n$ of what $P$ has sent. If $V$ sent 0, she checks that this square is $x$, i.e., $r^2\equiv x$ mod $n$. If $V$ sent 1 then she checks that this square is $y$, i.e., $s^2\equiv y$ mod $n$.

This cycle can be iterated many times to guarantee security and to convince $V$ that $P$ knows his private information $u$, which shows that $P$ validly possesses the identification number $I$, i.e., that $P$ is who he says he is. By passing this test over many cycles, $P$ has shown that he can produce either $r$ or $s$ upon request. Hence $P$ must know $u$ because in each cycle, he knows both $r$ and $s$, and $u\equiv rs$ mod $n$. 
Moreover, $V$ is unable to discover what $u$ is because that would require $V$ to be able to solve the square-root problem $x^2\equiv w$ mod $n$ \emph{without } knowing $p$ and $q$. This problem is known in cryptology circles as the \emph{quadratic residuosity problem}, and is regarded to be computationally intractable, hence essentially impossible to solve in any feasible length of time, when the modulus $n$ is the product of two very large unknown primes.

Because the bit chosen by $V$ is random, the probability that it is a 0 is $1/2$ and the probability that it is a 1 is $1/2$. If someone does not know $u$, the modular square root of $w$, then the probability that they will pass one iteration of the cycle is almost exactly $1/2$. If an impostor is attempting to deceive $V$ into believing that she, the impostor, is $P$, the probability of the impostor passing, say, 30 iterations of the cycle is hence approximately $1/2^{30}$, less than one in a billion. This makes it virtually impossible for $V$ to be deceived in this manner.

We now ask the following question: what has quadratic reciprocity got to do with all of this? We begin our answer to this question by recalling that in the initial steps of the Shamir zero-knowledge proof, the prover needs to find an integer $c$ such that the concatenation $w$ of $I$ with $c$ is a quadratic residue of $n=pq$, where $p$ and $q$ are very large primes. This can be done if and only if $\chi_p(w)=\chi_q(w)=1$, hence the prover must be able to compute Legendre symbols quickly and efficiently. As we have seen, if one can find sufficiently many factors of $w$ then the LQR can be used to perform this computation in the desired manner. Unfortunately, one of the outstanding, and very difficult, unsolved problems in computational number theory is the design of a computationally fast and efficient algorithm for factoring very large integers, and the ID number $I$, and also the integer $w$ in Shamir's algorithm, is often taken to be large for reasons of security. This difficulty precludes quadratic reciprocity from being used directly to compute Legendre symbols in Shamir's algorithm. On the other hand, fortunately, there is a very fast and efficient algorithm for computing Legendre symbols which avoids factoring, so much so that when using it, one can, using high-speed computers, of course, very quickly find the quadratic residues required in Shamir's zero-knowledge proof. We will now describe this algorithm for computing Legendre symbols, and it is in the verification of this algorithm that quadratic reciprocity will find its application.
\section{Jacobi Symbols}
The device on which our algorithm is based is a generalization of the Legendre symbol, due to Jacobi. We first define the \emph{Jacobi symbol} $\chi_1(m)$ to be 1 for all integers $m$. Now let $n>1$ be an odd integer, with prime factorization $n=p_1^{t_1}\cdots p_k^{t_k}$. If $m$ is a positive integer relatively prime to $n$, then the \emph{Jacobi symbol} $\chi_n(m)$ is defined as the product of Legendre symbols
\[
\chi_n(m)=\prod_{i=1}^k\ \chi_{p_i}(m)^{t_i}.\]
We emphasize here that this notation for the Jacobi symbol is not standard; we have chosen it to align with the character-theoretic notation we have used for Legendre symbols.

The Jacobi symbols satisfy exactly the same algebraic properties of the Legendre symbols, i.e., 

(a) if $a$ and $b$ are both relatively prime to $n$ and $a\equiv b$ mod $n$ then $\chi_n(a)=\chi_n(b)$;

(b) if $a$ and $b$ are both relatively prime to $n$ then $\chi_n(ab)=\chi_n(a)\chi_n(b)$.

\noindent It follows from (a) and (b) that if $n>1$ and we define the Jacobi symbol $\chi_n(m)$ to be zero whenever $\gcd(m, n)>1$ then $\chi_n$ is a real Dirichlet character of modulus $n$. Moreover the Jacobi symbols satisfy an exact analog of the first and second supplementary laws for the Legendre symbols:

(c) $\chi_n(-1)=(-1)^{\frac{1}{2}(n-1)}$;

(d) $\chi_n(2)=(-1)^{\frac{1}{8}(n^2-1)}$.

\noindent But that is not all! The Jacobi symbols also satisfy an exact analog of the Law of Quadratic Reciprocity, to wit,
\begin{thm}
$($Reciprocity Law for the Jacobi symbol $)$ If m and n are relatively prime odd positive integers then
\[
\chi_m(n)\chi_n(m)=(-1)^{\frac{1}{2}(m-1)\cdot \frac{1}{2}(n-1)}.\]
\end{thm}

Because they are necessary for the verification of the algorithm for the computation of Legendre symbols that we require, we will now prove properties (a), (b) and (d) and Theorem 4.13. As the verification of (a) and (b) are easy consequences of the definition of the Jacobi symbol and the analogous properties of the Legendre symbol, we can safely leave those details to the reader.

In order to verify (d), begin by letting $p_1^{t_1}\cdots p_m^{t_m}$ be the prime factorization of $n$. Then from Theorem 2.6 it follows that
\[
\chi_n(2)=\prod_{i=1}^m\ \chi_{p_i}(2)^{t_i}=(-1)^{\sigma},\]
where 
\[
\sigma=\sum_{i=1}^m\ \frac{t_i(p_i^2-1)}{8}\ .\]
We have that
\[
n^2=\prod_{i=1}^m\ \big(1+(p_i^2-1)\big)^{t_i}.\]
Because $p_i^2-1\equiv 0$ mod 8, for $i=1,\dots, m$, it follows that
\[
\big(1+(p_i^2-1)\big)^{t_i}\equiv 1+t_i(p_i^2-1)\ \textrm{mod}\ 64\]
and
\[
\big(1+t_i(p_i^2-1)\big)\big(1+t_j(p_j^2-1)\big)\equiv 1+t_i(p_i^2-1)+t_j(p_j^2-1)\ \textrm{mod}\ 64.\]
Hence
\[
n^2\equiv1+\sum_{i=1}^m\ t_i(p_i^2-1)\ \textrm{mod}\ 64,\]
which implies that
\[
\frac{n^2-1}{8}\equiv \sum_{i=1}^m\ \frac{t_i(p_i^2-1)}{8}=\sigma\ \textrm{mod}\ 8.\]
Therefore
\[
\chi_n(2)=(-1)^{\sigma}=(-1)^{\frac{1}{8}(n^2-1)}.\]

$\hspace{15cm}\ \textrm{QED}$

We begin the proof of Theorem 4.13 by  letting $p_1^{a_1}\cdots p_s^{a_s}$ and $q_1^{b_1}\cdots q_r^{b_r}$ be the prime factorizations of $m$ and $n$. Then
\[
\chi_n(m)=\prod_{i=1}^r\chi_{q_i}(m)^{b_i}=\prod_{i=1}^r\ \prod_{j=1}^s\chi_{q_i}(p_j)^{b_ia_j}\]
and
\[
\chi_m(n)=\prod_{j=1}^s\chi_{p_j}(n)^{a_j}=\prod_{j=1}^s\ \prod_{i=1}^r\chi_{p_j}(q_i)^{b_ia_j}.\]
Hence
\[
\chi_n(m)\chi_m(n)=\prod_{i=1}^r\ \prod_{j=1}^s\left[\chi_{q_i}(p_j)\chi_{p_j}(q_i)\right]^{a_jb_i}.\]
Because $m$ and $n$ are odd and relatively prime, all of the primes in the prime factorizations of $m$ and $n$ are odd and no prime factor of $m$ is a factor of $n$. The LQR thus implies that
\[
\chi_{q_i}(p_j)\chi_{p_j}(q_i)=(-1)^{\frac{1}{2}(p_j-1_)\frac{1}{2}(q_i-1)}.\]
Hence
\begin{equation*}
\chi_n(m)\chi_m(n)=\prod_{i=1}^r\ \prod_{j=1}^s(-1)^{a_j\frac{1}{2}(p_j-1)b_i\frac{1}{2}(q_i-1)}=(-1)^{\kappa},\tag{23}\]
where
\[
\kappa=\sum_{i=1}^r\sum_{j=1}^s\ \frac{a_j(p_j-1)}{2}\cdot \frac{b_i(q_i-1)}{2}.\]

We have that
\[
\sum_{i=1}^r\sum_{j=1}^s\ \frac{a_j(p_j-1)}{2}\cdot \frac{b_i(q_i-1)}{2}=\sum_{j=1}^s\ \frac{a_j(p_j-1)}{2}\sum_{i=1}^r\ \frac{b_i(q_i-1)}{2}.\]
Because
\[
m=\prod_{i=1}^s\ \big(1+(p_i-1)\big)^{a_i}\]
and $p_i-1$ is even, it follows that
\[
\big(1+(p_i-1)\big)^{a_i}\equiv 1+a_i(p_i-1)\ \textrm{mod}\ 4,\]
and
\[
\big(1+a_i(p_i-1)\big)\big(1+a_j(p_j-1)\big)\equiv 1+a_i(p_i-1)+a_j(p_j-1)\ \textrm{mod}\ 4.\]
Hence
\[
m\equiv 1+\sum_{i=1}^s\ a_i(p_i-1)\ \textrm{mod}\ 4,\] 
and so
\[
\sum_{i=1}^s\ \frac{a_i(p_i-1)}{2}\equiv \frac{m-1}{2}\ \textrm{mod}\ 2.\]
Similarly,
\[
\sum_{i=1}^r\ \frac{b_i(q_i-1)}{2}\equiv \frac{n-1}{2}\ \textrm{mod}\ 2.\]
Therefore,
\begin{equation*}
\kappa=\sum_{i=1}^r\sum_{j=1}^s\ \frac{a_j(p_j-1)}{2}\cdot \frac{b_i(q_i-1)}{2}\equiv \frac{m-1}{2}\cdot \frac{n-1}{2}\ \textrm{mod}\ 2.\tag{24}\]
It now follows from (23) and (24) that
\[
\chi_n(m)\chi_m(n)=(-1)^{\kappa}=(-1)^{\frac{1}{2}(m-1)\cdot \frac{1}{2}(n-1)}.\]

$\hspace{15cm}\ \textrm{QED}$
\section{An Algorithm for Fast Computation of Legendre Symbols}
The key ingredient of the algorithm for the computation of Legendre symbols that we want is a formula for the computation of certain Jacobi symbols. That formula uses data given in the form of two finite sequences of integers which are generated by a successive division and factorization procedure. In order to state that formula we start with two relatively prime positive integers $a$ and $b$ with $a>b$. We will generate two finite sequences of integers from $a$ and $b$ by using a modification of the Euclidean algorithm as follows: let $a=R_0$ and $b=R_1$. Using the division algorithm and then factoring out the highest power of 2 from the remainder, we obtain
\[
R_0=R_1q_1+2^{s_1}R_2,\]
where $\gcd(R_1, R_2)=1$ and $R_2$ is odd. Now successively apply the division algorithm as follows, factoring out the highest power of 2 from the remainders as you do so:
\begin{eqnarray*}
R_1&=&R_2q_2+2^{s_2}R_3\\
R_2&=&R_3q_3+2^{s_3}R_4\\   
&\vdots&\\
R_{n-2}&=&R_{n-1}q_{n-1}+2^{s_{n-1}}\cdot 1\\
R_n&=&1,\ s_n=0.
\end{eqnarray*}
Note that $R_i$ is an odd positive integer and $s_i$ is a nonnegative integer for $i=1,\dots,n$, and $\gcd(R_i, R_{i+1})=1$ for $i=0,\dots,n-1$. Because $R_{i+1}<R_i$ for each $i$, this division process will always terminate. The formula for the computation of the Jacobi symbols that is required can now be stated and proved:
\begin{prp}
If a and b be relatively prime positive integers such that $a>b$, b is odd, and $R_i$ and $s_i$, $i=1,\dots,n$, are the sequences of integers generated by the preceding algorithm, then
\[
\chi_b(a)=(-1)^{\sigma},\]
where
\[
\sigma=\sum_{i=1}^{n-1}\ \left(s_i\frac{R_i^2-1}{8}+\frac{(R_i-1)(R_{i+1}-1)}{4}\right).\]
\end{prp}

\emph{Proof}. From properties (a), (b), and (d) of the Jacobi symbol, it follows that
\begin{eqnarray*}
\chi_b(a)&=&\chi_{R_1}(R_0)=\chi_{R_1}(2^{s_1}R_2)\\
&=&\chi_{R_1}(2)^{s_1}\chi_{R_1}(R_2)\\
&=&(-1)^{s_1\cdot \frac{R_1^2-1}{8} }\chi_{R_1}(R_2),
\end{eqnarray*}
and it follows from Theorem 4.13 that
\[
\chi_{R_1}(R_2)=(-1)^{\frac{R_1-1}{2}\frac{R_2-1}{2}}\chi_{R_2}(R_1),\]
hence
\[
\chi_b(a)=(-1)^{\sigma_1}\chi_{R_2}(R_1),\]
where
\[
\sigma_1=s_1\frac{R_1^2-1}{8}+\frac{(R_1-1)(R_2-1)}{4}\ .\]
In the same manner, we obtain for $i=2\dots,n-1$,
\[
\chi_{R_i}(R_{i-1})=(-1)^{\sigma_i}\chi_{R_{i+1}}(R_i),\]
where
\[
\sigma_i=s_i\frac{R_i^2-1}{8}+\frac{(R_i-1)(R_{i+1}-1)}{4}\ .\]
When all of these equations are combined, the desired expression for $\chi_b(a)$ is produced. QED

The algorithm for the computation of Legendre symbols can now be described in a simple three-step procedure like so: let $p$ be an odd prime, $a$ a positive integer less than $p$; we wish to compute the Legendre symbol $\chi_p(a)$. 

\emph{Step} 1. Factor $a=2^sb$ where $b$ is odd (in Shamir's algorithm, this step can always be avoided by concatenating an odd integer to the integer $I$).

Theorem 2.6 implies that
\begin{equation*}
\chi_p(a)=\chi_p(2)^s\chi_p(b)=(-1)^{s\cdot \frac{p^2-1}{8} }\chi_p(b).\tag{25}\]
Now use Theorem 4.13 to obtain
\begin{equation*}
\chi_p(b)=(-1)^{\frac{1}{2}(p-1)\frac{1}{2}(b-1)}\chi_b(p).\tag{26}\]
Substitution of (26) into (25) yields

\emph{Step} 2. Write
\[
\chi_p(a)=(-1)^{\varepsilon}\chi_b(p),\]
where
\[
\varepsilon=\frac{s(p^2-1)}{8}+\frac{(p-1)(b-1)}{4}\ .\]

\emph{Step} 3. Use the formula from Proposition 4.14 to compute $\chi_b(p)$ and substitute that value into the formula for $\chi_p(a)$ in Step 2.

As an example, we use this algorithm to calculate $\chi_{311}(141)$ without factoring the argument 141. Because 141 is odd, Step 1 yields $s=0$, hence from Step 2 we obtain
\[
\chi_{311}(141)=\chi_{141}(311).\]
In Step 3, we need the sequence of divisions
\begin{eqnarray*}
311&=&141\cdot12+2^0\cdot29\\
141&=&29\cdot4+2^0\cdot25\\
29&=&25\cdot1+2^2\cdot1,
\end{eqnarray*}
and so the sequences that are required to apply Proposition 4.14 are $R_1=141, R_2=29, R_3=25, R_4=1$ and $s_1=0, s_2=0, s_3=2$. Hence from Step 2, we see that
\[
\chi_{311}(141)=(-1)^{\sigma},\]
where
\begin{eqnarray*}
\sigma&=&0\cdot\frac{141^2-1}{8}+0\cdot\frac{29^2-1}{8}+2\cdot\frac{25^2-1}{8}+\frac{(141-1)(29-1)}{4}+\frac{(29-1)(25-1)}{4}\\
&\equiv& 0\ \textrm{mod}\ 2,
\end{eqnarray*}
hence
\[
\chi_{311}(141)=1.\]

Of course in this simple example, we can obviously factor 141 completely and then use the LQR as before, but the whole point of the example is to calculate  $\chi_{311}(141)$ \emph{without any factoring}. In practical applications of quadratic residues in cryptology, such as Shamir's zero-knowledge proof, the arguments of Legendre symbols are frequently very large, and so complete factorization of the argument becomes computationally unfeasible. 

How can the efficiency of our algorithm for the calculation of Legendre symbols be measured when it is implemented for computation on modern high-speed computers? Integer calculations on a computer are done by using base-2 expansions of the integers, which are called \emph{bit strings}. A \emph{bit operation} is the addition, subtraction or multiplication of two bit strings of length 1, the division of a bit string of length 2 by a bit string of length 1 using the division algorithm, or the shifting of a bit string by one place. The computational efficiency of an algorithm is measured by its \emph{computational complexity}, which is an estimate of the number of bit operations that are needed to carry out the algorithm when it is programmed to run on a computer. Because our algorithm for the computation of $\chi_p(a)$ uses a variation of the Euclidean algorithm in Step 3, which accounts for most of the computational complexity, one can show that the algorithm requires only $O\big((\log_2a)^2\big)$ bit operations to compute  $\chi_p(a)$, which means that the algorithm is very fast and efficient. Thus one can very quickly determine the integer $w$ that is needed to implement Shamir's algorithm.

In addition to finding a quadratic residue $w$ of $n$, the initial steps in Shamir's algorithm also requires the determination of the square root of $w$ modulo $n$. The simple procedure that we described for computing this square root uses the powers $w^{\frac{1}{4}(p+1)}$ and $w^{\frac{1}{4}(q+1)}$ in an application of the Chinese remainder theorem, with the exponents of $w$ here being extremely large. This situation thus calls for a quick and efficient procedure for the computation of high-powered modular exponentiation, and so we will now present an algorithm which does that. 

The problem is to compute, for given positive integers $b, n$, and $N$ with $b<n$, the power $b^N$ mod $n$. We do this by first expressing the exponent $N$ in its base-2 expansion $(a_k a_{k-1}\dots a_1 a_0)_{\textrm{base}\ 2}$. Then compute the nonnegative minimal ordinary residues mod $n$ of $b, b^2,\dots,b^{2^k}$ by successively squaring and reducing mod $n$. The final step is to multiply together the minimal nonnegative ordinary residues of $b^{2^i}$ which correspond to $a_i=1$, reducing modulo $n$ after each multiplication. It can be shown that the nonnegative minimal ordinary residue of $b^N$ mod $n$ can be computed by this algorithm using only $O\big((\log_2 n)^2\log_2 N\big)$ bit operations.

The following example illustrates the calculations which are typically involved. We wish to compute $15^{402}$ mod 1607. The binary expansion of 402 is $110010010$. We calculate that
\begin{eqnarray*}
15&\equiv& 15\ \textrm{mod}\ 1607\\
15^2&\equiv& 225\ \textrm{mod}\ 1607\\
15^4&\equiv& 808\ \textrm{mod}\ 1607\\
15^8&\equiv& 422\ \textrm{mod}\ 1607\\
15^{16}&\equiv& 1314\ \textrm{mod}\ 1607\\
15^{32}&\equiv& 678\ \textrm{mod}\ 1607\\
15^{64}&\equiv& 82\ \textrm{mod}\ 1607\\
15^{128}&\equiv& 296\ \textrm{mod}\ 1607\\
15^{256}&\equiv& 838\ \textrm{mod}\ 1607.
\end{eqnarray*}
It follows that
\begin{eqnarray*}
15^{402}&=&15^{256+128+16+2}\\
&\equiv&838\cdot 296\cdot 1314\cdot 225\ \textrm{mod}\ 1607\\
&\equiv&570\cdot1314\cdot 225\ \textrm{mod}\ 1607\\
&\equiv&118\cdot 225\ \textrm{mod}\ 1607\\
&\equiv&838\ \textrm{mod}\ 1607.
\end{eqnarray*}

\chapter{The Zeta Function of an Algebraic Number Field and Some Applications}

At the end of section 6 of chapter 4, we left ourselves with the problem of determining the finite nonempty subsets $S$ of the positive integers such that for infinitely many primes $p$, $S$ is a set of non-residues of $p$. We observed there that if $S$ has this property then the product of all the elements in every subset of $S$ of odd cardinality is never a square. The object of this chapter is to prove the converse of this statement, i.e., we wish to prove Theorem 4.12. The proof of Theorem 4.12 that we present uses ideas that are closely related to the ones that Dirichlet used in his proof of Theorem 4.5, together with some technical improvements due to Hilbert. The key tool that we need is an analytic function attached to algebraic number fields, called the \emph{ zeta function} of the field. The definition of this function requires a significant amount of mathematical technology from the theory of algebraic numbers, and so in section 1 we begin with a discussion of the results from algebraic number theory that will be required, with Dedekind's Ideal Distribution Theorem as the final goal of this section. The zeta function of an algebraic number field is defined and studied in section 2; in particular, the Euler-Dedekind product formula for the zeta function is derived here. In section 3 a product formula for the zeta function of a quadratic number field that will be required in the proof of Theorem 4.12 is derived from the Euler-Dedekind product formula. The proof of Theorem 4.12, the principal object of this chapter, is carried out in section 4 and some results which are closely related to that theorem are also established there. In the interest of completeness, we prove in section 5 the Fundamental Theorem of Ideal Theory, Theorem 3.16 of Chapter 3, since it is used in an essential way in the derivation of the Euler-Dedekind product formula.

\section{Dedekind's Ideal Distribution Theorem}

We have already seen in sections 11 and 12 of Chapter 3 how the factorization of ideals in a quadratic number field can be used to prove the Law of Quadratic Reciprocity. The crucial fact on which that proof of quadratic reciprocity relies is the Fundamental Theorem of Ideal Theory (Theorem 3.16), the result which describes the fundamental algebraic structure of the ideals in the ring $R$ of algebraic integers in an algebraic number field $F$. As we mentioned in Chapter 3, the Fundamental Theorem of Ideal Theory is due to Richard Dedekind. In order to define and study the zeta function of $F$, we will need another very important theorem of Dedekind which provides a precise numerical measure of how the ideals of $R$ are distributed in $R$ according to the cardinality of the quotient rings of $R$ modulo the ideals. This result is often called Dedekind's Ideal Distribution Theorem, and the purpose of this section is to develop enough of the theory of ideals in $R$ so that we can state the Ideal Distribution Theorem precisely. All of this information will then be used in the next section to define the zeta function and establish the properties of the zeta function that we will need to prove Theorem 4.12.

Let $F$ denote an algebraic number field of degree $n$ that will remain fixed in the discussion until indicated otherwise, and let $R$ denote the ring of algebraic integers in $F$. In section 11 of chapter 3, we mentioned that every prime ideal of $R$ is maximal and that the cardinality of the quotient ring $R/I$ of $R$ is finite for all nonzero ideals $I$ of $R$.  Consequently, the ideals of $R$ are exceptionally ``large" subsets of $R$. We begin our discussion here by proving these facts as part of the following proposition.
\begin{prp}
$(i)$ An ideal of $R$ is prime if and only if it is maximal.

$(ii)$ If $I$ is a non-zero ideal of $R$ then the cardinality of the quotient ring $R/I$ is finite.

$(iii)$ If $I$ is a prime ideal of $R$ then there exists a rational prime $q \in \mathbb{Z}$ such that $I \cap \mathbb{Z}=q\mathbb{Z}$. In particular $q$ is the unique rational prime contained in $I$.

$(iv)$ If $I$ is a prime ideal of $R$ and $q$ is the rational prime in $I$ then $R/I$ is a finite field of characteristic $q$, hence there exists a unique positive integer $d$ such that $|R/I|=q^d$.
\end{prp}

\emph{Proof}. The proof of statements $(i)$ and $(ii)$ of Proposition 5.1 depend on the existence of an integral basis of $R$. A subset $\{\alpha_1,\dots,\alpha_k\}$  of $R$ is an \emph{integral basis of R} if for each $\alpha \in R$, there exists a $k$-tuple $(z_1,\dots,z_k)$ of integers, \emph{uniquely determined by $\alpha$}, such that
\[
\alpha=\sum_{i=1}^k z_i\alpha_i.\]
It is an immediate consequence of the definition that an integral basis $\{\alpha_1,\dots,\alpha_k\}$ is linearly independent over $\mathbb{Z}$, i.e., if $(z_1,\dots,z_k)$ is a $k$-tuple of integers such that $\sum_{i=1}^k z_i\alpha_i=0$ then $z_i=0$ for $i=1,\dots,k$. $R$ always has an integral basis (the interested reader may consult Hecke [27], section 22, Theorem 64, for a proof of this), and it is not difficult to prove that every integral basis of $R$ is a basis of $F$ as a vector space over $\mathbb{Q}$; consequently, all integral bases of $R$ contain exactly $n$ elements. 

Now for the proof of $(i)$. Let $I$ be a prime ideal of $R$: we need to prove that $I$ is a maximal ideal, i.e., we take an ideal $J$ of $R$ which properly contains $I$ and show that $J=R$.

Toward that end, let $\{\alpha_1,\dots,\alpha_n\}$ be an integral basis of $R$, and let $0 \not= \beta \in I$. If
\[
x^m+\sum_{i=0}^{m-1} z_i x^i \]
is the minimal polynomial of $\beta$ over $Q$ then $z_0 \not=0$ (otherwise, $\beta$ is the root of a nonzero polynomial over $Q$ of degree less that $m$) and
\[
z_0=-\beta^m-\sum_1^{m-1} z_i\beta^i \in I,\]
hence $\pm z_0 \in I$, and so $I$ contains a positive integer $a$. We claim that each element of $R$ can be expressed in the form
\[
a\gamma+\sum_1^n r_i \alpha_i,\] 
where $\gamma \in R,\  r_i \in [0, a-1], i=1,\dots,n.$

Assume this for now, and let $\alpha \in J \setminus I$. Then for each $k \in [1, \infty)$,
\[
\alpha^k=a\gamma_k+\sum_1^n r_{ik} \alpha_i,\ \gamma_k \in R,\  r_{ik} \in [0, a-1],\ i=1,\dots,n,\]
hence the sequence $(\alpha^k-a\gamma_k: k \in [1, \infty))$ has only finitely many values; consequently there exist positive integers $l<k$ such that
\[
\alpha^l-a \gamma_l=\alpha^k-a\gamma_k.\]
Hence
\[
\alpha^l(\alpha^{k-l}-1)=\alpha^k-\alpha^l=a(\gamma_k-\gamma_l) \in I\  (a \in I\ !).\]
Because $I$ is prime, either $\alpha^l \in I$ or $\alpha^{k-l}-1 \in I$. However, $\alpha^l \not \in I$ because $\alpha \not \in I$ and $I$ is prime. Hence
\[
 \alpha^{k-l}-1 \in I \subseteq J.\]
 But $k-l>0$ and $\alpha \in J$ (by the choice of $\alpha$), and so $-1 \in J$. As $J$ is an ideal, this implies that $J=R$ .

Our claim must now be verified. Let $\alpha \in R$, and find $z_i \in \mathbb{Z}$ such that 
\[
\alpha=\sum_{i=1}^n  z_i\alpha_i.\]
The division algorithm in $\mathbb{Z}$ implies that there exist $m_i \in \mathbb{Z}$, $r_i \in [1, a-1],\ i=1,\dots,n$, such that $z_i=m_ia+r_i,\ i=1,\dots,n$. Thus
\[
\alpha=a\sum_i m_i\alpha_i+\sum_ir_i\alpha_i=a\gamma+\sum_ir_i\alpha_i,\]
with $\gamma \in R$.  

We verify $(ii)$ next. Let $L \not=\{0\}$ be an ideal of $R$. We wish to show that $|R/L|$ is finite. A propos of that, choose $a \in L \cap \mathbb{Z}$ with $a>0$ (that such an $a$ exists follows from the previous  proof of statement $(i)$). Then $aR \subseteq L$, hence there is a surjection of $R/aR$ onto $R/L$, whence it suffices to show that $|R/aR|$ is finite.

We will in fact prove that $|R/aR|=a^n$. Consider for this the set
\[
S=\Big\{\sum_i z_i\alpha_i: z_i \in [0, a-1] \Big\}.\]
We show that $S$ is a set of coset representatives of $R/aR$; if this is true then clearly $|R/aR|=|S|=a^n$. Thus, let $\alpha=\sum_i z_i\alpha_i \in R$. 
Then there exist $m_i \in \mathbb{Z}$, $r_i \in [0, a-1],\ i=1,\dots,n$, such that $z_i=m_ia+r_i,\ i=1,\dots,n$. Hence 
\[
\alpha-\sum_i r_i\alpha_i=\Big(\sum_i m_i\Big)a \in aR\ \textrm{and}\ \sum_i r_i\alpha_i \in S, \]
and so each coset of $R/aR$ contains an element of $S$.

Let $\sum_i a_i\alpha_i, \sum_i a_i^{\prime}\alpha_i$ be elements of $S$ in the same coset. Then
\[
\sum_i (a_i-a_i^{\prime})\alpha_i=a\alpha,\ \textrm{for some}\  \alpha \in R.\]
Hence there exists $m_i \in \mathbb{Z}$ such that
\[
\sum_i (a_i-a_i^{\prime})\alpha_i=\sum_i m_ia\alpha_i,\]
and so the linear independence (over $\mathbb{Z}$) of $\{\alpha_1,\dots,\alpha_n\}$ implies that
\[
a_i-a_i^{\prime}=m_ia,\ i=1,\dots,n\]
i.e., $a$ divides $a_i-a_i^{\prime}$ in $\mathbb{Z}$. Because $|a_i-a_i^{\prime}|<a$ for all $i$, it follows that $a_i-a_i^{\prime}=0$ for all $i$. Hence each coset of $R/aR$ contains exactly one element of $S$. 

In order to verify $(iii)$, note first that the proof of statement $(i)$ implies that $I \cap \mathbb{Z} \not=\{0\}$ and $I \cap \mathbb{Z} \not=\mathbb{Z}$ because $1\not \in I$. Hence $I \cap \mathbb{Z}$ is a prime ideal of $\mathbb{Z}$, and is hence generated in $\mathbb{Z}$ by a unique prime number $q$.

Finally, we prove $(iv)$ by concluding from Proposition 5.1$(i)$ that $I$ is a maximal ideal of $R$: a standard result in elementary ring theory asserts that if $M$ is a maximal ideal in a commutative ring $A$ with identity then the quotient ring $A/M$ is a field (Hungerford [29], Theorem III.2.20), hence $R/I$ is a field, and is finite by Proposition 5.1$(ii)$.

To see that $R/I$ has characteristic $q$, note first that $I \cap \mathbb{Z}=q\mathbb{Z}$, and so there is a natural isomorphism of the field $\mathbb{Z}/q\mathbb{Z}$ into $R/I$ such that the identity in $\mathbb{Z}/q\mathbb{Z}$ is mapped onto the identity of $R/I$. Because $\mathbb{Z}/q\mathbb{Z}$ has characteristic $q$, it follows that if $\bar{1}$ is the identity in $R/I$ then $q\bar{1}=0$ in $R/I$, and $q$ is the least positive integer $n$ such that $n\bar{1}=0$ in $R/I$. Hence $R/I$ has characteristic $q$. $\hspace{11.8cm} \textrm{QED}$
 
\emph{Remark}. It is a consequence of Theorem 3.16 and Proposition 5.1 that $R$ contains infinitely many prime ideals.

It follows from Proposition 5.1$(ii)$ that if $I \not= \{0\}$ is an ideal of $R$ then $|R/I|$ is finite. We set
\[
N(I)=|R/I|,
\]
and call this the \emph{norm of I} (we defined the norm of an ideal in this way already in section 11 of Chapter 3 for ideals in a quadratic number field). The norm function $N$ on nonzero ideals is multiplicative with respect to the ideal product, i.e., we have
\begin{prp}
If $I$ and $J$ are $($not necessarily distinct$)$ nonzero ideals of R then
\[
N(IJ)=N(I)N(J).
\]
\end{prp}

\emph{Proof}. Hecke [27], section 27, Theorem 79. $\hspace{7.3cm} \textrm{QED}$

The multiplicativity of the norm function on ideals will play a crucial role in the derivation of  a very important product expansion formula for the zeta function that will  be done in the next section.

Now, let
\[
\mathcal{I}= \ \textrm{the set of all nonzero ideals of $R$}.
\]
If $n \in [1, \infty)$, let
\[
Z(n)=|\{I \in \mathcal{I}: N(I) \leq n\}|.
\]
The following proposition states a very important fact about the parameters $Z(n)$!

\begin{prp}
$Z(n)< +\infty$, for all $n \in [1, \infty)$.
\end{prp}

As a result of Proposition 5.3, $(Z(1), Z(2), Z(3)\dots)$ is a sequence of positive integers whose behavior determines how the ideals of $R$ are distributed throughout $R$ in accordance with the cardinality of the quotient rings of $R$. Useful information about the behavior of this sequence can hence be converted into useful information about the distribution of the ideals in $R$, and, as we shall see shortly, the Ideal Distribution  Theorem gives very useful information about the behavior of this sequence. We turn now to the

\emph{Proof of Proposition} 5.3. 

Perhaps the most elegant way to verify Proposition 5.3 is to make use of the ideal class group of $R$. We defined this group in section 11 of Chapter 3, and for the benefit of the reader, we will recall how that goes. First declare  that the ideals $I$ and $J$ of $R$ are equivalent if there exist nonzero elements $\alpha$ and $\beta$ of $R$ such that $\alpha I=\beta J$. This defines an equivalence relation on the set of all ideals of $R$, and the corresponding equivalence classes are the ideal classes of $R$. If we let $[I]$ denote the ideal class which contains the ideal $I$ then we define a multiplication on the set of ideal classes by declaring that the product of $[I]$ and $[J]$ is $[IJ]$. It can be shown that when endowed with this product (which is well-defined), the ideal classes of $R$ form an abelian group, called the ideal-class group of $R$. It is easy to see that the set of all principal ideals of $R$ is an ideal class, called the principal class, and one can prove that the principal class is the identity element of the ideal-class group. The ideal-class group is always finite, and the order of the ideal-class group of $R$ is called the class number of $R$. 

We begin the proof of Proposition 5.3 by letting $C$ be an ideal class of $R$ and for each $n \in [1, \infty)$, letting $\mathcal{Z}_C(n)$ denote the set
\[
\{I \in C \cap \mathcal{I}: N(I) \leq n\}.\]
We claim that $|\mathcal{Z}_C(n)|$ is finite. In order to verify this, let $J$ be a fixed nonzeo ideal in $C^{-1}$ (the inverse of $C$ in the ideal-class group), and let $0 \not=\alpha \in J$. Then there is a unique ideal $I$ such that $\alpha R=IJ$, and since $[I]=C[IJ]=C[\alpha R]=C$, it follows that $I \in C\cap \mathcal{I}$. Moreover, the map $\alpha R \rightarrow I$ is a bijection of the set of all nonzero principal ideals contained in $J$ onto $C\cap \mathcal{I}$. Proposition 5.2 implies that
\[
N(\alpha R)=N(I)N(J),\]
hence
\[
N(I) \leq n\ \textrm{if and only if}\ N(\alpha R) \leq nN(J).\]
Hence there is a bijection of  $\mathcal{Z}_C(n)$ onto the set
\[
\mathcal{J}=\{ \{0\} \not= \alpha R \subseteq J: N(\alpha R) \leq nN(J)\},\]
and so it suffices to show that $\mathcal{J}$ is a finite set.

That $|\mathcal{J}|$ is finite will follow if we prove that there is only a finite number of principal ideals of $R$ whose norms do not exceed a fixed constant. Suppose that this latter statement is false, i.e., there are infinitely many elements $\alpha_1, \alpha_2,\dots$ of $R$ such that the principal ideals $\alpha_i R, i=1, 2,\dots$ are distinct and $(N(\alpha_1 R), N(\alpha_2 R),\dots)$ is a bounded sequence. As all of the numbers $N(\alpha_i R)$ are positive integers, we may suppose with no loss of generality that $N(\alpha_i R)$ all have the same value $z$.

We now wish to locate $z$ in each ideal $\alpha_i R$. Toward that end, use the Primitive Element Theorem (Hecke [27], section 19, Theorem 52) to find $\theta \in F$, of degree $n$ over $\mathbb{Q}$, such that for each element $\nu$ of $F$, there is a unique polynomial $f \in \mathbb{Q}[x]$ such that $\nu=f(\theta)$ and the degree of $f$ does not exceed $n-1$. For each $i$, we hence find $f_i \in \mathbb{Q}[x]$ of degree no larger than $n-1$ and for which $\alpha_i=f_i(\theta)$. If $\theta_1,\dots,\theta_n$, with $\theta_1=\theta$, are the roots of the minimal polynomial of $\theta$ over $\mathbb{Q}$, then one can show that 
\[
N(\alpha_i R)=\Big|\prod_{k=1}^{n} f_i(\theta_k) \Big| \]
(Hecke [27], section 27, Theorem 76). Moreover, the degree  $d_i$ of $\alpha_i$ over $\mathbb{Q}$ divides $n$ in $\mathbb{Z}$, and if  $\alpha_i^{(1)},\dots,\alpha_i^{(d_i)}$, with $\alpha_i^{(1)}=\alpha_i$, denote the roots of the minimal polynomial of $\alpha_i$ over $Q$, then the numbers on the list $f_i(\theta_k), k=1,\dots,n,$ are obtained by repeating each $\alpha_i^{(j)}\  n/d_i$ times (Hecke [27], section 19, Theorem 54). If $c_0$ denotes the constant term of the minimal polynomial of $\alpha_i$ over $Q$, it follows that
\[
\prod_{k=1}^{n} f_i(\theta_k)=\Big(\prod_{k=1}^{d_i} \alpha_i^{(k)} \Big)^{n/d_i}=((-1)^{d_i}c_0)^{n/d_i} \in \mathbb{Z}. \] 
Because  $f_i(\theta_k)$ is an algebraic integer for all $i$ and $k$, it hence follows that
\[
\frac{z}{\alpha_i}=\pm\prod_{k=2}^n f_i(\theta_k)\in \mathcal{R} \cap F=R,\]
whence $z \in \alpha_i R$, for all $i$.

If we now let $\{\beta_1,\dots,\beta_n\}$ be an integral basis of $R$ then the claim in the proof of Proposition 5.1$(i)$ shows that for each $i$ there exists $\gamma_i \in R$ and $z_{ij} \in [0, z-1], j=1,\dots,n,$ such that
\[
\alpha_i=z\gamma_i+\sum_1^n z_{ij} \beta_j.\]
Because $z \in \alpha_i R$, it follows that
\[
\alpha_i R=zR+\Big(\sum_1^n z_{ij} \beta_j \Big)R,\ \textrm{for all}\ i.\]
However, the sum $\sum_1^n z_{ij} \beta_j$ can have only finitely many values; we conclude that the ideals $\alpha_i R, i=1, 2,\dots$ cannot all be distinct, contrary to their choice.

We now have what we need to easily prove that $Z(n)$ is finite. Let $C_1,\dots,C_h$ denote the distinct ideal classes of $R$. The set of all the ideals of $R$ is the (pairwise disjoint) union of the $C_i$'s hence $\{I \in \mathcal{I}: N(I) \leq n \}$ is the union of $\mathcal{J}_{C_1}(n),\dots, \mathcal{J}_{C_h}(n)$. Because each set $\mathcal{J}_{C_i}(n)$ is finite, so therefore is $|\{I \in \mathcal{I}: N(I) \leq n \}|=Z(n)$. $\hspace{6.4cm}\ \textrm{QED}$

We can now state the main result of this section:
\begin{thm}
$($Dedekind's Ideal Distribution Theorem$)$.
 The limit
 \[
 \lim_{n \rightarrow \infty} \frac{Z(n)}{n}= \lambda \]
 exists, is positive, and its value is given by the formula
 \[
 \lambda=\frac{2^{r+1} \pi^e \rho}{w \sqrt {|d|}} h,\]
 where
 \begin{eqnarray*}
 d&=&\ \textnormal{discriminant of $F$},\\
 e&=&\ \textnormal{$\frac{1}{2}$(number of complex embeddings of $F$ over $\mathbb{Q}$)},\\
 h&=&\ \textnormal{class number of $R$},\\
 r&=& \ \textnormal{unital rank of $R$},\\
 \rho&=&\ \textnormal{regulator of $F$},\\
 w&=&\ \textnormal{order of the group of roots of unity in $R$}.
 \end{eqnarray*}
\end{thm}
Thus the number of nonzero ideals of $R$ whose norms do not exceed $n$ is asymptotic to $\lambda n$ as $n \rightarrow +\infty$. 

The establishment of Theorem 5.4 calls for several results from the theory of algebraic numbers whose exposition would take us too far from what we wish to do here, so we omit the proof and instead refer the interested reader to Hecke [27], section 42, Theorem 122. Although we will make no further use of them, readers who are also interested in the definition of the discriminant of $F$ and the regulator of $F$, should see, respectively, the definition on p.73 and the definition on p.116 of Hecke [27].  We will define the parameter $e$ and the unital rank of $R$ in the two paragraphs after the next one.  The integers $d, e, h, r, w,$ and the real number $\rho$ are fundamental parameters associated with $F$ which govern many aspects of the arithmetic and algebraic structure of $F$ and $R$; Theorem 5.4 is a remarkable example of how these parameters work in concert to do that.  

Although the parameters which are used in the formula for the value of the limit  $\lambda=\lim_{n \rightarrow \infty} Z(n)/n$ are rather complicated to define for an arbitrary algebraic number field, they are much simpler to describe for a quadratic number field, so as to gain a better idea of how they determine the asymptotic behavior of the sequence $Z(1),Z(2),\dots$, we will take a closer look at what they are for quadratic fields. Thus, let $\mathbb{Q}(\sqrt m)$ be the quadratic number field determined by the square-free integer $m\not=0$ or 1. As we pointed out in section 11 of Chapter 3, the discriminant of $\mathbb{Q}(\sqrt m)$ is either $m$ or $4m$, if $m$ is, or respectively, is not, congruent to 1 mod 4. 

In order to calculate the parameter $e$ in Theorem 5.4, one needs to consider the \emph{embeddings} of an algebraic number field, i.e., the ring isomorphisms of the field into the set of complex numbers which fixes each element of $\mathbb{Q}$. An embedding is said to be \emph{real} if its range is a subset of the real numbers, otherwise, the embedding is said to be \emph{complex}. It can be shown that the number of embeddings is equal to the degree of the field and that the number of complex embeddings is even, and so $e$ is well-defined in Theorem 5.4. It follows that the quadratic field $\mathbb{Q}(\sqrt m)$ has precisely two embeddings: one is the trivial embedding which maps each element of $\mathbb{Q}(\sqrt m)$ to itself, and the other is the mapping on $\mathbb{Q}(\sqrt m)$ induced by the algebraic conjugate of $\sqrt m$ which sends the element $q+r\sqrt m$ for $(q, r)\in \mathbb{Q}\times \mathbb{Q}$ to the element $q-r\sqrt m$. It follows that if $m>0$ then there are no complex embeddings of 
$\mathbb{Q}(\sqrt m)$ and if $m<0$ then there are exactly 2 complex embeddings. Thus if $m>0$ then $e=0$ and if $m<0$ then $e=1$. 

We have already defined the class number $h$, and so we turn next to the unital rank $r$. This parameter is determined by the structure of the group of units in an algebraic number field. It can be shown that the group of units in the ring of algebraic integers $R$ in the algebraic number field $F$ is isomorphic to the direct sum of the finite cyclic group of roots of unity that are contained in $F$ and a free abelian group of finite rank $r$ (Hecke [27], section 34, Theorem 100). The rank $r$ of this free-abelian summand is by definition the \emph{unital rank} of $R$. When we now let $F=\mathbb{Q}(\sqrt m)$, it can be shown that if $m<0$ then the group of units of $\mathbb{Q}(\sqrt m)$  has no free-abelian summand, and so $r=0$ in this case. On the other hand, if $m>0$ then there is a unit $\varpi$ of $R$ in the group of units $U(R)$ such that $U(R)=\{\pm\varpi^n: n\in \mathbb{Z}\}$. If $\varpi$ is chosen to exceed 1 then it is uniquely determined as a generator of $U(R)$ in this way and is called the \emph{fundamental unit of R}. It follows that when $m>0$, the group of units of $R$ is isomorphic to the direct sum of the cyclic group of order 2 and the free abelian group $\mathbb{Z}$, hence the unital rank $r$ is 1 in this case.

The regulator $\rho$ of an algebraic number field $F$ is also determined by the group of units of $R$ by means of a rather complicated formula that uses a  determinant that is calculated from a basis of the free-abelian summand of the group of units. For a quadratic number field $\mathbb{Q}(\sqrt m)$ with $m<0$, whose group of units has no free-abelian summand, the regulator is taken to be 1, and if $m>0$ then the regulator of $\mathbb{Q}(\sqrt m)$ turns out to be $\log \varpi$, where $\varpi$ is the fundamental unit of $R$.

If $m>0$ then the group of roots of unity in  $\mathbb{Q}(\sqrt m)$ is simply $\{-1, 1\}$, and so the order $w$ of the group of roots of unity is 2. If $m<0$ then it can be shown that $w$ is 2 when $m<-4$, it is 6 when $m=-3$, and it is 4 when $m=-1$.

Taking all of this information into account, we see that for the quadratic number field $\mathbb{Q}(\sqrt m)$, the conclusion of Theorem 5.4 can be stated as follows: if $m>0$ and $\varpi$ is the fundamental unit in $R=\mathcal{R}\cap \mathbb{Q}(\sqrt m)$ then
\[
\lim_{n \rightarrow \infty} \frac{Z(n)}{n}=\left\{\begin{array}{rl}\displaystyle{\frac{2\log \varpi}{\sqrt m}}h,& \textrm{if  $m\equiv 1$ mod 4,}\\
\vspace{0.01cm}\\
\displaystyle{\frac{\log \varpi}{\sqrt m}h},& \textrm{if $m\not \equiv 1 $ mod 4,}\\\end{array}\right.\]
and if $m<0$ then
\[
\lim_{n \rightarrow \infty} \frac{Z(n)}{n}=\left\{\begin{array}{ll}\displaystyle{\frac{2\pi}{w\sqrt {|m|}}}h,\ \textrm{if  $m\equiv 1$ mod 4,}\\
\vspace{0.01cm}\\
\displaystyle{\frac{\pi}{w\sqrt {|m|}}h},\ \textrm{if $m\not \equiv 1 $ mod 4,}\\\end{array}\right.\]
where $w$ is 2 when $m<-4$, 6 when $m=-3$, and 4 when $m=-1$. The Ideal Distribution Theorem for quadratic number fields is in fact due to Dirichlet; after a careful study of Dirichlet's result, Dedekind generalized it to arbitrary algebraic number fields.

\section{The Zeta Function of an Algebraic Number Field}

We are now in a position to define and study the zeta function. Let $F$ be an algebraic number field of degree $n$ and let $R$ denote the ring of algebraic integers in $F$, as before.  Consider next the set $\mathcal{I}$ of all nonzero ideals of $R$. It is a consequence of Proposition 5.3 that $\mathcal{I}$ is countable, and so if $s \in \textbf{C}$ then the formal series
\begin{equation*}
\sum_{I \in \mathcal{I}} \frac{1}{N(I)^s} \tag{$*$}\]
is defined, relative to some fixed enumeration of $\mathcal{I}$. As we shall see, the zeta function of $F$ will be defined by this series. However, in order to do that precisely and rigorously, a careful examination of the convergence of this series must be done first. That is what we will do next.

If we let
\[
L(n)=|\{I \in \mathcal{I}: N(I) =n\}|,\ n\in [1, \infty),\]
then by formal rearrangement of its terms, we can write the series $(*)$ as
\begin{equation*}
\sum_{n=1}^{\infty}\frac{L(n)}{n^s}. \tag{$**$}\]
The series $(**)$ is a \emph{Dirichlet series}, i.e., a series of the form
\[
\sum_{n=1}^{\infty}\frac{a_n}{n^s},\]
where $(a_n)$ is a given sequence of complex numbers. The $L$-function of a Dirichlet character is another very important example of a Dirichlet series. 

We will determine the convergence of the series $(*)$ by studying the convergence of the Dirichlet series $(**)$. This will be done by way of the following proposition, which describes how a Dirichlet series converges.
\begin{prp}
Let $(a_n)$ be  sequence of complex numbers, let
\[
S(n)=\sum_{k=1}^n a_k,\]
and suppose that there exits $\sigma \geq 0, C>0$ such that
\[
\Big|\frac{S(n)}{n^{\sigma}} \Big| \leq C,\ \textrm{for all $n$ sufficiently large}.\]
Then the Dirichlet series 
\[
\sum_{n=1}^{\infty} \frac{a_n}{n^s}\] 
converges in the half-plane $\textnormal{Re}\ s> \sigma$ and uniformly in each closed and bounded subset of this half-plane. Moreover, if
\[
\lim_{n \rightarrow \infty} \frac{S(n)}{n}=d\]
then
\[
\lim_{s \rightarrow 1^+} (s-1)\sum_{n=1}^{\infty} \frac{a_n}{n^s} =d.\]
\end{prp}

\emph{Proof}(according to Hecke [27], section 42, Lemmas (a), (b), (c)). Let $m$ and $h$ be integers, with $m>0$ and $h \geq 0$, and let $K \subseteq \{ s: \textrm{Re}\ s> \sigma \}$ be a compact (closed and bounded) set. Then
\begin{eqnarray*}
\sum_{n=m}^{m+h} \frac{a_n}{n^s}&=&\sum_{n=m}^{m+h} \frac{S(n)-S(n-1)}{n^s}\\
&=&\frac{S(m+h)}{(m+h)^s}-\frac{S(m-1)}{m^s}+\sum_{n=m}^{m+h-1} S(n)\Big(\frac{1}{n^s}-\frac{1}{(n+1)^s} \Big)\\
&=&\frac{S(m+h)}{(m+h)^s}-\frac{S(m-1)}{m^s}+s\sum_{n=m}^{m+h-1} S(n)\int_n^{n+1}\frac{dx}{x^{s+1}}.
\end{eqnarray*}
If we now use the stipulated bound on the quotients $S(n)/n^{\sigma}$, it follows that
\begin{eqnarray*}
\Big|\sum_{n=m}^{m+h} \frac{a_n}{n^s} \Big|&\leq&\frac{2C}{m^{\textrm{Re}\ s-\sigma}}+C|s|\int_m^{\infty} \frac{dx}{x^{\textrm{Re}\ s- \sigma+1}}\\
&=&\frac{2C}{m^{\textrm{Re}\ s-\sigma}}+\frac{C|s|}{\textrm{Re}\ s-\sigma}\ \frac{1}{m^{\textrm{Re}\ s-\sigma}}.
\end{eqnarray*}
Because $K$ is a compact subset of Re $s>\sigma$, it is bounded and lies at a positive distance $\delta$ from Re $s=\sigma$, i.e., there is a positive constant $C^{\prime}$ such that
\[
\textrm{Re}\ s-\sigma \geq \delta\ \textrm{and}\ |s| \leq C^{\prime},\ \textrm{for all}\ s \in K.\]
Hence there is a positive constant $C^{\prime \prime}$, independent of $m$ and $h$, such that
\[
\Big|\sum_{n=m}^{m+h} \frac{a_n}{n^s} \Big| \leq C^{\prime \prime}\Big(1+\frac{1}{\delta} \Big)\frac{1}{ m^{\delta}}\ ,\ \textrm{for all}\ s \in K.\]
As $m$ and $h$ are chosen arbitrarily and $\delta$ depends on neither $m$ nor $h$, this estimate implies that the Dirichlet series converges uniformly on $K$, and as $K$ is also chosen arbitrarily, it follows that the series converges to a function continuous in Re $s>\sigma$.

We now assume that
\[
\lim_{n \rightarrow \infty} \frac{S(n)}{n}=d;\]
we wish to verify that
\[
\lim_{s \rightarrow 1^+} (s-1)\sum_{n=1}^{\infty} \frac{a_n}{n^s} =d.\]

From what we have just shown, it follows that the Dirichlet series now converges for $s>1$. Let
\[
S(n)=dn+\varepsilon_nn,\ \textrm{where}\ \lim_{n \rightarrow \infty} \varepsilon_n=0,\]
\[
\varphi(s)=\sum_{n=1}^{\infty} \frac{a_n}{n^s},\ s>1.\]
Then for $s>1$, we have that
\begin{eqnarray*}
|\varphi(s)-d\zeta(s)|&=&s\Big|\sum_{n=1}^{\infty} n\varepsilon_n\int_n^{n+1} \frac{dx}{x^{s+1}} \Big|\\
&<&s\sum_{n=1}^{\infty} |\varepsilon_n| \int_n^{n+1} \frac{dx}{x^s}.
\end{eqnarray*}

Let $\epsilon>0$, and choose an integer $N$ and a positive constant $A$ such that $|\varepsilon_n|< \epsilon$, for all $n \geq N$, and  $|\varepsilon_n| \leq A$, for all $n$. Then
\begin{eqnarray*}
|(s-1)\varphi(s)-d(s-1)\zeta(s)|&<&As(s-1)\sum_{n=1}^{N-1} \int_n^{n+1} \frac{dx}{x}+\epsilon s(s-1)+\sum_{n=N}^{\infty} \int_n^{n+1} \frac{dx}{x^s}\\
&=&As(s-1)\log N+\epsilon s(s-1) \int_N^{\infty} \frac{dx}{x^s}.
\end{eqnarray*}
Because the last expression has limit $\epsilon$ as $s \rightarrow 1$, it follows that
\[
\lim_{s \rightarrow 1^+} \big((s-1)\varphi(s)-d(s-1)\zeta(s)\big)=0.\]
We now claim that
\[
\lim_{s \rightarrow 1^+} (s-1)\zeta(s)=1;
\]
if this is so, then
\[
\lim_{s \rightarrow 1^+} (s-1)\varphi(s)=d,\]
as desired. This claim can be verified upon noting that
\[
\int_n^{n+1} \frac{dx}{x^s}<\frac{1}{n^s}<\int_{n-1}^n \frac{dx}{x^s}\ ,\ \textrm{for all}\ n \in [2, \infty)\ \textrm{and for all}\ s>1.\]
Hence
\[
\frac{1}{s-1}=\int_1^{\infty} \frac{dx}{x^s}<\sum_{n=1}^{\infty} \frac{1}{n^s}= \zeta(s)<1+\int_1^{\infty} \frac{dx}{x^s}=\frac{s}{s-1},\]
and so
\[
1<(s-1)\zeta(s)<s,\ \textrm{for all}\ s>1,\]
from which the claim follows immediately. $\hspace{8cm} \textrm{QED}$

Because each function $a_n/n^s$ is an entire function of $s$, a Dirichlet series which satisfies the hypotheses of Proposition 5.5 is a series of functions each term of which is analytic in Re $s> \sigma$ and which also  converges uniformly on every compact subset of Re $s>\sigma$. Hence the sum of the series is analytic in  Re $s>\sigma$.

We wish to apply Proposition 5.5 to the series $(**)$, and so we must study the behavior of the sequence 
\[
Z(n)=\sum_{k=1}^n L(k).\]
It is here that we make use of Theorem 5.4; it follows from that theorem that there is a positive constant $\lambda$ such that
 \[
 \lim_{n \rightarrow \infty} \frac{Z(n)}{n}= \lambda, \]
whence the sequence $(Z(n)/n)_{n=1}^{\infty}$ is bounded. Therefore the hypotheses of Proposition 5.5 are satisfied for $a_n=L(n)$ with $\sigma=1$, hence the series $(**)$ converges to a function analytic in Re $ s>1$.

We now let $s>1$. Because $L(n) \geq 0$ for all $n$, the convergence of $(**)$ is absolute for  $s>1$, hence we can rearrange the terms of $(**)$ in any order without changing its value. It follows that the value of the series 
\[
\sum_{I \in \mathcal{I}} \frac{1}{N(I)^s}\]
for $s>1$ is finite, is independent of the enumeration of $\mathcal{I}$ used to define the series, and is given by the value of the Dirichlet series $(**)$. 
\vspace{0.3cm}

\emph{Definition}. The \emph{$($Dedekind-Dirichlet$)$ zeta function of F} is the function $\zeta_F(s)$ defined for $s>1$ by
\[
\zeta_F(s)=\sum_{I \in \mathcal{I}} \frac{1}{N(I)^s}.\]

\emph{Remark}. One can show without difficulty that if $\sum_n a_n/n^s$ is a Dirichlet series which satisfies the hypotheses of Proposition 5.5 then $\sum_n a_n/n^s$ converges absolutely in Re $s>1+ \sigma$. If we apply this fact to the series $(**)$, it follows that $(**)$ converges absolutely in  Re $s> 2$. Hence the value of the series 
\[
\sum_{I \in \mathcal{I}} \frac{1}{N(I)^s}\]
for  Re $s>2$ is finite, is independent of the enumeration of $\mathcal{I}$ used to define the series, and is given by the value of the series $(**)$. Although we will make no use of this fact, it follows that the zeta function of $F$ can be defined by the series $(**)$ not only for $s>1$, but also for Re $s>1$, and when so defined, is analytic in that half-plane.

For emphasis, we record in the following proposition the observation that we made about the value of the zeta function of $F$ in the paragraph which immediately preceded its definition:
\begin{prp}
If
\[
L(n)=|\{I \in \mathcal{I}: N(I) =n\}|,\ n\in [1, \infty),\]
then
\[
\zeta_F(s)=\sum_{n=1}^{\infty}\frac{L(n)}{n^s}.\]
\end{prp}

\vspace{0.3cm}
For future reference, we also observe that Proposition 5.5 and Theorem 5.4 imply 
\begin{lem}
If $\zeta_F(s)$ is the zeta function of F and $\lambda$ is the positive constant in the conclusion of Theorem $5.4$ then
\[
\lim_{s \rightarrow 1^+} (s-1) \zeta_F(s)= \lambda.\]
\end{lem}

If $F=\mathbb{Q}$ then $R= \mathcal{R} \cap \mathbb{Q}=\mathbb{Z}$, hence the nonzero ideals of $R$ in this case are the principal ideals $n\mathbb{Z}, n \in [1, \infty)$. Then 
\[
N(n\mathbb{Z})=|\mathbb{Z}/n\mathbb{Z}|=n,\]
and so 
\[
\{I\in \mathcal{I}: |N(I)|=n \}= \{n\mathbb{Z} \}.\]
Hence the zeta function of $\mathbb{Q}$ is
\[
\zeta_{\mathbb{Q}}(s)=\sum_{n=1}^{\infty} \frac{1}{n^s},\]
the Riemann zeta function.

The next theorem gives a product formula for $\zeta_F(s)$ that is reminiscent of the product formula for the Dirichlet $L$-function of a Dirichlet character that we pointed out in section 4 of Chapter 4. It is a very useful tool for analyzing certain features of the behavior of $\zeta_F(s)$ and will play a key role in our proof of Theorem 4.12.
\begin{thm}
$($Euler-Dedekind product formula for $\zeta_F$$)$ Let $\mathcal{Q}$ denote the set of all prime ideals of R. Then
\begin{equation*}
\zeta_F(s)=\prod_{I \in \mathcal{Q}}\ \frac{1}{1-N(I)^{-s}} \ ,\ s>1. \tag{1}
\end{equation*}
\end{thm}

\emph{Proof}. Note that because a prime ideal $I$ of $R$ is proper, $N(I)>1$, and so each term of this product is defined for $s>1$. In order to prove the theorem we will need some standard facts about the convergence of infinite products, which we record in the following definitions and Proposition 5.9.
\vspace{0.3cm}

\emph{Definitions}. Let $(a_n)$ be a sequence of complex numbers such that $a_n \not= -1$, for all $n$. The infinite product
\[
\prod_1^{\infty}\ (1+a_n)\]
\emph{converges} if 
\[
\lim_{n \rightarrow \infty}\prod_1^n\ (1+a_k)\]
exists and is finite, and it \emph{converges absolutely} if 
\[\prod_1^{\infty}\ (1+|a_n|)\]
converges.

\vspace{0.3cm}
\begin{prp}
 $(i)$ $\prod_n (1+a_n)$ converges absolutely if and only if the series $\sum_n |a_n|$ converges.

$(ii)$ The limit of an absolutely convergent infinite product is not changed by any rearrangement of the factors. 
\end{prp}

\emph{Proof}. See Nevanlinna and Paatero [42], Sections 13.1, 13.2. $\hspace{4.2cm} \textrm{QED}   $

Returning to the proof of Theorem 5.8, we next consider the product on the right-hand side of (1). Because $N(I) \geq 2$ for all $I \in \mathcal{Q}$ it follows that for $s>1$,
\[
0<\frac{1}{1-N(I)^{-s}}-1=\frac{N(I)^{-s}}{1-N(I)^{-s}} \leq 2N(I)^{-s},\]
hence
\[
\sum_{I \in \mathcal{Q}} \Big(\frac{1}{1-N(I)^{-s}}-1 \Big) \leq 2\sum_{I \in \mathcal{Q}} N(I)^{-s}< +\infty \]
and so by Proposition 5.9, the product on the right-hand side of (1) converges absolutely for $s>1$ and its value is independent of the order of the factors.

The next step is to prove that this product converges to $\zeta_F(s)$ for $s>1$. Let
\[
\Pi(x)=\prod_{I \in \mathcal{Q}: N(I) \leq x}\ \frac{1}{1-N(I)^{-s}};\]
this product has only a finite number of factors by Proposition 5.3 and 
\[
\lim_{x \rightarrow +\infty} \Pi(x)=\prod_{I \in \mathcal{Q}}\ \frac{1}{1-N(I)^{-s}}.\]
We have that
\[
 \frac{1}{1-N(I)^{-s}}=\sum_{n=0}^{\infty} \frac{1}{N(I)^{ns}},\]
 hence $\Pi(x)$ is a finite product of absolutely convergent series, which we can hence multiply together and, in the resulting sum, rearrange terms in any order without altering the value of the sum. Proposition 5.2 implies that each term of this sum is either 1 or of the form
 \[
 N(I_1^{\alpha_1}\cdots I_r^{\alpha_r})^{-s},\]
 where $(\alpha_1,\dots,\alpha_r)$ is an $r$-tuple of positive integers, $I_i$ is a prime ideal for which $N(I_i)\leq x, i=1,\dots,r$, and all products of powers of prime ideals $I$ with $N(I) \leq x$ of this form occur exactly once. Hence
 \[
 \Pi(x)=1+\sum \frac{1}{N(I)^s},\]
 where the sum here is taken over all ideals $I$ of $R$ such that all prime ideal factors of $I$ have norm no greater than $x$. Now the Fundamental Theorem of Ideal Theory (Theorem 3.16) implies that all nonzero ideals of $R$ have a unique prime ideal factorization, hence 
 \[
 \zeta_F(s)-\Pi(x)=\sum \frac{1}{N(I)^s},\]
 where the sum here is taken over all ideals $I \not=\{0\}$ of $R$ such that at least one prime ideal factor of $I$ has norm greater than $x$. Hence this sum  does not exceed
 \[
 \sum_{n>x} \frac{L(n)}{n^s},\]
 and so
 \[
 \lim_{x \rightarrow +\infty} (\zeta_F(s)-\Pi(x))=\lim_{x \rightarrow +\infty} \sum_{n>x} \frac{L(n)}{n^s}=0.  \]
$\hspace{15.6cm} \textrm{QED}$

If $F=\mathbb{Q}$ then the prime ideals of $R=\mathbb{Z}$ are the principal ideals generated by the rational primes $q \in \mathbb{Z}$, and so it follows from
Theorem 5.8 that
\begin{equation*}
\zeta(s)=\prod_q\ \frac{1}{1-q^{-s}},\ s>1, \tag{2}
\] 
the Euler-product expansion of Riemann's zeta.

We are now going to use Theorem 5.8 to obtain a factorization of $\zeta_F$ over rational primes that is the analog of the product expansion (2) of the Riemann zeta function. In order to derive it,  we first recall from Proposition 5.1$(iii)$ and $(iv)$ that if $I$ is a prime ideal of $R$ then $I$ contains a unique rational prime $q$ and there is a unique positive integer $d$ such that $N(I)$ is $q^d$. The integer $d$ is called the \emph{degree of I} and we will denote it by deg $I$. We can now state and prove

\begin{thm}
If $\mathcal{Q}$ denotes the set of all prime ideals of R then the zeta function $ \zeta_F(s)$ of F has a product expansion given by
\begin{equation*}
\zeta_F(s)=\prod_{q\ \textnormal{a rational prime}} \Big(\prod_{I \in \mathcal{Q}: q \in I}\  \frac{1}{1-q^{-(\textnormal{deg}\ I)s}} \Big),\ s>1. \tag{3}
\]
\end{thm}

\emph{Proof}. If $n \in \mathbb{Z}$ then the ideal $nR$ is contained in a prime ideal of $R$ (Theorem 3.16) and so Proposition 5.1($iii$) implies that $\mathcal{Q}$ can be expressed as the pairwise disjoint union
\[
\bigcup_{q\ \textrm{a rational prime}} \{I \in \mathcal{Q}: q \in I\}.\]
Hence as a consequence of Theorem 5.8 and Proposition 5.9($ii$), we can rearrange the factors in (1) so as to derive the expansion (3) for $\zeta_F(s)$.$\hspace{6.8cm}\ \textrm{QED}$

The ideal $qR$ of $R$ is contained in only finitely many prime ideals (because of Theorem 3.16) and so each product inside the parentheses in (3) has only a finite number of factors; these finite products are called the \emph{elementary factors of} $\zeta_F$.

\section{The Zeta Function of a Quadratic Number Field} 

As has been the case frequently in much of our previous work, quadratic number fields provide interesting and important examples of various phenomena of great interest and importance in algebraic number theory, and zeta functions are no exception to this rule. In this section we will illustrate how the decomposition law for the rational primes in a quadratic number field, Proposition 3.17 from section 11 of chapter 3, and Theorem 5.10 can be used to derive a very useful product expansion for the zeta function of a quadratic number field. It is precisely this result that will be used to prove Theorem 4.12 in the next section.

For a square-free integer $m\not= 1$, let $F=\mathbb{Q}(\sqrt m), R=\mathcal{R} \cap F$. We recall for our convenience what the decomposition law for the rational primes in $R$ says. First, let $p$ be an odd prime. Then 

$(i)$ If $\chi_p(m)=1$ then $pR$ factors into the product of two distinct prime ideals, each of degree $1$.

$(ii)$ If $\chi_p(m)=0$ then $pR$ is the square of a prime ideal $I$, and the degree of $I$ is $1$.

$(iii)$ If $\chi_p(m)=-1$ then $pR$ is prime in $R$ of degree $2$.

\noindent The decomposition of the prime 2 in $R$ occurs as follows:

$(iv)$ If $m \equiv 1\ \textnormal{mod}\ 8$ then $2R$ factors into the product of two distinct prime ideals, each of degree $1$.

$(v)$ If $m \equiv$ 2 or 3 mod 4 then $2R$ is the square of a prime ideal $I$, and the degree of $I$ is 1.

$(vi)$ If $m\equiv$ 5 mod 8 then $2R$ is prime in $R$ of degree 2.

It follows from $(i)$-$(vi)$ that if $p$ is an odd prime in $\mathbb{Z}$ then the corresponding elementary factor of $\zeta_F$ is
\[
\frac{1}{(1-p^{-s})^2}\ ,\   \textrm{if}\ \chi_p(m)=1,\]
\[
\frac{1}{1-p^{-s}}\ ,\ \textrm{if}\ \chi_p(m)=0,\]
\[
\frac{1}{1-p^{-2s}}\ ,\ \textrm{if}\ \chi_p(m)=-1,\]
and the elementary factor corresponding to 2 is
\[
 \frac{1}{(1-2^{-s})^2}\ ,\   \textrm{if $m \equiv 1$  mod 8,}\]
\[
\frac{1}{1-2^{-s}}\ ,\ \textrm{if $m\equiv 2\ \textrm{or}\ 3\ \textnormal{mod}\ 4$ } ,\]
\[
\frac{1}{1-2^{-2s}}\ ,\ \textrm{if}\ m \equiv 5\  \textrm{mod}\ 8. \]
Observe next that each of the elementary factors corresponding to $p$ can be expressed as
\[
\frac{1}{1-p^{-s}} \frac{1}{1-\chi_p(d)p^{-s}}.\]
Hence from the product expansion (2) of the Riemann zeta function and the product expansion (3)  of $\zeta_F(s)$ we deduce
\begin{prp}
The zeta function of $\mathbb{Q}(\sqrt m)$ has the product expansion
\begin{equation*}
\zeta_{\mathbb{Q}(\sqrt m)}(s)=\theta(s)\zeta(s) \prod_p\ \frac{1}{1-\chi_p(m)p^{-s}},\ s>1, \tag{4}\]
where
\[
\theta(s)=\left\{\begin{array}{rl} \displaystyle \frac{1}{1-2^{-s}}\ ,& \textnormal{if $m \equiv 1$ \textrm{ mod 8},}\\
1\ ,& \textnormal{if $m\equiv 2\ \textnormal{or}\ 3\ \textnormal{mod}\ 4$ ,}\\ 
\displaystyle \frac{1}{1+2^{-s}}\ ,& \textnormal{if $m\equiv 5\ \textrm {mod }\ 8$.}\\\end{array}\right. 
\]
\end{prp}
We will use this factorization of $\zeta_{\mathbb{Q}(\sqrt m)}(s)$ to prove, in due course, the following lemma, the crucial fact that we will need to prove Theorem 4.12.
\begin{lem}
If $a \in \mathbb{Z}$ is not a square then
\[
\sum_p\ \chi_p(a)p^{-s}\]
remains bounded as $s \rightarrow 1^+$.
\end{lem}
 Note that Lemma 5.12 is very similar in form and spirit to the hypothesis of Lemma 4.7, which was a key step in Dirichlet's proof of Theorem 4.5. We will eventually see that this is no accident!

 \section{Proof of Theorem 4.12 and Related Results}

 We now have assembled all of the ingredients necessary for a proof of Theorem 4.12. As we have already verified the ``only if" implication in Theorem 4.12, we hence let $S$ be a nonempty finite subset of $[1, \infty)$ and suppose that for each subset $T$ of $S$ such that $|T|$ is odd, 
\[
\prod_{i \in T} i\ \textrm{is not a square}.\]
Let
\[
X=\{p: \chi_p \equiv -1\ \textrm{on}\ S\}.\]
We must prove that $X$ has infinite cardinality.

Consider the sum
\begin{equation*}
\Sigma(s)=\sum_{(p)} \Big(\prod_{i \in S} \big(1-\chi_p(i) \big) \Big)\cdot \frac{1}{p^s},\ s>1, \tag{5}
\]
where $(p)$ means that the summation is over all primes $p$ such that $p$ divides no element of $S$. Then
\[
\Sigma(s)=2^{|S|}\sum_{p \in X}\ \frac{1}{p^s},\ s>1, \]
hence  if we can show that
\begin{equation*}
\lim_{s \rightarrow 1^{+}} \Sigma(s)= +\infty, \tag{6}
\]
then the cardinality of $X$ will be infinite.

In order to get (6), we first calculate that
 \[
 \prod_{i \in S} \big(1-\chi_p(i) \big)=1+\sum_{\emptyset \not= T \subseteq S} (-1)^{|T|} \chi_p\Big(\prod_{i \in T} i \Big),\]
 substitute this into (5) and interchange the order of summation to obtain 
 \[
 \Sigma(s)=\sum_{(p)} \frac{1}{p^s}+\sum_{\emptyset \not= T \subseteq S} (-1)^{|T|} \Big( \sum_{(p)} \chi_p \Big(\prod_{i \in T} i \Big) \cdot \frac{1}{p^s} \Big).\]
Now divide $\{T: \emptyset \not= T \subseteq S \}$ into $U \cup V \cup W$, where 
\[
U=\Big\{\emptyset \not= T \subseteq S: |T|\ \textrm{is even and $\prod_{i \in T} i$ is a square} \Big\},
\]
\[
V=\Big\{\emptyset \not= T \subseteq S: |T|\ \textrm{is even and $\prod_{i \in T} i$ is not a square} \Big\}, 
\]
\[
W= \{ T \subseteq S: |T|\ \textrm{is odd} \}.\]
Then
\begin{eqnarray*}
\Sigma(s)&=&(1+|U|)\sum_{(p)} \frac{1}{p^s}\\
&+&\sum_{T \in V}\Big( \sum_{(p)} \chi_p \Big(\prod_{i \in T} i \Big) \cdot \frac{1}{p^s} \Big)\\
&-&\sum_{T \in W}\Big( \sum_{(p)} \chi_p \Big(\prod_{i \in T} i \Big) \cdot \frac{1}{p^s} \Big)\\
&=&\Sigma_1(s)+\Sigma_2(s)-\Sigma_3(s).
\end{eqnarray*}
Because the range of the summation here is over all but finitely many primes, Lemma 5.12, the definition of $V$ and the hypothesis on $S$ imply that $\Sigma_2(s)$ and $\Sigma_3(s)$ remain bounded as $s \rightarrow 1^{+}$, and so (6) will follow once we prove Lemma 5.12 and verify that
\begin{equation*}
\lim_{s \rightarrow 1^{+}}\sum_{(p)} \frac{1}{p^s}=+\infty. \tag{7}
\]

We check (7) first. Because the summation range in (7) is over all but finitely many primes, we need only show that
\begin{equation*}
\lim_{s \rightarrow 1^{+}}\sum_p \frac{1}{p^s}=+\infty. \tag{8}
\]
To see (8), recall from the proof of Proposition 5.5 that
\[
\lim_{s \rightarrow 1^{+}} (s-1) \zeta(s)=1,\]
hence
\begin{equation*}
 \lim_{s \rightarrow 1^{+}}\log \zeta(s)= \lim_{s \rightarrow 1^{+}} \log \frac{1}{s-1}+ \lim_{s \rightarrow 1^{+}} \log (s-1) \zeta(s)=+\infty. \tag{9}
 \]
 
 Now let $s>1$. The mean value theorem implies that
 \[
 |\log(1+x)| \leq 2|x|\ \textrm{for}\ |x| \leq \frac{1}{2},\]
 and so
 \[
 |\log (1-q^{-s})| \leq 2q^{-s},\ \textrm{for all}\ q \in P.\]
 Because $\sum_q q^{-s}<\sum_{n=1}^{\infty} n^{-s}< \infty$ it follows that the series
 \[
 \sum_q \log(1-q^{-s})\]
 is absolutely convergent. Hence
\begin{eqnarray*}
\log \zeta(s)&=&\log \Big(\prod_q \frac{1}{1-q^{-s}} \Big)\ (\textrm{from}\ (2))\\
&=&-\sum_q \log(1-q^{-s})\\
&=&\sum_q \frac{1}{q^{s}}+\sum_q \Big(-\log(1-q^{-s})-\frac{1}{q^{s}} \Big)\\
&=&\sum_q \frac{1}{q^{s}}+\sum_q \Big(\sum_{n \geq 2} \frac{1}{nq^{ns}} \Big),
\end{eqnarray*}
where we use the series expansion $\log(1-x)=-\sum_1^{\infty} x^n/n, |x|<1$, to obtain the last equation. Then
\begin{eqnarray*}
0<\sum_{n \geq 2} \frac{1}{nq^{ns}}&=&\frac{1}{q^{2s}} \Big(\sum_{n=0}^{\infty} \frac{1}{(n+2)q^{ns}} \Big)\\
&\leq&\frac{1}{q^{2s}}\sum_{n=0}^{\infty} q^{-ns}\\
&=&\frac{1}{q^{2s}} \frac{1}{1-q^{-s}}\\
&<&\frac{2}{q^2},\ \textrm{for all $q \geq 2$ and for all $s \geq 1$}.
\end{eqnarray*}
and so
\[
0<\sum_q \Big(\sum_{n \geq 2} \frac{1}{nq^{ns}} \Big)<2\sum_q \frac{1}{q^2}<+\infty\ \textrm{for all $s \geq1$.}
\]
It follows that
\[
\sum_q \frac{1}{q^s}= \log \zeta(s)+H(s),\ \textrm{$H(s)$ bounded on $s>1$},\]
hence this equation and (9) imply (8).

It remains only to prove Lemma 5.12. Let $d\not=1$ be a square-free integer. Then it is a consequence of the factorization (4) of $\zeta_F, F=\mathbb{Q}(\sqrt d)$ in Proposition 5.11 that
\[
\zeta_F(s)=\theta(s)\zeta(s) L(s),\ \textrm{where}\ L(s)=\prod_p \frac{1}{1-\chi_p(d)p^{-s}}.\]
By virtue of Lemma 5.7, 
\[
\lim_{s \rightarrow 1^{+}} (s-1)\zeta_F(s)= \lambda>0,\]
hence
\begin{eqnarray*}
\lim_{s \rightarrow 1^{+}} L(s)&=& \lim_{s \rightarrow 1^{+}} \frac{1}{\theta(s)} \frac{(s-1) \zeta_F(s)}{(s-1) \zeta(s)}\\
&=& \frac{\lambda}{\theta(1)}>0,
\end{eqnarray*}
and so
\begin{equation*}
\lim_{s \rightarrow 1^{+}} \log L(s)\ \textrm{is finite}. \tag{10}
\]
Now let $s>1$. Then
\begin{eqnarray*}
(11)\hspace{1cm} \ \log L(s)&=&-\sum_p \log(1-\chi_p(d)p^{-s})\\
&=&\sum_p \sum_{n=1}^{\infty} \frac{\chi_p(d)^n}{np^{ns}}\\
&=&\sum_p \chi_p(d)p^{-s}+\sum_p \sum_{n=2}^{\infty} \frac{\chi_p(d)^n}{np^{ns}}.
\end{eqnarray*}
Because
\[
\left|\sum_p \sum_{n=2}^{\infty} \frac{\chi_p(d)^n}{np^{ns}}\right| \leq \sum_p \sum_{n \geq 2} \frac{1}{np^{ns}},\]
the second term on the right-hand side of the last equation in (11) can be estimated as before to verify that it is bounded on $s>1$. Hence (10) and (11) imply that
\begin{equation*}
\sum_p \chi_p(d)p^{-s}\ \textrm{is bounded as $s \rightarrow 1^{+}$}. \tag{12}
\]
The integer $d$ here can be any integer $\not=1$ that is square-free, but every integer is the product of a square and a square-free integer, hence (12) remains valid if $d$ is replaced by any integer which is not a square.$\hspace{11.6cm} \textrm{QED}$ 

The technique used in the proof of Theorem 4.12 can also be used to obtain an interesting generalization of Basic Lemma 4.4 which answers the following question: if $S$ is a nonempty, finite subset of $[1, \infty)$ and $\varepsilon: S \rightarrow \{-1, 1\}$ is a given function, when does there exist infinitely many primes $p$ such that $\chi_p \equiv \varepsilon$ on $S$? There is a natural obstruction to $S$ having this property very similar to the obstruction that prevents the conclusion of Theorem 4.12 from being true for $S$. Suppose that there exists a subset $T \not= \emptyset$ of $S$ such that $\prod_{i \in T} i$ is a square. If we choose $i_0 \in T$ and define
\[
\varepsilon(i)=\left\{\begin{array}{rl} -1,& \textrm{if $i=i_0$,}\\
1,& \textrm{if $i \in S \setminus \{i_0\}$,} 
\\\end{array}\right. 
\]
then $\chi_p \not \equiv \varepsilon$ on $S$ for all sufficiently large $p$: otherwise there exits a $p$ exceeding all prime factors of the elements of $T$ such that 
\[
-1=\prod_{i \in T} \varepsilon(i)= \chi_p\Big(\prod_{i \in T} i \Big)=1.\]
By tweaking the proof of Theorem 4.12, we will show that this is the only obstruction to $S$ having this property.
\begin{thm}
Let $S$ be a nonempty finite subset of $[1, \infty)$. The following statements are equivalent: 

$(i)$ The product of all the elements in each nonempty subset of $S$ is not a square; 

$(ii)$ If  $\varepsilon: S \rightarrow \{-1, 1\}$ is a fixed but arbitrary function, then there exist infinitely many primes $p$ such that $\chi_p \equiv \varepsilon$ on $S$.
\end{thm}

\emph{Proof}. We have already observed that $(i)$ follows from $(ii)$, hence suppose that $S$ satisfies $(i)$ and let $\varepsilon: S \rightarrow \{-1, 1\}$ be a fixed function. Consider the sum\[
\Sigma_{\varepsilon}(s)=\sum_{(p)} \Big(\prod_{i \in S} \big(1+\varepsilon(i)\chi_p(i) \big) \Big)\cdot \frac{1}{p^s},\ s>1.\]
If 
\[
X_{\varepsilon}=\{p: \chi_p \equiv \varepsilon\ \textrm{on}\ S\} \]
then
\[
\Sigma_{\varepsilon}(s)=2^{|S|}\sum_{p \in X_{\varepsilon} }\ \frac{1}{p^s}\ .\]
Also,
\[
 \Sigma_{\varepsilon}(s)=\sum_{(p)} \frac{1}{p^s}+\sum_{\emptyset \not= T \subseteq S}\ \prod_{i \in T} \varepsilon(i)\ \Big( \sum_{(p)} \chi_p \Big(\prod_{i \in T} i \Big) \cdot \frac{1}{p^s} \Big).\]
Lemma 5.12 and the hypotheses on $S$ imply that the second term on the right-hand side of this equation is bounded as $s \rightarrow 1^{+}$, hence from (7) we conclude that
\[
\lim_{s \rightarrow 1^{+}} \Sigma_{\varepsilon}(s)=+\infty,\]
and so $X_{\varepsilon}$ is infinite. $\hspace{11.9cm} \textrm{QED}$
\vspace{0.0cm} 

\emph{Definition}. Any set $S$ satisfying statement $(ii)$ of Theorem 5.13 will be said to \emph{support all patterns}.

\vspace{0.3cm}
\emph{Remark}. The proof of Theorems 4.12 and 5.13 follows exactly the same strategy as Dirichlet's proof of Theorem 4.5. One wants to show that a set $X$ of primes with a certain property is infinite. Hence take $s>1$, attach a weight of $p^{-s}$ to each prime $p$ in $X$ and then attempt to prove that the weighted sum
\[
\sum_{p \in X}\ \frac{1}{p^s}
\]
of the elements of $X$ is unbounded as $s \rightarrow 1^{+}$. In order to achieve this (using ingenious methods!), one writes this weighted sum as $\sum_p\ 1/ p^s$ plus a term that is bounded as $s \rightarrow 1^{+}$. The similarity of all of these arguments is no accident; Theorem 5.13 is in fact also due to Dirichlet, and appeared in his great memoir [11], \emph{Recherches sur diverses applications de l'analyse infinit$\acute{\textrm{e}}$simal  $\grave{\textrm{a}}$ la th$\acute{\textrm{e}}$orie des nombres}, of 1839-40, which together with [10] founded modern analytic number theory. The proof of Theorem 5.13 given here is a variation on Dirichlet's original argument due to Hilbert [28], section 80, Theorem 111.

A straightforward modification of the proof of Theorem 4.9 can now be used to establish
\begin{thm}
 If  S is a nonempty, finite subset of $[1, \infty)$ such that  for all subsets $T$ of $S$ of odd cardinality, $\prod_{i \in T} i$ is not a square, $\mathcal{S}$ and $v: 2^S \rightarrow F^n$ are defined by $S$ as in the statement of Theorem $4.9$, and $d$ is the dimension of the linear span of $v(\mathcal{S})$ in $F^n$, then the density of the set $\{p: \chi_p \equiv -1\ \textrm{on}\ S\}$ is $2^{-d}$.
\end{thm}

If $p<q<r<s$ are distinct primes and we let, for example, $S_1=\{p, pq, qr, rs\}$ and $S_2=\{p, ps, pqr, pqrs\}$, then it follows from Theorem 5.14 and the row reduction of the incidence matrices of $S_1$ and $S_2$ that we performed in section 6 of Chapter 4 that the density of $\{p: \chi_p \equiv -1\ \textrm{on}\ S_1\}$ is $2^{-4}$ and the density of
$\{p: \chi_p \equiv -1\ \textrm{on}\ S_2\}$ is $2^{-3}$. As we pointed out in section 6 of chapter 4, a 2-dimensional subspace of $F^4$ contains only 3 nonzero vectors, and so if $S$ is a set of 4 nontrivial square-free integers such that $S$ is supported on 4 primes then the density of  $\{p: \chi_p \equiv -1\ \textrm{on}\ S\}$ cannot be $2^{-2}$. But it is also true that all of the vectors in a 2-dimensional subspace of $F^4$ must sum to 0 and so if $S$ is a set of 3 nontrivial square-free integers such that $S$ is supported on 4 primes then 
$\{p: \chi_p \equiv -1\ \textrm{on}\ S\}$ is in fact empty. In order to get a set $S$ from $p$, $q$, $r$, and $s$ such that the density of $\{p: \chi_p \equiv -1\ \textrm{on}\ S\}$ is $2^{-2}$, $S$ has to have 2 elements, and it follows easily from Theorem 5.14 that $S=\{pq,qrs\}$ is one of many examples for which the density of $\{p: \chi_p \equiv -1\ \textrm{on}\ S\}$ is $2^{-2}$.

A straightforward modification of the proof of Lemma 4.10 can also be used to establish
\begin{thm}
$($Filaseta and Richman, $[18]$, Theorem $2$$)$ If  S is a nonempty, finite subset of $[1, \infty)$ such that the product of all the elements in each nonempty subset of $S$ is not a square and $\varepsilon: S \rightarrow \{-1, 1\}$ is a fixed but arbitrary function, then the density of the set $\{p: \chi_p \equiv \varepsilon\ \textrm{on}\ S\}$ is $2^{-|S|}$.
\end{thm}

\section{Proof of the Fundamental Theorem of Ideal Theory}

Because the Fundamental Theorem of Ideal Theory was used at its full strength in the proof of the Euler-Dedekind product expansion of the zeta function (Theorem 5.8), and also because of the important role that it played (although not at full strength) in the results on the factorization of ideals in a quadratic number field from Chapter 3, we will present a proof of it in this final section of Chapter 5. Our account follows the outline given by O. Ore in [43].

Let $F$ be an algebraic number field of degree $n$ and let $R$ be the ring of algebraic integers in $F$.  We want to prove that every nonzero proper ideal of $R$ is a product of a  finite number of prime ideals and also that this factorization is unique up to the order of the prime-ideal factors. The strategy of our argument is to prove first that each nonzero proper ideal of $R$ \emph{contains} a finite product of prime ideals. We hence chose for each nonzero proper ideal $I$ a product of prime ideals with the smallest number of factors that is contained in $I$, and then by use of appropriate mathematical technology that we will develop, proceed by induction on this smallest number of prime-ideal factors to prove that $I$ is in fact \emph{equal} to a product of prime ideals.  Uniqueness will then follow by further use of the mathematical technology that we have at our disposal. We proceed to implement this strategy.

Let $I$ be an ideal of $R$, $\{0\}\not= I \not= R$. 
\begin{lem}
There exists a sequence of prime ideals $P_1,\dots,P_s$ of $R$ such that $I\subseteq P_i$, for all $i$ and $P_1\cdots P_s\subseteq I$.
\end{lem}

\emph{Proof}. If $I$ is prime then we are done, with $s=1$, hence suppose that $I$ is not prime. Then there exists a product $\beta \gamma$ of elements of $R$ which is in $I$ and $\beta\not \in I,\ \gamma \not \in I$. Let $\{\alpha_1,\dots,\alpha_n\} $ be an integral basis of $I$, and set
\[
J=(\alpha_1,\dots,\alpha_n,\beta),\ K=(\alpha_1,\dots,\alpha_n,\gamma).\]
Then
\[
JK\subseteq I,\ I\subsetneqq J,\ I\subsetneqq K.\]
If $J, K$ are both prime then we are done, with $s=2$. Otherwise apply this procedure to each nonprime ideal that occurs, and continue in this way as long as the procedure produces nonprime ideals. Note that after each step of the procedure,

$(i)$ the product of all the ideals obtained in that step is contained in $I$,

$(ii)$ $I$ is contained in each ideal obtained in that step, and

$(iii)$ each ideal obtained in that step is properly contained in an ideal from the immediately preceding step.

\emph{Claim}: this procedure terminates after finitely many steps.

If this is true then each ideal obtained in the final step is prime; otherwise the procedure would continue by applying it to a nonprime ideal. If $P_1,\dots,P_s$ are the prime ideals obtained in the final step then this sequence of ideals satisfies Lemma 5.16 by virtue of $(i)$ and $(ii)$ above. 

\emph{Proof of the claim}. Suppose this is false. Then $(ii)$ and $(iii)$ above imply that the procedure produces an infinite sequence of ideals $J_0, J_1,\dots, J_n,\dots$ such that $J_0=I$ and $J_i\subsetneqq J_{i+1}$, for all $i$. We will now prove that $I$ is contained in only finitely many ideals, hence no such sequence of ideals is possible.

The proof of Proposition 5.1$(i)$ implies that $I$ contains a positive rational integer $a$. We show that $a$ belongs to only finitely many ideals.

Suppose that $J$ is an ideal, with integral basis $\{\beta_1,\dots,\beta_n\}$, and $a \in J$. Then we also have that
\[
J=(\beta_1,\dots,\beta_n, a).\]
By the claim in the proof of Proposition 5.1$(i)$, for each $i$, there is $\gamma_i, \delta_i\in R$ such that $\beta_i=a\gamma_i+\delta_i$, and $\delta_i$ can take on only at most $an$ values. But then
\[
J=(a\gamma_1+\delta_1,\dots,a\gamma_n+\delta_s)=(\delta_1,\dots,\delta_n, a).\]
Because each $\delta_i$ assumes at most $an$ values, it follows that $J$ is one of only at most $an^2$ ideals. $\hspace{14.2cm}\textrm{QED}$

The statement of the next lemma requires the following definition:

\vspace{0.4cm}
\textit{Definition}. If $J$ is an ideal of $R$ then 
\[
J^{-1}=\{\alpha \in F: \alpha \beta \in R,\ \textrm{ for all}\ \beta \in J\}.\]
\begin{lem}
If P is a prime ideal of R then $P^{-1}$ contains an element of $F\setminus R$.
\end{lem}

\emph{Proof}. Let $x\in P$. Lemma 5.16 implies that $(x)$ contains a product $P_1\cdots P_s$ of prime ideals. Choose a product with the smallest number $s$ of factors.

Suppose that $s=1$. Then $P_1\subseteq (x)\subseteq P$. $P_1$ maximal (Proposition 5.1$(i)$) implies that $P=P_1=(x)$. Hence $1/x\in P^{-1}$. Also, $1/x\not \in R$; otherwise, $1=x\cdot 1/x \in P$, contrary to the fact that $P$ is proper.

Suppose that $s>1$. Then $P_1\cdots P_s \subseteq (x) \subseteq P$, and so the fact that $P$  is prime implies that $P$ contains a $P_i$, say $P_1$. $P_1$ maximal implies that $P=P_1$. $P_2\cdots P_s \nsubseteq (x)$ by minimality of $s$, hence there exits $\alpha \in P_2\cdots P_s$ such that $\alpha \not \in (x)$, and so $\alpha/x\not \in R$.

\emph{Claim}: $\alpha/x \in P^{-1}$. 

Let $\beta \in P$. We must prove that $\beta (\alpha/x)\in R$. To do that, observe that
\[
(\alpha)P \subseteq P_2\cdots P_sP=P_1\cdots P_s \subseteq (x),\]
and so there is a $\gamma \in R$ such that $\alpha \beta=x\gamma$, i.e., $\beta(\alpha/x)=\gamma$.$\hspace{4.6cm}\textrm{QED}$

The next lemma is the key technical tool that allows us to prove the Fundamental Theorem of Ideal Theory; it will be used to factor an ideal into a product of prime ideals and to show that this factorization is unique up to the order of the factors. In order to state it, we need to extend the definition of products of ideals to products of  arbitrary subsets of $R$ like so:

\vspace{0.4cm}
\textit{Definition}. If $S$ and $T$ are subsets of $R$ then the \emph{product ST of S and T} is the set consisting of all sums of the form $\displaystyle{\sum_i s_it_i}$, where $(s_i, t_i)\in S \times T$ for all $i$.

\vspace{0.4cm}
This product is clearly commutative and associative, and it agrees with the product defined before when $S$ and $T$ are ideals of $R$.
\begin{lem}
If P is a prime ideal of R and I is an ideal of R then $P^{-1}PI=I$.
\end{lem}

\emph{Proof}. It suffices to show that $P^{-1}P=(1)$. It is straightforward to show that $J=P^{-1}P$ is an ideal of $R$. As $1\in P^{-1}$, it follows that $P\subseteq J$ and so $P$ maximal implies that $P=J$ or $J=(1)$.

Suppose that $J=P$. Let $\{\alpha_1,\dots,\alpha_n\}$ be an integral basis of $P$, and use Lemma 5.17 to find $\gamma \in P^{-1},\ \gamma \not \in R$. Then $\gamma \alpha_i \in P$, for all $i$, and so
\[
\gamma \alpha_i=\sum_j a_{ij}\alpha_j,\ \textrm{where $a_{ij} \in \mathbb{Z}$ for all $i$, $j$}.\]
As a consequence of these equations, $\gamma$ is an eigenvalue of the matrix $[a_{ij}]$, hence it is a root of the characteristic polynomial of $[a_{ij}]$, and this characteristic polynomial is a monic polynomial in $\mathbb{Z}[x]$. As we showed in the proof of  Theorem 3.11, this implies that $\gamma$ is an algebraic integer, contrary to its choice. Hence $P\not=J$, and so $J=(1)$.$\hspace{5cm}\textrm{QED}$

The Fundamental Theorem of Ideal Theory is now a consequence of the next two lemmas.
\begin{lem}
Every nonzero proper ideal of R is a product of prime ideals.
\end{lem}

\emph{Proof}. Lemma 5.16 implies that every nonzero proper ideal of $R$ contains a product $P_1\cdots P_r$ of prime ideals, where we choose a product with the smallest number $r$ of factors. The argument now proceeds by induction on $r$.

Let $\{0\}\not= I \not= R$ be an ideal with $r=1$, i.e., $I$ contains a prime ideal $P$. $P$ maximal implies that $I=P$, and we are done.

Assume now that $r>1$ and every nonzero, proper ideal that contains a product of fewer than $r$ prime ideals is a product of prime ideals. 

Let $\{0\}\not= I \not= R$ be an ideal that contains a product $P_1\cdots P_r$ of prime ideals, with $r$ the smallest number of prime ideals with this property. Lemma 5.16 implies that $I$ is contained in a prime ideal $Q$. Hence  $P_1\cdots P_r\subseteq Q$, and so $Q$ contains a $P_i$, say $P_1$.  $P_1$ maximal implies that $Q=P_1$. Hence $I \subseteq P_1$. Then $IP_1^{-1}$ is an ideal of $R$; $I\subseteq IP_1^{-1}$ ($1\in P^{-1}$), and so  $IP_1^{-1} \not= \{0\}$. $IP_1^{-1} \not=R$; otherwise, $P_1\subseteq I$, hence $I=P_1$, contrary to the fact that $r>1$. Lemma 5.18 implies that
\[
P_2\cdots P_r=P_1^{-1}P_1\cdots P_r \subseteq IP_1^{-1},\]
hence by the induction hypothesis, $IP_1^{-1}$ is a product $P_1^{\prime}\cdots P_k^{\prime}$ of prime ideals, and so by Lemma 5.18 again,
\[
I=(IP_1^{-1})P_1=P_1^{\prime}\cdots P_k^{\prime}P_1\]
is a product of prime ideals.$\hspace{10.7cm}\textrm{QED}$
\begin{lem}
Factorization as a product of prime ideals is unique up to the order of the factors.
\end{lem}
\emph{Proof}. Suppose that $P_1\cdots P_r=Q_1\cdots Q_s$ are products of prime ideals, with $r\leq s$, say. $Q_1\cdots Q_s\subseteq Q_1$, hence $P_1\cdots P_r\subseteq Q_1$  and so the fact that $Q_1$ is a prime ideal and the maximality of the $P_i$'s imply, after reindexing one of the $P_i$'s, that $Q_1=P_1$. Then Lemma 5.18 implies that
\[
P_2\cdots P_r=P_1^{-1}P_1\cdots P_r=Q_1^{-1}Q_1\cdots Q_s=Q_2\cdots Q_s.\]
Continuing in this way, we deduce, upon reindexing of the  $P_i$'s, that $P_i=Q_i$, $i=1,\dots r$, and also, if $r<s$, that
\[
(1)=Q_{r+1}\cdots Q_s.\]
But this equation implies that $R=(1)\subseteq Q_{r+1}$, which is impossible as $Q_{r+1}$ is a proper ideal. Hence $r=s$.$\hspace{12.1cm}\textrm{QED}$

Dedekind's own proof of The Fundamental Theorem of Ideal Theory in [8], Chapter 4, section 25, is a model of clarity and insight which amply repays careful study. We strongly encourage the reader to take a look at it. 
\newpage
\afterpage{\null\newpage}
\thispagestyle{empty}
\chapter{Elementary Proofs}

After providing in section 1 of this chapter some motivation for the use of elementary methods in number theory, we present proofs of Theorems 4.12 and 5.13 in sections 3 and 2, respectively, which employ only Lemma 4.4 from Chapter 4 and linear algebra over the Galois field of order 2, thereby avoiding the use of zeta functions.

\section{Whither Elementary Proofs in Number Theory?}
 Dirichlet's incorporation of transcendental methods into number theory allowed him to establish many deep and far-reaching results in that subject. The importance of Dirichlet's work led to a very strong desire to understand it in as many different ways as possible, and this desire naturally motivated the search for different ways to prove his results.  As the years passed, particular attention was focused on removing any use of methods which relied on mathematical analysis, replacing them instead by ideas and techniques which deal with or stem directly from the fundamental structure of the integers, as this was sometimes viewed as being more suitable for the development of the most important results of the theory. The viewpoint that the preferred methods in number theory should be based only on fundamental properties of the integers was originally held evidently by none other than Leonard Euler (see the remarks by Gauss in [19], article 50 about Euler's proof of Fermat's little theorem), and much of the fundamental contributions to number theory by Euler, Lagrange, Legendre, Gauss, Dedekind, and many others can be seen as evidence of the value of that philosophy. The subject of \emph{elementary number theory}, i.e., the practice of number theory using methods which have their basis in the algebra and/or the geometry of the integers, and which, in particular, avoid the use of any of the infinite processes coming from analysis, has thus attained major importance. Indeed, among the results of twentieth-century number theory which generated exceptional excitement and interest is the discovery by Selberg [51], [52] and Erd$\ddot{\textrm{o}}$s [14] in 1949 of the long-sought elementary proofs of the Prime Number Theorem, Dirichlet's theorem on primes in arithmetic progression, and the Prime Number Theorem for primes in arithmetic progression.

The philosophical spirit of elementary number theory resonates with particular force in the mind of anyone who compares the way that we proved Theorems 4.2 and 4.3 to the way that we proved Theorems 4.12 and 5.13. The proof of the former two results are easy consequences of Lemma 4.4, which in turn depends on an elegant application of quadratic reciprocity and Dirichlet's theorem on primes in arithmetic progression. In contrast to that line of reasoning, our proof of Theorems 4.12 and 5.13 requires, by comparison, a rather sophisticated application of transcendental methods based on the Riemann zeta function and the zeta function of a quadratic number field. Because all of these results are very similar in content, this raises a natural question: can we give elementary proofs of Theorems 4.12 and 5.13 which, in particular, avoid the use of zeta functions and are more in line with the ideas used in the proof of Theorems 4.2 and 4.3? The answer: yes we can, and that will be done in this chapter by proving Theorems 4.12, in section 3, and 5.13, in section 2, using only Lemma 4.4 and linear algebra over $GF(2)$. Taking into account the fact that Dirichlet's theorem and the Prime Number Theorem for primes in arithmetic progression also have elementary proofs, the proofs that we have given of Theorems 4.2, 4.3, 4.9, 5.14, and 5.15 are already elementary.
\section{An Elementary Proof of Theorem 5.13}

We begin with Theorem 5.13: let $S$ be a nonempty finite subset of $[1, \infty)$ such that
\begin{equation*}
\textrm{for all $\emptyset \not= T \subseteq S,\ \prod_{i \in T} i$ is not a square.} \tag{1} \]
We wish to prove that for each function $\varepsilon: S\rightarrow \{-1, 1\}$, the cardinality of $\{p: \chi_p \equiv \varepsilon\ \textrm{on}\ S\}$ is infinite.

The first step in our reasoning is to reduce to the case in which every integer in $S$ is square-free. Recall that the square-free part $\sigma(z)$ of $z \in [1, \infty)$ is 
\[
\sigma(z)=\prod_{q \in \pi_{\textrm{odd}}(z)} q,\]
and observe that if $\emptyset \not= T \subseteq [1, \infty)$ is finite then
\[
\textrm{$\prod_{i \in T} i$ is not a square if and only if $\prod_{i \in T}\ \sigma(i)$ is not a square}.\]
(There is an integer $n$ such that $\prod_{i \in T} i= \prod_{i \in T}\ \sigma(i) \times n^2$, so the multiplicity $m$ of a prime factor $q$ of  $\prod_{i \in T} i$ in  $\prod_{i \in T} i$ is congruent mod 2 to the multiplicity $m^{\prime}$ of $q$ in $ \prod_{i \in T}\ \sigma(i)$ hence $m$ is odd if and only if $m^{\prime}$ is odd.) Also 
\[
\chi_p(z)=\chi_p(\sigma(z)),\ \textrm{ for all $p \notin \pi(z)$}.\]
Hence, upon replacing $S$ by the set formed from the integers $\sigma(z)$ for $z \in S$, we may suppose with no loss of generality that all elements of $S$ are square-free. Hence
\[
z= \prod_{q \in \pi(z)} q,\ z \in S,\]
$\pi(z) \not= \emptyset$, for all $z \in S (1 \notin S)$, and if $\{w, z\} \subseteq S$ then $\pi(w) \not= \pi(z)$.

The next step is to look
 for a purely combinatorial condition on the sets $\pi(z), z \in S$, that is equivalent to condition (1). The following notation will be helpful with regard to that: if $T \subseteq S$, let
\[
\Pi(T)=\bigcup_{i \in T}\ \pi(i),\]
\[
\mathcal{S}(T)=\{ \pi(i): i \in T\},\]
\[
p(T)= \prod_{i \in T} i,\]
and let
\[
\Pi=\bigcup_{i \in S}\ \pi(i),\]
\[
\mathcal{S}=\{ \pi(i): i \in S\}.\]
Now
\[
\Pi(T)=\ \textrm{ the set of all prime factors of}\ p(T)\]
and
\[
\textrm{the multiplicity in $p(T)$ of $q \in \Pi(T)$}= |\{X \in \mathcal{S}(T): q \in X\}|.\]
Hence
\begin{equation*}
p(T)\ \textrm{is not a square iff}\ \{q \in \Pi(T): |\{X \in \mathcal{S}(T): q \in X\}|\ \textrm{is odd}\ \} \not= \emptyset. \tag{2}\]

Condition (2) can be elegantly expressed by using the symmetric difference operation on sets. Recall that if $A$ and $B$ are sets then the \emph{symmetric difference $A \triangle B$ of A and B} is the set $(A\setminus B) \cup (B \setminus A)$. The symmetric difference operation is commutative and associative, hence if $A_1,\dots,A_k$ are distinct sets then the repeated symmetric difference 
\[
\triangle_i\ A_i=A_1\triangle \cdots \triangle A_k\]
is unambiguously defined. In fact, one can show that
\begin{equation*}
\bigtriangleup_i\ A_i=\Big\{a \in \bigcup_i\ A_i: |\{A_j: a \in A_j \}|\ \textrm{is odd} \Big\}. \tag{3}\]
Statements (2) and (3) imply that
\[
p(T)\ \textrm{is not a square if and only if $\triangle_{i \in T}\ \pi(i) \not= \emptyset.$}\]
Hence
\[
\textrm{condition (1) holds if and only if for all  nonempty subsets $T$ of $S , \triangle_{i \in T}\ \pi(i) \not= \emptyset$.}\]
 As the map $i \rightarrow \pi(i)$ is a bijection of $S$ onto $\mathcal{S}$, it follows that
 \begin{equation*}
 \textrm{condition (1) holds if and only if for all nonempty subsets $\mathcal{T}$ of $\mathcal{S}, \triangle_{T \in \mathcal{T}}\ T \not= \emptyset$.} \tag{4} \]
 Statement (4) is the combinatorial formulation of condition (1) that we want.
 
In order to express things more concisely, we recall now from section 5 of Chapter 5 that $S$ is said to support all patterns if for each function $\varepsilon: S \rightarrow \{-1, 1\}$, the set $\{p: \chi_p \equiv \varepsilon \ \textrm{on}\ S\}$ is infinite. Consequently from (4), in order to prove Theorem 5.13, we must show that
 \begin{equation*}
 \textrm{if $\triangle_{T \in \mathcal{T}} T \not= \emptyset$ for all $\emptyset \not= \mathcal{T} \subseteq \mathcal{S}$ then $S$ supports all patterns.} \tag{5} \]
 Hence we next look for a combinatorial condition on $\mathcal{S}$ which guarantees that $S$ supports all patterns. This is provided by
 \begin{lem}
 Suppose that $\mathcal{S}$ satisfies the following condition:
 \begin{quote}
$(6)$  for each nonempty subset $\mathcal{T}$ of $\mathcal{S}$, there exists a subset $N$ of $\Pi$ such that 
 \end{quote}
\begin{center} 
$\mathcal{T}=\{S \in \mathcal{S}: |N \cap S|\ \textnormal{ is odd} \}  $. 
\end{center}
Then $S$ supports all patterns.
 \end{lem}

\emph{Proof}. Let $\varepsilon$ be a function of $S$ into $\{-1, 1\}$. We must prove: $\{p: \chi_p \equiv \varepsilon\ \textrm{on}\ S \}$ is infinite. 

The map $\pi(i) \rightarrow \varepsilon(i), i \in S$ defines a function $\varepsilon^{\prime}$ of $\mathcal{S}$ into $\{-1, 1\}$. Let
\[
\mathcal{T}=(\varepsilon^{\prime})^{-1}(-1).\]
If $\mathcal{T}= \emptyset$ then $\varepsilon \equiv 1$, hence apply Theorem 4.3. Suppose that $\mathcal{T} \not= \emptyset$, and then find $N \subseteq \Pi$ such that $N$ satisfies the conclusion of (6) for this $\mathcal{T}$. Basic Lemma 4.4 implies that there are infinitely many primes $p$ for which
\begin{equation*}\
\{q \in \Pi: \chi_p(q)= -1\}=N. \tag{7}\]
Let $p$ be any one of these primes which divides no element of $S$.

We claim that $\chi_p \equiv \varepsilon$ on $S$. To verify this, note first that because of (7),
\[
\chi_p(i)=(-1)^{|N \cap \pi(i)|},\ \textrm{for all}\ i \in S.\]
Hence
\[
i \in S \cap\chi_p^{-1}(-1)\ \textrm{if and only if}\ |N\cap \pi(i)|\ \textrm{is odd}.\]
Since the conclusion of (6) holds for $N$ and $\mathcal{T}$, it follows that
\[
|N\cap \pi(i)|\ \textrm{is odd if and only if}\ \pi(i) \in \mathcal{T},\ \textrm{for all}\ i \in S.\]
The definition of $\varepsilon^{\prime}$ implies that
\[
\pi(i) \in \mathcal{T}\ \textrm{if and only if}\ i \in \varepsilon^{-1}(-1),\]
Hence
\[
S \cap\chi_p^{-1}(-1)=\varepsilon^{-1}(-1),\]
and so $\chi_p \equiv \varepsilon$ on $S$.$\hspace{12cm}      \textrm{QED}$

\emph{Remark}. The converse of Lemma 6.1 is valid.

In order to verify statement (5), and hence prove Theorem 5.13, it suffices by virtue of Lemma 6.1 to prove that if
\begin{equation*}
\textrm{$\triangle_{T \in \mathcal{T}} T \not= \emptyset$ for all $\emptyset \not= \mathcal{T} \subseteq \mathcal{S}$} \tag{8}
\]
then
\begin{equation*}
\textrm{for each $\emptyset \not= \mathcal{T} \subseteq \mathcal{S}$, there exits  $N \subseteq \Pi$ such that $\mathcal{T}=\{S \in \mathcal{S}: |N \cap S|\ \textnormal{ is odd} \}.  $ 
} \tag{9}\]
We have now completely removed residues and non-residues from the scene and have reduced everything to proving the following purely combinatorial statement about finite sets:
\begin{quote}
if $A$ is a nonempty finite set, $\emptyset \not= \mathcal{S} \subseteq 2^A \setminus \{\emptyset\}$, and $\mathcal{S}$ satisfies (8), then, with $\Pi$ replaced by $A$, $\mathcal{S}$ satisfies (9).
\end{quote} 

This can be done via linear algebra over $F=GF(2)$, by means of the same idea that we used in the proof of Lemma 4.11. We may suppose with no loss of generality that $A=[1, n]$ for some $n \in [1, \infty)$. Let 
\[
v: 2^A \rightarrow F^n\]
be the map defined in section 6 of Chapter 4. If $\mathcal{S}=\{S_1,\dots,S_m\}$, note that if $\emptyset \not= \mathcal{T} \subseteq \mathcal{S}$ then there is a bijection of the set of solutions over $F$ of the $m \times n$ system of linear equations 
\[
\sum_i \ v(T)(i)x_i=1,\ T \in \mathcal{T},\]
\[
\sum_i \ v(S)(i)x_i=0,\ S \in \mathcal{S} \setminus  \mathcal{T}, \]
onto the set 
\[
\{N \subseteq [1, n]: N\ \textrm{satisfies the conclusion of (9) (with $\Pi$ replaced by $A$) for $\mathcal{T}$}\} \]
given by
\[
(x_1,\dots,x_n) \rightarrow \{i: x_i=1\}.\]
Hence (9) holds with $\Pi$ replaced by $A$ if and only if the linear transformation of $F^n \rightarrow F^m$ with matrix
\[
B=\left( \begin{array} {ll} v(S_1)(1) \dots v(S_1)(n)\\
 \vdots \hspace{2.1cm}      \vdots\\
v(S_m)(1) \dots v(S_m)(n) \end{array} \right)
\]
is surjective, i.e., $B$ has rank $m$, i.e., the row vectors of $B$ are linearly independent over $F$.

We now show that
\begin {equation*}
\textrm{ the row vectors of $B$ are linearly independent over $F$ if and only if $\mathcal{S}$ satisfies $(8)$;} \tag{10} \]
this will prove Theorem 5.13 using only Lemma 4.4 and linear algebra over $F$!

If $w=(w_1,\dots,w_n) \in F^n$, recall that the \emph{support $\textnormal{supp}(w)$ of $w$} is the set 
\[
\textrm{supp}(w)=\{i: w_i=1\}.\]
It is easy to see that if $\emptyset \not= U \subseteq F^n$ then
\[
\textrm{supp}\Big(\sum_{w \in U} w\Big)= \triangle_{w \in U}\ \textrm{supp}(w), \]
and so
\begin{equation*}
\sum_{w \in U} w \not= 0\ \textrm{if and only if }\ \triangle_{w \in U}\ \textrm{supp}(w) \not= \emptyset. \tag{11} \]
Observe now that 
\begin{equation*}
\textrm{$U$ is linearly independent over $F$ if and only if for all $\emptyset \not= W \subseteq U,  \sum_{w \in W} w \not=0$.}  \tag{12}\]
Statement (10) is now a consequence of (11), (12), and the fact that
\[
\textrm{ supp$\big(v(T)\big)=T$, for all $T \in \mathcal{S}$.}\] 
$\hspace{15.5cm} \textrm{QED}$

\section{An Elementary Proof of Theorem 4.12}
Now for the proof of Theorem 4.12. Let $S$ be a nonempty finite subset of $[1, \infty)$ such that 
\begin{equation*}
p(T)\ \textrm{is not a square for all $T \subseteq S$ with $|T|$ odd.} \tag{13} \]
We need to prove that the set $\{p: \chi_p\equiv -1\ \textrm{on}\  S\}$ is infinite.

If we replace $S$ by the set $S^{\prime}$ of integers formed by the square-free parts of the elements of $S$ then (13) is true with $S$ replaced by $S^{\prime}$ hence we may again suppose with no loss of generality that all integers in $S$ are square-free. 

The argument now proceeds along the same line of reasoning that we used to prove Theorem 5.13. It follows as before that, with $\mathcal{S}= \{\pi(i): i \in S\}$, 
\begin{equation*}
\textrm{condition (13) holds if and only if $\triangle_{T \in \mathcal{T}} T \not= \emptyset$ for all $\mathcal{T} \subseteq \mathcal{S}$ with $|\mathcal{T}|$ odd.} \tag{14} \]
We then look for a combinatorial condition on $\mathcal{S}$ which implies that the set of primes
\[
\{p: \chi_p \equiv -1\ \textrm{on}\ S\} \]
is infinite, in analogy with Lemma 6.1. Such a condition is provided by
\begin{lem}
If there exists a subset $N$ of $\Pi= \bigcup_{i\in S} \pi(i)$ such that
\[
|N \cap \pi(i)|\  is\ odd\ for\ all\ i \in S,\]
then 
\[
\{p: \chi_p \equiv -1\ \textrm{on}\ S\} \]
is infinite.
\end{lem}

\emph{Proof}.  Let $N$ be a subset of $\Pi$ which satisfies the hypothesis of Lemma 6.2. As before, use Lemma 4.4 to find infinitely many primes $p$ such that
\[
\{q \in \Pi: \chi_p(q)= -1\}=N\ ; \]
then for all such $p$ which divides no element of $S$,
\[
\chi_p(i)=(-1)^{|N \cap \pi(i)|}=-1,\ \textrm{for all}\ i \in S.\]
$\hspace{15.4cm}\textrm{QED}$
 
The final step is to prove that if $A$ is a nonempty finite set, $\emptyset \not= \mathcal{S} \subseteq 2^A \setminus \{\emptyset\}$, and
\begin{equation*}
\textrm{$\triangle_{T \in \mathcal{T}} T \not= \emptyset$ for all $ \mathcal{T} \subseteq \mathcal{S}$ with $|\mathcal{T}|$ odd}, \tag{15}
\]
then there is a subset $N$ of $A$ such that
\[
|N\cap S|\ \textrm{is odd, for all}\ S \in \mathcal{S},\]
which can be done again by linear algebra over $F$.

We may take $A=[1, n]$, list the elements of $\mathcal{S}$ as $\mathcal{S}=\{S_1,\dots,S_m\}$ and then observe that, as in the proof just given of Theorem 5.13, there is a bijection of the set of solutions in $F^n$ of the system of equations
\[
\sum_i v(S_j)(i)x_i=1,\ j=1,\dots,m,\]
onto the set 
\[
\{N \subseteq [1, n]: |N \cap S|\ \textrm{is odd, for all}\ S \in \mathcal{S}\}.\] 
This system has a solution if and only if the matrices
\[
B=\left( \begin{array} {ll} v(S_1)(1) \dots v(S_1)(n)\\
 \vdots \hspace{2.1cm}      \vdots \\
v(S_m)(1) \dots v(S_m)(n) \end{array} \right)
\]
and
\[
B^{\prime}=\left( \begin{array} {ll} v(S_1)(1) \dots v(S_1)(n)\hspace{.4cm} 1\\
 \vdots \hspace{2cm}      \vdots\ \hspace{1.7cm} \vdots\\
v(S_m)(1) \dots v(S_m)(n)\hspace{.2cm} 1 \end{array} \right)
\]
have the same rank (over $F$), hence we must verify that if (15) holds then $B$ and $B^{\prime}$ have the same rank.

Assuming that (15) is valid, we let $v_1,\dots,v_m, v_1^{\prime},\dots,v_m^{\prime}$ denote the row vectors of $B$ and $B^{\prime}$, respectively. We will use (15) to prove that 
\begin{equation*}
\textrm{for all}\ \emptyset \not= T \subseteq [1, m],\ \sum_{i \in T} v_i=0\ \textrm{iff}\ \sum_{i \in T} v_i^{\prime}=0. \tag{16} \]
Statement (16) implies that if $\mathcal{L}$ (respectively, $\mathcal{L}^{\prime}$) is the set of all sets of linearly independent rows of $B$ (respectively, $B^{\prime}$) then the map $v_i \rightarrow v_i^{\prime}$ induces a bijection $\Lambda$ of $\mathcal{L}$ onto $\mathcal{L}^{\prime}$ such that 
\[
|\Lambda(L)|=|L|,\ \textrm{for all}\ L \in \mathcal{L},\]
and so
\[
\textrm{rank of $B$}= \max_{L \in \mathcal{L}} |L|=\max_{L \in \mathcal{L}^{\prime}} |L|=\textrm{rank of $B^{\prime}$}.\]

In order to verify (16), note first that if $\emptyset \not= T \subseteq [1, m]$ then
\begin{equation*}  
\textrm{$i$-th coordinate of $\sum_{j\in T} v_j=i$-th coordinate of $\sum_{j \in T} v_j^{\prime}, i=1,\dots,n$,} \tag{17} 
\]
and so if $\sum_{j \in T} v_j^{\prime}=0$ then $\sum_{j\in T} v_j=0$. Conversely, if $\sum_{j\in T} v_j=0$ then (15) implies that $|T|$ is even. Consequently, 
\[
(n+1)\textrm{-}\textrm{th coordinate of}\  \sum_{j \in T} v_j^{\prime}= |T| \cdot 1=0,\]
hence this equation and (17) imply that $\sum_{j \in T} v_j^{\prime}=0$. $\hspace{6.1cm} \textrm{QED}$ 
 
We close this chapter by discussing what happens if instead of subsets of $[1, \infty)$ we allow nonempty, finite subsets of $\mathbb{Z} \setminus \{0\}$ in the hypotheses of all of the theorems in Chapters 4 and 5. Theorem 4.2 remains valid if the positive integer in its hypothesis is replaced by a non-zero integer, and Theorems 4.3, 4.12, 5.13, and 5.15 remain valid with no change in their statements if the set $S$ in the hypotheses there is replaced by an arbitrary nonempty, finite subset of $\mathbb{Z} \setminus \{0\}$. In this more general situation, the integer $-1$ behaves like an additional prime, and once that is taken into account, all of our arguments, both elementary and non-elementary, can be modified without too much additional effort to verify these more general results. If the subset of $[1, \infty)$ in the hypotheses of Theorems 4.9 and 5.14 is replaced by a nonempty, finite subset $S$ of $\mathbb{Z} \setminus \{0\}$ and if the dimension $d$ is determined by $S$ as in the statements of those theorems, then the density of the sets in their conclusions is now either $2^{-d}$ or $2^{-(1+d)}$, with the latter value occurring if either $-1 \in S$ or the sets $\pi_{\textrm{odd}}(z), z \in S$, possess a certain combinatorial structure. However, the proof of this version of Theorems 4.9 and 5.14 proceeds along the same lines as the arguments that we have given, with only a few additional technical adjustments (see Wright [61], section 3 for the details). 

\newpage
\afterpage{\null\newpage}
\thispagestyle{empty}
\chapter{Dirichlet $L$-functions and the Distribution of Quadratic Residues}

In section 4 of Chapter 4, we saw how the non-vanishing at $s=1$ of the $L$-function $L(s, \chi)$ of a non-principal Dirichlet character $\chi$ played an essential role in the proof of Dirichlet's theorem on prime numbers in arithmetic progression (Theorem 4.5). In this chapter, the fact that $L(1, \chi)$ is not only nonzero, but \emph{positive}, when $\chi$ is real and non-principal, will be of central importance. The positivity of  $L(1, \chi)$ comes into play because we are interested in the following problem concerning the distribution of residues and non-residues of a prime $p$. Suppose that $I$ is an interval of the real line contained in the interval from 1 to $p$. Are there more residues of $p$ than non-residues in $I$, or are there more non-residues than residues, or is the number of residues and non-residues in $I$ the same? We will see that this question can be answered if we can determine if certain sums of values of the Legendre symbol of $p$ are positive, and it transpires that the positivity of the sum of these Legendre-symbol values, for certain primes $p$, are determined precisely by the positivity of $L(1, \chi)$ for certain Dirichlet characters $\chi$. We make all  of this precise in section 1, where the principal theorem of this chapter, Theorem 7.1, is stated and then used to  obtain some very interesting answers to our question about the distribution of residues and non-residues. In the next section, the proof of Theorem 7.1 is outlined; in particular we will see how the proof can be reduced to the verification of formulae, stated in Theorems 7.2, 7.3, and 7.4, which express the relevant Legendre-symbol sums in terms of the values of $L$-functions at $s=1$. Sections 3-7 are devoted to the proof Theorems 7.2-7.4. In section 3, the fact that $L(1, \chi)>0$ for real, non-principal Dirichlet  characters  is established, and sections 4-6 are devoted to discussing various results concerning Gauss sums, analytic functions of a complex variable, and Fourier series which are required for the arguments we take up in section 7. Because it plays such an important role in the results of this chapter, we prove in section 8 Dirichlet's fundamental Lemma 4.8 on the non-vanishing of $L(1, \chi)$ for real, non-principal characters. Motivated by the result on the convergence of Fourier series that is proved in section 6, we give yet another proof of quadratic reciprocity in section 9 that uses \emph{finite} Fourier series expansions.

\section{Positivity of Sums of Values of a Legendre Symbol}
Dirichlet [11] proved the following theorem in 1839: 
\begin{thm}

$(i)$ If $p \equiv 3\ \textnormal{mod}\ 4$ then
\[
\sum_{0<n<p/2}\ \chi_p(n)>0.\]

$(ii)$ If $p \equiv 1\ \textnormal{mod}\ 4$ then
\[
\sum_{0<n<p/4}\ \chi_p(n)>0.\]

$(iii)$ If $p>3$ then
\[
\sum_{0<n<p/3}\ \chi_p(n)>0.\]
\end{thm}
This result initiated a line of intense research on the positivity of various sums of values of a Dirichlet character, and of characters on more general groups, that continues unabated to the present day. The importance to us of Theorem 7.1 lies in its connection with the distribution of residues and non-residues of a prime $p$ throughout the interval $[1, p-1]$. In order to see how that goes, we consider a subinterval $I$ of $\mathbb{R}$ contained in $\{x\in \mathbb{R}: 0< x< p\}$, and, following Berndt [1], we define the \emph{ quadratic excess of I} to be the sum
\[
q(I)=\sum_{n \in I}\ \chi_p(n).\]
If $q(I)>0$ (respectively, $q(I)<0$) then the number of residues (respectively, non-residues) of $p$  inside $I$ \emph{exceeds} the number of non-residues (respectively, residues) of $p$ there, and if $q(I)=0$ then the number of residues and non-residues are the same. Hence Theorem 7.1 implies that if $p \equiv 3\ \textnormal{mod}\ 4$ then the number of residues inside the interval $(0, p/2)$  exceeds the number of non-residues there, or if $p \equiv 1\ \textnormal{mod}\ 4$ then the number of residues inside the interval $(0, p/4)$  exceeds the number of non-residues there, or if  $p>3$ then the number of residues inside the interval $(0, p/3)$  exceeds the number of non-residues there.

By taking Proposition 2.1 and Theorem 2.4 of Chapter 2 into account, we can say more. If $\{X_1,\dots,X_k\}$ is a set of pairwise disjoint subintervals of  $\{x\in \mathbb{R}: 0< x< p\}$ such that $[1, p-1]=\mathbb{Z} \cap \big( \bigcup_i X_i \big)$ then because of Proposition 2.1, we have that 
\begin{equation*}
\sum_i q(X_i)=0. \tag{1} \]
Now, using $(a, b)$ to denote the interval $\{x\in \mathbb{R}: a< x< b\}$, let
\[
I_1=(0, p/3),\ I_2=(p/3, 2p/3),\ I_3=(2p/3, p),\]
\[
J_1=(0, p/4),\ J_2=(p/4, p/2),\ J_3=(p/2, 3p/4),\ J_4=(3p/4, p).\]

Assume first that $p \equiv 3\ \textnormal{mod}\ 4$. Theorem 2.4 implies that $\chi_p(-1)=-1$, hence
\begin{eqnarray*}
(2)\hspace{2cm}     q(I_1)&=&\sum_{0<n<p/3}\ \chi_p(n)\\
&=&-\sum _{0<n<p/3}\ \chi_p(-n)\\
&=&-\sum _{0<n<p/3}\ \chi_p(p-n)\\
&=&-\sum _{2p/3<n<p}\ \chi_p(n)\\
&=&-q(I_3),
\end{eqnarray*}
and so by (1) and Theorem 7.1 ($iii$),
\[
q(I_2)=0\ \textrm{and}\ q(I_3)<0.\]
It follows that $(p/3, 2p/3)$ contains the same number of residues as non-residues of $p$ and the number of non-residues in $(2p/3, p)$ exceeds the number of residues there.

Assume next that $p \equiv 1\ \textnormal{mod}\ 4$. Theorem 2.4 implies that $\chi_p(-1)=1$ hence the minus signs in (2) can be dropped to conclude that 
\[
q(I_1)=q(I_3),\]
and so by (1) and Theorem 7.1$(iii)$ yet again,
\[
q(I_3)>0\ \textrm{and}\ q(I_2)=-q(I_1)-q(I_3)<0.\]
It follows that the number of non-residues of $p$ in $(p/3, 2p/3)$ exceeds the number of residues there and the number of residues of $p$ in $(2p/3, p)$ exceeds the number of non-residues there.

Similar arguments show that if $p \equiv 1\ \textnormal{mod}\ 4$ then
\begin{equation*}
q(J_1)=q(J_4),\ q(J_2)=q(J_3),\ q(J_1)=-q(J_3), \tag{3}\]
hence we conclude from (3) by way of Theorem 7.1$(ii)$ that the number of residues of $p$ in each of the intervals $(0, p/4)$ and $(3p/4, p)$ exceeds the number of non-residues  there and the number of non-residues in each of the intervals $(p/4, p/2)$ and $(p/2, 3p/4)$) exceeds the number of residues there.

Figures 1-4 below display graphically the distribution of the residues and non-residues in each of the cases that we have discussed. A $+$ above an interval indicates that the number of residues in that interval exceeds the number of non-residues there, a $-$ indicates that the number of non-residues exceeds the number of residues, and a $0$ indicates that the number of residues and non-residues are the same. We now turn to the proof of Theorem 7.1.
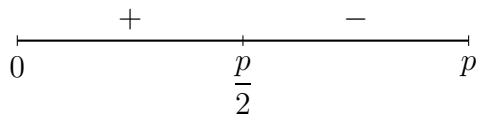
\begin{figure}[h]
  \centering
  \begin{tikzpicture}
    \draw[name path=left] (0, 2pt) -- (0, -2pt) node[anchor=north] {0};
    \draw[name path=right] (6, 2pt) -- (6, -2pt) node[anchor=north] {$p$};
    \draw[name path=middle] (3, 2pt) -- (3, -2pt) node[anchor=north] {$\displaystyle{\frac{p}{2}}$};
    \draw (1.5,0) node[above] {$+$};
    \draw (4.5,0) node[above] {$-$};
    \draw[name path=xaxis, thick] (0,0) -- (6,0);
  \end{tikzpicture}
  \caption{$p \equiv 3 \bmod 4$}
  \label{fig:7.1}
\end{figure}
\begin{figure}[h]
  \centering
  \begin{tikzpicture}
    \draw[name path=left] (0, 2pt) -- (0, -2pt) node[anchor=north] {0};
    \draw[name path=right] (6, 2pt) -- (6, -2pt) node[anchor=north] {$p$};
    \draw[name path=middle] (2, 2pt) -- (2, -2pt) node[anchor=north] {$\displaystyle{\frac{p}{3}}$};
    \draw[name path=middle] (4, 2pt) -- (4, -2pt) node[anchor=north] {$\displaystyle{\frac{2p}{3}}$};
    \draw (1,0) node[above] {$+$};
    \draw (3,0) node[above] {$-$};
    \draw (5,0) node[above] {$+$};
    \draw[name path=xaxis, thick] (0,0) -- (6,0);
  \end{tikzpicture}
  \caption{$p \equiv 1 \bmod 4$}
  \label{fig:7.2}
\end{figure}
\begin{figure}[h]
  \centering
  \begin{tikzpicture}
    \draw[name path=left] (0, 2pt) -- (0, -2pt) node[anchor=north] {0};
    \draw[name path=right] (6, 2pt) -- (6, -2pt) node[anchor=north] {$p$};
    \draw[name path=middle] (2, 2pt) -- (2, -2pt) node[anchor=north] {$\displaystyle{\frac{p}{3}}$};
    \draw[name path=middle] (4, 2pt) -- (4, -2pt) node[anchor=north] {$\displaystyle{\frac{2p}{3}}$};
    \draw (1,0) node[above] {$+$};
    \draw (3,0) node[above] {$0$};
    \draw (5,0) node[above] {$-$};
    \draw[name path=xaxis, thick] (0,0) -- (6,0);
  \end{tikzpicture}
  \caption{$p \equiv 3 \bmod 4$,\ $p>3$}
  \label{fig:7.3}
\end{figure}

\begin{figure}[h]
  \centering
  \begin{tikzpicture}
    \draw[name path=left] (0, 2pt) -- (0, -2pt) node[anchor=north] {0};
    \draw[name path=right] (8, 2pt) -- (8, -2pt) node[anchor=north] {$p$};
    \draw[name path=middle] (4, 2pt) -- (4, -2pt) node[anchor=north] {$\displaystyle{\frac{p}{2}}$};
    \draw[name path=right] (2, 2pt) -- (2, -2pt) node[anchor=north] {$\displaystyle{\frac{p}{4}}$};
    \draw[name path=middle] (6, 2pt) -- (6, -2pt) node[anchor=north] {$\displaystyle{\frac{3p}{4}}$};
    \draw (1,0) node[above] {$+$};
    \draw (3,0) node[above] {$-$};
    \draw (5,0) node[above] {$-$};
    \draw (7,0) node[above] {$+$};
    \draw[name path=xaxis, thick] (0,0) -- (8,0);
  \end{tikzpicture}
  \caption{$p \equiv 1 \bmod 4$}
  \label{fig:7.2}
\end{figure}

\section{Proof of Theorem 7.1: Outline of the Argument}

The proof of Theorem 7.1 depends on formulae for the quadratic excesses in the statement of that theorem which are given in terms of certain Dirichlet $L$-functions. Recall from section 4 of Chapter 4 that if $\chi$ is a Dirichlet character then the $L$-function of $\chi$ is defined by the Dirichlet series
\[
L(s, \chi)=\sum_{n=1}^{\infty} \frac{\chi(n)}{n^s},\ s \in \mathbf{C}.\]
The property of these $L$-functions that will be essential for our proof of Theorem 7.1 is the fact that when $\chi$ is real and non-principal, $L(1, \chi)>0$. It follows from Dirichlet's fundamental Lemma 4.8 in Chapter 4 that $L(1, \chi)\not=0$ for any non-principal Dirichlet character $\chi$, and we will show (among other things) in the next section that when $\chi$ is also  real, $L(1, \chi)\geq 0$, whence  $L(1, \chi)>0$. Consequently, if we can prove that each of the sums in Theorem 7.1 can be expressed as a positive multiple of the value at $s=1$ of the $L$-function of a real non-principal Dirichlet character, then Theorem 7.1 will follow from the positivity of that $L$-function value. That is the line of reasoning which we will follow to the proof of Theorem 7.1.

In order to carry out this argument, we therefore require formulae which express the quadratic excesses $q(0, p/2)$, $q(0, p/4)$, and $q(0, p/3)$ in terms of the value of $L$-functions at $s=1$.  The formula for $q(0, p/2)$ is given in
\begin{thm}
If $p \equiv 3\ \textnormal{mod}\ 4$ then
\[
q(0, p/2)=\frac{\sqrt p}{\pi} \big(2-\chi_p(2)\big)L(1, \chi_p).\]
\end{thm}
This theorem implies that statement $(i)$ of Theorem 7.1 is true.

In order to state the $L$-function formulae that will verify Theorem 7.1$(ii)$ and $(iii)$, we will need to make use of the fact that if $\chi_m$ and $\chi_n$ are Dirichlet characters of modulus $m$ and $n$, and if gcd$(m, n)=1$, then the point-wise product $\chi_m \chi_n$ is a Dirichlet character of modulus $mn$. This follows from the fact that if gcd$(m, n)=1$ then the Chinese remainder theorem implies that $U(mn)$ is isomorphic to the direct product $U(m) \times U(n)$, and so the point-wise product $\chi_m \chi_n$ clearly defines a homomorphism of $U(mn)$ into the circle group.

Our proof of Theorem 7.1 $(ii)$ will make use of the character $\chi_{4p}$ of modulus $4p$ given by point-wise multiplication of $\chi_p$ and the character $\chi_4$ of modulus 4 defined by 
\[
\chi_4(n)=\left\{\begin{array}{rl} (-1)^{(n-1)/2}\ ,& \textrm{$n$ odd,}\\
0\ ,& \textrm{$n$ even.} 
\\\end{array}\right. 
\]
Also, if $p>3$ then we let $\chi_{3p}$ denote the point-wise product of $\chi_3$ and $\chi_p$. It is clear that the characters $\chi_{3p}$ and $\chi_{4p}$ are real and non-principal.

Theorem 7.1$(ii)$ and $(iii)$ now follow, respectively, from
\begin{thm}
If $p \equiv 1\ \textnormal{mod}\ 4$ then
\[
q(0, p/4)=\frac{\sqrt p}{\pi} L(1, \chi_{4p}).\]
\end{thm}

\begin{thm}
Let $p>3$.

$(i)$ If $p \equiv 1\ \textnormal{mod}\ 4$ then
\[
q(0, p/3)=\frac{\sqrt{3p}}{2\pi} L(1, \chi_{3p}).\]

$(ii)$ If $p \equiv 3\ \textnormal{mod}\ 4$ then
\[
q(0, p/3)=\frac{\sqrt p}{2\pi}\big(3-\chi_{p}(3)\big) L(1, \chi_{p}).\]
\end{thm}

We have now reduced the proof of Theorem 7.1 to the proof of Theorems 7.2, 7.3, and 7.4. The proof of these results will require quite a bit of preliminary preparation. The facts about Dirichlet $L$-functions that will be needed are set forth in section 3. Our arguments will also require the calculation of a very useful Gauss sum, which will be carried out in section 4 (we used Gauss sums in section 9 of Chapter 3 to verify the Law of Quadratic Reciprocity). Our proof of Theorems 7.2 and 7.4$(ii)$ will follow a very nice argument of Bruce Berndt [1] which employs the theory of analytic functions of a complex variable, and so we will discuss the relevant facts from that subject in section 5. The convergence of Fourier series is the subject of section 6, required because we will use Fourier series to prove Theorems 7.3 and 7.4$(i)$, in the same spirit as Dirichlet's original proof of those theorems. Finally, we will bring all of this together for the proof of Theorems 7.2, 7.3, and 7.4 in section 7.

\section{Some Useful Facts About Dirichlet $L$-functions}

The information concerning $L$-functions of Dirichlet characters that we will need are recorded in
\begin{lem}
Let $\chi$ be a Dirichlet character $\textnormal{mod}\ m$.

$(i)$ If $\chi$ is non-principal then $L(s, \chi)$ is analytic in the half-plane $\textnormal{Re}\ s>0$.

$(ii)$ $L(s, \chi)$ has the absolutely convergent Euler-Dirichlet product expansion given by
\[
L(s, \chi)=\prod_q \frac{1}{1-\chi(q)q^{-s}},\ \textnormal{Re}\ s>1,\]
where the product is taken over all prime numbers $q$.

$(iii)$ If $\chi$ is real-valued and non-principal then $L(1, \chi)>0$.
\end{lem}
 
 \emph{Proof}. $(i)$ This will follow immediately from Proposition 5.5 after we prove that the sums
 \[
 \sum_{k=1}^n \chi(k) \]
 are uniformly bounded as a function of $n$. To see this, we claim first that
 \begin{quote}
(4) $\sum \chi(k)=0$, whenever this sum is taken over any complete system of ordinary residues mod $m$. 
\end{quote}
Assuming this is true, we take $n \in [1, \infty)$, write $n=r+lm,\ 0 \leq r<m$, and then calculate that
\begin{eqnarray*}
\sum_1^n \chi(k)&=&\sum_1^{lm-1}\chi(k)+\sum_{k=0}^r \chi(k+lm)\\
&=&\sum_{k=0}^r \chi(k+lm),\ \textrm{by (4)}\\
&=&\sum_{k=0}^r \chi(k),
\end{eqnarray*} 
hence\[
\Big|\sum_1^n \chi(k) \Big| \leq \sum_0^r |\chi(k)| \leq m-1. \]

In order to verify (4) use the fact that $\chi$ is periodic of period $m$ (Proposition 4.6) and the fact that $k$ in (4) runs through a complete set of ordinary residues mod $m$ to write
\[
\sum_k \chi(k)=\sum_{k \in U(m)} \chi(k),\]
so we need only show that this latter sum is 0.

 Because $\chi$ is non-principal, there is a $k_0 \in U(m)$ such that $\chi(k_0) \not= 1$. 
 The map $k \rightarrow kk_0$ is a bijection of $U(m)$ onto $U(m)$, hence
 \[
 \sum_{k \in U(m)} \chi(k)= \sum_{k \in U(m)} \chi(kk_0)=\chi(k_0)\sum_{k \in U(m)} \chi(k),\]
 hence
 \[
 (1-\chi(k_0)) \sum_{k \in U(m)} \chi(k)=0.\]
 As $1-\chi(k_0) \not= 0$, it follows that
 \[
 \sum_{k \in U(m)} \chi(k)=0.\]
 
$(ii)$ This product formula can be derived by appropriate modifications of our proof of Theorem 5.8, which verified the product formula for the zeta function of an algebraic number field. Note first that
\begin{eqnarray*}
\Big|\frac{1}{1-\chi(q)q^{-s}}-1\Big|&=&\Big|\frac{\chi (q)q^{-s}}{1-\chi(q)q^{-s}} \Big|\\
&\leq&\frac{q^{-\textrm{Re}\ s}}{1-q^{-\textrm{Re}\ s}}\\
&\leq&2q^{-\textrm{Re}\ s},\ \textrm{for all}\ q \geq 2,\ \textrm{Re}\ s>1,
\end{eqnarray*}
consequently Proposition 5.9 implies that the product in $(ii)$ is absolutely convergent for Re $s>1$. The proof of Theorem 5.8 can now be easily modified by replacing the set of prime ideals of $R$, the set of nonzero ideals of $R$, Proposition 5.2 and the Fundamental Theorem of Ideal Theory in that proof by, respectively, the set $P$ of all primes, the set $[1, \infty)$,  the complete multiplicativity of $\chi$, and the Fundamental Theorem of Arithmetic to obtain
\[
\sum_{n=1}^{\infty} \frac{\chi(n)}{n^s}=\prod_q \frac{1}{1-\chi(q)q^{-s}},\ \textnormal{for Re}\ s>1.\]

$(iii)$ If $\chi$ is real then every value of $\chi$ is 0 or $\pm1$, hence each factor in the Euler product expansion of $L(s, \chi)$ is positive for $s>1$. Consequently $L(s, \chi)$ is not less than 0, and so by the continuity of $L(s, \chi)$ on $s>0$ it follows that
\[
L(1, \chi)=\lim_{s \rightarrow 1^{+}} L(s, \chi) \geq 0.\]
But $L(1, \chi) \not= 0$, because of Dirichlet's fundamental Lemma 4.8, hence $L(1, \chi)>0$. $\hspace{0.8 cm} \textrm{QED}$   

\section{Calculation of a Gauss Sum}

In addition to $L$-functions, our derivation of the formulae in Theorems 7.2, 7.3, and 7.4 will also employ some very useful properties of Gauss sums. Recall from the proof of quadratic reciprocity in section 9 of Chapter 3 the Gauss sums
\[
G(n, p)=\sum_{j=0}^{p-1} \chi_p(j) \exp\Big(\frac{2 \pi inj}{p} \Big).\]
In that proof (Lemma 3.15 and Theorem 3.14), we showed that
\begin{equation*}
G(n, p)=\chi_p(n)G(1, p) \tag{5}\]
and that
\[
G(1, p)^2=\left\{\begin{array}{rl} p,& \textrm{if $p \equiv 1$ mod $4$,}\\
-p,& \textrm{if $p \equiv 3$ mod $4$.} 
\\\end{array}\right. 
\]
Determining the sign of $G(1, p)$ from this equation turns out to be a very difficult problem, and was solved by Gauss in 1805 after four long years of intense effort on his part. The plus sign is the correct one in both cases; we will present a very nice proof of this fact due to L. Kronecker, according to the account of it given in Ireland and Rosen [30], section 6.4.
\begin{thm}
\[
G(1, p)=\left\{\begin{array}{rl}\sqrt p\ ,& \textnormal{if}\ p \equiv 1\ \textnormal{mod}\ 4,\\
i \sqrt p\ ,& \textnormal{if}\ p \equiv 3\ \textnormal{mod}\ 4. 
\\\end{array}\right. 
\]
 \end{thm}
 \emph{Proof}. Let $\zeta=\exp(2\pi i/p)$. The argument proceeds through a series of claims and their verifications.
 
 \emph{Claim} 1.
 \[
 (-1)^{(p-1)/2}p=\prod_{k=1}^{(p-1)/2} (\zeta^{2k-1}-\zeta^{-2k+1})^2.\]
 
 \emph{Claim} 2.
 \[
\prod_{k=1}^{(p-1)/2} (\zeta^{2k-1}-\zeta^{-2k+1})= \left\{\begin{array}{rl}\sqrt p\ ,& \textnormal{if}\ p \equiv 1\ \textnormal{mod}\ 4,\\
i \sqrt p\ ,& \textnormal{if}\ p \equiv 3\ \textnormal{mod}\ 4. 
\\\end{array}\right. 
\]

Once that Claim 1 is verified, we deduce from Theorem 3.14 that
\[
G(1, p)=\varepsilon \prod_{k=1}^{(p-1)/2} (\zeta^{2k-1}-\zeta^{-2k+1}).
\]
where $\varepsilon=\pm 1$. The conclusion of Theorem 7.6 will then be at hand once we verify Claim 2 and prove that $\varepsilon=1$. Hence we make  

\emph{Claim} 3. $\varepsilon=1$.
 
To verify Claim 1, start with the factorization
\[
x^p-1=(x-1)\prod_{j=1}^{p-1}\ (x-\zeta^j).\]
Divide this equation by $x-1$ and set $x=1$ to derive that
\[
p=\prod_r\ (1-\zeta^r),\]
where this product is taken over any complete system of ordinary residues mod $p$. It is easy to see that the integers $\pm(4k-2), k=1,\dots,(p-1)/2,$ is such a system of residues, and so
\begin{eqnarray*}
 p&=&\prod_1^{(p-1)/2}\ (1-\zeta^{4k-2})\ \prod_1^{(p-1)/2}\ (1-\zeta^{-(4k-2)})\\
 &=&\prod_1^{(p-1)/2}\ (\zeta^{-(2k-1)}-\zeta^{2k-1})\ \prod_1^{(p-1)/2}\ (\zeta^{2k-1}-\zeta^{-(2k-1)})\\
 &=& (-1)^{(p-1)/2}\ \prod_1^{(p-1)/2}\ (\zeta^{2k-1}-\zeta^{-2k+1})^2. 
\end{eqnarray*} 
 
 Now for Claim 2. Claim 1 implies that
 \[
\Big(\prod_1^{(p-1)/2}\ (\zeta^{2k-1}-\zeta^{-2k+1})\Big)^2=(-1)^{(p-1)/2}p,\]
hence Claim 2 will follow from this equation once the sign of the product in Claim 2 is determined. That product is
\[
i^{(p-1)/2}\ \prod_1^{(p-1)/2} 2\sin \frac{(4k-2)\pi}{p}.\]
Observe now that for $k \in [1, (p-1)/2]$,
\[
\sin \frac{(4k-2)\pi}{p}<0\ \textrm{ iff}\  \frac{p+2}{4}<k \leq \frac{p-1}{2},\]
hence this product has precisely $(p-1)/2-[(p+2)/4]$ negative factors, and so the number of negative factors is either $(p-1)/4$ or $(p-3)/4$ if, respectively, $p \equiv 1$ or 3 mod 4. It is now easy to see from this that the product in Claim 2 is a positive number if $p \equiv1$ mod 4 or is $i \times$(a positive number) if $p \equiv 3$ mod 4.

In order to verify Claim 3, consider the polynomial 
\[
f(x)=\sum_{j=1}^{p-1} \chi_p(j)x^j-\varepsilon\prod_{k=1}^{(p-1)/2} (x^{2k-1}-x^{p-2k+1}).\]
Then
\[
f(\zeta)=G(1, p)-\varepsilon  \prod_1^{(p-1)/2} (\zeta^{2k-1}-\zeta^{-2k+1})=0\]
and
\[
f(1)=\sum_{j=1}^{p-1} \chi_p(j)=0.\]
Now the minimal polynomial of $\zeta$ over $\mathbb{Q}$ is $\sum_{k=0}^{p-1} x^k$, and so we conclude from the proof of Proposition 3.7 that $\sum_{k=0}^{p-1} x^k$ divides $f(x)$ in $\mathbb{Q}[x]$. As $x-1$ and $\sum_{k=0}^{p-1} x^k$ are both irreducible over $\mathbb{Q}$, they are relatively prime in $\mathbb{Q}[x]$. Because $x-1$ divides $f(x)$ in $\mathbb{Q}[x]$, it follows that $x^p-1=(x-1)( \sum_{k=0}^{p-1} x^k)$ must also divide $f(x)$ in $\mathbb{Q}[x]$. Hence there exists $h \in \mathbb{Q}[x]$ such that $f(x)=(x^p-1)h(x)$. Now replace $x$ by $e^z$ to obtain the equation
\[
\sum_{j=1}^{p-1} \chi_p(j)e^{jz}-\varepsilon\prod_{k=1}^{(p-1)/2} \big(e^{(2k-1)z}-e^{(p-2k+1)z}\big)=(e^{pz}-1)h(e^z).\]
Insert the power series expansion of $e^z$ into this equation and then deduce that the coefficient of $z^{(p-1)/2}$ on the left-hand side of the equation is
\[
\frac{1}{\big((p-1)/2\big)!} \sum_{j=1}^{p-1} \chi_p(j)j^{(p-1)/2}-\varepsilon\prod_{k=1}^{(p-1)/2}(4k-p-2),\]
while the coefficient of $z^{(p-1)/2}$ on the right-hand side is of the form $pA/B$, where $A$ and $B$ are integers and $\gcd(B, p)=1$. Now equate coefficients, multiply through by $B\big((p-1)/2\big)!$ and reduce mod $p$ to derive
\begin{eqnarray*}
\sum_{j=1}^{p-1} \chi_p(j)j^{(p-1)/2} &\equiv& \varepsilon\Big(\frac{p-1}{2}\Big)! \prod_{k=1}^{(p-1)/2} (4k-2)\\
&\equiv& \varepsilon \prod_{k=1}^{(p-1)/2} 2k\ \prod_{k=1}^{(p-1)/2} (2k-1)\\
&\equiv& \varepsilon(p-1)!\\
&\equiv&-\varepsilon\ \textrm{mod}\ p, 
\end{eqnarray*}
where the last congruence follows from Wilson's theorem. But then by Euler's criterion (Theorem 2.5),
\[
j^{(p-1)/2} \equiv \chi_p(j)\ \textrm{mod}\ p, \]
hence
\[
p-1=\sum_{j=1}^{p-1} \chi_p(j)^2 \equiv -\varepsilon\ \textrm{mod}\ p,\]
and so
\[
\varepsilon \equiv 1\ \textrm{mod}\ p.\]
Because $\varepsilon=\pm 1$, it follows that $\varepsilon=1$. $\hspace{8.8cm} \textrm{QED}$

\section{Some Useful Facts About Analytic Functions of a Complex Variable}

The proof of Theorems 7.2 and 7.4$(ii)$ that we will present uses an elegant application of contour integration from complex analysis due to Bruce Berndt. In this section we will discuss the requisite facts from that subject.

Let $\emptyset \not= U \subseteq \mathbf{C}$ be an open set. A function $f: U \rightarrow \mathbf{C}$ is \emph{analytic in U} if for each $z \in U$,
\[
\lim_{w \rightarrow z} \frac{f(w)-f(z)}{w-z}=f^{\prime}(z) \]
exists and is finite, i.e., $f$ has a \emph{complex derivative} at each point of $U$. A complex-valued function with domain $\textbf{C}$ is said to be \emph{entire} if it is analytic in $\textbf{C}$. We will use the following fundamental theorem about analytic functions in our proof of Lemma 4.8 for real Dirichlet characters that we will present in section 8:
\begin{thm}
$($Taylor-series expansion of analytic functions$)$ If f is analytic in U then the $n$-th order derivative $f^{(n)}(z)$ exists and is finite for all $z \in U$ and for all $n \in [1, \infty)$. Moreover, if $a \in U$ and $r>0$ is the distance of a to the boundary of U then
\[
f(z)=\sum_{n=0}^{\infty} \frac{f^{(n)}(a)}{n!} (z-a)^n\ ,\ |z-a|<r.
\]
\end{thm}

Theorem 7.7 highlights the remarkable regularity which all analytic functions possess: not only is an analytic function always infinitely differentiable, but it even has a convergent Taylor-series expansion in a neighborhood of each point in its domain. This is far from true for differentiable functions of a real variable.

Now let $I$ denote the closed unit interval on the real line, and let $\gamma: I \rightarrow U$ be a \emph{contour in U}, i.e., a continuous, piecewise-smooth function defined on $I$ with range in $U$. Let $\{\gamma\}$ denote the range of $\gamma$. If $g: \{\gamma\} \rightarrow \mathbf{C}$ is a function continuous on $\{\gamma\}$, $u=$ Re$(g)$, and $v=$ Im$(g)$, then the \emph{ contour integral of g along $\gamma$}, denoted by
\[
 \int_{\gamma} g(z)\ dz,
 \]
is defined by
\[
\oint_{\gamma} (u\ dx-v\ dy)+i \oint_{\gamma} (v\ dx+u\ dy),\]
where, from multi-variable calculus, $\oint_{\gamma}$ denotes standard line integration in the plane along $\gamma$ of real-valued functions continuous on $\{\gamma\}$. Since it would take us too far afield to give a detailed account of the properties of this integral, we instead refer to J.B. Conway [3], section IV.1 for that. We will need only the basic estimate
\begin{equation*}
\Big| \int_{\gamma} g(z)\ dz\Big| \leq \big(\max \big\{|g(z)|: z \in \{\gamma\}\big \}\big)(\textrm{length of}\ \gamma). \tag{6}\]
 
A contour $\gamma$ is \emph{closed} if $\gamma(0)=\gamma(1)$. The next theorem is one of the most important and most useful in all of complex analysis.
\begin{thm}
$($Cauchy's integral theorem $)$ If f is analytic in U and $\gamma$ is a closed contour in U which does not wind around any point in $\mathbf{C} \setminus U$ then 
\[
\int_{\gamma} f(z)\ dz=0. \]
\end{thm}

The next theorem provides a very useful formula for computing certain contour integrals of functions which are analytic outside of a finite set of points. In order to state it, some terminology needs to be defined, and so we will do that first.

A closed contour $\gamma$ is a \emph{Jordan contour} if $\gamma$ is an injective function on the set $I\setminus \{1\}$. Geometrically, this says that the path $\{\gamma\}$ does not cross itself (see figure 5). If $\gamma$ is a Jordan contour then $\gamma$ divides $\mathbf{C}$ into a pairwise disjoint union
\[
V \cup \{\gamma\} \cup W,
\]
where $V$ and $W$ are open sets and
\[
\textrm{the boundary of}\ V=\{\gamma\}=\ \textrm{the boundary of}\ W.\]
Suppose that as $t$ increases from 0 to 1, $\gamma(t)$ traverses $\{\gamma\}$ in the counterclockwise direction: we  then say that $\gamma$ is \emph{positively oriented}. If $\gamma$ is positively oriented then as $t$ increases from 0 to 1, for exactly one of the sets $V$ or $W$, $\gamma(t)$ winds around each of the points in that set exactly once. The set for which this occurs, either all of the points of $V$ or all of the points of $W$, is called the \emph{interior of} $\gamma$. The set $\mathbf{C} \setminus \big(\{\gamma\} \cup (\textrm{interior of}\ \gamma)\big)$ is the \emph{exterior of} $\gamma$. It can be shown that the interior of $\gamma$ is a bounded set and the exterior of $\gamma$ is unbounded. All of the facts in this paragraph are the contents of the Jordan Curve Theorem: for a proof, consult Dugundji [13], section XVII.5. 
\begin{figure}[h]
  \centering
  \begin{tikzpicture}[decoration={markings,mark=at position 0.25 with {\arrow[scale=2]{stealth}},mark=at position 0.5 with {\arrow[scale=2]{stealth}},mark=at position 0.75 with {\arrow[scale=2]{stealth}}}]
    \draw[<-,postaction=decorate] plot [smooth cycle] coordinates {(0,1)  (-3,2)  (-4,0)  (-4,-3)  (-2,-1.5) (0,-1) (8, -2)  (8,0) (6,4) (5,5) (3,4.5)    } 
    node at (0,0) {Interior of $\gamma$} 
    node at (-3,3) {Exterior of $\gamma$}
    node at (8.5,0) {$\{\gamma\}$};
  \end{tikzpicture}
  \caption{Geometry and topology of a positively oriented Jordan contour $\gamma$}
  \label{fig:10}
\end{figure}

A function $f$ has an \emph{isolated singularity} at a point $a$ if there is an $r>0$ such that $f$ is analytic in $0<|z-a|<r$, but $f^{\prime}(a)$ does not exist. An isolated singularity of $f$ at $a$ is a \emph{pole of order $m \in [1, \infty)$} if there exists $\delta>0$ and a function $g$ analytic in $|z-a|< \delta$ such that $g(a) \not= 0$ and
\[
f(z)=\frac{g(z)}{(z-a)^m},\ 0<|z-a|< \delta.\]
The \emph{residue of f} at this pole, denoted Res$(f, a)$, is the number
\[
\frac{g^{(m-1)}(a)}{(m-1)!}.\]
If the order of the pole at $a$ is 1 then it is called a \emph{simple pole}, and its residue there is
\[
g(a)=\lim_{z \rightarrow a}(z-a)f(z).\]
We can now state the result on the calculation of contour integrals that we need. 

\begin{thm}
$($The residue theorem$)$ Let $U$ be an open subset of $\mathbf{C}$, f a function analytic in $U$ except for poles located in $U$. If $\gamma$ is a positively oriented Jordan contour in U which does not wind around a point in $\mathbf{C} \setminus U$ and which does not pass through any of the poles of f, and if $a_1,\dots,a_n$ are the poles of f that are in the interior of $\gamma$, then
\[
\frac{1}{2\pi i}\int_{\gamma} f(z) dz=\sum_{k=1}^n \textnormal{Res}(f, a_k).\]
\end{thm}

For proof of Theorems 7.7, 7.8, and 7.9, consult, respectively, Conway [3], sections IV.2, IV.5, and V.2.

We will apply Theorems 7.8 and 7.9 in the following situation. Let $U$ be an open set, $h$ and $g$ functions analytic in $U$, and suppose that $a \in U$ is a zero of $g$, i.e., $g(a)=0$. Moreover suppose that $a$ is a \emph{simple zero}, i.e., $g^{\prime}(a) \not= 0$. Then $h/g$ has a simple pole at $a$ if and only if $h(a) \not= 0$, and if $h(a) \not=0$ then by way of L'Hospital's rule, 
\[
\textnormal{Res}(h/g, a)=\lim_{z \rightarrow a}\frac{(z-a)h(z)}{g(z)}=\frac{h(a)}{g^{\prime}(a)}\ .\]
Hence Theorems 7.8 and 7.9 imply 
\begin{lem}
Let $U$ be an open subset of $\mathbf{C}$, let $h$ and $g$ be analytic in U, and suppose g has only simple zeros in U. If $\gamma$ is a positively oriented Jordan contour in U which does not wind around a point in $\mathbf{C} \setminus U$ and does not pass through any of the zeros of g, and $a_1,\dots,a_n$ are the zeros of g in the interior of $\gamma$, then
\[
\frac{1}{2\pi i}\int_{\gamma}\ \frac{h(z)}{g(z)}\ dz=\sum_{k=1}^n\ \frac{h(a_k)}{g^{\prime}(a_k)}\ .\]
\end{lem}

In section 7, Theorems 7.2 and 7.4$(ii)$ will be deduced by integrating around rectangles a cleverly designed function analytic except for poles  and then applying Lemma 7.10.

\section{The Convergence of Fourier Series}
 Theorems 7.3 and 7.4$(i)$ will be deduced by appeals to certain facts concerning the convergence of Fourier Series. We therefore preface the proof proper with a brief discussion of Fourier series and their convergence. 

If $f$ is a real-valued function defined and integrable over $-\pi \leq x \leq  \pi$, then the \emph{Fourier series $S(f, x)$ of f} is the series defined by
\[
\frac{a_0}{2}+\sum_{n=1}^{\infty} (a_n \cos nx+b_n \sin nx), \]
where
\[
a_0=\frac{1}{ \pi} \int_{-\pi}^{ \pi} f(x) dx,\]
\[
a_n= \frac{1}{\pi} \int_{-\pi}^{ \pi} f(x) \cos nx\ dx,\]
\[
b_n= \frac{1}{ \pi} \int_{-\pi}^{ \pi} f(x) \sin nx\ dx,\ n=1,2,\dots;\]
$a_n$ and $b_n$ are called, respectively, the \emph{Fourier cosine and sine coefficients of f}.

Recall that a real-valued function $f$ defined on a closed and bounded interval $J=\{x: c\leq x\leq d\}$ of the real line is \emph{piecewise differentiable on J} if there is a finite partition of $\{x : c\leq x<d\}$ into subintervals such that for each subinterval $a \leq x<b$, there exists a function $g$ differentiable on $a \leq x \leq b$ such that $f\equiv g$ on  $a <x<b$. A function $f$ that is piecewise differentiable on $J$ is clearly piecewise continuous there, hence if  $c<x<d$ then the one-sided limits
\[
f_{\pm}(x)=\lim_{t \rightarrow x^{\pm}} f(t),\ \lim_{t\rightarrow c^+} f(t),\ \textrm{and}\ \lim_{t\rightarrow d^-} f(t) \] 
exist and are finite. It follows that if $f$ is defined on the entire real line, is periodic of period $2\pi$, and is piecewise differentiable on $-\pi \leq x \leq \pi$ then both one-sided limits of $f$ at any real number exist and are finite, and so the functions $f_{\pm}(x)=\lim_{t \rightarrow x^{\pm}} f(t)$ are both defined and real-valued on the entire real line. Figure 6 illustrates what the graph of a typical piecewise differentiable function looks like.

\begin{figure}[h]
  \centering
  \begin{tikzpicture}
    \draw[name path=xaxis, thick] (0,0) -- (12,0);
    \draw (0,0) -- (1,2) -- (2,0);
    \fill (2,0) circle (2pt);
    \draw (2,-1) circle (2pt) .. controls (4,-0.5) and (5, 2.5) .. (6,3) circle (2pt);
    \fill (6,2) circle (2pt);
    \draw (6,1) circle (2pt) -- (7, -1) circle (2pt);
    \fill (7,3) circle (2pt);
    \draw (7,3) .. controls (8,3) and (8,-2) .. (9,-1) -- (10,1) -- (12,1);
    
  \end{tikzpicture}
  \caption{A piecewise differentiable function}
  \label{fig:7.5}
\end{figure}
\begin{figure}[h]
  \centering
  \begin{tikzpicture}
    \draw[name path=xaxis, thick] (0,0) -- (12,0);
    \foreach \x/\y/\z in {1/2/3,3/4/5,5/6/7,7/8/9,9/10/11,11/12/13,13/14/15,15/16/17,17/18/19,19/20/21,21/22/23}
      \draw ({\x*.97^\x}, 0) -- ({\y*.97^\y}, 6/\x) -- ({\z*.97^\z}, 0);
      \draw (11.8,0.25) node {$\cdots$};
  \end{tikzpicture}
  \caption{A continuous, non-piecewise differentiable function}
 
  \label{fig:7.6}
\end{figure}

Piecewise differentiable functions exist in abundance and examples are very easy to come by; functions continuous but not piecewise differentiable on an interval are not difficult to construct either. Probably the simplest such example of the latter is to take a closed and bounded interval $J$ on the real line, and let $(a_n)_{n=1}^{\infty}$ be a strictly increasing sequence of elements of $J$ converging to the right-hand endpoint $d$ of $J$, with $a_1$ equal to the left-hand endpoint of $J$, say. On each closed interval with endpoints $a_ n$ and $a_{n+1}$ define the function $f_n$ which is 0 at $a_n $ and $a_{n+1}$, is $1/n$ at $(a_n+a_{n+1})/2$, and is linear and continuous on each of the closed intervals with left-hand endpoints $a_n$ and $(a_n+a_{n+1})/2$ and corresponding right-hand endpoints $(a_n+a_{n+1})/2$ and $a_{n+1}$, $n=1, 2,\dots$.  Then define $f$ on $J$ to equal $f_n$ on the closed interval with endpoints $a_ n$ and $a_{n+1}$ for $n=1, 2,\dots$, and set $f(d)=0$. The function $f$ is continuous on $J$, it is not differentiable at each point $a_n, n=2, 3,\dots$, and because $a_n\rightarrow d$, $f$ is not piecewise differentiable on $J$. In Figure 7, we indicate what the graph of a continuous, non-piecewise differentiable function constructed along these lines would look like.

We will use the following basic theorem on the convergence of Fourier series, a variant of which was first proved by Dirichlet [9] in 1829.
\begin{thm}
If f is defined on the real line $\mathbb{R}$, is periodic of period $2 \pi$, and is piecewise differentiable on $-\pi \leq x \leq \pi$, then the Fourier series $S(f, x)$ of f converges to
\[
\frac{f_+(x)+f_-(x)}{2},\ x\in \mathbb{R}.\]
In particular, if f is continuous at x then $S(f, x)$ converges to $f(x)$.
\end{thm}

\emph{Proof}. Let
\[
S_n(x)=\frac{a_0}{2}+\sum_{k=1}^{n}\  (a_k \cos kx+b_n \sin kx), \]
denote the $n$-th partial sum of the Fourier series of $f$. The key idea of this argument, due to Dirichlet, and used more or less in all convergence proofs of Fourier series, is to first express $S_n(x)$ in an integral form that is more amenable to an analysis of the convergence involved. Using the definition of the Fourier cosine and sine coefficients of $f$, we thus calculate that
\begin{eqnarray*}
S_n(x)&=&\frac{1}{\pi} \int_{-\pi}^{\pi} f(t)\Big(\frac{1}{2}+\sum_{k=1}^{n}\  (\cos kx \cos kt+ \sin kx \sin kt)\Big)dt\\
&=&\frac{1}{\pi} \int_{-\pi}^{\pi} f(t)\Big( \frac{1}{2}+\sum_{k=1}^n \cos k(x-t) \Big)dt. 
\end{eqnarray*}
Using the trigonometric identity
\[
\frac{1}{2}+\sum_{k=1}^n \cos k\theta=\frac{\sin\big(n+\frac{1}{2} \big)\theta}{2 \sin\big(\frac{\theta}{2}\big)},\]
it follows that
\[
S_n(x)=\frac{1}{\pi} \int_{-\pi}^{\pi} f(t) D_n(x-t)dt,\]
where
\[
D_n(\theta)=\frac{\sin\big(n+\frac{1}{2} \big)\theta}{2 \sin\big(\frac{\theta}{2}\big)}\]
is the \emph{Dirichlet kernel} of $S_n(x)$ (at $\theta=k\pi$, $k$ an even integer, we define $D_n(\theta)$ to be $n+\frac{1}{2}$, so as to make $D_n$ a function continuous on $\mathbb{R}$). Using the facts that $f$ and $D_n$ are of period $2\pi$ and $D_n$ is an even function, we can rewrite the integral formula for $S_n$ as
\[
S_n(x)=\frac{1}{\pi} \int_0^{\pi} \big(f(x+t)+f(x-t)\big)D_n(t)dt\ ,\ x\in \mathbb{R}.\]
If we now let $f \equiv 1$ in this equation and check that for this $f$, $S_n \equiv 1$, we find that
\[
1=\frac{2}{\pi}  \int_0^{\pi} D_n(t)dt.\]
After multiplying this equation by $\frac{1}{2}(f_{+}(x)+f_-(x))$ and then subtracting the equation resulting from that from the equation given by the above integral formula for $S_n$, it follows that
\begin{equation*}
S_n(x)-\frac{f_{+}(x)+f_-(x)}{2}=\frac{1}{\pi} \int_0^{\pi}\big(f(x+t)-f_+(x)+f(x-t)-f_-(x) \big) D_n(t)dt. \tag{6} 
\]

Now let
\[
\Xi(x, t)=\frac{f(x+t)-f_+(x)+f(x-t)-f_-(x)}{2 \sin\displaystyle \left(\frac{t}{2}\right)}\ ,\ 0<t \leq \pi.
\]
With an eye toward defining $\Xi(x, \cdot)$ at $t=0$ so as to make $\Xi(x, \cdot)$ right-continuous there, we study the behavior  of $\Xi(x,t)$ as $t \rightarrow 0^+$.  To that end, first rewrite $\Xi(x, t)$ as
\[
\Xi(x, t)=\left(\frac{f(x+t)-f_+(x)}{t}+\frac{f(x-t)-f_-(x)}{t} \right)\cdot \frac{t}{2 \sin \displaystyle \left(\frac{t}{2}\right)}\ , \ 0<t\leq \pi.\]
Because $f$ is periodic of period $2\pi$ and $f$ is piecewise differentiable on $-\pi \leq \xi \leq \pi$, there exists subintervals $a \leq \xi<b$, $b \leq \xi<c$ of the real line and functions $g$ and $h$ differentiable on $a \leq \xi \leq b$ and  $b \leq \xi \leq c$, respectively, such that $b \leq x<c$ and $f(\xi)$ equals, respectively, $g(\xi)$ or $h(\xi)$ whenever, respectively, $a<\xi<b$ or $b<\xi<c$.  A moment's reflection now confirms that
\[
\lim_{t \rightarrow 0^+}\frac{f(x+t)-f_+(x)}{t}=h^{\prime}(x),\]
\[
\lim_{t \rightarrow 0^+}\frac{f(x-t)-f_-(x)}{t}=\left\{\begin{array}{rl}
-h^{\prime}(x)\ ,& \textnormal{if}\ x>b,\\
-g^{\prime}(b)\ ,& \textnormal{if}\ x=b, 
\\\end{array}\right. 
\]
and so we conclude that $\lim_{t \rightarrow 0^+}\Xi(x, t)$ exists and is finite. If we take $\Xi(x, 0)$ to be this finite limit, then $\Xi(x, \cdot)$ is defined and piecewise continuous on $0 \leq t \leq \pi$. 

It follows that the functions
\[
\Xi(x, t) \sin \Big(n+\frac{1}{2}\Big)t,\ 0\leq t \leq \pi,\]
and
\[
\big(f(x+t)-f_+(x)+f(x-t)-f_-(x) \big) D_n(t),\ 0\leq t \leq \pi,\]
are both piecewise continuous on $0\leq t \leq \pi$ and agree on $0<t\leq \pi$. The latter function can hence be replaced by the former function in the integrand of the integral on the right-hand side of (6) to obtain the equation
\[
S_n(x)-\frac{f_{+}(x)+f_-(x)}{2}=\frac{1}{\pi} \int_0^{\pi}\Xi(x, t) \sin \Big(n+\frac{1}{2}\Big)t\ dt.\]
The conclusion of Theorem 7.11 will now follow if we prove that
\[
\lim_{n \rightarrow +\infty} \frac{1}{\pi} \int_0^{\pi} \Xi(x, t) \sin \Big(n+\frac{1}{2}\Big)t\ dt=0.\]

In order to do that, use the formula for the sine of a sum to write
\[
 \int_0^{\pi} \Xi(x, t) \sin \Big(n+\frac{1}{2}\Big)t\ dt=\int_{-\pi}^{\pi} \alpha(t) \sin nt\ dt+ \int_{-\pi}^{\pi} \beta(t) \cos nt\ dt,\]
where
\[
\alpha(t)=\left\{\begin{array}{rl}
0\ ,& \textnormal{if}\ -\pi \leq t<0,\\
\Xi(x, t) \cos\displaystyle \left(\frac{t}{2}\right)\ ,& \textnormal{if}\ 0 \leq t \leq \pi, 
\\\end{array}\right. 
\]
\[
\beta(t)=\left\{\begin{array}{rl}
0\ ,& \textnormal{if}\ -\pi \leq t<0,\\
\Xi(x, t) \sin \displaystyle \left(\frac{t}{2}\right)\ ,& \textnormal{if}\ 0 \leq t \leq \pi. 
\\\end{array}\right. 
\]
Because $\alpha$ and $\beta$ are functions piecewise continuous  on $-\pi\leq t\leq \pi$, our proof will be done upon verifying that if a function $\psi$ is piecewise continuous on  $-\pi\leq t\leq \pi$ and if $a_n$ and $b_n$ are the Fourier cosine and sine coefficients of $\psi$ then
\[
\lim_n a_n=0=\lim_n b_n\]
(This very important fact is known as the \emph{Riemann-Lebesgue lemma}). In order to see that, note that the set of functions $\left\{\frac{1}{\sqrt{2\pi}}\right\}\cup \left\{\frac{1}{\sqrt{\pi}}\cos nt: n \in [1, \infty)\right\} \cup \left\{\frac{1}{\sqrt{\pi}}\sin nt: n \in [1, \infty)\right\}$ is orthonormal with respect to the inner product defined by integration over the interval $-\pi\leq t\leq \pi$, hence a straightforward calculation using this fact shows that if $\sigma_n$ denotes the $n$-th partial sum of the Fourier series of $\psi$ then
\[
0 \leq \frac{1}{\pi} \int_{-\pi}^{\pi} (\psi-\sigma_n)^2\ dx=\frac{1}{\pi} \int_{-\pi}^{\pi} \psi^2\ dx-\Big(\frac{a_0^2}{2}+\sum_{k=1}^n(a_k^2+b_k^2)\Big), \]
and so
\[
\frac{a_0^2}{2}+\sum_{k=1}^n(a_k^2+b_k^2)\leq \frac{1}{\pi} \int_{-\pi}^{\pi} \psi^2\ dx<+\infty,\ \textrm{for all}\ n \in [1,\infty)\]
(this is \emph{Bessel's inequality}). Hence the series
\[
\frac{a_0^2}{2}+\sum_{n=1}^{\infty}(a_n^2+b_n^2)\]
converges, and so $a_n$ and $b_n$ both tend to 0 as $n \rightarrow +\infty$. $\hspace{5.5cm} \textrm{QED}$

\emph{Remarks}.  

$(1)$ Another very useful class of real-valued functions for which the conclusion of Theorem 7.11 is also valid is the functions $f$ that are defined on the whole real line, periodic of period $2\pi$, and are of \emph{bounded variation} on $-\pi \leq x \leq \pi$. This means that the supremum of the sums
\[
\sum_{i=1}^m |f(x_i)-f(x_{i-1})| \]
as $\{-\pi=x_0<x_1<\dots<x_m=\pi\}$ varies over all divisions of the interval $-\pi \leq x \leq \pi$ by a finite number of points $x_0,\dots,x_m$ is finite. A result from elementary real analysis asserts that if $f$ is of bounded variation on $-\pi \leq x \leq \pi$ then $f$ is the difference of two functions both of which are non-decreasing on $-\pi \leq x \leq \pi$, and so if $f$ is also defined on the entire real line and is periodic of period $2\pi$ then the one-sided limits $f_{\pm}(x)$ exist and are finite for all $x$. That Theorem 7.11 is valid for all functions of bounded variation on $-\pi \leq x \leq \pi$ is in fact what Dirichlet proved in his landmark paper [9]. This version of Theorem 7.11 also works in our proof of Theorems 7.3 and 7.4 \emph{infra}; we have proved Theorem 7.11 for piecewise differentiable functions because the argument which covers that situation is a bit more elementary than the one which suffices for functions of bounded variation. For a proof of the latter theorem, the interested reader should consult Zygmund [65], Theorem II.8.1. However, note well: a function that is piecewise differentiable need not be of bounded variation and a function of bounded variation is not necessarily piecewise differentiable.

$(2)$ It transpires that if $f$ is a complex-valued function defined on the integers which is periodic in the sense that for some integer $m>1$, $f(a)=f(b)$ whenever $a\equiv b$ mod $m$, then $f$ can be expanded in terms of a \emph{finite Fourier series} in complete analogy with the expansion into infinite Fourier series that we have just discussed. When this finite Fourier series expansion is applied to Gauss sums, another proof of the Law of Quadratic Reciprocity results, as we will see in section 9 below.

\section{Proof of Theorems 7.2, 7.3, and 7.4}

We begin this section with the proof of Theorem 7.2. Let $p\equiv 3$ mod 4: we must prove that 
\[
q(0, p/2)=\frac{\sqrt p}{\pi} \big(2-\chi_p(2)\big)L(1, \chi_p).\]

Toward that end, consider the functions $F(z)$ and $f(z)$ defined by
\[
F(z)=\sum_{0<j<p/2} \chi_p(j) \cos \left( \left(1-\frac{4j}{p} \right) \pi z\right),\]
\[
f(z)=\frac{\pi F(z)}{z \cos (\pi z)}\ .\]
We will prove Theorem 7.2 by integrating $f(z)$ around rectangles and then applying Lemma 7.10.

Note first that the numerator and denominator of $f$ are entire functions, then that the zeros of the denominator of $f$ occur at $z=0, z_n=(2n-1)/2,\ n \in \mathbb{Z}$, and that they are all simple. In order to apply Lemma 7.10 to  $f$, we therefore need to calculate
\[\frac{\pi F(z)}{\frac{d}{dz}(z \cos \pi z)}\ \ \textrm{at}\ z=0,\ z_n,\ n \in \mathbb{Z}.\]
At $z=0$ this is
\begin{equation*}
\pi F(0)=\pi \sum_{0<j<p/2} \chi_p(j)=\pi q(0, p/2), \tag{7}\]
and at $z=z_n$, it is
\[
(-1)^n\frac{F(z_n)}{z_n}.\]
We claim that
\begin{equation*}
(-1)^n\frac{F(z_n)}{z_n}=-\frac{\sqrt p}{2n-1}\ \chi_p(2n-1),\ n \in \mathbb{Z}. \tag{8}
\end{equation*}
In order to check this, we will first use the elementary identity
\begin{equation*}
\cos z=\frac{e^{iz}+e^{-iz}}{2} \tag{9}
\]
to calculate $F(z_n)$ as a Gauss sum. Toward that end, let $\alpha_j=1-(4j/p)$; then
\begin{eqnarray*}
\exp \left(i\frac{2n-1}{2} \alpha_j \pi \right)&=&\exp \left(i\frac{2n-1}{2} \pi \right) \exp \left(-i\frac{2 \pi j(2n-1)}{p} \right)\\
&=&(-1)^{n+1}i \exp \left(-i\frac{2 \pi j(2n-1)}{p} \right),
\end{eqnarray*}
and similarly
\[\
\exp \left(-i\frac{2n-1}{2} \alpha_j \pi \right)=(-1)^ni \exp \left(i\frac{2 \pi j(2n-1)}{p} \right),\]
Hence from (9) we deduce that
\begin{eqnarray*}
F(z_n)&=&\frac{(-1)^{n+1}i}{2} \sum_{0<j<p/2} \chi_p(j) \exp \left(-\frac{2 \pi ij(2n-1)}{p} \right)\\
&+&\frac{(-1)^ni}{2} \sum_{0<j<p/2} \chi_p(j) \exp \left(\frac{2 \pi ij(2n-1)}{p} \right).
\end{eqnarray*}
Observe now that the exponential factors here are periodic of period $p$ in the variable $j$ and, as $p\equiv 3\ \textrm{mod}\ 4, \chi_p(-1)=-1$. We can hence shift the summation in the first term on the right-hand side of this equation to express that term as
\[
\frac{(-1)^ni}{2} \sum_{p/2<j<p} \chi_p(j) \exp \left(\frac{2 \pi ij(2n-1)}{p} \right),\]
hence 
\begin{equation*}
 F(z_n)=\frac{(-1)^ni}{2} \sum_{0<j<p} \chi_p(j) \exp \left(\frac{2 \pi ij(2n-1)}{p} \right)=\frac{(-1)^ni}{2}\ G(2n-1, p).\tag{10}
\end{equation*}
Hence (10), (5), and Theorem 7.6 imply
\begin{eqnarray*}
(-1)^n\frac{F(z_n)}{z_n}&=&\frac{i}{2z_n}\ G(2n-1, p)\\
&=&\frac{i}{2n-1}\ \chi_p(2n-1)\ G(1, p)\\
&=&-\frac{\sqrt p}{2n-1}\ \chi_p(2n-1).
\end{eqnarray*}
This verifies (8).

Now for the contour around which we will integrate $f$. Let $\gamma_N$ denote the positively oriented rectangle centered at the origin, with horizontal side length $4pN$ and vertical side length $2 \sqrt N$, where $N$ is a fixed positive integer. $\gamma_N$ is clearly a Jordan contour, and the zeros of $z \cos \pi z$ inside $\gamma_N$ are 0 and $z_n, n \in [-pN+1, pN]$. Hence (7), (8), and Lemma 7.10 imply that
\begin{equation*}
\frac{1}{2\pi i}\int_{\gamma_N} f(z) dz=\pi q(0, p/2)- \sqrt p \sum_{n=-pN+1}^{pN} \frac{\chi_p(2n-1)}{2n-1}\ . \tag{11}\]
Because $\chi_p(-1)=-1$,
\[
 \frac{\chi_p(k)}{k}= \frac{\chi_p(-k)}{-k},\ \textrm{for all}\ k \in Z\setminus \{0\},\]
hence 
\begin{equation*}
\sum_{n=-pN+1}^{pN} \frac{\chi_p(2n-1)}{2n-1}=2 \sum_{n=1}^{pN}\ \frac{\chi_p(2n-1)}{2n-1}\ . \tag{12}\]

We claim that
\begin{equation*}
\lim_{N \rightarrow \infty} \frac{1}{2 \pi i} \int_{\gamma_N}\ f(z)\ dz=0. \tag{13}\]
Assuming this for a moment, we deduce from (11), (12) and (13) that
\begin{equation*}
q(0, p/2)=\frac{2 \sqrt p}{\pi}\ \lim_{N \rightarrow \infty}  \sum_{n=1}^{pN}\ \frac{\chi_p(2n-1)}{2n-1}\ . \tag{14}\]
In order to evaluate the limit on the right-hand side of (14), note that for each integer $M>1$,
\[
\frac{\chi_p(2)}{2} \sum_1^{M-1}\ \frac{\chi_p(k)}{k}=\sum_1^{M-1}\ \frac{\chi_p(2k)}{2k}\ ,\]
hence
\[
\sum_{1}^{2M-1}\ \frac{\chi_p(k)}{k}-\frac{\chi_p(2)}{2} \sum_1^{M-1}\ \frac{\chi_p(k)}{k}=\sum_{1}^{M} \frac{\chi_p(2n-1)}{2n-1}\ .\]
Letting $M \rightarrow \infty$ in this equation, we obtain 
\begin{eqnarray*}
\lim_{M \rightarrow \infty} \sum_{1}^{M} \frac{\chi_p(2n-1)}{2n-1}&=&\left(1-\frac{\chi_p(2)}{2} \right) \sum_1^{\infty}\ \frac{\chi_p(k)}{k}\\
&=&\left(1-\frac{\chi_p(2)}{2} \right) L(1, \chi_p).
\end{eqnarray*}
Hence from (14) it follows that
\[
q(0, p/2)=\frac{\sqrt p}{\pi} \big(2-\chi_p(2)\big) L(1, \chi_p),\]
the conclusion of Theorem 7.2.

We now need only to verify (13). This requires appropriate estimates of $f$ along the sides of $\gamma_N$. Consider first the function
\[
g(z)=\frac{\cos(\alpha \pi z)}{\cos(\pi z)},\ \alpha=1- \frac{4j}{p},\]
coming from a term of $F(z)/\cos \pi z$. Using (9), we calculate that for $z=x+iy$,
\[
|g(z)|^2=h(z)e^{2 \pi (\alpha-1)|y|},\ \textrm{where}\]
\[
h(z)=\frac{e^{-4 \pi \alpha |y|}+2e^{-2 \pi (\alpha-1)|y|} \cos 2x+1}{e^{-4 \pi |y|}+2e^{-2 \pi |y|} \cos 2x+1}\ .\]
We have
\[
\alpha-1 \leq -4/p,\ \textrm{for all}\ \alpha,\]
\[
h(z) <4/(1/2)=8,\ \textrm{for all}\ |y| \geq 1,\]
and so
\[
|g(z)|<2 \sqrt 2\ e^{-(4\pi/p)|y|},\ \textrm{for all}\ |y| \geq 1.\]
Hence
\begin{equation*}
\left| \frac{F(z)}{\cos (\pi z)} \right|<p \sqrt 2\  e^{-(4\pi/p)|y|},\ \textrm{for all}\ |y| \geq 1. \tag{15} \]
From (15) it follows that 
\begin{equation*}
|f(z)|<\frac{p \sqrt 2\  e^{-(4\pi/p) \sqrt N}}{\sqrt N},\  \textrm{for all $z$ on the horizontal sides $H_N$ of $\gamma_N$}. \tag{16}
\end{equation*}
By (15), $F(z)/ \cos (\pi z)$ is bounded on the vertical line Re $z=2p$. But  $F(z)/ \cos (\pi z)$ is periodic of period $2p$, hence there is a constant $C$, independent of $N$, such that 
\[
\left| \frac{F(z)}{\cos (\pi z)} \right| \leq C,\ \textrm{for all $z$ on the vertical sides $V_N$ of $\gamma_N$}.\]
Hence
\begin{equation*}
|f(z)| \leq \frac{C}{2pN},\ \textrm{for all $z$ on the vertical sides $V_N$ of $\gamma_N$} . \tag{17}\]
The estimates (6), (16), and (17) now imply that
\begin{eqnarray*}
\left|\int_{\gamma_N}\ f(z)\ dz \right|&\leq& \left|\int_{H_N}\ f(z)\ dz \right|+\left|\int_{V_N}\ f(z)\ dz \right|\\
&\leq&\frac{p \sqrt 2\  e^{-(4\pi/p) \sqrt N}}{\sqrt N} \cdot 8pN+ \frac{C}{2pN} \cdot 4 \sqrt N\\
&\rightarrow& 0,\ \textrm{as}\ N \rightarrow \infty.
\end{eqnarray*}
$\hspace{15.6cm} \textrm{QED}$

Now for the proof of Theorem 7.3. Here $p \equiv 1\ \textrm{mod}\ 4$ and we must show that 
\[
q(0, p/4)=\frac{\sqrt p}{\pi} L(1, \chi_{4p}).\]

The proof we give is based on the convergence of Fourier series and is very much in the same spirit as Dirichlet's original argument. Let $f$ be the function defined on $\mathbb{R}$ which is
\[
1,\ \textrm{for}\ 0 \leq x< \pi/2,\ 3 \pi/2< x \leq 2 \pi,\]
\[
0, \ \textrm{for}\ x=\pi/2,\ 3\pi/2,\]
\[
-1,\ \textrm{for}\ \pi/2<x <3\pi/2,\]
 and is periodic of period $2 \pi$. Clearly $f$ is piecewise differentiable on $-\pi \leq x \leq \pi$, hence calculation of the Fourier series of $f$ and Theorem 7.11 imply that
\begin{equation*}
f(x)=-\frac{4}{\pi} \sum_{n=1}^{\infty}\ \frac{(-1)^n}{2n-1} \cos (2n-1)x,\ -\infty<x<+\infty. \tag{18}\]

Next, let $\chi=\chi_{4p}=\chi_4 \chi_p$. Multiply the equation of Gauss sums
\[
G(2n-1, \chi_p)=\chi_p(2n-1)G(1, p),\]
from (5), by
\[
\frac{(-1)^n}{2n-1}\]
to obtain
\begin{equation*}
  \frac{(-1)^n \chi_p(2n-1)}{2n-1}G(1, p)=\sum_{j=1}^{p-1} \chi_p(j)  \frac{(-1)^n}{2n-1} \exp\left(\frac{2 \pi i(2n-1)j}{p} \right). \tag{19}\]
By virtue of Theorem 7.6,
\[
G(1, p)=\sqrt p,\]
and so, upon taking the real part of (19), we arrive at
\begin{equation*}
\sqrt p\  \frac{(-1)^n \chi_p(2n-1)}{2n-1}=\sum_{j=1}^{p-1} \chi_p(j)  \frac{(-1)^n}{2n-1} \cos \left( (2n-1) \cdot \frac{2 \pi j}{p} \right). \tag{20}\]
The definition of $\chi_4$ implies that
\[
\chi(k)=0,\ k\ \textrm{even},\]
\[
\chi(2n-1)=(-1)^{n+1} \chi_p(2n-1),\]
hence
\begin{equation*}
 \sum_{n=1}^{\infty} \frac{(-1)^n \chi_p(2n-1)}{2n-1}=-\sum_{k=1}^{\infty} \frac{\chi(k)}{k}=- L(1, \chi). \tag{21}\]
On the other hand, we have from (18) that
\begin{equation*}
-\frac{\pi}{4} f\left(\frac{2 \pi j}{p}\right)=\sum_{n=1}^{\infty} \frac{(-1)^n}{2n-1} \cos \left( (2n-1) \cdot \frac{2 \pi j}{p} \right),\ j=1,\dots,p-1. \tag{22}\]
Consequently, we can sum (20) from $n=1$ to $\infty$, interchange the order of summation on the right-hand side of the equation that results from that, and then use (21) and (22) to deduce that
\begin{equation*}
\sqrt p\ L(1, \chi)=\frac{\pi}{4} \sum_{j=1}^{p-1} f\left(\frac{2 \pi j}{p}\right) \chi_p(j). \tag{23} \]

The final step is to evaluate the right-hand side of (23). Note that
\[
0<j<\frac{p}{4}\ \textrm{iff}\ 0<\frac{2 \pi j}{p}<\frac{\pi}{2},
\]
\[
\frac{p}{4}<j<\frac{p}{2}\ \textrm{iff}\ \frac{\pi}{2}<\frac{2 \pi j}{p}<\pi,\]
\[
\frac{p}{2}<j<\frac{3p}{4}\ \textrm{iff}\ \pi<\frac{2 \pi j}{p}<\frac{3\pi}{2},\]
\[
\frac{3p}{4}<j<p\ \textrm{iff}\ \frac{3\pi}{2}<\frac{2 \pi j}{p}<2 \pi.\]
Hence, according to the definition of $f$,
\[
\textrm{right-hand side of (23)}=\frac{\pi}{4}\big(q(0, p/4)-q(p/4, p/2)-q(p/2, 3p/4)+q(3p/4, p) \big).\]
But by way of (3), 
\[
q(0, p/4)=q(3p/4, p),\]
\[
q(p/4, p/2)=-q(0, p/4),\]
\[
q(p/2, 3p/4)=-q(0, p/4),\]
and so
\[
\textrm{right-hand side of (23)}=\pi q(0, p/4),\]
whence
\[
q(0, p/4)=\frac{\sqrt p}{\pi} L(1, \chi).\]
$\hspace{15.6cm} \textrm{QED}$

The proof of Theorem 7.4 naturally divides into the verification of each of the statements $(i)$ and $(ii)$, and so we will verify each of these in turn

Begin with statement $(i)$. We have here that $p \equiv 1\ \textrm{mod}\ 4$, we want to verify that
\[
q(0, p/3)=\frac{\sqrt{3p}}{2\pi} L(1, \chi_{3p}),\]
and we will use Fourier series once more. Let $f$ be the function that is
\[
1,\ \textrm{for}\ 0 \leq x<2\pi /3,\ 4\pi /3< x \leq 2\pi,\]
\[
1/2,\  \textrm{for}\ x=2\pi /3,\ 4\pi /3,\]
\[
0,\ \textrm{for}\ 2\pi /3<x<4\pi /3,\]
and is periodic of period $2\pi$. Calculation of the Fourier series of $f$ and Theorem 7.11 imply that
\[
f(x)=\frac{2}{3}+\frac{\sqrt 3}{\pi}\sum_{n=1}^{\infty} \frac{a_n}{n} \cos nx,\ -\infty<x<+\infty,\]
where\[
a_n=\left\{\begin{array}{rl}0,& \textrm{if $3$ divides $n$,}\\
1,& \textrm{if $n \equiv 1\ \textrm{mod}\ 3$,}\\
-1,& \textrm{if $n \equiv 2\ \textrm{mod}\ 3$ .}\\\end{array}\right.
\]
Observe now that
\[
a_n=\chi_3(n),\ \textrm{for all}\ n,\]
and so
\begin{equation*}
f(x)=\frac{2}{3}+\frac{\sqrt 3}{\pi}\sum_{n=1}^{\infty} \chi_3(n) \frac{ \cos nx}{n},\ -\infty<x<+\infty. \tag{24} \]
Now multiply both sides of 
\[
G(n, \chi_p)=\chi_p(n) G(1, p)
\]
by
\[
\frac{\sqrt 3}{\pi n}\ \chi_3(n),\]
equate real parts in the equation which results, and then use Theorem 7.6, (24), and summation of the resulting terms from $n=1$ to $\infty$ as was done in the proof of Theorem 7.3 to obtain
\[
\frac{\sqrt{3p}}{\pi} L(1, \chi_{3p})=\frac{\sqrt 3}{\pi} G(1, p) \sum_{n=1}^{\infty} \frac{\chi_3(n) \chi_p(n)}{n}=\sum_{j=1}^{p-1}\left(f\left(\frac{2\pi j}{p} \right)-\frac{2}{3} \right) \chi_p(j).\]
Because $\sum_1^{p-1} \chi_p(j)=0$, the sum on the right is
\begin{eqnarray*}
\sum_{j=1}^{p-1}f\left(\frac{2\pi j}{p} \right) \chi_p(j)&=&\sum_{0<j<p/3} \chi_p(j)+\sum_{2p/3<j<p} \chi_p(j),\ \textrm{by definition of $f$},\\
&=&2\sum_{0<j<p/3} \chi_p(j),\ \textrm{because}\ \chi_p(-1)=1,\\
&=&2q(0, p/3).
\end{eqnarray*}
Hence
\[
q(0, p/3)=\frac{\sqrt{3p}}{2\pi} L(1, \chi_{3p}),\]
which is the conclusion of Theorem 7.4$(i)$.

The verification of statement $(ii)$ of Theorem 7.4 follows by either contour integration or the method of Fourier series along the same lines of argument that we have used before. We will outline the main ideas in the contour-integration proof and leave the rest of the details (and the proof via Fourier series) as an instructive exercise for the interested reader.

Let
\[
f(z)=\frac{\pi F(z)}{z\sin \pi\left(z+\displaystyle{\frac{1}{3}}\right)},\]
where
\[
F(z)=2i\sum_{0<j<p/3}\ \chi_p(j)\sin\left(\pi z+\displaystyle{\frac{\pi}{3}}-\displaystyle{\frac{6\pi jz}{p}}\right)+e^{-3\pi iz}\sum_{p/3<j<2p/3}\ \chi_p(j)e^{6\pi ijz/p}.\]
We must calculate
\[
\frac{\pi F(z)}{\displaystyle{\frac{d}{dz}}\left(z\sin \pi\left(z+\displaystyle{\frac{1}{3}}\right)\right)}\]
at $z=0, n-\frac{1}{3}, n\in \mathbb{Z}$. At $z=0$, we obtain
\[
\frac{\pi F(0)}{\sin\left(\displaystyle{\frac{\pi}{3}}\right)}=2\pi i\ q(0, p/3),\] 
and at $z=n-\frac{1}{3}$, we find that the value is
\[
\frac{3(-1)^n}{3n-1}F\left(n-\displaystyle{\frac{1}{3}}\right)=-\frac{3}{3n-1}G(3n-1,p)=-\frac{3}{3n-1}\chi_p(3n-1)G(1, p).\]

We now integrate $f(z)$ over a suitable rectangle $\gamma_N$ as in the proof of Theorem 7.2, apply Lemma 7.10 using the values of $\displaystyle{\frac{\pi F(z)}{\frac{d}{dz}\big(z\sin \pi\big(z+\frac{1}{3})\big)\big)}}$ that we have calculated, and then let $N\rightarrow +\infty$ as before to deduce that
\begin{eqnarray*}
0&=&2\pi i\ q(0, p/3)-3G(1, p)\sum_{n=-\infty}^{\infty}\ \frac{\chi_p(3n-1)}{3n-1}\\
&=&2\pi i\ q(0, p/3)-3G(1, p)\left(\sum_{n=1}^{\infty}\ \frac{\chi_p(3n-1)}{3n-1}+\sum_{n=0}^{\infty}\ \frac{\chi_p(3n+1)}{3n+1}\right)\\
&=&2\pi i\ q(0, p/3)-3G(1, p)\left(\sum_{n=1}^{\infty}\ \frac{\chi_p(n)}{n}-\sum_{n=1}^{\infty}\ \frac{\chi_p(3n)}{3n}\right),
\end{eqnarray*}
from which Theorem 7.4$(ii)$ follows easily after another application of Theorem 7.6. QED

\emph{Remarks}

$(1)$  Berndt's paper [1] is well worth studying; in it, he establishes many other results on positivity and negativity of the quadratic excess over various intervals: for example if $p \equiv 11,19\ \textrm{mod}\ 40$ then $q(0, p/10)>0$ and if $p \equiv 5\ \textrm{mod}\ 24$ then $q(3p/8, 5p/12)<0$. He also gives a very interesting discussion of the history of this problem with numerous pertinent references to the literature. 

$(2)$ Because the statements in Theorem 7.1 are so important in the theory of quadratic residues, elementary proofs of them would be of great interest. However, despite numerous efforts by many people during the intervening 175 years, those proofs continue to remain elusive.

\section{An Elegant Proof of Lemma 4.8 for Real Dirichlet Characters}
Because of the crucial role that it has played in the work done in this chapter, we will now prove Lemma 4.8 for real, non-principal Dirichlet characters $\chi$, i.e., we will show that if $\chi(\mathbb{Z})=[-1, 1]$ then $L(1, \chi) \not= 0$. The proof that we will present is due to de la Valle$\acute{\textrm{e}}$ Poussin [45] and is one of the most elegant arguments available for this. Following  Davenport [6], pp. 32-34, we start by recalling some well-known facts about analytic continuation of Riemann's zeta.

Following long tradition in these matters, we let $s=\sigma+it$ denote a complex variable. Proposition 5.5 implies that $\zeta(s)$ is analytic in $\sigma>1$; we want to show that $\zeta$ can be extended to a function analytic in $\sigma>0$ except for a simple pole at $s=1$. In order to do that, let $\sigma>1$ and then write 
\begin{eqnarray*}
\zeta(s)=\sum_{n=1}^{\infty}n^{-s}&=&\sum_{n=1}^{\infty} n(n^{-s}-(n+1)^{-s})\\
&=&s\sum_{n=1}^{\infty}n \int_n^{n+1} x^{-(s+1)} dx\\
&=&s\int_1^{\infty} [x]x^{-(s+1)} dx,
\end{eqnarray*}
where $[x]$ denotes the greatest integer which does not exceed $x$. Now let $[x]=x-(x)$, so that $(x)$ denotes the fractional part of $x$. This gives
\begin{equation*}
\zeta(s)=\frac{s}{s-1}-s\int_1^{\infty} (x)x^{-(s+1)}dx,\ \sigma>1. \tag{25} \]
The integral on the right is absolutely convergent for $\sigma>0$, uniformly convergent for $\sigma \geq \epsilon>0$, and all Riemann sums of the integrand are entire functions of $s$, hence this integral defines a function analytic in $\sigma>0$. Consequently the right-hand side of (25) extends $\zeta(s)$ to a function analytic in $\sigma>0$ except for a simple pole at $s=1$. It hence follows that
\begin{equation*}
\lim_{s \rightarrow 1^{+}} \zeta(s)= +\infty. \tag{26} 
\]
 
Next we observe that the proof of the Euler-Dedekind product expansion of the zeta function of an algebraic number field $F$ given in Theorem 5.8 can be easily modified to show that that product expansion is valid for all $\sigma>1$. If we hence take the number field $F$ in that theorem to be $\mathbb{Q}$, we deduce that $\zeta$ has the Euler-product expansion
 \[
 \zeta(s)=\prod_q (1-q^{-s})^{-1},\ \sigma>1.\]
 We also have from the estimate in the proof of (7) in Chapter 5 that the series
 \[
 \sum_q \log(1+q^{-\sigma})\]
 is absolutely convergent for $\sigma>1$.  Hence 
 \[
 |\zeta(s)| \geq  \prod_q (1+q^{-\sigma})^{-1}= \exp\left(-\sum_q \log(1+q^{-\sigma}) \right)>0,\ \sigma>1,\]
and so $\zeta(s)$ never vanishes in $\sigma>1$.

Now let $\chi$ be a real, non-principal Dirichlet character, and suppose by way of contradiction that $L(1, \chi)=0$. Because $L(s, \chi)$ is analytic in $\sigma>0$ (Lemma 7.5$(i)$) and $\zeta$ has a simple pole at $s=1$ as its only singularity in $\sigma>0$, it follows that
\[
L(s, \chi) \zeta(s)\ \textrm{is analytic in}\ \sigma>0.\]
Because $\zeta(2s) \not= 0$ in $\sigma>1/2$, the function
\[
\psi(s)=\frac{L(s, \chi)\zeta(s)}{\zeta(2s)}\]
is analytic in $\sigma>1/2$. Equation (26) implies that $\lim_{s \rightarrow \frac{1}{2}^{+}} \zeta(2s)= +\infty$, hence
\begin{equation*}
\lim_{s \rightarrow \frac{1}{2}^{+}}\psi(s)=0. \tag{27}\]
For $\sigma>1, \psi$ has the Euler product expansion
\[
\psi(s)=\prod_q \frac{(1-\chi(q)q^{-s})^{-1}(1-q^{-s})^{-1}}{(1-q^{-2s})^{-1}}\ .\]
Let $m=$ the modulus of $\chi$. $\chi(q)=0$ if and only if $q$ divides $m$, and the factor of the Euler product corresponding to such $q$ is
\[
1+q^{-s}.\]
If $\chi(q)=-1$ then the factor corresponding to $q$ is 
\[
\frac{(1+q^{-s})^{-1}(1-q^{-s})^{-1}}{(1-q^{-2s})^{-1}}=1.\]
Hence
\begin{equation*}
\frac{\psi(s)}{\displaystyle{\prod_{q|m} (1+q^{-s})}}=\prod_{q: \chi(q)=1} \frac{1+q^{-s}}{1-q^{-s}}\ ,\ \sigma>1.\tag{28}\]
(We note incidentally that $X=\{q:\chi(q)=1\}$ must be infinite; otherwise
\[
\psi(s)=\prod_{q|m} (1+q^{-s}) \prod_{q \in X} \frac{1+q^{-s}}{1-q^{-s}} \]
and this product has only a \emph{finite} number of factors, hence $\lim_{s \rightarrow \frac{1}{2}^{+}}\psi(s)>0$, contrary to (27)). 

Next let
\[
\phi(s)=\frac{\psi(s)}{\displaystyle{\prod_{q|m} (1+q^{-s})}}\ .\]
As the denominator here is nonzero in $\sigma>0$, $\phi(s)$ is analytic in $\sigma>1/2$, and (27) implies that
\begin{equation*}
\lim_{s \rightarrow \frac{1}{2}^{+}} \phi(s)=0. \tag{29}\]
We will now show that the product expansion (28) of $\phi$ implies that 
\begin{equation*}
\phi(s)>1\ \textrm{ for}\ \frac{1}{2}<s<2. \tag{30}\]
This contradicts (29) and so Lemma 4.8 follows for real non-principal characters. 

In order to verify (30), observe that
\[
 \frac{1+q^{-s}}{1-q^{-s}}=1+2\sum_{n=1}^{\infty} q^{-ns}\ , \sigma>1,\]
  hence we can use (28) to express $\phi(s)$ as a Dirichlet series
 \[
 \phi(s)= \sum_{n=1}^{\infty} \frac{a_n}{n^s}\ , \sigma>1,
 \]
 where the coefficients $a_n$ are calculated like so: $a_1=1$, and if $n \geq 2$ then
\[
a_n=\left\{\begin{array}{rl}2^{|\pi(n)|}\ ,& \textrm{if} \  \pi(n) \subseteq \{q: \chi(q)=1\},\\
0\ ,& \textrm{otherwise}. 
\\\end{array}\right. 
\]
In particular, $a_n \geq 0$, for all $n$.

Because $\phi$ is analytic in $\sigma>\frac{1}{2}$, it is a consequence of Theorem 7.7 that $\phi$ has a Taylor series expansion centered at 2 with radius of convergence at least $\frac{3}{2}$, i.e.,
\[
\phi(s)=\sum_{m=0}^{\infty} \frac{\phi^{(m)}(2)}{m!} (s-2)^m ,\ |s-2|< \frac{3}{2}.\] 
We can calculate $\phi^{(m)}(2)$ by term-by-term differentiation of the Dirichlet series: this series is locally uniformly convergent in $\sigma>1$ and so we can apply the theorem which asserts that a series of functions analytic in an open  set $U$ and locally uniformly convergent there has a sum that is analytic in $U$ and the derivative can be calculated by term-by-term differentiation of the series. The result is
\[
\phi^{(m)}(2)=(-1)^m \sum_{n=1}^{\infty} \frac{a_n(\log n)^m}{n^2}=(-1)^mb_m\ , b_m \geq 0.\]
Hence
\[
\phi(s)=\sum_{m=0}^{\infty} \frac{b_m}{m!} (2-s)^m,\ |s-2|<\frac{3}{2}.\]
If $\frac{1}{2}<s<2$ then all terms of this series are non-negative, hence $\phi(s)\geq \phi(2)>1$ for $\frac{1}{2}<s<2. \hspace{14.5cm} \textrm{QED}$

\section{A Proof of Quadratic Reciprocity via Finite Fourier Series}

In this section we will apply a finite Fourier series expansion to powers of the Gauss sums
\[
\sum_{n=0}^{p-1}\ \chi_p(n)\zeta^{nm},\ \zeta=\exp(2\pi i/p),\]
to derive once more the Law of Quadratic Reciprocity. This proof is due to H. Rademacher, and we follow his account of it from [46].

The Fourier series expansion to which we are referring is given in the following lemma:
\begin{lem}
If $m>1$ is an integer, $F(t)$ is a $($complex-valued$)$ function defined on $\mathbb{Z}$ of period m, and $\zeta_m=\exp(2\pi i/m)$, then
\[
F(t)=\sum_{n=0}^{m-1}\ a(n)\zeta_m^{nt},\]
where
\[
a(n)=\frac{1}{m}\sum_{t=0}^{m-1}\ F(t)\zeta_m^{-nt}.\]
\end{lem}

In order to see why the expansion of $F$ in Lemma 7.12 can be viewed as a discrete Fourier series, we need to consider the complex version of the real-valued Fourier series that was defined and studied in section 6. To that end, for a \emph{complex-valued} function $f(x)$ defined on the interval $-\pi\leq x \leq \pi$ of the real line, define the (complex) \emph{Fourier series of f} as
\begin{equation*}
\sum_{n=-\infty}^{\infty}\ f_n e^{inx},\tag{$*$}\]
where
\[
f_n=\frac{1}{2\pi}\int_{-\pi}^{\pi}\ f(x)e^{-inx},\ n=0, \pm1, \pm2,\dots,\]
are the (complex) \emph{Fourier coefficients of f}. When $f$ is real-valued, the Fourier series $S(f)$ that we defined in section 6 can be recovered from the series $(*)$ by grouping the terms with indices $-n$ and $n$ together and formally rewriting the series as
\[
f_0+\sum_{n=1}^{\infty} (f_ne^{inx}+f_{-n}e^{-inx}).\]
If the Fourier series of $f$ converges to $f(x)$ then we can write
\begin{equation*}
f(x)=\sum_{n=-\infty}^{\infty}\ f_n e^{inx}.\tag{31}\]

We now proceed to discretize this continuous picture. First, replace the continuous variable $x$ by the discrete variable $t=0, 1, 2,\dots,m-1$. Second, substitute the exponential function $\zeta_m^t=\exp(2\pi it/m)$ and the function $F(t)$ in Lemma 7.12 for the exponential function $e^{ix}$ and the function $f(x)$, respectively. The analog of the Fourier coefficient $f_n$, which is the mean value of the function $f(x)e^{-inx}$ over the interval $-\pi\leq x \leq \pi$, is then the mean $a(n)$ of the values at $t=0, 1, 2,\dots,m-1$ of the function $ F(t)\zeta_m^{-nt}$. It follows that the discrete version of (31), in other words, a \emph{finite} Fourier series expansion of $F(t)$, is precisely the conclusion of Lemma 7.12. 

\emph{Proof of Lemma} 7.12.

By direct substitution using the stated formula for $a(n)$, we compute that
\begin{eqnarray*}
\sum_{n=0}^{m-1}\ a(n)\zeta_m^{nt}&=&\frac{1}{m}\sum_{n=0}^{m-1}\left(\sum_{s=0}^{m-1}\ F(s)\zeta_m^{-ns}\right)\zeta_m^{nt}\\
&=&\frac{1}{m}\sum_{s=0}^{m-1}\ F(s)\sum_{n=0}^{m-1}\ \zeta_m^{n(t-s)}\\
&=&F(t),
\end{eqnarray*}
where the last line follows by the same calculation that we used to derive equation (21) in section 9 of Chapter 3.$\hspace{11.5cm}\ \textrm{QED}$

N.B. We will refer to the coefficients $a(n)$ in the conclusion of Lemma 7.12 as the \emph{Fourier coefficients of F}. We also note that it follows from the fact that $F(t)$ and $\zeta_m^{\pm nt}$ are periodic of period $m$ that the sums in the formulae for $F(t)$ and $a(n)$ in Lemma 7.12 can be taken over any complete set of ordinary residues modulo $m$ without a change in their values; we will use this observation without further reference in the rest of this section.

Let $p$ and $q$ be distinct odd primes, with $\zeta=\exp(2\pi i/p)$. In order to make the writing less cumbersome, we change the notation slightly as follows: for each $t\in \mathbb{Z}$,
let
\[
G(\zeta^t)=\sum_{n=0}^{p-1}\ \chi_p(n)\zeta^{tn}.\]

Recall from Theorem 3.14 that
\[
G(\zeta)^2=(-1)^{\frac{1}{2}(p-1)}p,\]
and so 
\[
G(\zeta)^{q-1}=\big(G(\zeta)^2\big)^{\frac{1}{2}(q-1)}=(-1)^{\frac{1}{2}(p-1)\frac{1}{2}(q-1)}p^{\frac{1}{2}(q-1)}.\]
We conclude from Euler's criterion (Theorem 2.5) that
\begin{equation*}
G(\zeta)^{q-1}\equiv (-1)^{\frac{1}{2}(p-1)\frac{1}{2}(q-1)}\chi_q(p)\ \textrm{mod}\ q.\tag{32}\]
The LQR will follow from (32) if we can prove that 
\begin{equation*}
G(\zeta)^{q-1}\equiv \chi_p(q)\ \textrm{mod}\ q,\tag{33}\]
because then (32) and (33) will imply that
\[
(-1)^{\frac{1}{2}(p-1)\frac{1}{2}(q-1)}\chi_q(p)\equiv \chi_p(q)\ \textrm{mod}\ q,\]
and this congruence must in fact be an equality since the difference of both sides must be either 0 or  $\pm 2$ and also divisible by the odd prime $q$. In order to deduce the LQR, we must therefore verify (33).

This verification will be done by expanding the function $G(\zeta^t)^q$ as a finite Fourier series by way of Lemma 7.12 with $m=p$. Thus, we calculate that
\begin{equation*}
G(\zeta^t)^q=\sum_{u\ \textrm{mod}\ p}\ a_q(u)\zeta^{ut},\tag{34}\]
with Fourier coefficients
\begin{eqnarray*}
a_q(u)&=&\frac{1}{p}\sum_{v\ \textrm{mod}\ p}\ G(\zeta^v)^q\zeta^{-vu}\\
&=&\frac{1}{p}\sum_{v\ \textrm{mod}\ p}\sum_{m_1\ \textrm{mod}\ p}\ \chi_p(m_1)\zeta^{vm_1}\cdots\sum_{m_q\ \textrm{mod}\ p}\ \chi_p(m_q)\zeta^{vm_q}\zeta^{-vu}\\
&=&\frac{1}{p}\sum_{v\ \textrm{mod}\ p}\sum_{\substack{v, m_1,\dots,m_q\\ \textrm{mod}\ p} }\ \chi_p(m_1m_2\cdots m_q)\zeta^{v(m_1+m_2+\cdots +m_q-u)}\\
&=&\frac{1}{p}\sum_{\substack{v, m_1,\dots,m_q\\ \textrm{mod}\ p} }\ \chi_p(m_1m_2\cdots m_q)\sum_{v\ \textrm{mod}\ p}\zeta^{v(m_1+m_2+\cdots +m_q-u)}.
\end{eqnarray*}
Thus
\begin{equation*}
a_q(u)=\sum_{\substack{m_j\ \textrm{mod}\ p\\m_1+m_2+\cdots +m_q\equiv u\ \textrm{mod}\ p}}\chi_p(m_1m_2\cdots m_q).\tag{35}\]

We will now use the equation
\begin{equation*}
G(\zeta^t)=\chi_p(t)G(\zeta),\ t\in \mathbb{Z},\tag{36}\]
from Lemma 3.15 to calculate $a_q(u)$ in a different way. Because $q$ is odd, it follows from (36) that
\[
G(\zeta^v)^q=\chi_p(v)G(\zeta)^q.\]
Hence
\begin{eqnarray*}
a_q(u)&=&\frac{1}{p}\ G(\zeta)^q\sum_{v\ \textrm{mod}\ p}\chi_p(v)\zeta^{-vu}\\
&=&\frac{1}{p}\ G(\zeta)^qG(\zeta^{-u})\\
&=&\frac{1}{p}\ G(\zeta)^q\chi_p(u)G(\zeta^{-1}),
\end{eqnarray*}
where in the last line, we have applied (36) again, this time with $\zeta$ replaced by $\zeta^{-1}$. Comparing the last formula with its own special case $u=1$, we find that
\begin{equation*}
a_q(u)=\chi_p(u)a_q(1).\tag{37}\]

Next, insert (37) into (34) to deduce that
\[
G(\zeta^t)^q=a_q(1)\sum_{u\ \textrm{mod}\ p}\ a_q(u)\zeta^{ut}=a_q(1)G(\zeta^t).\]
Because $G(\zeta)\not=0$ it follows from this equation and another application of (37) that
\begin{equation*}
G(\zeta)^{q-1}=a_q(1)=\chi_p(q)a_q(q).\tag{38}\]

We will use (35) and (38) to deduce (33); that will complete this proof of the LQR. Toward that end, we substitute (35) into (38) to obtain
\begin{equation*}
G(\zeta)^{q-1}=\chi_p(q)\sum_{\substack{m_j\ \textrm{mod}\ p\\m_1+m_2+\cdots +m_q\equiv q\ \textrm{mod}\ p}}\chi_p(m_1m_2\cdots m_q).\tag{39}\]
This equation must be reduced modulo $q$. In order to do that, we first examine the possibilities for the set of summation variables $m_1\dots,m_q$. Suppose that 
\[
m_1\equiv m_2\equiv \dots \equiv m_q\ \textrm{mod}\ p.\]
This requires that 
\begin{equation*}
m_1+m_2+\cdots +m_q\equiv qm_j\equiv q\ \textrm{mod}\ p,\]
hence $m_j\equiv 1$ mod $p$, yielding only the one summand
\begin{equation*}
\chi_p(1)=1.\tag{40}\]
All other solutions $m_1,\dots, m_q$ of 
\[
m_1+m_2+\cdots +m_q\equiv q\ \textrm{mod}\ p\]
must contain non-congruent integers. If such a solution is cyclically permuted, then another solution of this type will be obtained. Indeed, a cyclic permutation of 
\[
m_1, m_2,\dots m_q\]
can be expressed as
\[
m_{1+s}, m_{2+s},\dots,m_{q+s},\]
for some $s$, $1\leq s<q$, with the subscripts here being taken modulo $q$. If this set of solutions is of the same kind as the previous one, then
\[
m_j\equiv m_{j+s}\ \textrm{mod}\ p.\]
Hence, upon successively setting $j=s, 2s,\dots,(q-1)s$, we obtain 
\[
m_s\equiv m_{2s}\equiv\cdots \equiv m_{qs}\ \textrm{mod}\ p,\]
where the subscripts here form a complete set of ordinary residues mod $q$, and this puts us in the case that we have already considered. It follows that those solutions of $m_1+m_2+\cdots m_q\equiv q\ \textrm{mod}\ p$ in which non-congruent numbers appear produce the term
\[
q\chi_p(m_1m_2\cdots m_q)\]
and these terms sum to a number congruent to 0 mod $q$. Consequently,  when the sum 
\[
\sum_{\substack{m_j\ \textrm{mod}\ p\\m_1+m_2+\cdots +m_q\equiv q\ \textrm{mod}\ p}}\chi_p(m_1m_2\cdots m_q)\]
is reduced modulo $q$, the only nonzero term is the single term (40). Therefore by (39),
\[
G(\zeta)^{q-1}\equiv \chi_p(q)\ \textrm{mod}\ q,\]
which verifies (33), and hence also the LQR.$\hspace{7.5cm}\ \textrm{QED}$ 

In analogy with the Fourier-series expansion of a function $f(x)$ of a real variable $x$, the functions $\zeta_m^{nt}$ in Lemma 7.12 can be viewed as the harmonics of the function $F(t)$, with the amplitude of the harmonics given by the Fourier coefficients $a(n)$. The proof of quadratic reciprocity that we have presented here can thus be interpreted as showing that the LQR comes from the fact that, modulo the prime $q$, the dominant harmonic of $G(\zeta^t)^q$ is $\zeta^{qt}$ with amplitude congruent to 1 modulo $q$. 

\chapter{Dirichlet's Class-Number Formula}

Although De la Valle$\acute{\textrm{e}}$ Poussin's proof of Lemma 4.8 is elegant and efficient, it fails to explain exactly why the values of $L$-functions at $s=1$ are positive, and consequently the actual reason why the sums in Theorem 7.1 turn out to be positive is not yet clear. Often in number theory, integers which occur in interesting situations are in fact positive because they count something, and it transpires that this is the case in Theorem 7.1. The sums there in fact count equivalence classes of binary quadratic forms, or, to say the same thing in a different way, ideal classes in quadratic number fields. In order to see this, we will now  present a second proof of Lemma 4.8 that uses Dirichlet's famous class-number formula, which formula calculates the value of $L(1, \chi)$ when $\chi$ is a real (primitive) Dirichlet character in terms of the number of certain equivalence classes of quadratic forms. As we did in our first proof of Lemma 4.8, we will follow the exposition for this as set forth in Davenport [6].

We begin things in this chapter by using some structure theory of Dirichlet characters in the first section to reduce the problem of proving that $L(1, \chi)\not=0$ for an arbitrary real non-principal Dirichlet character $\chi$ to proving that that is true for an arbitrary real \emph{primitive} character.  This highlights the importance of real primitive characters in our discussion, and so the facts concerning the structure of those characters that will be required   are recorded in section 2. In order to precisely state what the class-number formula asserts, we need to explain what a class number is, and that is the subject of sections 3 and 4. Section 3 recalls the fundamental equivalence relation defined on the set of all primitive and irreducible quadratic forms with a given discriminant that we first discussed in section 12 of Chapter 3, and section 4 uses the equivalence classes coming from this equivalence relation and some information on the representation of integers by quadratic forms to define the class number. Dirichlet's class-number formula is stated precisely and proved in section 5. Dirichlet's formula calculates the value at $s=1$ of the $L$-function of a primitive Dirichlet character as the product of a canonically determined positive constant and the class number of an appropriate set of quadratic forms, thereby providing the definitive explanation of the positivity of that value of the $L$-function. In section 7, we reformulate the class-number formula in terms of the class number determined by the ideal classes in quadratic number fields, and we also give class-number formulae for the sums of the Legendre symbols in Theorem 7.1. As an extra bonus from the theory discussed in section 1, we give in the last section of this chapter our seventh and final proof of the Law of Quadratic Reciprocity.

\section{Some Structure Theory for Dirichlet Characters}

We begin by discussing some useful facts regarding the structure of Dirichlet characters. Let $\chi$ be a Dirichlet character of modulus $b$, and let $d$ be a positive divisor of $b$. The number $d$ is an \emph{induced modulus of} $\chi$ if $\chi(n)=1$ whenever $\gcd(n,b)=1$ and $n \equiv 1$ mod $d$. In other words, $d$ is an induced modulus of $\chi$ if $\chi$ acts like a character with modulus $d$ on the integers in an ordinary residue class of 1 mod $d$ which are relatively prime to $b$. One can show without too much difficulty that a positive divisor $d$ of $b$ is an induced modulus of $\chi$ if and only if there exists a Dirichlet character $\xi$ of modulus $d$ such that
\begin{equation*}
\chi(n)=\chi_1(n)\xi(n), \textrm{for all}\ n \in \mathbb{Z}, \tag{1}
\]
where $\chi_1$ denotes the principal character of modulus $b$. Hence the induced moduli of $\chi$ are precisely the moduli of Dirichlet characters which ``induce"  the character  $\chi$ in the sense of (1). It is a straightforward consequence of (1) that $\chi$ is non-principal if and only if $\xi$ is non-principal. Because $d$ is a factor of $b$, it also follows from (1) and the Dirichet product formula for $L$-functions (Lemma 7.5$(ii)$) that the $L$-function of $\chi$ is a factor of the $L$-function of $\xi$.

The modulus $b$ is clearly an induced modulus of $\chi$, and if $b$ is the only positive divisor of $b$ which is an induced modulus, $\chi$ is said to be \emph{primitive}. If $d$ is the smallest positive divisor of $b$ that is an induced modulus of $\chi$ i.e., $d$ is the \emph{conductor of} $\chi$, then one can show that the factor $\xi$ in (1) is a primitive character modulo $d$. Taking the character $\xi$ in (1) to be the character modulo the conductor of $\chi$, it follows from (1) that $\chi$ is real if and only if $\xi$ is real, and so from this fact and our observations in the previous paragraph, we conclude that in order to show that the value at $s=1$ of the $L$-function of an arbitrary real and non-principal Dirichlet character is nonzero, and hence positive by the proof of Lemma 7.5$(iii)$, it suffices to prove that this is true for all real \emph{primitive} Dirichlet characters. It hence behooves us to take a closer look at the structure of those characters, and that is what we will do in the next section.

\section{The Structure of Real Primitive Dirichlet Characters}

The structure of a real primitive Dirichlet character is of a very particular type; in this section we will describe that structure precisely. It is here that we will also begin to see how the classical theory of quadratic forms enters the picture. We begin by noting that every Legendre symbol is real, and because the modulus is prime, they are also primitive. However, there are also real primitive characters of modulus 4 and 8, and so we will describe those next.

For the modulus 4, there is only one non-principal character, defined by
\[
\chi_4(n)=\left\{
\begin{array}{rl}1,& \textrm{if $n\equiv $ 1 mod 4,}\\
-1,& \textrm{if $n\equiv -1$ mod 4,}\\\end{array}
\right.
\]
and $\chi_4$ is obviously real and primitive. For the modulus 8, there are only two real primitive characters, defined by
\[
\chi_8(n)=\left\{\begin{array}{rl}1,& \textrm{if $n\equiv \pm 1$ mod 8,}\\
-1,& \textrm{if $n\equiv \pm 3$ mod 8.}\\\end{array}\right.
\]
and 
\[
\chi_{-8}(n)=\left\{\begin{array}{rl}1,& \textrm{if $n\equiv$ 1 or 3 mod 8,}\\
-1,& \textrm{if $n\equiv  -1$ or $-3$ mod 8.}\\\end{array}\right.
\]
In fact, we have that $\chi_{-8}=\chi_4 \chi_8$.

It can be shown that the prime-power moduli for which real primitive characters exist are: any odd prime $p$, with corresponding character $\chi_p$ (the Legendre symbol of $p$), 4, with corresponding character $\chi_4$, and 8, with corresponding characters $\chi_{\pm 8}$. Following Davenport, we will call the moduli 4, 8, $p, p>2$, the \emph{basic moduli} and the corresponding real primitive characters the \emph{basic characters}. More generally, the following theorem asserts that the moduli for which real primitive characters exist are determined by the products of basic moduli and the real primitive characters are determined by the products of the basic characters which correspond to the basic moduli. For a proof of the theorem, we refer the interested reader to Davenport [6], pp. 38-40.
\begin{thm}
Let
\[
\mathcal{M}=\{4, 8\}\cup \{p\in P: p>2\}\]
denote the set of basic moduli. If $b>1$ is an integer then b is the modulus of a real primitive Dirichlet character if and only if b is the product of relatively prime factors from $\mathcal{M}$. If b is such a modulus and B is the set of basic moduli into which b factors, then the real primitive Dirichlet characters of modulus b are given precisely by the product (or products)
\[
\prod_{n\in B}\ \chi_n,\]
where $\chi_n$ is a basic character of modulus n.
\end{thm}
\vspace{0.1cm}

It follows from Theorem 8.1 that if $b$ is a modulus as specified in that theorem then there is exactly one real primitive Dirichlet character of modulus $b$, unless 8 is a factor of $b$, in which case there are exactly two such characters.  Real primitive characters can also be used to give yet another proof of the Law of Quadratic Reciprocity, as we will see in section 7.

If $\chi$ is a real primitive Dirichlet character of modulus $b$ then we let $d_{\chi}$ denote the parameter
\[
\chi(-1)b.\]
It is a consequence of Theorem 8.1 that as $\chi$ varies over all real primitive characters, $d_{\chi}$ varies over all integers, positive or negative, that are either

(a) square-free and congruent to 1 mod 4, or

(b) of the form $4n$, where $n$ is square-free and $n\equiv$ 2 or 3 mod 4

\noindent (Davenport [6], pp. 40-41).

The range of the parameter $d_{\chi}$ marks the first occurrence of what we will eventually see as an intimate connection between real primitive Dirichlet characters and the classical theory of quadratic forms. If
\[
ax^2+bxy+cy^2\]
is a quadratic form with integer coefficients $a, b, c$, then the discriminant of the form is the familiar algebraic invariant $b^2-4ac$. In this theory, the forms which are not a product of linear factors are the primary objects of study, and in this case, the discriminant is not a perfect square, and is hence a non-square integer which is congruent to 0 or 1 mod 4. A discriminant  $b^2-4ac$ is said to be a fundamental discriminant if $\gcd(a, b, c)=1$. It is then not difficult to see that the set of fundamental discriminants consists precisely of the integers which satisfy the conditions (a) and (b) above, i.e., the set of fundamental discriminants coincides with the range of $d_{\chi}$ as $\chi$ varies over all real primitive Dirichlet characters.

The range of $d_{\chi}$ can also be used to make an important connection between Dirichlet characters and quadratic number fields. As we saw in section 11 of Chapter 3, a quadratic number field is generated over $\mathbb{Q}$ by the square root of a square-free integer $m$, and the discriminant of the quadratic field is defined as $m$ or $4m$ if $m$ is, respectively, congruent to 1 mod 4 or congruent to 2 or 3 mod 4. Consequently, the range of $d_{\chi}$ also coincides with the  discriminants of the set of all quadratic number fields.
  
The proof of Dirichlet's class-number formula that we will present makes use of a canonical representation of the real primitive Dirichlet characters that is parametrized by the set of fundamental discriminates coming from the theory of quadratic forms. In order to explain what this parameterization is, let $d$ denote a fundamental discriminant, which we understand at this point as simply an integer that satisfies either of the conditions (a) or (b) described above. It can be shown (Landau [35], pp. 221-222) that $d$ has a unique factorization into a product of relatively prime factors taken from the set of \emph{primary discriminants}
\[
-4,\ 8,\ -8,\ (-1)^{\frac{1}{2}(p-1)}p,\ \textrm{$p$ an odd prime};\]
we call this factorization the \emph{primary factorization of} $d$. To each primary discriminant, we associate one of the basic characters introduced in this section: to the integer $-4$ we associate the character $\chi_4$, to 8 and $-8$ the character $\chi_8$ and $\chi_{-8}$, respectively, and to each integer $(-1)^{\frac{1}{2}(p-1)}p$, we associate the Legendre symbol $\chi_p$. If $D$ is the set of primary discriminants in the primary factorization of $d$, we let $\chi(d)$ denote the real primitive character of modulus $|d|$ given by the product
\[
\chi(d)=\prod_{n\in D}\ \chi_n,\]
where $\chi_n$ is the character that we have associated with $n$ above. We will call $\chi(d)$ \emph{the character determined by d}. A proof of the following theorem can be found in Davenport [6], pp. 38-41:
\begin{thm}
The map $d\rightarrow \chi(d)$ is a bijection of the set of fundamental discriminants onto the set of all real primitive Dirichlet characters, and its inverse map is given by $\chi \rightarrow \chi(-1)(\textnormal{modulus of }\chi)$.
\end{thm}

The correspondence given in Theorem 8.2 is a special case of the correspondence  from $\mathbb{Z}\setminus \{0\}$ into the set of all real Dirichlet characters defined by the \emph{Kronecker symbol}, a generalization of the Legendre and Jacobi symbols, which is a very useful device in the study of real Dirichlet characters. We will make no further use of the Kronecker symbol in these notes; for its definition and basic properties, the interested reader is referred to either Landau [35], Definition 20 and Theorems 96-101 or Cohen [2], section 2.2.2.

We now have at our disposal all of the information about real primitive Dirichlet characters that we need in order to state and prove Dirichlet's class-number formula. In addition to that information, we will also require some results from the classical theory of quadratic forms, to which we turn next.

\section{Elements of the Theory of Quadratic Forms}

Let $d$ be a fundamental discriminant, i.e., an integer which is not a square and which is either square-free and congruent to 1 mod 4 or of the form $4n$, where $n$ is square-free and is congruent to either 2 or 3 mod 4. The classes to which Dirichlet's class-number formula originally referred are the equivalence classes determined by the basic equivalence relation on quadratic forms defined by modular substitutions, which we discussed in section 12 of Chapter 3. For the reader's convenience, we will recall the essential features of that equivalence relation.

We begin with some convenient notation. If $(a, b, c)$ is an ordered triple of integers then $[a, b, c]$ will denote the quadratic form $ax^2+bxy+cy^2$ and if $d$ is a fundamental discriminant  then $\mathcal{Q}(d)$ will denote the set of all irreducible and primitive quadratic forms of discriminant $d$, i.e., the set of all quadratic forms $[a, b, c]$ of discriminant $d$ which do not factor into  linear forms with integer coefficients and for which $\gcd(a,b, c)=1$. In particular, if $[a, b, c]\in \mathcal{Q}(d)$ then $acd\not=0$. On the set $\mathcal{Q}(d)$ we will declare that two forms $q(x, y)=ax^2+bxy+cy^2$ and  $q_1(X, Y)=a_1X^2+b_1XY+c_1Y^2$ in $\mathcal{Q}(d)$ are equivalent if there is a linear transformation defined by
\[
x=\alpha X+\beta Y,\ y=\gamma X+\delta Y,\]
where $\alpha, \beta, \gamma,$ and $\delta$ are integers satisfying $\alpha \delta-\beta \gamma=1$, such that 
\[
q(\alpha X+\beta Y,\ \gamma X+\delta Y)=q_1(X, Y).\]
These transformations are called modular substitutions, and each modular substitution maps $\mathcal{Q}(d)$ bijectively onto $\mathcal{Q}(d)$. As we pointed out in section 12 of Chapter 3, it follows from  classical results of Lagrange that the number of equivalence classes is finite. There is always at least one form in $\mathcal{Q}(d)$, called the principal form, defined by
\[
x^2-\frac{1}{4}dy^2,\ \textrm{if $d\equiv 0$ mod 4,}\]
or
\[
x^2+xy-\frac{1}{4}(d-1)y^2,\ \textrm{if $d\equiv 1$ mod 4,}\]
hence the number of equivalence classes is a positive integer.

An important parameter which enters into Dirichlet's formula for the class number is determined by the \emph{automorphs} of a form in $\mathcal{Q}(d)$. These automorphs are the modular substitutions which leave a given form invariant. There are always two automorphs for every form: the trivial substitution $x=X,\ y=Y$ and the negative of the trivial substitution,  $x=-X,\ y=-Y$. If $d\leq -5$ then there are no other automorphs. If $d=-3$ or $-4$ then there is only one equivalence class of forms, represented by the principal form. If $d=-3$, this is the form $x^2+xy+y^2$, with the additional automorphs $x=-Y,\ y=X+Y$ and $x=X+Y,\ y=-X$ and their negatives. If $d=-4$ then the principal form is $x^2+y^2$, and this form has, in addition to the ones already mentioned, the automorph $x=Y,\ y=X$ and its negative. We denote the number of automorphs by $w$, so that
\begin{equation*}
w=\left\{
\begin{array}{ll}2,\ \textrm{if $d<-4$,}\\
4,\ \textrm{if $d=-4$,}\\
6, \ \textrm{if $d=-3$.}\\\end{array}
\right. \tag{2}
\]
On the other hand, if $d>0$ then each form of discriminant $d$ has infinitely many automorphs, which are determined by the integral solutions $(t, u)$ of the associated \emph{Pell's equation}
\[
t^2-du^2=4.\]
This equation has the trivial solution $t=\pm 2, u=0$, and if $(t_0, u_0)$ is the solution for which $t_0>0$ and $u_0$ is positive and is as small as possible, the so-called \emph{minimal positive solution} of Pell's equation, then all nontrivial solutions $(t, u)$ are generated from the minimal positive solution via the equation
\[
\frac{1}{2}\big(t+u\sqrt d\ \big)=\pm\left[\frac{1}{2}\big(t_0+u_0\sqrt d\ \big)\right]^n,\]
where $n$ varies over all positive and negative integers. For a form $[a, b, c]$ of discriminant $d$, it can be shown that all automorphs of the form are given by
\[
x=\frac{1}{2}(t-bu)X-cuY,\
y=auX+\frac{1}{2}(t+bu)Y,\]
with the trivial automorphs determined by the trivial solutions of Pell's equation. For a proof of the results in this paragraph, the interested reader should consult Landau [35], Theorems 111 and 202.

\section{Representation of Integers by Quadratic Forms and the Class Number}
 
Let $n$ be an integer, $q(x, y)$ a quadratic form. We say that $n$ is represented by $q(x, y)$ if there exist integers $x$ and $y$ such that $n=q(x, y)$, and two representations $q(x, y)$ and $q(X, Y)$ of $n$ are \emph{distinct} if the ordered pairs $(x, y)$ and $(X, Y)$ are distinct. In this section we will be interested in the number of distinct ways a positive integer can be represented by quadratic forms in $\mathcal{Q}(d)$. The problem of determining what integers are represented by a given quadratic form, which is naturally closely related to the number of such representations, was one of the main motivations for the theory of quadratic forms which Gauss developed in Section V of the \emph{Disquisitiones}. The number of representations of a positive integer will also be a crucial idea in the proof of the class-number formula. 

As we mentioned in section 12 of Chapter 3, the manner in which a given quadratic form $q=[a, b, c]$ in $\mathcal{Q}(d)$ represents integers depends on the sign of the discriminant $d$. If $d>0$ then $q$ represents both positive and negative integers. If $d<0$ and $a>0$ then $q$ represents no negative integers, and $q(x, y)$ represents 0 only if $x=y=0$. If $d<0$ and $a<0$ then $q$ represents no positive integers, and $q$ represents 0 only if $x=y=0$. Hence forms with positive discriminant are called indefinite and forms with negative discriminant are called positive or negative definite if $a$ is, respectively, positive or negative. A straightforward calculation shows that any unimodular substitution transforms a positive definite form into a positive definite form, and multiplication of a form by $-1$ transforms a positive definite form into a negative definite form. When $d<0$, it hence follows that the equivalence classes of $\mathcal{Q}(d)$ which are determined by modular substitutions are divided evenly into a positive number of classes consisting of positive definite forms and a positive number of classes consisting of negative definite forms. When $d<0$, we define the \emph{class number of d} to be the number of equivalence classes of positive definite forms, and we denote it by $h(d)$. When $d>0$, we define the \emph{class number of d} to be the number of all equivalence classes in $\mathcal{Q}(d)$, and we also denote it by $h(d)$. As we shall see, this is the number to which Dirichlet's class-number formula refers and in the context of understanding why the value of $L$-functions at $s=1$ is positive, it is the most important parameter in that formula.

When $d>0$ then each form $q(x, y)$ is indefinite, and so we can chose a positive integer $k$ such that $k=q(x_0, y_0)$ with $\gcd(x_0, y_0)=1$. It can then be shown that $q(x, y)$ is equivalent to a form with $k$ as the coefficient of its $x^2$-term (Landau [35], Theorem 201), and so we can choose a form $[a, b, c]$ from each equivalence class with $a>0$. When $d<0$, we can chose such a form from each equivalence class of positive definite forms. We formalize this procedure by defining a \emph{representative system of forms} as a set of representatives $[a, b, c]$, one from each equivalence class (positive definite if $d<0$), such that each representative has $a>0$. 

When $d<0$ then, as per the convention that we will now follow, all forms are positive definite, and so the number of representations of a positive integer $n$ is finite. We hence denote by $R(n)$ the total number of (distinct) representations of $n$ by all forms from a representative system. We observe that $R(n)$ does not depend on the representative system use to define it.

On the other hand, when $d>0$, each representation of the positive integer $n$ by a fixed form gives rise to infinitely many representations of $n$ by that form, obtained by applying the automorphs of the form. We are going to circumvent this difficulty by restricting the representations of $n$ to ones which satisfy a particular condition which we will specify next, and which is designed to canonically select a \emph{finite} set of representations of $n$ by each quadratic form.

Let $q=[a, b, c]$ be a quadratic form of discriminant $d>0$. Let $(t_0, u_0)$ be the minimal positive solution of the associated Pell's equation 
\[
t^2-du^2=4,\]
and let
\[
\varepsilon=\frac{t_0+u_0\sqrt d}{2}.\]
If $n$ is a positive integer then the representation $n=q(x, y)$ of $n$ by $q$ is said to be \emph{primary} if 
\[
2ax+(b-\sqrt d)y>0,\]
and
\[
1\leq \frac{2ax+(b+\sqrt d)y)}{2ax+(b-\sqrt d)y}<\varepsilon^2.\]
It can then be shown that the number of primary representations of $n$ by $q$ is finite (Landau [35], Theorem 203 and the remark after this theorem). When $d>0$, we hence denote by $R(n)$ the total number of primary representations of $n$ by all forms from a representative system, and we also note that, as before, $R(n)$ is independent of the representative system of forms used to define it.

One of the fundamental results in the classical theory of quadratic forms gives a formula for the calculation of $R(n)$ by means of Dirichlet characters, and this theorem provides the link between quadratic forms and Dirichlet characters on which the derivation of the class-number formula is based. It goes like so: for a proof, see Landau [35], Theorem 204:
\begin{thm}
Let d be a fundamental discriminate and let $\chi(d)$ be the character determined by d via Theorem $8.2$. If n is a positive integer that is relatively prime to d and $R(n)$ denotes the total number of representations of n $($primary if $d>0)$ by a representative system of quadratic forms in $\mathcal{Q}(d)$, then the value of $R(n)$ is given by the formula \[
R(n)=w\sum_{m|n}\ \chi(d)(m),\]
where w is given by $(1)$ if $d<0$, and $w=1$ if $d>0$.
\end{thm}

\noindent As Davenport summarizes, this theorem is deduced by expressing $R(n)$ in terms of the number of solutions of the congruence $x^2\equiv d$ mod $4n$ and then evaluating that sum using appropriate characters defined by the Legendre symbol. With Theorem 8.3 in hand, we now have all the ingredients required for the statement and proof of the class-number formula.

\section{The Class-Number Formula}

If $\chi$ is a real primitive Dirichlet character of modulus $b$ and if $d=\chi(-1)b$ then it follows from Theorem 8.2 that $d$ is a fundamental discriminant. Dirichlet's class-number formula calculates the value at $s=1$ of the $L$-function of $\chi$ as the product of a certain positive parameter and the class number of $d$. In this sense, we can hence interpret this result to say that $L(1, \chi)$ counts the equivalence classes of quadratic forms in $\mathcal{Q}(d)$, thereby providing the definitive reason why $L(1, \chi)$ is positive. We now state and prove Dirichlet's remarkable formula.
\begin{thm}
$($Dirichlet's Class-Number Formula$)$. Let $\chi$ be a real primitive Dirichlet character of modulus b, let $d=\chi(-1)b$, and let $h(d)$ denote the class number of d. 

$(i)$ Suppose that $d<0$. Then
\[
L(1, \chi)=\frac{2\pi}{w\sqrt{|d|}}h(d),\]
where w is determined by $(2)$.

$(ii)$ Suppose that $d>0$. Let $(t_0, u_0)$ be the minimal positive solution of $t^2-du^2=4$, and let $\varepsilon=\frac{1}{2}(t_0+u_0\sqrt d\ )>1$. Then
\[
L(1, \chi)=\frac{\log \varepsilon}{\sqrt d}h(d).\]
\end{thm}

\emph{Proof}. Let $\chi$, $d$, and $h(d)$ be as given in the statement of Theorem 8.4. Starting with a positive integer $n$ relatively prime to $d$ and the formula 
\begin{equation*}
R(n)=w\sum_{m|n}\ \chi(m) \tag{3}\]
for $R(n)$ as given by Theorem 8.3, the key idea of this argument is to calculate the  asymptotic average value of $R(n)$ as $n\rightarrow +\infty$ in two different ways and compare the results. The first way is to use (3) and the fact that $\chi$ is non-principal to deduce that the  asymptotic average value of $R(n)$ as $n\rightarrow +\infty$ is $wL(1, \chi)$. The second way is to directly use the definition of $R(n)$ and an integer-lattice count in the plane to express the asymptotic average value of $R(n)$ as the product of a positive parameter and the class number of $d$. Setting the two expressions of the asymptotic average equal to each other produces the class-number formula that we desire. 

In order to implement the first way, we begin by using the formula (3) to write the sum
\begin{equation*}
\frac{1}{w}\sum_{\substack{n=1\\ \gcd(n,\  d)=1}}^N R(n)=\sum_{\substack{m_1m_2\leq N\\ \gcd(m_1m_2,\  d)=1}} \chi(m_1). \tag{4}\]
We then split the sum on the right-hand side of this equation into
\[
\sum_{m_1\leq \sqrt N}\ \chi(m_1)\Big(\sum_{\substack{m_2\leq N/m_1\\ \gcd(m_2,\ d)=1}} 1\Big)+\sum_{\substack{m_2<\sqrt N\\ \gcd(m_2,\ d)=1}}\Big(\sum_{\sqrt N<m_1\leq N/m_2}\ \chi(m_1)\Big),\]
because the sum on the left is summed over all pairs $m_1, m_2$ such that $m_1\leq \sqrt N$ and the sum on the right ranges over all pairs for which $m_1>\sqrt N$. The inner sum in the sum on the left is
\[
\frac{N}{m_1}\frac{\varphi(|d|)}{|d|}+ \textit{O}\big(\varphi(|d|)\big),\]
where $\varphi$ is Euler's totient, hence the double sum on the left is estimated by
\[
N\frac{\varphi(|d|)}{|d|}\sum_{m_1\leq \sqrt N}\ \frac{\chi(m_1)}{m_1}+\textit{O}(\sqrt N).\]
We now exploit the fact that $\chi$ is non-principal, to wit, the sum of the values $\chi(m)$ as $m$ varies throughout any finite interval is bounded (as was shown in the proof of Proposition 7.5$(i)$). Hence the double sum on the right is \textit{O}$(\sqrt N)$. Inserting these estimates into (3) and then multiplying the resulting equation by $w/N$, we obtain the estimate
\[
\frac{1}{N}\sum_{\substack{n=1\\ \gcd(n,\  d)=1}}^N R(n)=w\frac{\varphi(|d|)}{|d|}\sum_{m\leq \sqrt N}\ \frac{\chi(m)}{m}+\textit{O}(N^{-\frac{1}{2}}).\]
The sum on the right-hand side of this equation, when extended from $m=1$ to $+\infty$, has a remainder that, when estimated by a summation-by-parts argument, is
\[
\sum_{m>\sqrt N}\ \frac{\chi(m)}{m}=\textit{O}(N^{-\frac{1}{2}}).\]
It follows that
\begin{equation*}
\lim_{N\rightarrow +\infty}\frac{1}{N}\sum_{\substack{n=1\\ \gcd(n,\  d)=1}}^N R(n)=w\frac{\varphi(|d|)}{|d|}\sum_{m=1}^{\infty}\ \frac{\chi(m)}{m}=w\frac{\varphi(|d|)}{|d|}L(1, \chi).\tag{5}\] 
As Davenport points out, the quotient $\varphi(|d|)/|d|$ measures the density of the integers that are relatively prime to $|d|$, and so this equation asserts that the asymptotic average of $R(n)$ with respect to $n$ is $wL(1, \chi)$.

We now want to calculate the limit
\[
\lim_{N\rightarrow +\infty}\frac{1}{N}\sum_{\substack{n=1\\ \gcd(n,\  d)=1}}^N R(n)\]
in the second way, i.e.,  by means of the definition of $R(n)$ as the total number of representations of $n$ by quadratic forms in $\mathcal{Q}(d)$ from a representative system. Hence we let $R(n, f)$ denote the number of representations of $n$ (primary when $d>0$) by a particular form $f\in \mathcal{Q}(d)$, so that, by definition,
\[
R(n)=\sum_{f}\ R(n, f),\]
where the sum here is taken over all quadratic forms $f$ from a representative system. Note that the number of terms in this sum is therefore $h(d)$. We now wish to calculate the limit 
\begin{equation*}
\lim_{N\rightarrow +\infty}\frac{1}{N}\sum_{\substack{n=1\\ \gcd(n,\  d)=1}}^N R(n, f),\tag{6}\]
and it will transpire that this limit has a value $\kappa(d)$ independent of $f$. Hence
\begin{equation*}\
\lim_{N\rightarrow +\infty}\frac{1}{N}\sum_{\substack{n=1\\ \gcd(n,\  d)=1}}^N R(n)=\kappa(d)h(d).\tag{7}\]
It follows from (5) and (7) that
\begin{equation*}
L(1, \chi)=\frac{|d|}{w\varphi(|d|)}\kappa(d)h(d),\tag{8}\]
and so in order to complete the proof of Theorem 8.4, we must calculate $\kappa(d)$.

Suppose first that $d<0$. Let $f=[a, b, c]$. The sum
\[
\sum_{\substack{n=1\\ \gcd(n,\  d)=1}}^N R(n, f)\]
is the number of ordered pairs of integers $(x, y)$ satisfying
\[
0<ax^2+bxy+cy^2\leq N,\ \gcd(ax^2+bxy+cy^2, d)=1.\]
As $x$ and $y$ each run through a complete set of ordinary residues mod $|d|$, there are exactly  $|d|\varphi(|d|)$ of the numbers $ax^2+bxy+cy^2$ that are relatively prime to $d$ (Landau [35], Theorem 206). It therefore suffices to fix a pair of integers $(x_0, y_0)$ for which 
\[
\gcd(ax_0^2+bx_0y_0+cy_0^2, d)=1,\]
and consider the number of pairs of integers $(x, y)$ which satisfy
\[
ax^2+bxy+cy^2\leq N,\ x\equiv x_0,\ y\equiv y_0\ \textrm{mod $|d|$}.\]
The first inequality asserts that the point $(x, y)$ is in an ellipse centered at the origin and passing through the points $(\pm\sqrt{N/a}, 0)$ and $(0, \pm\sqrt{N/c})$, with the ellipse expanding uniformly as $N\rightarrow+\infty$. Figures 1 and 2 exhibit two typical elliptical regions which arise in this manner.

\begin{figure}[h]
  \centering
  \begin{tikzpicture}
    \draw[name path=xaxis-minus, thick, <-] (-5,0) -- (0,0);
    \draw[name path=xaxis-plus, thick, ->] (0,0) -- (5,0);
    \draw[name path=yaxis-minus, thick, <-] (0,-4) -- (0,0);
    \draw[name path=yaxis-plus, thick, ->] (0,0) -- (0,4);

    \draw[pattern=north east lines, name path=ellipse] (0,0) ellipse [rotate=-30, x radius = 3, y radius=1.5];

    \fill[name intersections={of=ellipse and yaxis-plus, by=beta}] 
    (beta) circle (2pt) node[anchor=south west] {$\beta$};
    \fill[name intersections={of=ellipse and yaxis-minus, by=beta-minus}] 
    (beta-minus) circle (2pt) node[anchor=north east] {$-\beta$};

    \fill[name intersections={of=ellipse and xaxis-plus, by=alpha}] 
    (alpha) circle (2pt) node[anchor=south west] {$\alpha$};
    \fill[name intersections={of=ellipse and xaxis-minus, by=alpha-minus}] 
    (alpha-minus) circle (2pt) node[anchor=north east] {$-\alpha$};

    \draw (3,4) node[anchor=north] {$\alpha = \sqrt{\frac{N}{a}}$};
    \draw (3,3) node[anchor=north] {$\beta = \sqrt{\frac{N}{c}}$};
    
  \end{tikzpicture}
  \caption{Elliptical region , $\displaystyle{\frac{b}{a-c}} < 0$}
  \label{fig:8.1}
\end{figure}

\begin{figure}[h]
  \centering
  \begin{tikzpicture}
    \draw[name path=xaxis-minus, thick, <-] (-5,0) -- (0,0);
    \draw[name path=xaxis-plus, thick, ->] (0,0) -- (5,0);
    \draw[name path=yaxis-minus, thick, <-] (0,-4) -- (0,0);
    \draw[name path=yaxis-plus, thick, ->] (0,0) -- (0,4);

    \draw[pattern=north east lines, name path=ellipse] (0,0) ellipse [rotate=30, x radius = 3, y radius=1.5];

    \fill[name intersections={of=ellipse and yaxis-plus, by=beta}] 
    (beta) circle (2pt) node[anchor=south east] {$\beta$};
    \fill[name intersections={of=ellipse and yaxis-minus, by=beta-minus}] 
    (beta-minus) circle (2pt) node[anchor=north west] {$-\beta$};

    \fill[name intersections={of=ellipse and xaxis-plus, by=alpha}] 
    (alpha) circle (2pt) node[anchor=north west] {$\alpha$};
    \fill[name intersections={of=ellipse and xaxis-minus, by=alpha-minus}] 
    (alpha-minus) circle (2pt) node[anchor=south east] {$-\alpha$};

    \draw (3,4) node[anchor=north] {$\alpha = \sqrt{\frac{N}{a}}$};
    \draw (3,3) node[anchor=north] {$\beta = \sqrt{\frac{N}{c}}$};
    
  \end{tikzpicture}
  \caption{Elliptical region , $\displaystyle{\frac{b}{a-c}} > 0$}
  \label{fig:8.2}
\end{figure}

\newpage 
 
Using the fact that the area enclosed by an ellipse is $\pi$ times the product of the lengths of the semi-axes, we find that the area of the ellipse is
\[
\frac{2\pi}{\sqrt{4ac-b^2}}N=\frac{2\pi}{\sqrt{|d|}}N.\]
By dividing the plane into squares of side length $|d|$ centered at the points of the integer lattice, it can be shown without difficulty that as $N\rightarrow +\infty$, the number of points in the integer lattice and lying inside the ellipse is asymptotic to 
\[
\frac{1}{|d|^2}\frac{2\pi}{\sqrt{|d|}}N.\]
This must now be multiplied by $|d|\varphi(|d|)$ in order to account for the number of points $(x_0, y_0)$. It follows that $\kappa(d)$, the value of the limit (5), is
\[
\frac{\varphi(|d|)}{|d|}\frac{2\pi}{\sqrt{|d|}}.\]
Substituting this value of $\kappa(d)$ into (8) yields the conclusion $(i)$ of Theorem 8.4.

Now let $d>0$. The lattice-point count here is a bit more involved than the one that was done for $d<0$ because we need to count primary representations. If we set
\[
\theta=\frac{-b+\sqrt d}{2a},\ \theta^{\prime}=\frac{-b-\sqrt d}{2a},\ \varepsilon=\frac{1}{2}(t_0+u_0\sqrt d\ )>1,\]
where $(t_0, u_0)$ is the positive minimal solution of Pell's equation $t^2-du^2=4$ as before, then an integer pair $(x, y)$ determines a primary representation of an integer if
\[
x-\theta y>0\ \textrm{and}\ 1\leq \frac{x-\theta^{\prime}y}{x-\theta y}<\varepsilon^2.\]
Arguing as we did in the case of negative $d$, for a fixed pair of integers $(x_0, y_0)$ for which
\[
\gcd(ax_0^2+bx_0y_0+cy_0^2, d)=1,\]
we need to count the number of integer points $(x, y)$ such that
\[
ax^2+bxy+cy^2\leq N,\  x-\theta y>0,\  1\leq \frac{x-\theta^{\prime}y}{x-\theta y}<\varepsilon^2,\]
and
\[
x\equiv x_0,\  y\equiv y_0\ \textrm{mod}\ d.\]

The set of conditions
\[
ax^2+bxy+cy^2\leq N,\  x-\theta y>0,\  1\leq \frac{x-\theta^{\prime}y}{x-\theta y}<\varepsilon^2,\]
determines a hyperbolic sector in the upper half-plane of two possible types: one is bounded by the nonnegative $x$-axis, the branch of the hyperbola $ax^2+bxy+cy^2=N$ passing through the point $\sqrt{N/a}$ on the $x$-axis, and the ray emanating from the origin and lying along the line $\nu x=(\nu \theta+1)y,\ \nu=\displaystyle{\frac{a(\varepsilon^2-1)}{\sqrt d}}$, and the other is bounded by the nonpositive $x$-axis, the branch of the hyperbola $ax^2+bxy+cy^2=N$ passing through the point $-\sqrt{N/a}$ on the $x$-axis, and the ray emanating from the origin and lying along the line $\nu x=(\nu \theta+1)y$. Figures 3, 4, and 5 illustrate typical hyperbolic sectors which arise in the first way; we invite the reader to provide figures for the hyperbolic sectors arising in the second way.

\begin{figure}[h]
  \centering
  \begin{tikzpicture}[scale=1]
    \draw[name path=xaxis, thick, ->] (-0.5,0) -- (5,0);
    \draw[name path=yaxis, thick, ->] (0,-0.5) -- (0,5);
    \draw[name path=hyper, very thick , <->] (-1,4) .. controls (0,1) and (1,0) .. (4,-1);
    \draw[name path=arrow, very thick, ->] (0,0) -- (1,4);
    
    \fill[name intersections={of=hyper and xaxis, by=alpha}] 
    (alpha) circle (2pt) node[anchor=north] {$\alpha$};

    \fill[name intersections={of=hyper and yaxis, by=beta}] 
    (beta) circle (2pt) node[anchor=west] {$\beta$};

    \path[name intersections={of=hyper and arrow, by=gamma}];

    \draw (2,3) node[anchor=north] {$\alpha = \sqrt{\frac{N}{a}}$};
    \draw (2,2) node[anchor=north] {$\beta = \sqrt{\frac{N}{c}}$};

    \begin{scope}
    \clip (0,0) -- (gamma) -- (alpha);
    \draw[pattern=north east lines] (-1,4) .. controls (0,1) and (1,0) .. (4,-1) -- (0,0) -- (gamma);
    \end{scope}
  \end{tikzpicture}
  \caption{Hyperbolic sector, $c > 0$, $\nu\theta + 1 > 0$}
  \label{fig:a}
\end{figure}

\begin{figure}[h]
  \centering
  \begin{tikzpicture}[scale=1]
    \draw[name path=xaxis, thick, ->] (-1.5,0) -- (5,0);
    \draw[name path=yaxis, thick, ->] (0,-0.5) -- (0,5);
    \draw[name path=hyper, very thick , <->] (-1,4) .. controls (0,1) and (1,0) .. (4,-1);
    \draw[name path=arrow, very thick, ->] (0,0) -- (-0.9,5);
    
    \fill[name intersections={of=hyper and xaxis, by=alpha}] 
    (alpha) circle (2pt) node[anchor=north] {$\alpha$};

    \fill[name intersections={of=hyper and yaxis, by=beta}] 
    (beta) circle (2pt) node[anchor=west] {$\beta$};

    \path[name intersections={of=hyper and arrow, by=gamma}];

    \draw (2,3) node[anchor=north] {$\alpha = \sqrt{\frac{N}{a}}$};
    \draw (2,2) node[anchor=north] {$\beta = \sqrt{\frac{N}{c}}$};

    \begin{scope}
    \clip (0,0) -- (gamma) -- (alpha);
    \draw[pattern=north east lines] (-1,4) .. controls (0,1) and (1,0) .. (4,-1) -- (0,0) -- (gamma);
    \end{scope}
  \end{tikzpicture}
  \caption{Hyperbolic sector, $c > 0$, $\nu\theta + 1 < 0$}
  \label{fig:a}
\end{figure}

\begin{figure}[h]
  \centering
  \begin{tikzpicture}[scale=1]
    \draw[name path=xaxis, thick, ->] (-1.5,0) -- (5,0);
    \draw[name path=yaxis, thick, ->] (0,-4) -- (0,4);
    \draw[name path=hyper, very thick , <->] (4,4) .. controls (0,1) and (0,-1) .. (4,-4);
    \draw[name path=arrow, very thick, ->] (0,0) -- (4,3.5);
    
    \fill[name intersections={of=hyper and xaxis, by=alpha}] 
    (alpha) circle (2pt) node[anchor=north west] {$\alpha$};

    \path[name intersections={of=hyper and arrow, by=gamma}];

    \draw (4,2) node[anchor=north] {$\alpha = \sqrt{\frac{N}{a}}$};

    \begin{scope}
    \clip (0,0) -- (gamma) -- (alpha);
    \draw[pattern=north east lines] (4,4) .. controls (0,1) and (0,-1) .. (4,-4) -- (0,0) -- (gamma);
    \end{scope}
  \end{tikzpicture}
  \caption{Hyperbolic sector, $c < 0$}
  \label{fig:8.3}
\end{figure}

\newpage

We need to calculate the area of this sector, which we will do by using the change of variables
\[
\xi=x-\theta y,\  \eta=x-\theta^{\prime}y.\]
Because of the fact that
\[
ax^2+bxy+cy^2=a(x-\theta y)(x-\theta^{\prime}y),\]
it follows that the hyperbolic sector is mapped onto the sector in the $\eta, \xi$ plane given by
\[
\eta \xi\leq \frac{N}{a},\  \xi>0,\  \xi\leq \eta<\varepsilon^2\xi,\]
or equivalently,
\[
0<\xi\leq \sqrt{\frac{N}{a}},\  \xi \leq \eta<\min\left(\varepsilon^2\xi,\ \frac{N}{a\xi}\right).\]
The Jacobian of the change of variables is
\[
\frac{\sqrt d}{a}\ ,\]
hence, upon setting $\xi_1=\varepsilon^{-1}\sqrt{N/a}$, the area of the sector in the $x, y$ plane is
\[
\frac{a}{\sqrt d}\left(\int_0^{\xi_1}\ (\varepsilon^2\xi-\xi) d\xi+\int_{\xi_1}^{\varepsilon\xi_1}\ \left(\frac{N}{a\xi}-\xi\right)d\xi \right)=\frac{N\log \varepsilon}{\sqrt d}.\]
Following the argument that we outlined when $d<0$, we conclude that the number of integer points inside the hyperbolic sector is asymptotic to
\[
\frac{1}{d^2}\frac{\log \varepsilon}{\sqrt d}N\]
as $N\rightarrow +\infty$. When this is multiplied by $d\varphi(d)$, in order to account for the number of choices of the pair $(x_0, y_0)$, the value of $\kappa(d)$ is now
\[
\frac{\varphi(d)}{d}\frac{\log \varepsilon}{\sqrt d},\]
hence we obtain the conclusion  $(ii)$ of Theorem 8.4 from (8) using this value of $\kappa(d)$. QED

\section{The Class-Number Formula and the Class Number of Quadratic Fields}

As we alluded to in the introduction to this chapter, the class-number formula can also be used to count ideal classes in quadratic number fields in exactly the same way that it counts equivalence classes of quadratic forms, and this provides another good explanation of the positivity of $L$-functions at $s=1$. In this section, we will discuss in more detail exactly how this goes.

Recall from section 11 of Chapter 3 that if $K$ is an algebraic number field and $R$ is the ring of algebraic integers in $K$, then the ideals $I$ and $J$ of $R$ are said to be equivalent if there exists nonzero elements $\alpha$ and $\beta$ of $R$ such that $\alpha I=\beta J$.  The number of equivalence classes of ideals with respect to this equivalence relation is finite, and that number is called the class number of $R$.

As we saw in sections 11 and 12 of Chapter 3, there is a very close connection between the theory of ideals in quadratic number fields and the classical theory of quadratic forms. For the reader's convenience, we will recapitulate the principal features of that connection here.

Let $d$ be a fundamental discriminant. Then $d$ is either odd and square-free, or $d$ is divisible by 4 and $d/4$ is square-free. If we let $m=d$ in the former case and $m=d/4$ in the latter case, then the quadratic number field $F$ generated over $\mathbb{Q}$ by $\sqrt m$ has discriminant $d$. In order to relate the ideal theory of $F$ to quadratic forms of discriminant $d$, as we did in section 11 of  Chapter 3, we will need to recall the definition of the the norm of an element of $F$. One starts by taking an integral basis $\{\omega_1, \omega_2\}$ of $F$, letting $k\in F$, and then expressing $k$ uniquely as $k=r\omega_1+s\omega_2$, for some ordered pair $(r, s)\in \mathbb{Q}\times \mathbb{Q}$. We then define the norm of $k$ as the number defined by $N(k)=(r\omega_1+s\omega_2)(r\omega_1^{\prime}+s\omega_2^{\prime})$,
where $\{\omega_1^{\prime}, \omega_2^{\prime}\}$ denotes the algebraic conjugates of $\{\omega_1, \omega_2\}$ over $\mathbb{Q}$. This definition of $N(k)$ does not depend on the integral basis used to define it, hence selecting either the standard integral basis $\{1, \frac{1}{2}(1+\sqrt m)\}$ or $\{1, \sqrt m\}$, depending on whether $m$ is or, respectively, is not congruent to 1 mod 4 (section 11, Chapter 3), a simple calculation shows that if $d<0$ then $N(k)>0$ whenever $k\not= 0$. 

If $I$ is an ideal in $R=\mathcal{R}\cap F$ then the norm mapping $N$ maps $I$ into $\mathbb{Z}$, and if $\{\alpha, \beta\}$ is an integral basis of $I$ and $\xi=x\alpha+y\beta,\ (x, y)\in \mathbb{Z}\times \mathbb{Z}$, then we also have that
\[
N(\xi)=(x\alpha+y\beta)(x\alpha^{\prime}+y\beta^{\prime}).\]
The right-hand side of this equation defines a quadratic form $q(x, y)$ with integer coefficients. When we recall that the norm $N(I)$ of $I$ is defined as the cardinality of the finite quotient ring $R/I$ (section 1 of Chapter 5), it can be shown that all coefficients of $q(x, y)$ are divisible by $N(I)$, and that if we let
\[
\frac{N(\xi)}{N(I)}=ax^2+bxy+cy^2,\]
then this form is in $\mathcal{Q}(d)$. The equivalence class of $\mathcal{Q}(d)$ which contains this form can be shown to be independent of the integral basis of $I$ which is used to define the form. 

As we discussed in section 12 of Chapter 3, the map which sends an ideal $I$ of $R$ to the quadratic form  $N(\xi)/N(I)=ax^2+bxy+cy^2$ induces a bijective correspondence between the equivalence classes of ideals of $R$ in the narrow sense and the forms in a representative system of discriminant $d$ (Theorem 3.22). This correspondence will now be used to express the class number of $d$ in terms of the class number of $R$. 

Recall from section 12 of Chapter 3 that if $I$ and $J$ are ideals of $R$ then $I$ is said to be equivalent to $J$ in the narrow sense if there exists $k\in F$ such that $N(k)>0$ and $I=kJ$. If $d<0$ then there is no difference between ordinary equivalence of ideals and equivalence in the narrow sense because the values of $N$ are all positive on $F\setminus \{0\}$. When $d>0$ then the difference between these two equivalence relations is mediated by the units in $R$. If $R$ has a unit of norm $-1$ then there is also no difference between the two equivalence relations. On the other hand, if $d>0$ and there is no unit with a negative norm then it can be shown that each ideal class in the ordinary sense in the union of exactly two ideal classes in the narrow sense (section 12, Chapter 3). 

It follows that if if we denote the class number of $R$ by $h_1(d)$ and recall that $h(d)$ denotes the class number of $d$ as determined from the classes of forms in $\mathcal{Q}(d)$, then $h(d)=h_1(d)$ whenever either $d<0$ or $d>0$ and $R$ has a unit of norm $-1$, and $h(d)=2h_1(d)$, otherwise. 

The parameters $w$ and $\varepsilon=\frac{1}{2}(t_0+u_0\sqrt d\ )$ which occur in the class-number formula can be calculated by means of the units in $R$. When $d<0$, the parameter $w$ is defined by the  values as stipulated for it by (2), and those values are also the number of roots of unity in $R$ . As for $\varepsilon$, which enters the class-number formula when $d>0$, its value can be calculated by the fundamental unit in the group of units $U(R)$ of $R$. Recall from section 12 of Chapter 3 that when $d>0$, there is a unit $\varpi$ of $R$ such that
\[
U(R)=\{\pm\varpi^n: n\in \mathbb{Z}\}.\]
If $\varpi$ is chosen to exceed 1 then it is uniquely determined as a generator of $U(R)$ in this sense and is called the {fundamental unit of $R$. If $\varpi$ is the fundamental unit, then it can be shown that $\varepsilon=\varpi$ if $N(\varpi)=1$ and $\varepsilon=\varpi^2$ if $N(\varpi)=-1$.

Combining all of these results leads to the conclusion that
\[
h(d)=h_1(d),\ \textrm{if $d<0$},\]
and
\[
h(d)\log \varepsilon=2h_1(d)\log \varpi,\ \textrm{if $d>0$}.\]
When these equations are combined with Dirichlet's class-number formula it follows that if $\chi$ is a real primitive Dirichlet character of modulus $b$ and $d=\chi(-1)b$, then
\[
L(1, \chi)=\frac{2\pi}{w\sqrt{|d|}} h_1(d),\ \textrm{if $d<0$},\]
and
\[
L(1, \chi)=\frac{2\log \varpi}{\sqrt{d}} h_1(d),\ \textrm{if $d>0$}.\]

Using the first of these formulae for $L(1, \chi)$, we can now calculate the quadratic excesses in Theorems 7.2, 7.3, and 7.4 in terms of class numbers of quadratic fields. If $p\equiv 3$ mod 8 then
\[
q(0, p/2)=3 h_1(-p),\]
if $p\equiv 7$ mod 8 then
\[
q(0, p/2)=h_1(-p),\]
if $p\equiv 1$ mod 4 then
\[
q(0, p/4)=\frac{1}{2}h_1(-p),\]
if $p>3$ and $p\equiv 1$ mod 4 then
\[
q(0, p/3)=\frac{1}{2}h_1(-3p),\]
if $p>3$ and $p\equiv 7$ mod 12 then
\[
q(0, p/3)=2 h_1(-p),\]
and if $p>3$ and $p\equiv 11$ mod 12 then
\[
q(0, p/3)=h_1(-p).\]
A noteworthy consequence of these equations is that $h_1(-p)$ is even if $p\equiv 1$ mod 4 and $h_1(-3p)$ is even if $p>3$ and $p\equiv 1$ mod 4. As a matter of fact, it is an interesting open problem to determine if similar sums of this type which are positive are always linear combinations of class numbers; for further information on sums and congruences for $L$-functions, consult  Urbanowicz and Williams [56] (I thank an anonymous referee for calling my attention to this problem and this reference).

\vspace{0.4cm}
\section{A Character-Theoretic Proof of Quadratic Reciprocity}

\vspace{0.4cm}
In this section, we will present our final proof of quadratic reciprocity. The main idea of this approach is to define Gauss sums for a general Dirichlet character and then use computations with these sums and the structure theory for Dirichlet characters that was discussed in section 1 to derive a formula which calculates the value of a real primitive character $\chi$ at odd primes $p$ as an appropriate value of the Legendre symbol $\chi_p$. If we then set $\chi$ equal to $\chi_q$ for an odd prime $q$ distinct from $p$ into that formula, the LQR will immediately result from that choice of $\chi$. This argument reveals how quadratic reciprocity is caused by the fact that the Legendre symbols are real primitive characters which can reproduce the value at the primes of an arbitrary real primitive character. We follow the exposition as recorded in Cohen [2], sections 2.1 and 2.2.

Let $\chi$ be a Dirichlet character of modulus $b$, let $a\in \mathbb{Z}$, set $\zeta_b=\exp(2\pi i/b)$, and define the \emph{Gauss sum of $\chi$ at a} by
\[
G(\chi, a)=\sum_{x\ \textrm{mod}\ b}\ \chi(x) \zeta_b^{ax},\]
where the sum is taken over all $x$ in a complete set of ordinary residues mod $b$. It is clear, because of the periodicity of $\chi(x)$ and $ \zeta_b^{ax}$ in $x$, that this sum does not depend on the set of ordinary residues mod $b$ used to define it. We first encountered Gauss sums in Chapter 3; in particular, the proof of Lemma 3.15 can be easily modified to establish
\begin{lem}
If a and b are relatively prime then
\[
G(\chi, a)=\overline{\chi(a)}G(\chi, 1),\]
where the bar over $\chi(a)$ denotes the complex conjugate.\end{lem}

The next lemma records a useful condition on $a$ which implies that $G(\chi, a)=0$.
\begin{lem}
Let $d=\gcd(a, b)$ and suppose that b/d is not an induced modulus of $\chi$. Then $G(\chi, a)=0$.
\end{lem}

\emph{Proof}. Because $b/d$ is not an induced modulus, there is an integer $m$ such that $m\equiv 1$ mod $b/d$, $\gcd(m, b)=1$, and $\chi(m)\not= 1$. Upon letting $m^{-1}$ denote the inverse of $m$ mod $b$, we have that
\[
\chi(m)G(\chi, a)=\sum_{x\ \textrm{mod}\ b}\ \chi(mx) \zeta_b^{ax}=\sum_{y\ \textrm{mod}\ b}\ \chi(y) \zeta_b^{ay m^{-1}}.\]
Because $m\equiv 1$ mod $b/d$, it follows that
\[
ay m^{-1}=\frac{a}{d}dy m^{-1}\equiv ay\ \textrm{mod}\ b,\]
hence
\[
\chi(m)G(\chi, a)=\sum_{y\ \textrm{mod}\ b}\ \chi(y) \zeta_b^{ay}=G(\chi, a),\]
so that $G(\chi, a)=0$, as $\chi(m)\not=1$. $\hspace{9cm}\ \textrm{QED}$

The next result shows that when $\chi$ is primitive, Lemma 8.5 holds for all integers $a$,  and not just the integers that are relatively prime to $b$.

\begin{lem}
If $\chi$ is primitive then
\[
G(\chi, a)=\overline{\chi(a)}G(\chi, 1),\ \textrm {for all }\ a\in \mathbb{Z}.\]
\end{lem}

\emph{Proof}. We need only check this equation if $d=\gcd(a, b)>1$. But then $\chi(a)=0$, and $b/d$ cannot be induced modulus of $\chi$, hence $G(\chi, a)=0$ by Lemma 8.6.$\hspace{2.6cm}\ \textrm{QED}$

If $\chi$ is a primitive character then we can use Lemma 8.7 to calculate $|G(\chi, 1)|$ as follows:
\begin{eqnarray*}
|G(\chi, 1)|&=&\left(\overline{G(\chi, 1)}G(\chi, 1)\right)^{1/2}\\
&=&\left(\sum_{a\ \textrm{mod}\ b}\ \overline{\chi(a)}G(\chi, 1)\zeta_b^{-a}\right)^{1/2} \\
&=&\left(\sum_{a\ \textrm{mod}\ b}\ G(\chi, a)\zeta_b^{-a}\right)^{1/2},\ \textrm{by Lemma 8.7}\\
&=&\left(\sum_{x=1}^b\ \chi(x)\left(\sum_{a=1}^b\ \zeta_b^{a(x-1)}\right)\right)^{1/2}\\
&=&\sqrt b,
\end{eqnarray*}
where the last line follows because the inner sum here is a geometric series, which is 0 if $b$ does not divide $x-1$, i.e., if $x\not= 1$, and is $b$ if $x=1$. This calculation now has the following corollary, which is a direct generalization of Theorem 3.14.
\begin{cor}
If $\chi$ is a primitive character of modulus b then
\[
G(\chi, 1)G(\overline{\chi}, 1)=\chi(-1)b.\]
In particular, if $\chi$ is also real then
\[
G(\chi, 1)^2=\chi(-1)b.\]
\end{cor}

\emph{Proof}. From Lemma 8.5 (applied to $\overline{\chi}$), we have that
\[
\overline{G(\chi, 1)}=G(\overline{\chi}, -1)=\chi(-1)G(\overline{\chi}, 1).\]
Now multiply this equation by $\chi(-1)G(\chi, 1)$ and then calculate that
\[
\chi(-1)b=\chi(-1)G(\chi, 1)\overline{G(\chi, 1)}=\chi(-1)^2G(\chi, 1)G(\overline{\chi}, 1)=G(\chi, 1)G(\overline{\chi}, 1).\]
$\hspace{15.5cm}\ \textrm{QED}$

The next result shows that the value of a real primitive character at the odd prime $p$ can be calculated by the Legendre symbol $\chi_p$ of $p$; it is a direct generalization of equation (17) in Chapter 3. Note also that its proof follows the same ideas used in the third proof of the LQR that we gave in Chapter 3.
\begin{thm}
If $\chi$ is a real primitive character of modulus b then for all odd primes p,
\[
\chi(p)=\chi_p\big(\chi(-1)b\big).\] 
\end{thm}

\emph{Proof}. Since both sides of the equation to be verified are 0 when $p$ divides $b$, we may assume that $p$ is an odd prime not dividing $b$. Let $\mathcal{R}$ denote the ring of algebraic integers. Since the quotient ring $\mathcal{R}/p\mathcal{R}$ has characteristic $p$, the map which sends an element of that ring to its $p$-th power is additive, and we also have that $\zeta_b\in \mathcal{R}$. Hence
\[
G(\chi, 1)^p\equiv \sum_{x\ \textrm{mod}\ b}\ \chi(x)^p \zeta_b^{px}\ \textrm{mod}\ p\mathcal{R}.\]
As $\chi$ is real and $p$ is odd, $\chi(x)^p=\chi(x)$, and so it follows from Lemma 8.5 that 
\[
G(\chi, 1)^p\equiv \chi(p)G(\chi, 1)\ \textrm{mod}\ p\mathcal{R}.\]
On the other hand, $\chi$ is real and primitive, hence by Corollary 8.8, $G(\chi, 1)^2=\chi(-1)b$, and so multiplication of this congruence by $G(\chi, 1)$ yields the congruence 
\[
\chi(-1)b\big((\chi(-1)b)^{\frac{1}{2}(p-1)}-\chi(p)\big)\equiv 0 \ \textrm{mod}\ p\mathcal{R}.\]
 If we now multiply this congruence by the inverse of $\chi(-1)b$ in $\mathbb{Z}/p\mathbb{Z}$ to obtain 
 \[
 (\chi(-1)b)^{\frac{1}{2}(p-1)}\equiv \chi(p)\ \textrm{mod}\ p\mathcal{R}\]
and then use the congruence
\[
 (\chi(-1)b)^{\frac{1}{2}(p-1)}\equiv \chi_p(\chi(-1)b)\ \textrm{mod}\ p\mathbb{Z}\]
 that we get from Euler's criterion, it follows that
 \[
 \chi(p)\equiv \chi_p(\chi(-1)b)\ \textrm{mod}\ p\mathcal{R}.\]
Hence
\[
\frac{\chi(p)-\chi_p(\chi(-1)b)}{p}\in \mathcal{R}\cap \mathbb{Q}=\mathbb{Z}.\]
Because the numerator of this rational integer is either 0 or $\pm2$ and $p$ is odd, it follows that the numerator must be 0, whence the conclusion of the theorem. $\hspace{4.7cm}\ \textrm{QED}$

The LQR can now be deduced immediately from Theorem 8.9 and Theorem 2.4 of Chapter 2. Let $p$ and $q$ be distinct odd primes. We take $\chi =\chi_q$ in Theorem 8.9 and then apply Theorem 2.4 twice to calculate that
\begin{eqnarray*}
\chi_q(p)&=&\chi_p(\chi_q(-1)q)\\
&=&\chi_p\big((-1)^{\frac{1}{2}(q-1)}q\big)\\
&=&\chi_p(-1)^{\frac{1}{2}(q-1)}\chi_p(q)\\
&=&(-1)^{\frac{1}{2}(p-1)\frac{1}{2}(q-1)}\chi_p(q).
\end{eqnarray*}

\newpage
\afterpage{\null\newpage}
\thispagestyle{empty}

\chapter{Quadratic Residues and Non-Residues in Arithmetic Progression}

The distribution problem for residues and non-residues has been intensively studied for 175 years using a rich variety of formulations and techniques. The work done in Chapter 7 gave a window through which we viewed one of these formulations and also saw a very important technique used to study it. Another problem that has been studied almost as long and just as intensely is concerned with the arithmetic structure of residues and non-residues. In this chapter, we will sample one aspect of that very important problem by studying when residues and non-residues form very long sequences in arithmetic progression. The first major advance in that problem came in 1939 when Harold Davenport proved the existence of residues and non-residues which form arbitrarily long sets of consecutive integers. As an introduction to the circle of ideas on which the work of this chapter is based, we briefly discuss Davenport's results and the technique that he used to obtain them in section 1. Davenport's approach uses another application of the Dirichlet-Hilbert trick, which we used in the proofs of Theorems 4.12 and 5.13 presented in Chapter 5, together with an ingenious estimate of the absolute value of certain Legendre-symbol sums with polynomial values in their arguments. Davenport's technique is quite flexible, and so we will adapt it in order to detect long sets of residues and non-residues in arithmetic progression. In section 2, we will formulate our results precisely as a series of four problems which will eventually be solved in sections 4 and 10. This will require the estimation of the sums of values of Legendre symbols with polynomial arguments a la Davenport, which estimates we will derive in section 3 by making use of a very important result of Andr$\acute{\textrm{e}}$ Weil concerning the number of rational points on a nonsingular algebraic curve over $\mathbb{Z}/p\mathbb{Z}$. In addition to these estimates, we will also need to calculate a term which will be shown to determine the asymptotic behavior of the number of sets of residues or non-residues which form long sequences of arithmetic progressions, and this calculation will be performed in sections 6-9. Here we will see how techniques from combinatorial number theory are applied to study residues and non-residues. In section 11, an interesting class of examples will be presented, and we will use it to illustrate exactly how the results obtained in section 10, together with some results of section 11, combine to describe asymptotically how many sets there are of residues or non-residues which form long arithmetic progressions. Finally, the last section of this chapter discusses a result which, in certain interesting situations, calculates the asymptotic density of the set of primes which have residues and non-residues which form long sets of specified arithmetic progressions.

\section{Long Sets of Consecutive Residues and Non-Residues}

The following question began to attract interest in the early 1900's: if $s$ is a fixed positive integer and $p$ is sufficiently large, does there exist an $n \in[1, \infty)$ such that $\{n, n+1,\dots,n+s-1\}$ is a set of residues (respectively, non-residues) of $p$ inside $[1, p-1]$, i.e., for all sufficiently large primes $p$, does $[1, p-1]$ contain arbitrarily long sets of consecutive residues, (respectively, non-residues) of $p$? For $s=2, 3 ,4,$ and 5, various authors showed that the answer is yes; in fact it was shown that if $R_s(p)$ (respectively, $N_s(p)$) denotes the number of sets of $s$ consecutive residues (respectively, non-residues) of $p$ inside $[1, p-1]$ then as $p \rightarrow +\infty$,
\begin{equation*}
R_s(p) \sim 2^{-s}p \sim N_s(p),\ \textrm{for}\  s=2, 3, 4,\ \textrm{and}\ 5. \tag{1}\]
This shows in particular that for $s=2, 3, 4,$ and 5, not only are $R_s(p)$ and $N_s(p)$ both positive, but as $p \rightarrow +\infty$, they both tend to $+\infty$. Based on this  evidence and extensive numerical calculations, the speculation was that (1) in fact is valid without any restriction on $s$, and in 1939, Harold Davenport [5] proved that this is indeed the case. 

Davenport established the validity of (1) in general by yet another application of the Dirichlet-Hilbert trick that was used in the proof of Theorems 4.12 and 5.13. Let $\mathbb{F}_p$ denote the field $\mathbb{Z}/p\mathbb{Z}$ of $p$ elements. Then $U(p)$ can be viewed as the group of nonzero elements of $\mathbb{F}_p$, and if $\varepsilon \in \{-1, 1\}$ then, a la Dirichlet-Hilbert, the sum
\[
2^{-s} \sum_{x=1}^{p-s} \prod_{i=0}^{s-1} \big(1+\varepsilon \chi_p(x+i)\big) \]
is $R_s(p)$ (respectively, $N_s(p)$) when $\varepsilon=1$(respectively, $\varepsilon=-1$). Davenport rewrote this sum as
\begin{equation*}
2^{-s}(p-s)+2^{-s}\sum_{\emptyset \not= T \subseteq [0, s-1]} \varepsilon^{|T|} \left( \sum_{x=1}^{p-s} \chi_p\left(\prod_{i \in T} (x+i)\right) \right), \tag{2} \]
and then proceeded to estimate the size of the second term of the expression (2). This term is a sum of terms of the form
\[
\pm\sum_{x=1}^{p-s} \chi_p\big(f(x)\big),\]
where $f$ is a monic polynomial of degree at most $s$ over $\mathbb{F}_p$ with distinct roots in $\mathbb{F}_p$. Using results from the theory of certain $L$-functions due to Hasse, Davenport found absolute constants $C>0$ and $0<\sigma<1$ such that
\[
\left|\sum_{x=1}^{p-s} \chi_p\big(f(x)\big) \right| \leq Csp^{-\sigma},\ \textrm{for all $p$ large enough}. \]
This estimate, the heart of Davenport's argument, implies that the modulus of the second term in (2) does not exceed $Csp^{\sigma}$, and so
\[
|R_s(p)-2^{-s}(p-s)| \leq Csp^{\sigma},\ \textrm{for all $p$ large enough}.\]
Hence
\begin{eqnarray*} 
\left|\frac{R_s(p)}{2^{-s}p}-1 \right| &\leq& \frac{s}{p}+Cs2^sp^{\sigma-1}\\
&\rightarrow& 0\ \textrm{as}\ p \rightarrow +\infty.
\end{eqnarray*}
The same argument also works for $N_s(p)$

It transpires that Davenport's technique is quite flexible and can be used to investigate the occurrence of residues and non-residues with specific arithmetical properties. We are going to use it to detect arbitrarily long arithmetic progressions of residues and non-residues of a prime.

\section{Long Sets of Residues and Non-Residues in Arithmetic Progression}
 Our point of departure from Davenport's work is to notice that the sequence $\{x,x+1,\dots, x+s-1\}$ of $s$ consecutive positive integers is an instance of the sequence$\{x,x+b,\dots,x+b(s-1)\}$, an arithmetic progression of length $s$ and common difference $b$, with $b=1$. Thus, if $(b,s)\in [1,\infty)\times[1,\infty)$, and we set
 \[
 AP(b;s)=\Big\{ \{n+ib: i\in [0,s-1]\}: n\in [1,\infty)\Big\},
 \]
the family of all arithmetic progressions of length $s$ and common difference $b$, it is natural to inquire about the asymptotics as $p\rightarrow +\infty$ of the number of elements of $AP(b;s)$ that are sets of quadratic residues (respectively, non-residues) of $p$ that occur inside $[1,p-1]$. We also consider the following related question: if $a\in [0,\infty)$, set
\[
AP(a,b;s)=\Big\{ \{a+b(n+i):i\in [0,s-1]\}: n\in [1,\infty)\Big\},
\]
the family of all arithmetic progressions of length $s$ taken from a fixed arithmetic progression 
\[
AP(a, b)= \{a+bn: n\in [1,\infty)\}.
\]
We then ask for the asymptotics of the number of elements of $AP(a,b;s)$ that  are sets of quadratic residues (respectively, non-residues) of $p$ that occur inside $[1,p-1]$. 
Solutions of these problems will provide interesting insights into how often quadratic residues and non-residues appear as arbitrarily long arithmetic progressions.

 We will in fact consider the following generalization of these questions. For each $m\in [1,\infty)$, let
 \[
 \textbf{a}=(a_1,\dots,a_m)\  \textrm{and}\ \textbf{b}=(b_1,\dots,b_m)
 \] 
 be $m$-tuples of nonnegative integers such that $(a_i,b_i)\not=(a_j,b_j)$, for all $i\not=j$. Let $s\in [1, \infty)$. When the $b_i$'$s$ are distinct and positive, we set
 \[
 AP(\textbf{b}; s)=\Big\{\bigcup_{j=1}^{m} \{n+ib_j:i\in [0,s-1]\}:n\in [1, \infty)\Big\},
 \]
 and when the $b_i$'s are all positive (but not necessarily distinct), we set
 \[
 AP(\textbf{a},\textbf{b};s)=\Big\{\bigcup_{j=1}^{m} \{a_j+b_j(n+i):i\in [0,s-1]\}:n\in [1,\infty)\Big\}.
 \] 
 
 The elements of $AP(\textbf{b}; s)$ are formed by taking an $n\in [1, \infty)$, then selecting an arithmetic progression of length $s$ with initial term $n$ and common difference $b_i$ for each $i=1,\dots,m$, and then taking the union of the arithmetic progressions so chosen. Elements of $AP(\textbf{a},\textbf{b};s)$ are obtained by taking an $n\in [1, \infty)$, choosing from the arithmetic progression$ \{a_i+b_im: m\in [1,\infty)\}$ the arithmetic progression with initial term $a_i+b_in$ and length $s$ for each $i=1,\dots,m$, and then forming the union of the arithmetic progressions chosen in that way. If $m=1$ then we recover our original sets $AP(b;s)$ and $AP(a,b;s)$. We now pose
 \begin{quote}
 \textbf{Problem 1} (respectively, \textbf{Problem 2}): determine the asymptotics as $p\rightarrow+\infty$ of the number of elements of $AP(\textbf{b};s)$ (respectively, $AP(\textbf{a}, \textbf{b};s)$) that are sets of quadratic residues of $p$ inside $[1,p-1]$.
 \end{quote}
 We also pose as \textbf{Problem 3} and \textbf{Problem 4} the problems which result when the phrase ``quadratic residues" in the statements of Problems 1 and 2 is replaced by the phrase ``quadratic non-residues".

As we saw in section 1, the main step of Davenport's solution to the problem of finding long sets of consecutive residues and non-residues was finding good estimates of sums of the form
\[
\sum_{x=0}^{p-1}\ \chi_p\big(f(x)\big)\]
for certain polynomials $f(x)\in \mathbb{F}_p[x]$.
In order to solve Problems 1-4, we will use a variant of Davenport's reasoning, tailored to detect long sets of residues and non-residues in arithmetic progression. Our techniques will also require appropriate estimation of these sums. Estimates very much like Davenport's will suffice to solve Problems 1 and 3, but, for technical reasons, they are not sufficient to solve Problems 2 and 4. For those problems, we will need to obtain good estimates, which do not depend on $N$, for sums of the form
\[
\sum_{x=0}^N\ \chi_p\big(f(x)\big),\]
where $N$ can be any integer in $[0, p-1]$. As preamble to our solution of Problems 1-4, we show in the following section how some results of A. Weil can be used to efficiently and elegantly derive the estimates that we require. 
\section{Weil Sums and their Estimation}

In order to solve Problems 1-4, we will require estimates of sums of the form
\begin{equation*}
\sum_{x=1}^N \chi_p\big(f(x)\big), \tag{*} \]
where $f$ is a polynomial in $\mathbb{F}_p[x]$ and $N$ is a fixed integer in $[1, p-1]$.

Suppose first that $N=p-1$. In this case there is an elegant way to calculate the sum $(*)$ in terms of the number of rational points on an algebraic curve over $\mathbb{F}_p$.

If $F$ is a field, $\overline{F}$ is an algebraic closure of $F$, and $g(x, y)$ is a polynomial in two variables with coefficients in $F$, then the set of points
\[
C=\{(x, y) \in \overline{F} \times \overline{F}: g(x, y)=0 \} \]
is an \emph{algebraic curve over F}. A point $(x, y) \in C$ is a \emph{rational point of C over F} if $(x, y) \in F \times F$. If $F$ is finite then the set of rational points on an algebraic curve over $F$ is evidently finite, and so the determination of the cardinality of the set of rational points is an interesting and very important problem in combinatorial number theory.  In 1948, A. Weil's great treatise [57] on the geometry of algebraic curves over finite fields was published, which contained, among many other results of fundamental importance, an upper estimate of the number of rational points in terms of $\sqrt{|F|}$ and certain geometric parameters associated with an algebraic  curve. The Weil bound has turned  out to be very important for various problems in number theory; in particular, we will now show how it can be employed to obtain good estimates of the sums $(*)$ when $N=p-1$.

Let $f \in \mathbb{F}_p[x]$ and consider the algebraic curve $C$ over $\mathbb{F}_p$ defined by the polynomial
\[
y^2-f(x).\]
We will calculate the so-called \emph{complete Weil sum}\[
\sum_{x=0}^{p-1} \chi_p\big(f(x)\big)\]
in terms of the number of rational points of $C$ over $\mathbb{F}_p$.

Let $\mathcal{R}(p)$ denote the set of rational points of $C$, i.e.,
\[
\mathcal{R}(p)=\{(x, y) \in \mathbb{F}_p \times \mathbb{F}_p: y^2=f(x)\},\]
and let
\[
S_0=\{x \in \mathbb{F}_p:  f(x)=0\},\]
\[
S_{+}=\{x \in \mathbb{F}_p\setminus S_0: \chi_p(f(x))=1\},\]
\[
S_{-}=\{x \in \mathbb{F}_p \setminus S_0: \chi_p(f(x))=-1\}.\]
If $x \in S_{+}$ then there are exactly two solutions $\pm y_0 \not= 0$ of $y^2=f(x)$ in $\mathbb{F}_p$, hence $(x, \pm y_0) \in \mathcal{R}(p)$. Conversely, if $(x, y) \in \mathcal{R}(p)$ and $y \not= 0$ then $0 \not= y^2=f(x)$, hence $x \in S_{+}$ and $y=\pm y_0$. We conclude that 
\begin{equation*}
\left|\mathcal{R}(p) \right|= \left|S_0 \right|+2 \left|S_{+} \right|. \tag{2} \]
Because $\mathbb{F}_p$ is the pairwise disjoint union of $S_0, S_{+}$, and $S_{-}$,
\begin{equation*}
\left|S_0 \right|+\left| S_{+} \right|+ \left| S_{-} \right|=p. \tag{3} \] 
Observe now that
\begin{equation*}
\sum_{x=0}^{p-1} \chi_p\big(f(x)\big)=\left|S_{+} \right|-\left|S_{-} \right|. \tag{4}\]
Equations (2), (3), and (4) imply
\begin{eqnarray*}
\left|\mathcal{R}(p) \right|&=& \left|S_0 \right|+ \left|S_{+} \right|+\left|S_{-} \right|+\sum_{x=0}^{p-1} \chi_p\big(f(x)\big)\\
&=&p+\sum_{x=0}^{p-1} \chi_p\big(f(x)\big),
\end{eqnarray*}
i.e.,
\begin{equation*}\
\sum_{x=0}^{p-1} \chi_p\big(f(x)\big)=\left|\mathcal{R}(p) \right|-p. \tag{5}\]

We are ready to apply Weil's estimate of $\left|\mathcal{R}(p) \right|$. In this case, Weil ([57], Corollaire IV.3)  proved that if $y^2-f(x)$ is \emph{non-singular over $\mathbb{F}_p$}, which means essentially that $f$ is monic of degree at least 1 and there does not exist a polynomial $g \in \mathbb{F}_p[x]$ such that $f=g^2$, then
\begin{equation*}
\left|\mathcal{R}(p) \right|=1+p-r(p),\ \textrm{where}\ 1 \leq r(p)<d \sqrt p,\ d=\ \textrm{degree of $f$} \tag{6}\] 
(for an elementary proof of (6), see Schmidt [50], Theorem 2.2C). If $f \in \mathbb{F}_p[x]$ is monic with distinct roots in $\mathbb{F}_p$ then $f$ cannot be the  square of a polynomial over $\mathbb{F}_p$, and so $y^2-f(x)$ is non-singular over $\mathbb{F}_p$. Hence (5) and (6) imply
\begin{thm}
$($complete Weil-sum estimate$)$ If $f \in \mathbb{F}_p[x]$ is monic of degree $d \geq1$ and $f$ has distinct roots in $\mathbb{F}_p$ then
\[
\Big|\sum_{x=0}^{p-1} \chi_p\big(f(x)\big) \Big| < d \sqrt p.\]
\end{thm}

As its definition makes clear, a Weil sum is nothing more than the sum of a certain sequence of 0's and $\pm1$'s. The content of Theorem 9.1 (and also Theorem 9.2 to follow) is that for certain polynomials $f\in \mathbb{F}_p[x]$, a remarkable cancellation occurs in the terms of $\sum_{x=0}^{p-1} \chi_p\big(f(x)\big)$ so that the absolute value of this sum, ostensibly as large as $p$, is in fact less that $d\sqrt p$, where $d$ is the degree of $f$. 

The work of Weil in [57] is another seminal development in modern number theory. There Weil used methods from algebraic geometry to study number-theoretic properties of curves, thereby founding the subject of \emph{arithmetic algebraic geometry}.  This not only introduced important new techniques in both number theory and geometry, but it also led to the formulation of innovative strategies for attacking a wide variety of problems which until then had been intractable. Certainly one of the most spectacular examples of that is the proof of Fermat's Last Theorem by Andrew Wiles [60] in 1995 (with an able assist from Richard Taylor [55]), which employed arithmetic algebraic geometry as one of its crucial tools. 

We now turn to the problem of estimating the sums $(*)$ when $N< p-1$. An \emph{incomplete Weil sum} is a sum of the form
\begin{equation*}
\sum_{x=M}^N \chi_p\big(f(x)\big),\ \tag{**}\]
where $f \in \mathbb{F}_p[x]$, and either $0 \leq M \leq N<p-1$ or $0<M \leq N \leq p-1$. Our solution of Problems 2 and 4 will require an estimate of incomplete Weil sums similar to the estimate of complete Weil sums provided by Theorem 9.1, and also independent of the parameters $M$ and $N$. When $f(x)=x$, Polya proved in 1918 that 
\[
\Big|\sum_{x=M}^N \chi_p(x) \Big| \leq \sqrt p \log p,\]
and Vinogradov in the same year showed that if $\chi$ is a non-principal Dirichlet character mod $m$ then
\[
\Big|\sum_{x=M}^N \chi(x) \Big| \leq 6\sqrt m \log m.\]
Assuming the Generalized Riemann Hypothesis, in 1977 Montgomery and Vaughn improved this to
\[
\Big|\sum_{x=M}^N \chi(x) \Big| \leq C\sqrt m \log \log m.\]
By an earlier result of Paley (which holds without assuming GRH), this estimate, except for the choice of the constant $C$, is best possible. It follows that  an estimate of $(**)$ that is independent of $M$ and $N$ will most likely behave more or less like (an absolute constant)$\times \sqrt p \log p$. In fact, we will prove
\begin{thm}
$($incomplete Weil-sum estimate$)$ There exists $p_0>0$ such that the following statement is true: if $p \geq p_0$, if $f \in \mathbb{F}_p[x]$ is monic of degree $d \geq 1$ with distinct roots in $\mathbb{F}_p$, and $N \in [0, p-1]$, then
\[
\Big|\sum_{x=0}^N \chi_p\big(f(x)\big) \Big| \leq d(1+ \log p) \sqrt p.\]
\end{thm}

Our proof of Theorem 9.2 will make use of certain homomorphisms of the additive group of $\mathbb{F}_p$ into the circle group, defined like so. Let 
\[
e_p(\theta)=\exp \left(\frac{2\pi i\theta}{p} \right).\]
If $n \in \mathbb{Z}$ then we set
\[
\psi(m)=e_p(mn),\ m \in \mathbb{Z}.\]
Because $\psi(m)=\psi(m^{\prime})$ whenever $m \equiv m^{\prime}\ \textrm{mod}\ p, \psi$ defines a homomorphism of the additive group of $\mathbb{F}_p$ into the circle group, hence $\psi$ is called an \emph{additive character} mod $p$.

Now for each $n \in  \mathbb{Z}, \zeta=e_p(n)$ is a $p$-th root of unity, i.e., $\zeta^p=1$, and from the factorization
\[
(1-\zeta)\Big(\sum_{k=0}^{p-1} \zeta^k \Big)=1-\zeta^p=0
\]
we see that
\[
\sum_{k=0}^{p-1} \zeta^k =0,\]
unless $\zeta=1$. Applying this with $\zeta=e_p(n-a)$, we obtain
\begin{equation*}
\frac{1}{p}\sum_{x=0}^{p-1} e_p(-ax) e_p(nx)=\left\{\begin{array}{rl}1 ,& \textrm{if} \  n \equiv a\ \textrm{mod}\ p,\\
0 ,& \textrm{otherwise}, 
\\\end{array}\right. \tag{7} \] 
the so-called \emph{orthogonality relations} of the additive characters. These relations are quite similar to the orthogonality relations satisfied by Dirichlet characters (section 4 of Chapter 4), the latter of which Dirichlet used to prove Lemma 4.7, on his way to the proof of Theorem 4.5.

\emph{Proof of Theorem} 9.2.

 Let $f \in \mathbb{F}_p[x]$ be monic of degree $d \geq 1$, with distinct roots in $\mathbb{F}_p$, let $N \in [1, p-1]$ and set
\[
S(N)=\sum_{x=1}^N \chi_p\big(f(x)\big).\]
The strategy of this argument is to use the orthogonality relations of the additive characters to express $S(N)$ as a sum of terms $\lambda(x)S(x), x=0,1,\dots,p-1$, where $\lambda(x)$ is a sum of additive characters and $S(x)$ is a sum that is a ``twisted" or ``hybrid" version of a complete Weil sum. Appropriate estimates of these terms are then made to obtain the conclusion of Theorem 9.2.

We first decompose $S(N)$ like so: 
\begin{eqnarray*}
S(N)&=&\sum_{k=1}^N \sum_{j=0}^{p-1} \delta_{jk} \chi_p\big(f(j)\big),\ \delta_{jk}=\left\{\begin{array}{rl}1\ ,& \textrm{if} \  j=k,\\
0\ ,& \textrm{if}\ j \not= k. 
\\\end{array}\right. \\
&=&\sum_{k=1}^N \sum_{j=0}^{p-1} \chi_p\big(f(j)\big)\left(\frac{1}{p}\sum_{x=0}^{p-1} e_p(xk) e_p(-xj) \right),\ \textrm{by}\ (7)\\
&=&\frac{1}{p} \sum_{x=0}^{p-1} \left(\sum_{k=1}^N e_p(xk) \right) \sum_{j=0}^{p-1} \chi_p\big(f(j)\big)e_p(-xj)\\
&=&\frac{1}{p} \sum_{x=0}^{p-1} \lambda(x) S(x),
\end{eqnarray*}
where
\[
\lambda(x)=\sum_{k=1}^N e_p(xk),\  S(x)=\sum_{k=0}^{p-1} e_p(-xk)\chi_p\big(f(k)\big).\]

The next step is to estimate $|\lambda(x)|$ and $|S(x)|, x=0,1,\dots,p-1$. To get a useful estimate of $|\lambda(x)|$, use the trigonometric identities
\[
\sum_{k=1}^N \cos k\theta=\frac{\sin \left( \left(N+\frac{1}{2} \right) \theta \right)- \sin \left(\frac{\theta}{2}\right)}{2\sin\left(\frac{\theta}{2}\right)},\]
\[
\sum_{k=1}^N \sin k\theta=\frac{\cos \left(\frac{\theta}{2}\right)-\cos \left( \left(N+\frac{1}{2} \right) \theta \right) }{2\sin\left(\frac{\theta}{2}\right)},\]
to calculate that
\[
|\lambda(x)|=\left|\frac{\sin\left(N\pi x/p\right)}{\sin\left(\pi x/p \right)} \right|.\]
Now use the estimate
\[
\frac{2|\theta|}{\pi} \leq |\sin \theta|\ ,\ |\theta| \leq \frac{\pi}{2}, \]
to get
\begin{equation*}
|\lambda(x)| \leq \frac{p}{2|x|}\ ,\ 0<|x|<\frac{p}{2}. \tag{8}\]
The sums $\lambda(x)$ and $S(x)$ are periodic in $x$ of period $p$, hence
\begin{equation*}
S(N)=\frac{1}{p} \sum_{|x|< p/2} \lambda(x) S(x). \tag{9}\]
Note that $\lambda(0)=N$, hence (8), (9) imply that
\[
\left|S(N)-\frac{N}{p}
S(0) \right| \leq \frac{1}{2}\sum_{0<|x|<p/2} |x|^{-1} |S(x)|.\]

An estimate of each sum $S(x)$ is now required. These are so-called hybrid or mixed Weil sums, and consist of terms $e_p(-xy)\chi_p\big(f(y)\big), y=0,1,\dots,p-1$, which are the terms of the complete Weil sum $\sum_{y=0}^{p-1} \chi_p\big(f(y)\big)$ that are ``twisted" by the multiplier $e_p(-xy)$. As Perel'muter [44] proved in 1963  by means of the arithmetic algebraic geometry of Weil (see also Schmidt [50], Theorem 2.2G for an elementary proof), this twisting causes no problems, i.e., we have the estimate
\[
|S(x)| \leq d \sqrt p,\ \textrm{for all}\ x \in \mathbb{Z}.\]
Hence
\begin{eqnarray*}
|S(N)| &\leq& \frac{N}{p} |S(0)|+ \frac{1}{2}\sum_{0<|x|<p/2} |x|^{-1} |S(x)|\\
&\leq&d\sqrt p \Big(1+ \sum_{1 \leq n <p/2} \frac{1}{n} \Big).
\end{eqnarray*}
Because 
\[
\lim_{p\rightarrow +\infty} \Big(\gamma+ \log \left[\frac{p}{2} \right]-\sum_{1 \leq n<p/2} \frac{1}{n} \Big)=0,\]
where $\gamma=$ 0.57721$\dots$ is Euler's constant, we are done. $\hspace{5.5cm} \textrm{QED}$

\section{Solution of Problems 1 and 3}

Now that we have Theorem 9.1 at our disposal, Problems 1 and 3 can be solved, i.e., the asymptotic behavior, as the prime $p\rightarrow +\infty$, of the number of elements of 
\[
AP(\textbf{b}; s)=\Big\{\bigcup_{j=1}^{k} \{n+ib_j:i\in [0,s-1]\}:n\in [1, \infty)\Big\}\]
that are sets of residues (respectively, non-residues) of $p$ inside $[1, p-1]$ can be determined. We begin with some terminology and notation that will allow us to state our results precisely and concisely. Let $W=\{z_1,\dots,z_r\}$ be a nonempty, finite subset of $[0,\infty)$ with its elements indexed in increasing order $z_i<z_j$ for $ i<j$. 
We let
\[
\mathcal{S}(W)=\big\{\{n+z_i:i\in [1,r]\}:n\in [1,\infty)\big\},
\]
the set of all shifts of $W$ to the right by a positive integer. Let $\varepsilon$ be a choice of signs for $[1,r]$, i.e., a function from $[1,r]$ into $\{-1,1\}$. If $S=\{n+z_i:i\in [1,r]\}$ is an element of $\mathcal{S}(W)$, we will say that the pair $(S,\varepsilon)$ is a \emph{residue pattern of p} if
\[
\chi_p(n+z_i)=\varepsilon(i),\  \textrm{for all}\  i\in [1,r].
\]
The set $\mathcal{S}(W)$ has the \emph{universal pattern property} if there exists $p_0>0$ such that for all $p\geq p_0$ and for all choices of signs $\varepsilon$ for $[1,r]$, there is a set $S\in \mathcal{S}(W)\cap 2^{[1,p-1]}$ such that $(S,\varepsilon)$ is a residue pattern of $p$. $\mathcal{S}(W)$ hence has the universal pattern property if and only if for all $p$ sufficiently large, $\mathcal{S}(W)$ contains a set that exhibits any fixed but arbitrary pattern of quadratic residues and non-residues of $p$. This property is inspired directly by Davenport's work: using this terminology, we can state the result of [5, Corollary of Theorem 5] for quadratic residues as asserting that if $s\in [1,\infty)$ then $\mathcal{S}([0,s-1])$ has the universal pattern property, and moreover, for any choice of signs $\varepsilon$ for $[1,s]$, the cardinality of the set
\[
\{S\in \mathcal{S}([0,s-1])\cap 2^{[1,p-1]}: (S,\varepsilon)\  \textrm{is a residue pattern of}\ p\}
\]    
is asymptotic to $2^{-s}p$ as $p\rightarrow +\infty$. Note that if $\varepsilon$ is the choice of signs that is either identically 1 or identically $-1$ on $[1,s]$, then we recover the results that were discussed in section 1 of this chapter.

Suppose now that there exists nontrivial gaps between elements of $W$, i.e., $z_{i+1}-z_i\geq 2$ for at least one $i\in [1,r-1]$. It is then natural to search for elements $S$ of $\mathcal{S}(W)$ such that the quadratic residues (respectively, non-residues) of $p$ inside $[\min{S},\max{S}]$ consists precisely of the elements of $S$, so that $S$ acts as the ``support" of quadratic residues or non-residues of $p$ inside the minimal interval of consecutive integers containing $S$. We formalize this idea by declaring $S$ to be a \emph{residue} (respectively,  \emph{non-residue}) \emph{support set of p} if $S=(\textrm{the set of all residues of $p$ inside $[1, p-1]$})\cap [\min S,\max S]$ (respectively, $S=(\textrm{the set of all non-residues of $p$ inside $[1, p-1]$})\cap [\min S,\max S]$). We then define $\mathcal{S}(W)$ to have the \emph{residue} (respectively, \emph{non-residue}) \emph{support property} if there exist $p_0>0$ such that for all $p\geq p_0$, there is a set $S\in \mathcal{S}(W)\cap 2^{[1,p-1]}$ such that $S$ is a residue (respectively, non-residue) support set of $p$. 

We now use Davenport's method to establish the following proposition, which generalizes [5, Corollary of Theorem 5] for quadratic residues.

\begin{prp}
\label{prp1}
If $W$ is any nonempty, finite subset of $[0,\infty)$, then $\mathcal{S}(W)$ has the universal pattern property and both the residue and non-residue support properties. Moreover, if $\varepsilon$ is a choice of signs for $[1,|W|]$,
\[
c_\varepsilon(W)(p)=\left| \{S\in \mathcal{S}(W)\cap 2^{[1,p-1]}: (S,\varepsilon)\ \textrm{is a residue pattern of}\  p\}\right|, and
\]
\[
c_\sigma(W)(p)=\left| \{S\in \mathcal{S}(W)\cap 2^{[1,p-1]}: S\ \textrm{is a residue (respectively, non-residue) support set of}\  p\}\right|,
\]
then as $p\rightarrow +\infty$,
\[
 c_\varepsilon(W)(p)\sim 2^{-|W|}p\ and\ c_\sigma(W)(p)\sim 2^{-(1+\max W-\min W)}p.
 \]
 \end{prp} 
 
\textit{Proof}. Suppose that the asserted asymptotics of $c_\varepsilon(W)(p)$ has been established for all nonempty, finite subsets $W$ of $[0,\infty)$. Then the asserted asymptotics for $c_\sigma(W)(p)$ can be deduced from that by means of the following trick. Let $W\subseteq [0,\infty)$ be nonempty and finite. Define the choice of signs $\varepsilon$ for [min $W$, max $W$] to be 1 on $W$ and $-1$ on [min $W$, max $W$]$\setminus$$W$. Now for each $p$, let
\[
\mathcal{S}(p)=\{S\in \mathcal{S}(W)\cap 2^{[1,p-1]}: S\ \textrm{is a residue support set of}\  p\},
\]
\[ 
\mathcal{R}(p)=\{S\in \mathcal{S}([\min W,\max W])\cap 2^{[1,p-1]}: (S,\varepsilon)\ \textrm{is a residue pattern of}\  p\}.
\]    
If to each $E\in \mathcal{R}(p)$ (respectively, $F\in \mathcal{S}(p)$), we assign the set $f(E)=E\cap$(set of all residues of $p$ inside $[1, p-1]$) (respectively, $g(F)=[\min F,\max F]$), then $f$ (respectively, $g$) maps $\mathcal{R}(p)$ (respectively, $\mathcal{S}(p)$) injectively into $\mathcal{S}(p)$ (respectively, $\mathcal{R}(p)$). Hence $\mathcal{R}(p)$ and $\mathcal{S}(p)$ have the same cardinality.
Because of our assumption concerning the asymptotics of $c_\varepsilon([\min W,\max W])(p)$, it follows that as $p\rightarrow +\infty$,
\[
c_\sigma(W)(p)=\left| \mathcal{S}(p)\right|=\left| \mathcal{R}(p)\right|\sim 2^{-\left| [\min W,\max W]\right|}\ p=2^{-(1+\max W-\min W)}\ p.
\]
This establishes the conclusion of the proposition with regard to residue support sets, and the conclusion with regard to non-residue support sets follows by repeating the same reasoning after $\varepsilon$ is replaced by $-\varepsilon$.

If $\varepsilon$ is now an arbitrary choice of signs for$[1,|W|]$, it hence suffices to deduce the asserted asymptotics of $c_\varepsilon(W)(p)$. Letting $r(p)=p-\max W-1$, we have for all $p$ sufficiently large that
\[
c_\varepsilon(W)(p)=2^{-|W|}\sum_{x=1}^{r(p)}\ \prod_{i=1}^{|W|} \Big(1+\varepsilon(i)\chi_p(x+z_i)\Big).
\]
This sum can hence be rewritten as
\[
2^{-|W|}r(p)+2^{-|W|}\sum_{\emptyset\ \neq\ T\ \subseteq\ [1,|W|]}\  \prod_{i\in T} \varepsilon(i)\Big(\sum_{x=1}^{r(p)} \chi_p\Big(\prod_{i\in T} (x+z_i)\Big)\Big).
\] 
The asserted asymptotics for $c_\varepsilon(W)(p)$ now follows from an application of Theorem 9.1 to the Weil sums in the second term of this expression. $\hspace{6.9cm} \textrm{QED}$                                 

Now, let $(k, s)\in [1, \infty)\times [1, \infty),\{b_1,\dots, b_k\}\subseteq [1, \infty)$ and $\textbf{b}=(b_1,\dots,b_k)$. We will apply Proposition 9.3 to the family of sets defined by 
 \[
 AP(\textbf{b}; s)=\Big\{\bigcup_{j=1}^{k} \{n+ib_j:i\in [0,s-1]\}:n\in [1, \infty)\Big\};
 \]
we need only to observe that
\[ 
 AP(\textbf{b}; s)=\mathcal{S}\Big(\bigcup_{j=1}^{k}\ \{ib_j: i\in [0, s-1]\}\Big),
\]
 for then the following theorem is an immediate consequence of Proposition 9.3. In particular, if the choice of signs $\varepsilon$ in the theorem is taken to be either identically 1 or identically $-1$, we obtain the solution of Problems 1 and 3.
\begin{thm}
\label{thm3}
$($Wright $[62]$, Theorem $2.3)$ If $(k, s)\in [1, \infty)\times [1, \infty),\{b_1,\dots, b_k\}\subseteq [1, \infty)$ and $\textbf{b}=(b_1,\dots,b_k)$, then $ AP(\textbf{b}; s)$ has the universal pattern property and both the residue and non-residue support properties. Moreover, if $b=\max \{b_1,\dots,b_k\}$,
\[
\gamma=\Big|\bigcup_{j=1}^{k}\ \{ib_j: i\in [0, s-1]\}\Big|,
\]
$\varepsilon$ is a choice of signs for $[1, \gamma]$,
\[
c_\varepsilon(p)=|\{S\in AP(\textbf{b}; s)\cap 2^{[1, p-1]}: (S,\varepsilon)\textit{ is a residue pattern of p}\}|, and
\]
\[
c_\sigma(p)=|\{S\in AP(\textbf{b}; s)\cap 2^{[1, p-1]}: S\ \textit{is a residue (respectively, non-residue) support set of p}\}|,
\]
then as $p\rightarrow+\infty$,
\[
c_\varepsilon(p)\sim 2^{-\gamma}p\  and\ c_\sigma(p)\sim 2^{-(1+b(s-1))}p.
\]

\end{thm}

As an example of Theorem 9.4 in action, take  $k=5,\ s=6$, and $\textbf{b}=(b_1, b_2,  b_3,  b_4, b_5)=(1, 2, 3, 5, 7)$. Then
\[
1+b(s-1)=1+7\cdot5=36\]
and
\begin{eqnarray*}
\gamma&=&\Big|\bigcup_{j=1}^{5}\ \{ib_j: i\in [0, 5]\}\Big|\\
&=&\big|\{0\} \cup \{1,\ 2,\ 3,\ 5, \ 7\}\cup \{2,\ 4,\ 6,\ 10,\ 14\} \cup \{3,\ 6,\ 9,\ 15,\ 21\}\\
&\cup& \{4,\ 8,\ 12,\ 20,\ 28\}\cup \{5,\ 10,\ 15,\ 25,\ 35\}\big|\\
&=&\big|\{0,1, 2, 3, 4, 5, 6,7,8, 9, 10, 12, 14, 15, 20, 21, 25, 28, 35\}\big|\\
&=&19.
\end{eqnarray*}
We have that
\[
AP(\textbf{b}; 6)=\big\{\{n+z: z\in \{0,1, 2, 3, 4, 5, 6,7,8, 9, 10, 12, 14, 15, 20, 21, 25, 28, 35\}: n\in [1, \infty)\big\},
\]
and so if $\varepsilon$ is a choice of signs for $[1, 19]$ then Theorem 9.4 implies that as $p\rightarrow+\infty$,
\[
c_{\varepsilon}(p)\sim 2^{-19}p\ \textrm{and}\ c_{\sigma}(p)\sim 2^{-36}p.\]

\section{Solution of Problems 2 and 4: Introduction}

 Let  $(m, s)\in [1, \infty)\times [1, \infty)$, let $\textbf{a}=(a_1,\dots, a_m)$, (respectively, $\textbf{b}=(b_1,\dots, b_m)$) be an $m$-tuple of nonnegative (respectively, positive) integers such that $(a_i, b_i)\not= (a_j, b_j)$ for $i\not= j$, let  $(\textbf{a}, \textbf{b})$ denote the $2m$-tuple $(a_1,\dots, a_m, b_1,\dots, b_m)$ (we will call $(\textbf{a}, \textbf{b})$ a \emph{standard $2m$-tuple}) , and recall from section 2 that
 \[
 AP(\textbf{a},\textbf{b};s)=\Big\{\bigcup_{j=1}^{m} \{a_j+b_j(n+i):i\in [0,s-1]\}:n\in [1,\infty)\Big\}.
 \] 
Problems 2 and 4 ask for the asymptotic behavior as $p\rightarrow +\infty$ of the number of elements of $AP(\textbf{a},\textbf{b};s)\cap 2^{[1, p-1]}$ which are sets of residues (respectively, non-residues) of $p$. Because of certain arithmetical interactions which can take place between the elements of the sets in $AP(\textbf{a},\textbf{b};s)$, the asymptotic behavior of this sequence is somewhat more complicated than what occurs for $AP(\textbf{b}; s)$ as per Theorem 9.4.

In order to explain the situation, we set
\begin{equation*}
q_\varepsilon(p)=|\{A\in AP(\textbf{a}, \textbf{b}; s)\cap 2^{[1, p-1]}: \chi_p(a)=\varepsilon, \ \textrm{for all}\ a\in A\}|
\end{equation*} 
and note that the value of $q_\varepsilon(p)$ for $\varepsilon=1$ (respectively, $\varepsilon=-1$) counts the number of elements of $AP(\textbf{a}, \textbf{b}; s)$ that are sets of residues (respectively, non-residues) of $p$ which are located inside $[1, p-1]$. As we mentioned before, it  will transpire that the asymptotic behavior of $q_\varepsilon(p)$ depends on certain arithmetic interactions that can take place between the elements of $AP(\textbf{a}, \textbf{b}; s)$. In order to see how this goes, first consider the set $B$ of \emph{distinct} values of the coordinates of $\textbf{b}$. If we declare the coordinate $a_i$ of $\textbf{a}$ and the coordinate $b_i$ of $\textbf{b}$ to \emph{correspond} to each other, then for each $b\in B$, we let $A(b)$ denote the set of all coordinates of $\textbf{a}$ whose corresponding coordinate of $\textbf{b}$ is $b$. We then relabel the elements of $B$ as $b_1,\dots,b_k$, say, and for each $i\in [1, k]$, set  
\[
 S_i=\bigcup_{a\in A(b_i)} \{a+b_ij: j\in [0, s-1]\}, 
 \]
 and then let
 \[
 \alpha=\sum_i\ |S_i|,\   b=\max \{b_1,\dots, b_k\}.
 \]
 
 Next, suppose that
 \begin{quote}
 $(***)$  if $(i, j)\in [1, k]\times [1, k]$ with $i\not= j$ and $(x, y)\in A(b_i)\times A(b_j)$, then either $b_ib_j$ does not divide $yb_i-xb_j$ or  $b_ib_j$ divides $yb_i-xb_j$ with a quotient that exceeds $s-1$ in modulus.
 \end{quote}
Then we will show in section 10 that as $p\rightarrow +\infty,\  q_\varepsilon(p)$ is asymptotic to $(b\cdot 2^\alpha)^{-1}p$. On the other hand, if the assumption $(***)$ does not hold then we will also show in section 10 that the asymptotic behavior of $q_\varepsilon(p)$ falls into two distinct regimes, with each regime determined in a certain manner by the integral quotients
\begin{equation*}
\frac{yb_i-xb_j}{b_ib_j},\  (x, y)\in  A(b_i)\times A(b_j),\tag{$\natural$}
\]
whose moduli do not exceed $s-1$. More precisely, these quotients determine a positive integer $e<\alpha$ and a collection $\mathcal{S}$ of nonempty subsets of $[1, k]$ such that each element of  $\mathcal{S}$ has even cardinality and for which the following two alternatives hold:

$(i)$ if $\prod_{i\in S} b_i$ is a square for all $S\in \mathcal{S}$, then as $p\rightarrow +\infty,\  q_\varepsilon(p)$ is asymptotic to $(b\cdot 2^{\alpha-e})^{-1}p$, or

$(ii)$ if there is an $S\in \mathcal{S}$ such that $\prod_{i\in S} b_i$ is not a square, then there exist two disjoint, infinite sets of primes $\Pi_+$ and $\Pi_-$  whose union contains all but finitely many of the primes and such that $ q_\varepsilon(p)=0$ for all $p\in \Pi_-$, while as $p\rightarrow +\infty$ inside $\Pi_+$, $q_\varepsilon(p)$ is asymptotic to $(b\cdot 2^{\alpha-e})^{-1}p$. Thus we see that when $(***)$ does not hold and $p\rightarrow +\infty$, either $ q_\varepsilon(p)$ is asymptotic to $(b\cdot 2^{\alpha-e})^{-1}p$ or $q_\varepsilon(p)$ asymptotically oscillates infinitely often between 0 and $(b\cdot 2^{\alpha-e})^{-1}p$.

In light of what we have just discussed, it will come as no surprise that the solution of Problems 2 and 4 for  $AP(\textbf{a},\textbf{b};s)$ involves  a bit more effort than the solution of Problems 1 and 3 for  $AP(\textbf{b}; s)$. In order to analyze the asymptotic behavior of $q_\varepsilon(p)$, we follow the same strategy as before: using an appropriate sum of products involving $\chi_p,\  q_\varepsilon(p)$ is expressed as a sum of a dominant term and a remainder. If the dominant term is a non-constant linear function of $p$ and the remainder term does not exceed an absolute constant $\times \sqrt p\log p$, then the asymptotic behavior of $ q_\varepsilon(p)$ will be in hand.    
 
We in fact will implement this strategy when the set  $AP(\textbf{a}, \textbf{b}; s)$ in the definition of $q_\varepsilon(p)$ is replaced by a slightly more general set; for a precise statement of what we establish, see Theorem 9.9 in section 10. Also in section 10, we then deduce  the solution of Problems 2 and 4 from this more general result, where, in particular, we indicate more precisely the manner in which the integral quotients ($\natural$) whose moduli do not exceed $s-1$ determine the parameter $e$ and collection of sets $\mathcal{S}$ discussed above.

\section{Preliminary Estimate of $q_{\varepsilon}(p)$}
We begin the analysis of $q_{\varepsilon}(p)$ by taking a closer look at the structure of  $AP(\textbf{a},\textbf{b};s)$.  Let $\mathcal{J}$ denote the set of all subsets $J$ of $[1, m]$ that are of maximal cardinality with respect to the property that $b_j$  is equal to a fixed integer $b_J$ for all $j\in J$. We note that $\{J: J\in \mathcal{J}\}$ is a partition of $[1, m]$ and that $b_J\not= b_{J'}$ whenever $\{J, J'\}\subseteq \mathcal{J}$. Because $(a_i, b_i)\not= (a_j, b_j)$ whenever $i\not=j$, it follows that if $J\in \mathcal{J}$ then the integers $a_j$ for $j\in J$ are all distinct. Let
 \[
 S_J=\bigcup_{j\in J}\ \{a_j+b_Ji: i \in [0, s-1]\},\ J\in \mathcal{J}.
 \]
Then
 \begin{equation*}
\bigcup_{j=1}^{m}\ \{a_j+b_j(n+i): i\in [0, s-1]\}=\bigcup_{J\in \mathcal{J}}\ b_Jn+S_J\ ,\ \textrm{for all}\ n\in [1, \infty).\tag{11}
  \end{equation*} 

It follows that $AP(\textbf{a},\textbf{b};s)$ is a special case of the following more general situation. Let $k\in [1, \infty)$, let $B=\{b_1,\dots,b_k\}$ be a set of positive integers, and let $\textbf{S}=(S_1,\dots,S_k)$ be a $k$-tuple of finite, nonempty subsets of $[0, \infty)$. By way of analogy with the expression of the elements of $AP(\textbf{a}, \textbf{b}; s)$ according to (11), we will denote by $AP(B, \textbf{S})$ the collection of sets defined by
\[
\Big\{\bigcup_{i=1}^k\ b_in+S_i: n\in [1, \infty)\Big\}.
\] 
We are interested in the number of elements of $AP(B, \textbf{S})$ that are sets of quadratic residues or, respectively, quadratic non-residues of a prime $p$, and so if $\varepsilon\in \{-1, 1\}$, we replace $AP(\textbf{a},\textbf{b};s)$ by  $AP(B, \textbf{S})$ in the definition of $q_{\varepsilon}(p)$ and retain the same notation, so that now
\[
 q_\varepsilon(p)=|\{A\in AP(B, \textbf{S})\cap 2^{[1, p-1]}: \chi_p(a)=\varepsilon, \ \textrm{for all}\ a\in A\}|.
 \]
The problem is to find an asymptotic formula for $q_\varepsilon(p)$ as $p\rightarrow +\infty$.

We can now begin to implement the strategy for determining the asymptotic behavior of $q_{\varepsilon}(p)$ as set forth in section 5. Our goal here is to find an initial estimate of $q_{\varepsilon}(p)$ in terms of an expression that we will eventually denote by $\Sigma_4(p)$ such that $q_{\varepsilon}(p)-\Sigma_4(p)=O(\sqrt p\log p)$ as $p\rightarrow +\infty$. In sections 7, 8, and 9, we will then prove that  $\Sigma_4(p)$ is a non-constant linear function of $p$ for enough primes $p$ so that the precise asymptotic behavior of $q_{\varepsilon}(p)$ is captured.

Toward that end, begin by noticing that there is a positive constant $C$, depending only on $B$ and $ \textbf{S}$, such that for all $n\geq C$,
 \begin{equation*}
\textrm{the sets $b_in+S_i, i\in [1, k]$, are pairwise disjoint, and}\tag{12}
 \end{equation*} 
\begin{equation*}
\textrm{$\bigcup_{i=1}^k\ b_in+S_i$ is uniquely determined by $n$}.\tag{13} 
 \end{equation*} 
Because of (12) and (13), if
\[
\alpha=\sum_i|S_i|\ \textrm{and}\  r(p)=\min_i\left[\frac{p-1-\textrm{max}\ S_i}{b_i}\right],
\]
then the sum
\[
2^{-\alpha}\sum_{x=1}^{r(p)}\ \prod_{i=1}^k\ \prod_{j\in S_i}\ \big(1+\varepsilon\chi_p(b_ix+j)\big)
\]
differs from $q_\varepsilon(p)$ by at most $O(1)$, hence, as per the strategy as outlined in section 5, this sum can be used to determine the asymptotics of $q_\varepsilon(p)$. 

Apropos of that strategy, let
\[
\mathcal{T}=\bigcup_{i=1}^k\ \{(i, j): j\in S_i\},
\] 
and then rewrite the above sum as
 \begin{equation*}
2^{-\alpha}r(p)+2^{-\alpha}\sum_{\emptyset\not= T\subseteq \mathcal{T}}\ \varepsilon^{|T|}\prod_{i=1}^k\ \chi_p(b_i)^{|\{j: (i, j)\in T\}|}\sum_{x=1}^{r(p)}\ \chi_p\Big(\prod_{(i,j)\in T}\ (x+\bar{b_i}j)\Big),\tag{14}
 \end{equation*}
 where $\bar{b_i}$ denotes the inverse of $b_i$ modulo $p$, which clearly exists for all $p$ sufficiently large. Our intent now is to estimate the modulus of certain summands in the second term of (14) by means of Theorem 9.2. 

Let $\Sigma(p)$ denote the second term of the sum in (14). In order to carry out the intended estimate, we must first remove from $\Sigma(p)$ the terms to which Theorem 9.2 cannot be applied. Toward that end, let
\begin{quote}
$E(p)=\{\emptyset\not= T\subseteq \mathcal{T}$: the distinct elements, modulo $p$, in the list $\bar{b_i}j, (i, j)\in T$, each occurs an \emph{even} number of times$\}$.
\end{quote} 
We then split $\Sigma(p)$ into the sum $\Sigma_1(p)$  of terms taken over the elements of $E(p)$ and the sum $\Sigma_2(p)= \Sigma(p)- \Sigma_1(p)$. The sum $\Sigma_2(p)$ has no more than  
$2^{\alpha}-1$ terms each of the form
\[
\pm 2^{-\alpha}\sum_{x=1}^{r(p)} \chi_p\Big(\prod_{(i,j)\in T} (x+\bar{b_i}j)\Big),\ \emptyset \not= T\in 2^{\mathcal{T}}\setminus E(p).
\] 
Since $\emptyset \not= T\notin E(p)$, the polynomial in $x$  in this term at which $\chi_p$ is evaluated can be reduced to a product of at least one and no more than $\alpha$ distinct monic linear factors in $x$ over $\mathbb{F}_p$, and so the sum in each of the above terms of $\Sigma_2(p)$ is an incomplete Weil sum to which Theorem 9.2 can be applied. It therefore follows from that theorem that 
\[
\Sigma_2(p)=O(\sqrt p\log p)\ \textrm{as}\  p\rightarrow +\infty.
\]
We must now estimate
\[
\Sigma_3(p)=2^{-\alpha}r(p)+\Sigma_1(p),
\]
and, as we shall see, it is precisely this term that will produce the dominant term which determines the asymptotic behavior of $q_\varepsilon(p)$.

Since each element of $E(p)$ has even cardinality,
\[
\Sigma_1(p)=2^{-\alpha}\sum_{T \in E(p)}\ \prod_{i=1}^{k}\chi_p(b_i)^{|\{j:(i,j)\in T\}|}\ \sum_{x=1}^{r(p)}\ \chi_p\Big(\prod_{(i,j)\in T}(x+\bar{b_i}j)\Big).
\] 
We now examine the sum over $x\in [1, r(p)]$ on the right-hand side of this equation. Because $T \in E(p)$, each term in this sum is either 0 or 1, and a term is 0 precisely when the value of $x$ in that term agrees with the minimal nonnegative ordinary residue mod $p$ of $- \bar{b_i}j$, for some element $(i, j)$ of $T$. However, there are at most $\alpha/2$ of these values at which $x$ can agree for each $T\in E(p)$ and so it follows that $\Sigma_3(p)$ differs by at most $O(1)$ from 
\[
\Sigma_4(p)=2^{-\alpha}r(p)\Big(1+\sum_{T\in E(p)}\ \prod_{i=1}^{k}\ \chi_p(b_i)^{|\{j:(i,j)\in T\}|}\Big).
\]
Consequently,
\begin{equation*}
\textrm{for all $p$ sufficiently large},\  q_\varepsilon(p)-\Sigma_4(p)=O(\sqrt p\log p),\tag{15}
\end{equation*}
and so it suffices to calculate $\Sigma_4(p)$ in order to determine the asymptotics of $q_\varepsilon(p)$.

\section{Calculation of $\Sigma_4(p)$: Preliminaries}
The calculation of $\Sigma_4(p)$ requires a careful study of $E(p)$. In order to pin this set down a bit more firmly, we make use of the equivalence relation $\approx$ defined on
\[
 \mathcal{T}=\bigcup_{i=1}^k\ \{(i, j): j\in S_i\}
\] 
 as follows: if $((i, j), (l, m))\in  \mathcal{T}\times \mathcal{T} $ then $(i, j)\approx (l, m)$ if $b_lj=b_im$. For all $p$ sufficiently large, $(i, j)\approx (l, m)$  if and only if $\bar{b_i}j\equiv \bar{b_l}m$ mod $p$, and so if we let $\mathcal{E}(A)$ denote the set of all nonempty subsets of even cardinality of a finite set $A$, then
\begin{quote}
for all $p$ sufficiently large, $E(p)$ consists of all subsets $T$ of $\mathcal{T}$ such that there exists a nonempty subset $\mathcal{S}$ of equivalence classes of $\approx$ and elements $E_S\in \mathcal{E}(S)$ for $S\in \mathcal{S}$ such that
\begin{equation*}
T=\bigcup_{S\in  \mathcal{S}}\ E_S.\tag{16}
\end{equation*}
\end{quote}
In particular, it follows that for all $p$ large enough, $E(p)$ does not depend on $p$, hence from now on, we delete the ``$p$" from the notation for this set.

The description of $E$ given by (16) mandates that we determine the equivalence classes of the equivalence relation $\approx$. In order to do that in a precise and concise manner, it will be convenient to use the following notation: if $b\in [1, \infty)$ and $S\subseteq [0, \infty)$, we let $b^{-1}S$ denote the set of all rational numbers of the form $z/b$, where $z$ is an element of $S$. We next let 
\[
\mathcal{K}=\Big\{\emptyset\not= K\subseteq [1, k]: \bigcap_{i\in K}\ b_i^{-1}S_i\not= \emptyset\Big\}.
\]
If $K\in \mathcal{K}$ then we set
\[
T(K)=\Big(\bigcap_{i\in K}\ b_i^{-1}S_i\Big)\cap \Big(\bigcap_{i\in [1, k]\setminus K}\ ( \mathbb{Q}\setminus b_i^{-1}S_i)\Big).
\]
 Let
\[
\mathcal{K}_{\max}=\{K\in \mathcal{K}: T(K)\not= \emptyset\}.
\]
 Using Proposition 1.4, it is then straightforward to verify that the equivalence classes of $\approx$ consist precisely of all sets of the form
\[
\{(i, tb_i): i\in K\},
\]
where $K\in \mathcal{K}_{\max}$ and $t\in T(K)$. 

Observe next that if the set
\[
\Big\{\{(i, tb_i): i\in K\}: K\in \mathcal{K}, t\in \bigcap_{i\in K}\ b_i^{-1}S_i\Big\}
\]
is ordered by inclusion then the equivalence classes of $\approx$ are the maximal elements of this set. Hence 
 $T(K)\cap T(K')=\emptyset$ whenever $\{K, K'\}\subseteq \mathcal{K}_{\max}$. Consequently, if $(K, K')\in  \mathcal{K}_{\max}\times \mathcal{K}_{\max},\emptyset \not= \sigma\subseteq K,\emptyset \not= \sigma'\subseteq K', t\in T(K)$, and $t'\in T(K')$, then $\{(i, tb_i): i\in \sigma\}$ and $ \{(i, t'b_i): i\in \sigma'\}$ are each contained in distinct equivalence classes of $\approx$ if and only if $t\not= t'$ . The following lemma is now an immediate consequence of (16) and the structure just obtained for the equivalence classes of $\approx$.
\begin{lem}
\label{lem4}
If $T\in E$ then there exists a nonempty subset $\mathcal{S}$ of $\mathcal{K}_{\max}$ , a nonempty subset $\Sigma(S)$ of $\mathcal{E}(S)$ for each $S\in \mathcal{S}$ and a nonempty subset $T(\sigma, S)$ of $T(S)$ for each $\sigma\in \Sigma(S)$ and $S\in \mathcal{S}$ such that 
\begin{equation*}
\textrm{the family of sets}\  \Big\{ T(\sigma, S): \sigma\in \Sigma(S),\  S\in \mathcal{S}\Big\}\  \textrm{is pairwise disjoint, and}
\end{equation*}
\begin{equation*}
T=\bigcup_{S\in \mathcal{S}}\ \Big[\bigcup_{\sigma\in \Sigma(S)}\Big(\bigcup_{t\in T(\sigma, S)}\ \{(i, tb_i): i\in \sigma\}\Big)\Big].
 \end{equation*}
\end{lem}

We have now determined via Lemma 9.5 the structure of the elements of $E$ precisely enough for effective use in the calculation of $\Sigma_4(p)$. However, if we already know that $q_\varepsilon(p)=0$, the value of $\Sigma_4(p)$ is obviated in our argument. It would hence be very useful to have a way to mediate between the primes $p$ for which $q_\varepsilon(p)= 0$ and the primes $p$ for which  $q_\varepsilon(p)\not= 0$. We will now define and study a gadget which does that.

\section {The $(B, \textbf{S})$-signature of a Prime}

Denote by $\Lambda(\mathcal{K})$ the set
\[
\bigcup_{K\in \mathcal{K}_{\max}} \mathcal{E}(K).
\]
Then $\Lambda(\mathcal{K})$ is empty if and only if every element of  $\mathcal{K}_{\max}$ is a singleton.

Suppose that $\Lambda(\mathcal{K})$ is not empty. We will say that $p$ is an \emph{allowable prime} if no element of $B$ has $p$ as a factor. If  $p$ is an allowable prime, then the $(B, \textbf{S})$-\emph{signature of p} is defined to be the multi-set of $\pm 1$'s given by
\[
\Big\{\chi_p\Big(\prod_{i\in I}\  b_i\Big): I\in \Lambda(\mathcal{K})\Big\}.
\]

\noindent We declare the signature of $p$ to be \emph{positive} if all of its elements are 1, and \emph{non-positive} otherwise. Let
\begin{quote}
$\Pi_+(B, \textbf{S})$ (respectively, $ \Pi_-(B, \textbf{S})$) denote the set of all allowable primes $p$ such that the $(B, \textbf{S})$-signature of $p$ is positive (respectively, non-positive).
\end{quote}
We can now prove the following two lemmas: the first records some important information about the signature, and the second implies that we need only calculate $\Sigma_4(p)$ for the primes $p$ in $\Pi_+(B, \textbf{S})$.
\begin{lem}
\label{lem2}
$(i)$ The set $\Pi_+(B, \textbf{S})$ consists precisely of all allowable primes $p$ for which each of the sets
\begin{equation*}
\{b_i: i\in I\},\  I \in \Lambda(\mathcal{K}) ,\tag{$\sharp$}
\end{equation*}
is either a set of residues of $p$ or a set of non-residues of $p$. In particular, $\Pi_+(B, \textbf{S})$ is always an infinite set.

$(ii)$ The set $\Pi_-(B, \textbf{S})$ consists precisely of all allowable primes $p$ for which at least one of the sets $(\sharp)$ contains a residue of $p$ and a non-residue of $p$, $\Pi_-(B, \textbf{S})$ is always either empty or infinite, and $\Pi_-(B, \textbf{S})$ is empty if and only if for all $I\in \Lambda(\mathcal{K}),\  \prod_{i\in I}\ b_i$ is a square.
\end{lem} 
\emph{Proof}. Suppose that $p$ is an allowable prime such that each of the sets  $(\sharp)$ is either a set of residues of $p$ or a set of non-residues of $p$. Then
\[
\chi_p\Big(\prod_{i\in I}\ b_i\Big)=1
\]
whenever $I\in  \Lambda(\mathcal{K})$ because $|I|$ is even, i.e., $p\in \Pi_+(B, \textbf{S})$. On the other hand, let $p\in \Pi_+(B, \textbf{S})$ and let $I=\{i_1,\dots,i_n\} \in  \Lambda(\mathcal{K})$. Then because $p\in \Pi_+(B, \textbf{S})$,
\[
\chi_p(b_{i_j}b_{i_{j+1}})=1,\ j\in [1, n-1],
\]
and these equations imply that $\{b_i:i\in I\}$ is either a set of residues of $p$ or a set of non-residues of $p$. This verifies the first statement in ($i$), and the second statement follows from the fact (Theorem 4.3) that there are infinitely many primes $p$ such that $B$ is a set of residues of $p$.

Statement ($ii$) of the lemma follows from ($i$), the definition of $\Pi_-(B, \textbf{S})$, and the fact (Theorem 4.2) that a positive integer is a residue of all but finitely many primes if and only if it is a square. $\hspace{13.1cm} \textrm{QED}$

It is a consequence of the following lemma that we need only calculate $\Sigma_4(p)$ for the primes $p$ which are in $\Pi_+(B, \textbf{S})$. As we will see in the next section, this greatly simplifies that calculation.
\begin{lem}
\label{lem3}
 If $p\in \Pi_-(B, \textbf{S})$ then $q_\varepsilon(p)=0$.
 \end{lem}
\emph{Proof}. If $p\in \Pi_-(B, \textbf{S})$ then there is an $I\in  \Lambda(\mathcal{K})$ such that
\[
\chi_p\Big(\prod_{i\in I} b_i\Big)=-1.
\]
Because $I$ is nonempty and of even cardinality, there exists $\{m, n\}\subseteq I$ such that
\begin{equation*}
\chi_p(b_mb_n)=-1.\tag{17}
\end{equation*}
Because $\{m, n\}$ is contained in an element of  $\mathcal{K}_{\max}$, it follows that $b_m^{-1}S_m\cap b_n^{-1}S_n\not= \emptyset$, and so we find a non-negative rational number $r$ such that
\begin{equation*}
rb_m\in S_m\ \textrm{and}\ rb_n\in S_n.\tag{18}
\end{equation*}

By way of contradiction, suppose that $q_\varepsilon(p)\not= 0$. Then there exists a $z\in [1, \infty)$ such that $b_mz+S_m$ and $b_nz+S_n$ are both contained in $[1, p-1]$ and 
\begin{equation*}
\chi_p(b_mz+u)=\chi_p(b_nz+v),\ \textrm{for all}\  u\in S_m\ \textrm{and for all}\ v\in S_n.\tag{19}
\end{equation*}
If $d$ is the greatest common divisor of $b_m$ and $b_n$ then there is a non-negative integer $t$ such that $r=t/d$. Hence by (18) and (19),
\begin{eqnarray*}
\chi_p(b_m/d)\chi_p(dz+t)&=&\chi_p(b_mz+rb_m)\\
&=&\chi_p(b_nz+rb_n)\\
&=&\chi_p(b_n/d)\chi_p(dz+t).
\end{eqnarray*}
However, $dz+t\in [1, p-1]$ and so $\chi_p(dz+t)\not= 0$. Hence
\[
\chi_p(b_m/d)=\chi_p(b_n/d),
\]
 and this value of $\chi_p$, as well as $\chi_p(d)$, is nonzero because $d, b_m/d$, and $b_n/d$ are all elements of $[1, p-1]$. But then
\[
\chi_p(b_mb_n)=\chi_p(d^2)\chi_p(b_m/d)\chi_p(b_n/d)=1,
\]
contrary to (17).   $\hspace{12.6cm} \textrm{QED}$

\section{Calculation of $\Sigma_4(p)$: Conclusion}
With Lemmas 9.5 and 9.7 in hand, we now calculate the sum $\Sigma_4(p)$ that arose in (15).
By virtue of Lemma 9.7, we need only calculate $\Sigma_4(p)$ for $p\in \Pi_+(B, \textbf{S})$, hence let $p$ be an allowable prime for which
\begin{equation*}
\chi_p\Big(\prod_{i\in I} b_i\Big)=1,\ \textrm{for all}\ I\in  \Lambda(\mathcal{K}).\tag{20}
\end{equation*}

We first recall that
\begin{equation*}
\Sigma_4(p)=2^{-\alpha}r(p)\Big(1+\sum_{T\in E}\ \prod_{i=1}^{k}\ \chi_p(b_i)^{|\{j: (i, j)\in T\}|}\Big),\tag{21}
\end{equation*}
where
\[
 r(p)=\min_i\left[\frac{p-1-\max S_i}{b_i}\right],\]
and so we must evaluate the products over $T\in E$ which determine the summands of the third factor on the right-hand side of (21). Toward that end, let $T\in E$ and use Lemma 9.5 to find a nonempty subset $ \mathcal{S}$ of $\mathcal{K}_{\max}$, a nonempty subset $\Sigma(S)$ of $\mathcal{E}(S)$ for each $S\in \mathcal{S}$ and a nonempty subset $T(\sigma, S)$ of $T(S)$ for each $\sigma\in \Sigma(S)$ and $S\in \mathcal{S}$ such that 
\begin{equation*}
\textrm{the  sets}\ T(\sigma, S), \sigma\in \Sigma(S),\  S\in \mathcal{S},\  \textrm{are pairwise disjoint, and}
\end{equation*}
\begin{equation*}
T=\bigcup_{S\in \mathcal{S}}\ \Big[\bigcup_{\sigma\in \Sigma(S)}\Big(\bigcup_{t\in T(\sigma, S)}\ \{(n, tb_n): n\in \sigma\}\Big)\Big].
 \end{equation*}
Then
\[
\{j: (i, j)\in T\}=\bigcup_{S\in \mathcal{S}}\ \Big(\bigcup_{\sigma\in \Sigma(S): i\in \sigma}\{tb_i: t\in T(\sigma, S)\}\Big)
\]
and this union is pairwise disjoint. Hence
\[
|\{j: (i, j)\in T\}|=\sum_{S\in  \mathcal{S}}\sum_{\sigma\in \Sigma(S): i\in \sigma}\ |T(\sigma, S)|.
\]
Thus from this equation and (20) we find that
\begin{eqnarray*}
\prod_{i=1}^{k}\ \chi_p(b_i)^{|\{j: (i, j)\in T\}|}&=&\prod_{i\in \cup_{S\in \mathcal{S}}\cup_{\sigma\in \Sigma(S)}\ \sigma}\ \chi_p(b_i)^{\sum_{S\in  \mathcal{S}}\sum_{\sigma\in \Sigma(S): i\in \sigma} |T(\sigma, S)|}\\
&=&\prod_{S\in  \mathcal{S}}\Big(\prod_{\sigma\in \Sigma(S)}\ \Big(\chi_p\Big(\prod_{i\in \sigma} b_i\Big)\Big)^{|T(\sigma, S)|}\Big)\\
&=&1.
\end{eqnarray*}
Hence
\begin{equation*}
\sum_{T\in E}\ \prod_{i=1}^{k}\ \chi_p(b_i)^{|\{j: (i, j)\in T\}|}=|E|,\tag{22}
\]
and so we must count the elements of $E$. In order to do that, note first that the pairwise disjoint decomposition (16) of an element $T$ of $E$ is uniquely determined by $T$, and, obviously, uniquely determines $T$. Hence if $\mathcal{D}$ denotes the set of all equivalence classes of $\approx$ of cardinality at least 2 then
\begin{eqnarray*}
|E|&=&\sum_{\emptyset\not= \mathcal{S}\subseteq \mathcal{D}}\ \prod_{S\in \mathcal{S}}\ |\mathcal{E}(S)|\\
&=&-1+\prod_{D\in \mathcal{D}} (1+|\mathcal{E}(D)|)\\
&=&-1+ \prod_{D\in \mathcal{D}} 2^{|D|-1}\\
&=&-1+2^{-|\mathcal{D}|}\cdot 2^{\sum_{D\in \mathcal{D}}|D|}. 
\end{eqnarray*}
However, $\mathcal{D}$ consists of all sets of the form
\[
\{(i, tb_i): i\in K\}
\]
where $K\in \mathcal{K}_{\max}, |K|\geq 2$, and $t\in T(K)$. Hence  
\[
|\mathcal{D}|=\sum_{K\in  \mathcal{K}_{\max}: |K|\geq 2}\ |T(K)|,
\]
\[
\sum_{D\in \mathcal{D}}|D|=\sum_{K\in  \mathcal{K}_{\max}: |K|\geq 2}\ |K||T(K)|,
\]
and so if we set
\[
e=\sum_{K\in  \mathcal{K}_{\max}}\ |T(K)|(|K|-1),
\] 
then
\begin{equation*}
|E|=2^e-1.\tag{23}
\]

\noindent Equations (21), (22), and (23) now imply
\begin{lem}
\label{lem4}
If
\[
\alpha=\sum_i|S_i|,\ e=\sum_{K\in  \mathcal{K}_{\max}}\ |T(K)|(|K|-1),\ \textrm{and}\  r(p)=\min_i\left[\frac{p-1-\max S_i}{b_i}\right],
\] 
then
\[
\Sigma_4(p)=2^{e-\alpha}r(p),\ \textrm{for all}\  p\in \Pi_+(B, \textbf{S}).
\]
\end{lem}

If we set $\displaystyle{b=\max_i\{b_i\}}$ then it follows from Lemma 9.8 that as $p\rightarrow +\infty$ inside $\Pi_+(B, \textbf{S})$,
\[
\Sigma_4(p)\sim (b\cdot 2^{\alpha-e})^{-1}p.\]
When we insert this asymptotic approximation of $\Sigma_4(p)$ into the estimate (15), and then recall Lemma 9.7, we see that $ (b\cdot 2^{\alpha-e})^{-1}p$ is a linear function of $p$ which should work to determine the asymptotic behavior of $q_{\varepsilon}(p)$. We will now show in the next section that it does work in exactly that way.

\section{Solution of Problems 2 and 4: Conclusion}

 All of the ingredients are now assembled for a proof of the following theorem, which determines the asymptotic behavior of $q_\varepsilon(p)$.

\begin{thm}
\label{thm5}
$($Wright $[62]$, Theorem $6.1)$ Let $\varepsilon\in \{-1, 1\}, k\in [1, \infty)$, and let $B=\{b_1,\dots,b_k\}$ be a set of positive integers and $\textbf{S}=(S_1,\dots,S_k)$ a $k$-tuple of finite, nonempty subsets of $[0, \infty)$. If $ \mathcal{K}_{\max}$ is the set of  subsets of $[1, k]$ defined by $B$ and $\textbf{S}$ as in section $7$, 
let 
\[
\Lambda(\mathcal{K})=\bigcup_{K\in \mathcal{K}_{\max}} \mathcal{E}(K),
\]
\[
\alpha=\sum_i|S_i|,\ b=\max_i\{b_i\},\  e=\sum_{K\in  \mathcal{K}_{\max}}\ |T(K)|(|K|-1),\ and
\]
\[
 q_\varepsilon(p)=|\{A\in AP(B, \textbf{S})\cap 2^{[1, p-1]}: \chi_p(a)=\varepsilon, \ \textrm{for all}\ a\in A\}|.
\]
$(i)$ If the sets $b_1^{-1}S_1,\dots,b_k^{-1}S_k$ are pairwise disjoint  then
\[
q_\varepsilon(p)\sim (b\cdot 2^\alpha)^{-1}p\ as\ p\rightarrow +\infty.
\]

\noindent$(ii)$ If the sets $b_1^{-1}S_1,\dots,b_k^{-1}S_k$ are not pairwise disjoint then

$(a)$ the parameter $e$ is positive and less than $\alpha$;

$(b)$ if $\prod_{i\in I} b_i$ is a square for all $I\in \Lambda(\mathcal{K})$ then
\[
q_\varepsilon(p)\sim (b\cdot 2^{\alpha-e})^{-1}p\ as\ p\rightarrow +\infty;
\]

$(c)$ if there exists $I\in \Lambda(\mathcal{K})$ such that $\prod_{i\in I} b_i$ is not a square then

$(\alpha)$ the set $\Pi_+(B, \textbf{S})$ of primes with positive $(B, \textbf{S})$-signature and the set $\Pi_-(B, \textbf{S})$ of primes with non-positive $(B, \textbf{S})$-signature are both infinite,

$(\beta)$ $q_\varepsilon(p)=0$ for all $p$ in $\Pi_-(B, \textbf{S})$, and

$(\gamma)$ as $p\rightarrow +\infty$ inside $\Pi_+(B, \textbf{S})$,
\[
q_\varepsilon(p)\sim (b\cdot2^{\alpha-e})^{-1}p\ .
\]
 \end{thm}
 \emph{Proof}. If the sets $b_1^{-1}S_1,\dots,b_k^{-1}S_k$ are pairwise disjoint  then every element of $\mathcal{K}_{\max}$ is a singleton set, hence all of the equivalence classes of the equivalence relation $\approx$ defined above on $\bigcup_{i=1}^k\ \{(i, j): j\in S_i\}
$ by the set $B$ are singletons. It follows that the set $E$ which is summed over in (21) is empty and so  
\begin{equation*}
\Sigma_4(p)=2^{-\alpha}r(p),\ \textrm{for all}\ p\ \textrm{sufficiently large}.\tag{24}
\]
Upon recalling that
\[
r(p)=\min_i\left[\frac{p-1-\textrm{max}\ S_i}{b_i}\right],
\]
and then noting that as $p\rightarrow +\infty$, $r(p)\sim p/b$, the conclusion of $(i)$ is an immediate consequence of (15) and (24).

Suppose that the sets $b_1^{-1}S_1,\dots,b_k^{-1}S_k$ are not pairwise disjoint. Then $\Lambda(\mathcal{K})$ is not empty and so conclusion $(a)$ is an obvious consequence of the definition of $e$. If $\prod_{i\in I} b_i$ is a square for all $I\in \Lambda(\mathcal{K})$ then it follows from its definition that $\Pi_+(B, \textbf{S})$ contains all but finitely many primes, and so $(b)$ is an immediate consequence of (15) and Lemma 9.8. On the other hand, if there exists $I\in \Lambda(\mathcal{K})$ such that $\prod_{i\in I} b_i$ is not a square then $(\alpha)$ follows from Lemma 9.6, $(\beta)$ follows from Lemma 9.7, and $(\gamma)$  is an immediate consequence of (15) and Lemma 9.8. $\hspace{14.77cm} \textrm{QED}$  

Theorem 9.9 shows that the elements of $\Lambda(\mathcal{K})$ contribute to the formation of quadratic residues and non-residues inside $AP(B, \textbf{S})$. If no such elements exist then $q_\varepsilon(p)$ has the expected minimal asymptotic approximation $(b\cdot 2^{\alpha})^{-1}p$ as $p\rightarrow +\infty$. In the presence of elements of  $\Lambda(\mathcal{K})$, the parameter $e$ is positive and less than $\alpha$, the  asymptotic size of $q_\varepsilon(p)$ is increased by a factor of $2^e$, and whenever $\Pi_-(B, \textbf{S})$ is empty, $q_\varepsilon(p)$ is asymptotic to $(b\cdot 2^{\alpha-e})^{-1}p$ as $p\rightarrow +\infty$. However, the most interesting behavior occurs when $\Pi_-(B, \textbf{S})$ is not empty; in that case, as $p\rightarrow +\infty, q_\varepsilon(p)$ asymptotically oscillates infinitely often between 0 and $(b\cdot 2^{\alpha-e})^{-1}p$. 

\emph{Remark}. If we observe that the cardinality of the set
\[
\bigcup_{i=1}^k\ b_i^{-1}S_i
\] 
is equal to the number of equivalence classes of the equivalence relation $\approx$ that was defined on the set
\[
\mathcal{T}=\bigcup_{i=1}^k\ \{(i, j): j\in S_i\},
\]
then it follows that
\[
\big|\bigcup_{i=1}^k\ b_i^{-1}S_i\big|=\sum_{K\in  \mathcal{K}_{\max}}\ |T(K)|.
\]
But we also have that
\[
\alpha=|\mathcal{T}|=\sum_{K\in  \mathcal{K}_{\max}}\ |T(K)||K|.
\]
Consequently, the exponents in the power of $1/2$ that occur in the asymptotic approximation to $q_{\varepsilon}(p)$ in Theorem 9.9 are in fact all equal to the cardinality of $\bigcup_{i=1}^k\ b_i^{-1}S_i$.

Theorem 9.9 will now be applied to the situation of primary interest to us here, namely to the family of sets  $AP(\textbf{a}, \textbf{b}; s)$ determined by a standard $2m$-tuple $(\textbf{a}, \textbf{b})$. In this case, the decomposition (11) of the sets in $AP(\textbf{a}, \textbf{b}; s)$ shows that there is a set $B=\{b_1,\dots,b_k\}$ of positive integers (the set of \emph{distinct} values of the coordinates of $\textbf{b}$), a $k$-tuple $(m_1,\dots,m_k)$ of positive integers such that $m=\sum_i m_i$, and sets
\[
A_i=\{a_{i1},\dots,a_{im_i}\}
\]
of non-negative integers, all uniquely determined by  $(\textbf{a}, \textbf{b})$, such that if we let
\begin{equation*}
S_i=\bigcup_{j=1}^{m_i}\ \{a_{ij}+b_il: l\in [0, s-1]\},\  i\in [1, k],\tag{25}
\]  
and set
\[
\textbf{S}=(S_1,\dots,S_k)
\]
then
\[
AP(\textbf{a}, \textbf{b}; s)= AP(B, \textbf{S}).
\]
 It follows that
\[
b_i^{-1}S_i=\bigcup_{q\in b_i^{-1}A_i}\ \{q+j: j\in [0, s-1]\},\ i\in [1, k].
\]
These sets then determine the subsets of $[1, k]$ that constitute 
\[
\mathcal{K}=\{\emptyset\not= K\subseteq [1, k]: \bigcap_{i\in K}\ b_i^{-1}S_i\not= \emptyset\}\}
\]
and hence also the elements of $\mathcal{K}_{\max}$, according to the recipe given in section 7.
The sets in $\mathcal{K}_{\max}$, together with the parameters
\[
\alpha=\sum_i|S_i|,\ b=\max_i\{b_i\},\  \textrm{and}\  e=\sum_{K\in  \mathcal{K}_{\max}}\ |T(K)|(|K|-1),
\]
when used as specified in Theorem 9.9, then determine precisely the asymptotic behavior of the sequence $q_\varepsilon(p)$ that is defined upon replacement of $AP(B, \textbf{S})$ by $AP(\textbf{a}, \textbf{b}; s)$ in the statement of Theorem 9.9, thereby solving Problems 2 and 4. In particular, the sets $b_1^{-1}S_1,\dots,b_k^{-1}S_k$ are pairwise disjoint if and only if 
 \begin{quote}
 $(26)$  if $(i, j)\in [1, k]\times [1, k]$ with $i\not= j$ and $(x, y)\in A_i\times A_j$, then either $b_ib_j$ does not divide $yb_i-xb_j$ or  $b_ib_j$ divides $yb_i-xb_j$ with a quotient that exceeds $s-1$ in modulus.
 \end{quote}
Hence the conclusion of statement $(i)$ of Theorem 9.9  holds for $AP(\textbf{a}, \textbf{b}; s)$ when condition (26) is satisfied, while the conclusions of statement $(ii)$ of Theorem 9.9 hold for $AP(\textbf{a}, \textbf{b}; s)$ whenever condition (26) is not satisfied. In the following section we will present several examples which illustrate how Theorem 9.9 works in practice to determine the asymptotic behavior of $q_{\varepsilon}(p)$. We will see there, in particular, that for each integer $m\in [2, \infty)$ and for each of the hypotheses in the statement of Theorem 9.9, there exists infinitely many standard $2m$-tuples $(\textbf{a}, \textbf{b})$ which satisfy that hypothesis. 

\section{An Interesting Class of Examples}

In order to apply Theorem 9.9 to a standard $2m$-tuple $(\textbf{a}, \textbf{b})$, we need to calculate the parameters $\alpha$ and $e$, the set $\Lambda(\mathcal{K})$, and the associated signatures of the allowable primes. In general, this can be somewhat complicated, but there is a class of standard $2m$-tuples for which these computations can be carried out by means of easily applied algebraic and geometric formulae, which we will discuss next.

Let $k \in [2, \infty)$. We will say that a standard $2k$-tuple $(\textbf{a}, \textbf{b})$ of integers is \emph{admissible} if it satisfies the following two conditions:
\begin{equation*}
\textrm{the coordinates of}\ \textbf{b}\ \textrm {are distinct, and}\tag{27},
\]
\begin{equation*}
a_ib_j-a_jb_i\not= 0\ \textrm{for}\ i\not= j\tag{28}.
\end{equation*}
If $s\in [1, \infty)$ and $(\textbf{a}, \textbf{b})$ is admissible then it follows trivially from (27) that
\[
S_i=\{a_i+b_ij: j\in [0, s-1]\},\ i\in [1, k], \]
hence
\[
|S_i|=s,\ i \in [1, k],\]
and so the parameter $\alpha$ in the statement of Theorem 9.9 for $AP(\textbf{a}, \textbf{b}; s)$ is $ks$.

We turn next to the calculation of the parameter $e$. Let $q_i=a_i/b_i, i \in [1, k]$;  (28) implies that the $q_i$'s are distinct, and without loss of generality, we suppose that the coordinates of $\textbf{a}$ and $\textbf{b}$ are indexed so that $q_i<q_{i+1}$ for each $i\in [1, k-1]$. Let $\textbf{R}$ denote the set of all subsets $R$ of $\{q_1,\dots,q_k\}$ such that $|R|\geq 2$ and $R$ is maximal relative to the property that $w-z$ is an integer for all $(w, z)\in R\times R$. We note that $\textbf{R}$ is just the set of all equivalence classes of cardinality at least 2 of the equivalence relation $\sim$ defined on the set $\{q_1,\dots,q_k\}$ by declaring that $q_i \sim q_j$ if $q_i-q_j \in \mathbb{Z}$. After linearly ordering the elements of each $R\in \textbf{R}$, we let $D(R)$ denote the $(|R|-1)$-tuple of positive integers whose coordinates are the distances between consecutive elements of $R$. Then if $M_R(s)$ denotes the multi-set formed by the coordinates of $D(R)$ which do not exceed $s-1$, it can be shown that 
\begin{equation*}
e=\sum_{R\in \textbf{R}}\ \sum_{r\in M_R(s)} \ (s-r)\tag{29}
\]
(see Wright [62], section 8). We note in particular that $e=0$ if and only if the set $\{R\in \textbf{R}: M_R(s)\not= \emptyset\}$ is empty and that this occurs if and only if the sets $b_i^{-1}S_i, i \in [1, k]$, are pairwise disjoint. Formula (29) shows that $e$ can be calculated solely by means of information obtained directly and straightforwardly from the set $\{q_1,\dots,q_k\}$.   

In order to calculate the signature of allowable primes, the set $\Lambda(\mathcal{K})$ must be computed. There is an elegant geometric formula for this computation that is based on the concept of what we will call an overlap diagram, and so those diagrams will be described first.

Let $(n, s)\in [1,\infty)\times [1,\infty)$ and let $ \textbf{g}=(g(1),\dots,g(n))$ be an $n$-tuple of positive integers. We use \textbf{g} to construct the following array of points. In the plane, place $s$ points horizontally one unit apart, and label the $j$-th point as $(1, j-1)$ for each $j\in [1, s]$. This is \emph{row $1$}. Suppose that row $i$ has been defined. One unit vertically down and $g(i)$ units horizontally to the right of the first point in row $i$, place $s$ points horizontally one unit apart, and label the $j$-th point as $(i+1, j-1)$ for each $j\in [1, s]$. This is $\emph{row}\ i+1$.
The array of points so formed by these $n+1$ rows is called the \emph{overlap diagram of} \textbf{g}, the sequence \textbf{g} is called the \emph{gap sequence} of the overlap diagram, and a nonempty set that is formed by the intersection of the diagram with a vertical line is called a \emph{column} of the diagram. N.B. We do not distinguish between the different possible positions in the plane which the overlap diagram may occupy. A typical example with $n=3, s=8$, and gap sequence (3, 2, 2) looks like
\vspace{0.5cm}
\begin{center}
\begin{tabular}{ccccccccccccccc}
 $\bigcdot$&$\bigcdot$&$\bigcdot$&$\bigcdot$&$\bigcdot$&$\bigcdot$&$\bigcdot$&$\bigcdot$&&&&&&&\\

&&&$\bigcdot$&$\bigcdot$&$\bigcdot$&$\bigcdot$&$\bigcdot$&$\bigcdot$&$\bigcdot$&$\bigcdot$&&\\
&&&&&$\bigcdot$&$\bigcdot$&$\bigcdot$&$\bigcdot$&$\bigcdot$&$\bigcdot$&$\bigcdot$&$\bigcdot$\\
&&&&&&&$\bigcdot$&$\bigcdot$&$\bigcdot$&$\bigcdot$&$\bigcdot$&$\bigcdot$&$\bigcdot$&$\bigcdot$\     \\
\end{tabular}
\end{center}

\vspace{0.2cm} 
\begin{center}
 An overlap diagram
\end{center}

\vspace{0.5cm}
 We need to describe how and where rows overlap in an overlap diagram. Begin by first noticing that if $(g(1),\dots,g(n))$ is the gap sequence, then row $i$ overlaps row $j$ for $i<j$ if and only if
 \[
 \sum_{r=i}^{j-1} g(r)\leq s-1;
 \]                                                         
in particular, row $i$ overlaps row $i+1$ if and only if $g(i)\leq s-1$. Now let $\mathcal{G}$ denote the set of all subsets $G$ of $[1, n]$ such that $G$ is a nonempty set of consecutive integers maximal with respect to the property that $g(i)\leq s-1$ for all $i\in G$. If $\mathcal{G}$ is empty then $g(i)\geq s$ for all $i\in [1, n]$, and so there is no overlap of rows in the diagram. Otherwise there exists $m\in[1, 1+[(n-1)/2]]$ and strictly increasing sequences $(l_1,\dots,l_m)$ and $(M_1,\dots,M_m)$ of positive integers, uniquely determined by the gap sequence of the diagram, such that $l_i\leq M_i$ for all $i\in [1, m], 1+M_i\leq l_{i+1}$ if $i\in [1, m-1]$, and
\[
\mathcal{G}=\{[l_i, M_i]: i\in [1, m]\}.
\]
In fact, $l_{i+1}> 1+M_i$ if $i\in [1, m-1]$, lest the maximality of the elements of $\mathcal{G}$ be violated. It follows that the intervals of integers $[l_i, 1+M_i], i\in [1, m]$, are pairwise disjoint. 

The set $\mathcal{G}$ can now be used to locate the overlap between rows in the overlap diagram like so: for $i\in [1, m]$, let
\[
B_i=[l_i, 1+M_i],
\]
and set
\[
\mathcal{B}_i=\textrm{the set of all points in the overlap diagram whose labels are in}\ B_i\times [0, s-1].
\]
We refer to $\mathcal{B}_i$ as the \emph{i-th block} of the overlap diagram; thus the blocks of the diagram are precisely the regions in the diagram in which rows overlap. 

We will now use the elements of $\textbf{R}$ to construct a series of overlap diagrams. Let $R$ be an element of $\textbf{R}$ such that $D(R)$ has at least one coordinate that does not exceed $s-1$. Next, consider the nonempty and pairwise disjoint family of all subsets $V$ of $R$ such that $|V|\geq 2$ and $V$ is maximal with respect to the property that the distance between consecutive elements of $V$ does not exceed $s-1$. List the elements of $V$ in increasing order and then for each $i\in [1, |V|-1]$ let $q_V(i)$ denote the distance between the $i$-th element and the $(i+1)$-th element on that list. N.B. $q_V(i) \in [1, \infty)$, for all $i\in [1, |V|-1]$. Finally, let $\mathcal{D}(V)$ denote the overlap diagram of the $(|V|-1)$-tuple $(q_V(i): i\in [1, |V|-1])$. Because $q_V(i)\leq s-1$ for all $i\in [1, |V|-1]$, $\mathcal{D}(V)$ consists of a single block.  

Using a suitable positive integer $m$, we index all of the sets $V$ that arise from all of the elements of $\textbf{R}$ in the previous construction as $V_1,\dots, V_m$ and then define the \emph{quotient diagram} of $(\textbf{a}, \textbf{b})$ to be the $m$-tuple of overlap diagrams $(\mathcal{D}(V_n): n\in [1, m])$. We will refer to the diagrams $\mathcal{D}(V_n)$ as the \emph{blocks} of the quotient diagram.

The quotient diagram $\mathcal{D}$ of $(\textbf{a}, \textbf{b})$ will now be used to calculate the set $\Lambda(\mathcal{K})$ determined by $(\textbf{a}, \textbf{b})$ and hence the associated signature of an allowable prime. In order to see how this goes, we will need to make use of a certain labeling of the points of  $\mathcal{D}$  which we describe next. Let $V_1,\dots,V_m$ be the subsets of $\{q_1,\dots, q_k\}$ that determine the sequence of overlap diagrams $\mathcal{D}(V_1),\dots,\mathcal{D}(V_m)$ which constitute $\mathcal{D}$, and then find the subset $J_n$ of $[1, k]$ such that $V_n=\{q_j: j\in J_n\}$, with $j\in J_n$ listed in increasing order (note that this ordering of $J_n$ also linearly orders $q_j, j \in J_n$).The overlap diagram $\mathcal{D}(V_n)$ consists of $|J_n|$ rows, with each row containing $s$ points. If $i\in [1, |J_n|]$ is taken in increasing order then there is a unique element $j$ of $J_n$ such that the $i$-th element of $V_n$ is $q_j$. Proceeding from left to right in each row, we now take $l\in [1, s]$ and label the $l$-th point of row $i$ in $\mathcal{D}(V_n)$ as $(j, l-1)$. N.B. This labeling of the points of $\mathcal{D}(V_n)$ \emph{does not} necessarily coincide with the labeling of the points of an overlap diagram that was used before to define the blocks of the diagram.

Next let  $C$ denote a column of  one of the diagrams $\mathcal{D}(V_n)$ which constitute  $\mathcal{D}$. We identify  $C$  with the subset of $[1, k]\times [0, s-1]$ defined by  
\begin{equation*}
\{(i, j)\in [1, k]\times [0, s-1]: (i, j)\ \textrm{is the label of a point in}\ C\}, \tag{30}
\]
let $\mathcal{C}_n$ denote the set of all subsets of  $[1, k]\times [0, s-1]$ which arise from all such identifications, and then set $\mathcal{C}=\bigcup_n\ \mathcal{C}_n$. If $\theta$ denotes the projection of $[1, k]\times [0, s-1]$ onto $[1, k]$ then one can show (Wright [63], Lemma 2.5) that $K\in \mathcal{K}_{\max}$ if and only if there exists a $T\in \mathcal{C}$ such that $K=\theta(T)$, and so
\begin{equation*}
\Lambda(\mathcal{K})=\bigcup_{T\in \mathcal{C}}\ \mathcal{E}(\theta(T)).\tag{31}
\]
When this formula for $\Lambda(\mathcal{K})$ is now combined with (29), it follows that all of the data required for an application of Theorem 9.9 can be easily read off directly from the set $\{q_1, \dots ,q_k\}$ and the quotient diagram of $(\textbf{a}, \textbf{b})$.

At this juncture, some concrete examples which illustrate the mathematical technology that we have introduced are in order. But before we get to those, recall that if $(\textbf{a}, \textbf{b})$ is an admissible $2k$-tuple, $B$ is the set formed by the coordinates $b_1,\dots,b_k$ of $\textbf{b}, S_i=\{a_i+b_ij: j\in [0, s-1]\}$, where $a_i$ is the $i$-th coordinate of $\textbf{a},\  i\in [1, k],$ and $\textbf{S}$ is the $k$-tuple of sets $(S_1,\dots,S_k)$, then the pair $(B, \textbf{S})$ determines by way of Theorem 9.9 the asymptotic behavior of $|\{A\in AP(\textbf{a}, \textbf{b}; s)\cap 2^{[1, p-1]}: \chi_p(a)=\varepsilon, \ \textrm{for all}\ a\in A\}|,\ \varepsilon \in \{-1, 1\}$. Hence for this pair, we use the more specific notation $\Pi_{\pm}(\textbf{a}, \textbf{b})$ for the sets $\Pi_{\pm}(B, \textbf{S})$ in the statement of Theorem 9.9.

Now for the examples. We start with a simple example which illustrates how the parameter $e$ and the set $\Lambda(\mathcal{K})$ are calculated from (29) and (31). Suppose that $s=5$ and the quotient diagram of the admissible 8-tuple $(\textbf{a}, \textbf{b})$ consists of the single overlap diagram located in the plane as follows:
\begin{figure}[h]
  \centering
  \begin{tikzpicture}
    \draw[name path=xaxis, thick] (0,0) -- (12,0);
    \foreach \q/\y/\b in {1/4/1,2/3/3,3/2/5,4/1/8} {
      \draw (\b,0) circle (1pt) node[anchor=north] {$q_{\q}$};
        \foreach \x in {0,1,2,3,4} {
          \fill ({\b+\x},\y) circle (2pt);
        }
    }
  \end{tikzpicture}
  \caption{Location of the quotient diagram}
  \label{fig:9.1}
\end{figure}

\newpage

\noindent Then $k=4, q_2=q_1+2, q_3=q_2+2$, and $q_4=q_3+3$, hence the set $\textbf{R}$ consists of the single set $R=\{q_1, q_2, q_3, q_4\}$, whence $D(R)=(2, 2, 3)$ and $M_{R}(5)=\{2, 2, 3\}$. It therefore follows from (29) that
\[
e=(5-2)+(5-2)+(5-3)=8.\]
Consequently,\[
\alpha-e=4\cdot5-8=12.\]
We have that
\[
b_i^{-1}S_i=\{q_i+j: j\in [0, 4]\},\ i=1, 2, 3, 4,\]
and so 12 is also the cardinality of the union $\bigcup_{i=1}^4 b_i^{-1}S_i$. 

Turning to the calculation of $\Lambda(\mathcal{K})$ by means of (31), note first that the quotient diagram of $(\textbf{a}, \textbf{b})$ consists of a single block $\mathcal{D}(V)$ with $V=\{q_1, q_2, q_3, q_4\}$. The set  of indices of the elements of $V$ is $J=\{1, 2, 3, 4\}$, and so the points of $\mathcal{D}(V)$ are labeled as indicated in Figure 2:

\begin{figure}[h]
  \centering
  \begin{tikzpicture}
    \draw[name path=xaxis, thick] (0,0) -- (12,0);
    \foreach \q/\y/\b in {1/4/1,2/3/3,3/2/5,4/1/8} {
      \draw (\b,0) circle (1pt) node[anchor=north] {$q_{\q}$};
        \foreach \x in {0,1,2,3,4} {
          \fill ({\b+\x},\y) circle (2pt) node[anchor=south] {$(\q,\x)$};
        }
    }
  \end{tikzpicture}
  \caption{Labeled points of the quotient diagram}
  \label{fig:9.1}
\end{figure}
\newpage

\noindent The columns in $\mathcal{C}$, identified as subsets of $\{1, 2, 3, 4\}\times \{0, 1, 2, 3, 4\}$ are hence 
\begin{quote}
$\{(1, 0)\}, \{(1, 1)\}$, $\{(1, 2), (2, 0)\}, \{(1, 3), (2, 1)\}, \{(1, 4), (2, 2), (3, 0)\}, \{(2, 3)$, $(3, 1)\}, \{(2, 4), (3, 2)\}, \{(3, 3), (4, 0)\}$, $\{(3, 4), (4, 1)\}, \{(4, 2)\}, \{(4, 3)\}$ and 
$\{(4, 4)\}$. 
\end{quote}
From this, we find that the sets $\theta(T), T\in \mathcal{C},$ are 
\[
\{1\}, \{1, 2\}, \{1, 2, 3\}, \{2, 3\}, \{3, 4\}, \textrm{ and}\ \{4\},\] 
and so $\Lambda(\mathcal{K})$, according to (31), consists of the sets 
\[
\{1, 2\}, \{1, 3\}, \{2, 3\}, \textrm{and} \{3, 4\}. \]
Consequently, the $(\textbf{a}, \textbf{b})$-signature of an allowable prime $p$ is 
\[
\{\chi_p(b_1b_2), \chi_p(b_1b_3), \chi_p(b_2b_3), \chi_p(b_3b_4)\}.\]

Next, let $m\in [1, +\infty)$ and for each $n\in [1, m]$, let $\mathcal{D}(n)$ be a fixed but arbitrary overlap diagram with $k_n$ rows, $k_n\geq 2$, and gap sequence $(d(i, n): i\in [1, k_n-1])$, with no gap exceeding $s-1$. Let $k_0=0, k=\sum_{n=0}^m k_n$. We will now exhibit infinitely many admissible $2k$-tuples $(\textbf{a}, \textbf{b})$ whose quotient diagram is $\Delta=(\mathcal{D}(n): n\in [1, m])$. This is done by taking the $(k-1)$-tuple $(d_1,\dots, d_{k-1})$ in the following lemma to be 
\[
d_i=\left\{\begin{array}{rl}d\Big(i-\sum_0^n k_j, n+1\Big),& \textrm{if}\ n\in [0, m-1]\ \textrm{and}\  i\in \Big[1+\sum_0^n k_j, -1+\sum_0^{n+1} k_j\Big],\\
s,& \textrm{elsewhere},\end{array}\right.
\]
and then letting $(\textbf{a}, \textbf{b})$ be any $2k$-tuple obtained from the construction in the lemma.
\begin{lem}
\label{lem1}
For $k\in [2, \infty)$, let $(d_1,\dots, d_{k-1})$ be a $(k-1)$-tuple of positive integers. Define $k$-tuples $(a_1,\dots, a_k), (b_1,\dots, b_k)$ of positive integers inductively as follows: let $(a_1, b_1)$ be arbitrary, and if $i>1$ and $(a_i, b_i)$ has been defined, choose $t_i\in [2, \infty)$ and set
\[
a_{i+1}=t_i(a_i+d_ib_i),\ \ b_{i+1}=t_ib_i.
\]
Then
\[
\frac{a_i}{b_i}-\frac{a_j}{b_j}=\sum_{r=j}^{i-1}d_r,\ \textrm{for all}\  i>j.
\]
\end{lem}

\emph{Proof}. This is a straightforward calculation using the recursive definition of the $k$-tuples $(a_1,\dots, a_k)$ and $ (b_1,\dots, b_k)$. $\hspace{10.5cm} \textrm{QED}$

We can use Lemma 9.10 to also find infinitely many admissible $2k$-tuples $(\textbf{a}, \textbf{b})$  with quotient diagram $\Delta$ and  such that the set $\Pi_-(\textbf{a}, \textbf{b})$ is empty. To do this, simply choose the integer $b_1$ and all subsequent $t_i$'s used in the above construction from Lemma 9.10 to be squares. This shows that there are infinitely many admissible $2k$-tuples with a specified quotient diagram which satisfy the hypotheses of Theorem 9.9($ii$)($b$). On the other hand, if $b_1$ and all the subsequent $t_i$'s are instead chosen to be distinct primes, it follows that the $2k$-tuples determined in this way all have quotient diagram $\Delta$ and each have $\Pi_-(\textbf{a}, \textbf{b})$ of infinite cardinality, and so there are infinitely many admissible $2k$-tuples with specified quotient diagram which satisfy the hypotheses of Theorem 9.9($ii$)($c$). We also note that if all of the coordinates of $(d_1,\dots, d_{k-1})$ in Lemma 9.10 are chosen to exceed $s-1$ then we obtain infinitely many admissible $2k$-tuples which satisfy the hypothesis of Theorem 9.9($i$). 

For example, suppose we want to find infinitely many admissible 16-tuples $(\textbf{a}, \textbf{b})$ whose quotient diagram consists of two copies of the overlap diagram in Figure 1, and for which $\Pi_-(\textbf{a}, \textbf{b})$ is empty.  We note first that the gap sequence of the overlap diagram is $(2, 2, 3)$, then take the 7-tuple in Lemma 9.10 to be (2, 2, 3, 5, 2, 2, 3), select $n\in [2, \infty)$, set $t_i=n^2$  and $a_1=b_1=1$ in the recursive formulae for $a_{i+1}$ and $b_{i+1}, i=1, 2, 3, 4, 5, 6, 7$, to obtain 
\[
(\textbf{a}, \textbf{b})=(1, 3n^2, 5n^4, 8n^6, 13n^8, 15n^{10}, 17n^{12}, 20n^{14}, 1, n^2, n^4,n^6, n^8, n^{10}, n^{12}, n^{14}).\] 
If we also wish to find infinitely many admissible 16-tuples $(\textbf{a}, \textbf{b})$  with this quotient diagram, but for which $\Pi_-(\textbf{a}, \textbf{b})$ is infinite, then in this same recipe, let $\{p_1, p_2, p_3, p_4, p_5, p_6, p_7\}$ be any set of 7 primes and take $t_i=p_i$ for $i=1, 2, 3, 4, 5, 6, 7$ to obtain 
\begin{quote}
$(\textbf{a}, \textbf{b})=\big(1, 3p_1, 5\prod_{i=1}^2 p_i, 8\prod_{i=1}^3 p_i, 13\prod_{i=1}^4 p_i, 15\prod_{i=1}^5 p_i,17\prod_{i=1}^6 p_i, 20\prod_{i=1}^7 p_i,$
\end{quote}
\begin{quote}
\hspace{3.2cm} $1, p_1, \prod_{i=1}^2 p_i, \prod_{i=1}^3 p_i, \prod_{i=1}^4 p_i, \prod_{i=1}^5 p_i, \prod_{i=1}^6 p_i, \prod_{i=1}^7 p_i\big)$. 
\end{quote}
Finally, to find infinitely many admissible 16-tuples  $(\textbf{a}, \textbf{b})$ which satisfy the hypothesis of Theorem 9.9($i$), take the 7-tuple in Lemma 9.10 to be (5, 5, 5, 5, 5, 5, 5) and $t_i=n$ for $i=1, 2, 3, 4, 5, 6, 7$ for $n\in [2, \infty)$ to obtain 
\[
(\textbf{a}, \textbf{b})=(1, 6n, 11n^2, 16n^3, 21n^4, 26n^5, 31n^6, 36n^7, 1, n, n^2, n^3, n^4, n^5, n^6, n^7).\]

With this cornucopia of examples in hand, for $\varepsilon\in \{-1, 1\}$, we let $q_\varepsilon(p)$ denote the cardinality of the set
\begin{equation*}
\{A\in AP(\textbf{a}, \textbf{b}; s)\cap 2^{[1, p-1]}: \chi_p(a)=\varepsilon, \ \textrm{for all}\ a\in A\},
\end{equation*} 
where $(\textbf{a}, \textbf{b})$ is admissible. We will now use the quotient diagram of $(\textbf{a}, \textbf{b})$, formulae (29), (31), and Theorem 9.9 to study how $(\textbf{a}, \textbf{b})$ determines the asymptotic behavior of $q_{\varepsilon}(p)$ in specific situations. We will illustrate how things work when $k=2$ and 3, and for when ``minimal" or ``maximal" overlap is present in the quotient diagram of $(\textbf{a}, \textbf{b})$.     

When $k=2$, there is only at most a single overlap of rows in the quotient diagram of $(\textbf{a}, \textbf{b})$, and if, e.g., $a_1b_2-a_2b_1=qb_1b_2$ with $0<q\leq s-1$, then the quotient diagram looks like
\newpage
\begin{center}
\begin{tabular}{ccccccccccccc}
$\bigcdot$&$\bigcdot$&$\bigcdot$&$\bigcdot$&$\bigcdot$&$\bigcdot$&$\bigcdot$&$\bigcdot$\ &&&&&\\
$\leftarrow$&$q$&$\rightarrow$&$\bigcdot$&$\bigcdot$&$\bigcdot$&$\bigcdot$&$\bigcdot$&\ $\bigcdot$&\ $\bigcdot$&\ $\bigcdot$\  &\  \\
\end{tabular}
\end{center}

\vspace{0.2cm} 
\begin{center}
F\textsc{igure} 3. A quotient diagram for $k=2$
\end{center}

\vspace{0.5cm}
\noindent  We have that $\alpha=2s$ and, because of (29), $e=s-q$. Formula  (31) shows that the signature of $p$ is $\{\chi_p(b_1b_2)\}$, and so we conclude from Theorem 9.9 that when $b_1b_2$ is a square,
\[
q_\varepsilon(p)\sim (b\cdot 2^{s+q})^{-1}p,\ \textrm{as}\ p\rightarrow +\infty,
\]
and when $b_1b_2$ is not a square, $\Pi_+(\textbf{a}, \textbf{b})$ is the set of all allowable primes $p$ such that $\{b_1, b_2\}$ is either a set of  residues of $p$ or a set of  non-residues of $p$, $\Pi_-(\textbf{a}, \textbf{b})$ is the set of all allowable primes $p$ such that $\{b_1, b_2\}$ contains a residue of $p$ and a non-residue of $p$,  
\[
q_\varepsilon(p)=0,\ \textrm{for all}\ p\ \textrm{in}\ \Pi_-(\textbf{a}, \textbf{b}),
\]
and as $p\rightarrow +\infty$ inside $\Pi_+(\textbf{a}, \textbf{b})$,
\[
q_\varepsilon(p)\sim (b\cdot 2^{s+q})^{-1}p.
\]

 When $k=3$ there are exactly three types of overlap possible in the quotient diagram of $(\textbf{a}, \textbf{b})$, determined, e.g., when either

 $(i)$ exactly one,

 $(ii)$ exactly two, or

 $(iii)$ exactly three
 
 \noindent of $b_1b_2$, $b_2b_3$, and $b_1b_3$ divide, respectively, $a_2b_1-a_1b_2, a_3b_2-a_2b_3$, and $a_3b_1-a_1b_3$ with positive quotients not exceeding $s-1$.

In case $(i)$, with $a_2b_1-a_1b_2=qb_1b_2$, say, the block in the quotient diagram of $(\textbf{a}, \textbf{b})$ is formed by a single overlap between rows 1 and 2, and this block looks exactly like the overlap diagram that was displayed for $k=2$ above. It follows that the conclusions from (29), (31), and Theorem 9.9 in case $(i)$ read exactly like the conclusions in the $k=2$ case described before, except that the exponent of the power of $1/ 2$ in the coefficient of $p$ in the asymptotic approximation is now $2s+q$ rather than $s+q$.   

In case $(ii)$, with $a_2b_1-a_1b_2=qb_1b_2$ and $a_3b_2-a_2b_3=rb_2b_3$, say, the block in the quotient diagram is formed by an overlap between rows 1 and 2 and an overlap between rows 2 and 3, but no overlap between rows 1 and 3. Hence the diagram looks like 
\vspace{0.5cm}
\begin{center}
\begin{tabular}{cccccccccccccccccc}
$\bigcdot$\ \ &$\bigcdot$\ \ &$\bigcdot$\ \ &$\bigcdot$\  \ & $\bigcdot$\ \ &$\bigcdot$\ \  &$\bigcdot$\ \ &$\bigcdot$&& &&&&&&&&\\
$\leftarrow$&$q$&$\rightarrow$& $\bigcdot$  &$\bigcdot$\  &$\bigcdot$\  & $\bigcdot$&$\bigcdot$&$\bigcdot$&\  $\bigcdot$&\  $\bigcdot$\  &\  \\
&&&$\leftarrow$&&$r$&&$\rightarrow$ & $\bigcdot$\  &\  \ $\bigcdot$\  &\  \ $\bigcdot$\ \ &\  $ \bigcdot$\   &\ \ $\bigcdot$\ \ \ &$\bigcdot$\ \ \ &$\bigcdot$\ \ \ &$\bigcdot$\  &\  \\
\end{tabular}
\end{center}

\vspace{0.2cm} 
\begin{center}
F\textsc{igure} 4. A quotient diagram for $k=3$
\end{center}

\vspace{0.5cm}
\noindent Here $\alpha=3s$, and, because of (29) and (31), $e=2s-q-r$ and the signature of $p$ is $\{\chi_p(b_1b_2), \chi_p(b_2b_3)\}$. We hence conclude from Theorem 9.9 that if $b_1b_2$ and $b_2b_3$ are both squares then
\begin{equation*}
q_\varepsilon(p)\sim (b\cdot 2^{s+q+r})^{-1}p\ \textrm{as}\ p\rightarrow +\infty.\tag{32}
\end{equation*}  
On the other hand, if either $b_1b_2$ or $b_2b_3$ is not a square then $\Pi_+(\textbf{a}, \textbf{b}) $ consists of all allowable primes $p$ such that $\{b_1, b_2, b_3\}$  is either a set of residues of $p$ or a set of non-residues of $p$, $\Pi_-(\textbf{a}, \textbf{b})$ consists of all allowable primes $p$ such that  $\{b_1, b_2, b_3\}$ contains a residue of $p$ and a non-residue of $p$,
\begin{equation*}
q_\varepsilon(p)=0,\  \textrm{for all}\ p\in \Pi_-(\textbf{a}, \textbf{b}),\ \textrm{and}\tag{33}
\end{equation*} 
\begin{equation*}
q_\varepsilon(p)\sim (b\cdot 2^{s+q+r})^{-1}p\ \textrm{as}\ p\rightarrow +\infty\ \textrm{inside}\ \Pi_+(\textbf{a}, \textbf{b}).\tag{34}
\end{equation*} 

In case $(iii)$, with the quotients $q$ and $r$ determined as in case $(ii)$, and, in addition, $a_3b_1-a_1b_3=tb_1b_3$, say, the block in the quotient diagram is now formed by an overlap between each pair of rows, and so the diagram looks like
\vspace{0.5cm}
\begin{center}
\begin{tabular}{cccccccccccccccccc}
$\bigcdot$&$\bigcdot$&$\bigcdot$&$\bigcdot$\ &\ $\bigcdot$\ &\ $\bigcdot$\ &\ $\bigcdot$\ &\ $\bigcdot$&& &&&&&&&&\\
$\leftarrow$&$q$&$\rightarrow$&$\bigcdot$\ &\ $\bigcdot$\ &\ $\bigcdot$\ &\ $\bigcdot$\ &\ $\bigcdot$\ &\ $\bigcdot$\ &\ $\bigcdot$&$\bigcdot$&  \\
&&&$\leftarrow$&$r$&$\rightarrow$\ &\ $\bigcdot$\ &\ $\bigcdot$\ &\ $\bigcdot$\ &\ \ $\bigcdot$\ \ &\  $\bigcdot$\ \ &\ $\bigcdot$\ \ &\ $\bigcdot$\ \ &\ $\bigcdot$\  &\  \\
\end{tabular}
\end{center}  

\vspace{0.2cm} 
\begin{center}
F\textsc{igure} 5. Another quotient diagram for $k=3$
\end{center}

\vspace{0.5cm}

\noindent It follows that $\alpha= 3s, e=2s-q-r$, and the signature of $p$ is $\{\chi_p(b_1b_2), \chi_p(b_1b_3), \chi_p(b_2b_3)\}$. In this case, the asymptotic approximation (32) holds whenever $b_1b_2, b_1b_3$, and $b_2b_3$ are all squares, and when at least one of these integers is not a square, $\Pi_+(\textbf{a}, \textbf{b})$ and $\Pi_-(\textbf{a}, \textbf{b})$ are determined by $\{b_1, b_2, b_3\}$ as before and (33) and (34) are valid.

\emph{Minimal overlap}. Here we take the quotient diagram to consist of a single block with gap sequence $(s-1, s-1,\dots,s-1)$, so that the overlap between rows is as small as possible: a typical quotient diagram for $s=4$ and $k=5$ looks like
\vspace{0.5cm}
\begin{center}
\begin{tabular}{cccccccccccccccc}
$\bigcdot$&$\bigcdot$&$\bigcdot$&$\bigcdot$&&&&&&&&&&&&\\
&&&$\bigcdot$&$\bigcdot$&$\bigcdot$&$\bigcdot$&&&&&&&&&\\
&&&&&&$\bigcdot$&$\bigcdot$&$\bigcdot$&$\bigcdot$&&&&&&\\
&&&&&&&&&$\bigcdot$&$\bigcdot$&$\bigcdot$&$\bigcdot$&\\
&&&&&&&&&&&&$\bigcdot$&$\bigcdot$&$\bigcdot$&$\bigcdot$\ \   \\
\end{tabular}
\end{center}

\vspace{0.2cm} 
\begin{center}
F\textsc{igure} 6. A quotient diagram with minimal overlap
\end{center}

\vspace{0.5cm}
\noindent Here $\alpha=ks$, $e=k-1$, and the signature of $p$ is $\{\chi_p(b_ib_{i+1}): i\in [1, k-1]\}$. Hence via Theorem 9.9 , if $b_ib_{i+1}, i\in [1, k-1]$, are all squares then
\begin{equation*}
q_\varepsilon(p)\sim (b\cdot 2^{1+k(s-1)})^{-1}p\ \textrm{as}\ p\rightarrow +\infty,
\end{equation*}  
and if at least one of those products is not a square, then $\Pi_+(\textbf{a}, \textbf{b}) $ consists of all allowable primes $p$ such that $\{b_1, \dots, b_k\}$  is either a set of residues of $p$ or a set of non-residues of $p$, $\Pi_-(\textbf{a}, \textbf{b})$ consists of all allowable primes $p$ such that  $\{b_1, \dots, b_k\}$ contains a residue of $p$ and a non-residue of $p$,
\begin{equation*}
q_\varepsilon(p)=0,\  \textrm{for all}\ p\in \Pi_-(\textbf{a}, \textbf{b}),\ \textrm{and}\tag{35}
\end{equation*} 
\begin{equation*}
q_\varepsilon(p)\sim (b\cdot 2^{1+k(s-1)})^{-1}p\ \textrm{as}\ p\rightarrow +\infty\ \textrm{inside}\ \Pi_+(\textbf{a}, \textbf{b}).
\end{equation*} 

\emph{Maximal overlap} ($k\geq 3$). Here we take the quotient diagram to consist of a single block with gap sequence $(1,1,\dots,1)$ and $k=s$, so that the overlap between each pair of rows is as large as possible: the diagrams for $k=3, 4,$ and 5 look like
\vspace{0.5cm}
\begin{center}
\begin{tabular}{ccccccccccccccccccccccc}
$\bigcdot$&$\bigcdot$&$\bigcdot$&&&&$\bigcdot$&$\bigcdot$&$\bigcdot$&$\bigcdot$&&&&&$\bigcdot$&$\bigcdot$&$\bigcdot$&$\bigcdot$&$\bigcdot$&&&&\\
&$\bigcdot$&$\bigcdot$&$\bigcdot$&&&&$\bigcdot$&$\bigcdot$&$\bigcdot$&$\bigcdot$&&&&&$\bigcdot$&$\bigcdot$&$\bigcdot$&$\bigcdot$&$\bigcdot$&&&\\
&&$\bigcdot$&$\bigcdot$&$\bigcdot$&&&&$\bigcdot$&$\bigcdot$&$\bigcdot$&$\bigcdot$&&&&&$\bigcdot$&$\bigcdot$&$\bigcdot$&$\bigcdot$&$\bigcdot$&&\\
&&&&&&&&&$\bigcdot$&$\bigcdot$&$\bigcdot$&$\bigcdot$&&&&&$\bigcdot$&$\bigcdot$&$\bigcdot$&$\bigcdot$&$\bigcdot$ \\
&&&&&&&&&&&&&&&&&&$\bigcdot$&$\bigcdot$&$\bigcdot$&$\bigcdot$&$\bigcdot$\ \  \\
\end{tabular}
\end{center}

\vspace{0.2cm} 
\begin{center}
F\textsc{igure} 7. Quotient diagrams with maximal overlap
\end{center}

\vspace{0.5cm}
\noindent We have in this case that $\alpha=k^2$, $e=(k-1)^2$, and the signature of $p$ is
\[
\Big\{\chi_p\Big(\prod_{i\in I} b_i\Big): I\in \mathcal{E}([1, k])\ \Big\}.
\]
 Hence if $\prod_{i\in I} b_i$ is a square for all $I\in \mathcal{E}([1, k])$ then
\begin{equation*}
q_\varepsilon(p)\sim (b\cdot 2^{2k-1})^{-1}p\ \textrm{as}\ p\rightarrow +\infty,
\end{equation*}  
and if one of these products is not a square then $\Pi_+(\textbf{a}, \textbf{b})$ and $ \Pi_-(\textbf{a}, \textbf{b})$ are determined by  $\{b_1, \dots, b_k\}$ as before, (35) holds, and 
\begin{equation*}
q_\varepsilon(p)\sim (b\cdot 2^{2k-1})^{-1}p\ \textrm{as}\ p\rightarrow +\infty\ \textrm{inside}\ \Pi_+(\textbf{a}, \textbf{b}).
\end{equation*}

It follows from our discussion after the proof of Theorem 9.9 that an increase in the number of overlaps between rows in the quotient diagram of $(\textbf{a}, \textbf{b})$ leads to an increase in the asymptotic number of elements of $AP(\textbf{a}, \textbf{b}; s)\cap 2^{[1, p-1]}$ that are sets of residues or non-residues of $p$, and these examples now verify that principle quantitatively.  In order to make this explicit, note first that Lemma 9.10 can be used to generate examples in which the $(k-1)$-tuple $(d_1,\dots,d_{k-1})$ varies arbitrarily, while at the same time $b=\max\{b_1,\dots,b_k\}$ always takes the same value. Hence we may assume in the discussion to follow that the value of $b$ is constant in each set of examples, and so the only parameter that is relevant when comparing asymptotic approximations to $q_{\varepsilon}(p)$ is the exponent of the power of $1/2$ in the coefficient of that approximation. When $k=2$, there is either no overlap between rows or exactly 1 overlap; in the former case, the exponent in the power of $1/ 2$ that occurs in the asymptotic approximation to $q_\varepsilon(p)$ is $2s$ and in the latter case this exponent is less than $2s$. When $k=3$ there are 0, 1, 2, or 3 possible overlaps between rows, with the last three possibilities occurring, respectively, in cases ($i$), ($ii)$, and $(iii)$ above. It follows that $q<s$ in case  $(i)$, $q+r\geq s$ in case $(ii)$ and $q+r<s$ in case $(iii)$. Hence the exponent in the power of $1/ 2$ that occurs in the asymptotic approximation to $q_\varepsilon(p)$ is $3s$ when no overlap occurs, is greater than $2s$ and less than $3s$ in case $(i)$, is at least $2s$ and less than $3s$ in case ($ii)$, and is less than $2s$ in case ($iii)$. If we also take $k=s$ when there is minimal overlap in the quotient diagram and compare that to what happens when there is maximal overlap there, we see that the exponent in the power of  $1/ 2$ that occurs in the asymptotic approximation of $q_\varepsilon(p)$ is quadratic in $k$, i.e., $ k^2-k+1$, in the former case, but only linear in $k$, i.e., $ 2k-1$, in the latter case.  

\section{The Asymptotic Density of $\Pi_+(\textbf{a}, \textbf{b})$}
Suppose that $(\textbf{a}, \textbf{b})$ is a standard $2k$-tuple and assume that there exists an $I \in \Lambda(\mathcal{K})$ such that $\prod_{i \in I} b_i$ is not a square. Then, in accordance with Theorem 9.9, the sets $\Pi_+(\textbf{a}, \textbf{b})$ and $\Pi_-(\textbf{a}, \textbf{b})$  are both infinite, and so it is of interest to calculate their asymptotic density. Because  $\Pi_+(\textbf{a}, \textbf{b})$ and $\Pi_-(\textbf{a}, \textbf{b})$ are disjoint sets with only finitely many primes outside of their union, it follows that
\[
\textrm{the density of}\  \Pi_+(\textbf{a}, \textbf{b}) +\ \textrm{the density of}\ \Pi_-(\textbf{a}, \textbf{b})=1,\]
so it suffices to calculate only the density of $\Pi_+(\textbf{a}, \textbf{b})$.

In order to keep the technicalities from becoming too complicated, we will describe this calculation for the following special case: assume that 
\begin{quote}
(36)  $(\textbf{a}, \textbf{b})$ is admissible, the square-free parts $\sigma_i=\sigma(b_i)$ of the coordinates $b_i$ of $\textbf{b}$ are distinct and for each nonempty subset of $T$ of $[1, k], \prod_{i \in T} \sigma_i$ is not a square.
\end{quote}
This condition is satisfied, for example, if 
\begin{equation*}
\textrm{$b_i$ is square-free for all $i$ and $\pi(b_i)$ is a proper subset of $\pi(b_{i+1})$, for all $i \in [1, k-1]$}. \tag{37}\] 
Moreover for each $k \in [2, \infty)$, Lemma 9.10 can be used to construct infinitely many admissible $2k$-tuples with a fixed but arbitrary quotient diagram which satisfy (37).

Let $(\mathcal{D}(V_1),\dots,\mathcal{D}(V_m))$ be the quotient diagram of $(\textbf{a}, \textbf{b})$ and let $D_i$ be the subset of $[1, k]$ such that $V_i=\{q_j: j \in D_i\}, i \in [1, m]$;
as the sets $V_1,\dots,V_m$  are pairwise disjoint, so also are the sets $D_1,\dots,D_m$ .

Now, let $\mathcal{C}_i$ denote the set of columns of the overlap diagram $\mathcal{D}(V_i)$, realized as  subsets of $[1, k] \times [0, s-1]$ as per the identification given by (30), and let
\[
\Lambda_i(\mathcal{K})=\bigcup _{C \in \mathcal{C}_i} \mathcal{E}(\theta(C)).\]
Then
\begin{equation*}
\bigcup_{I \in \Lambda_i(\mathcal{K})} I=\bigcup_{C \in \mathcal{C}_i} \theta(C)=D_i,\  i \in [1, m], \tag{38}
\]
and so it follows from the pairwise disjointness of the $D_i$'s, these equations, and (31) that 
\begin{equation*}
\Lambda(\mathcal{K})=\bigcup_i   \Lambda_i(\mathcal{K}),\ \textrm{and this union is pairwise disjoint}. \tag{39}\]
 
Next, for each $I \in \Lambda(\mathcal{K})$ let
\[
S(I)=\{\sigma_i: i \in I\},\]
and then set
\[
\mathcal{M}_1=\{I \in \Lambda(\mathcal{K}): 1 \in S(I)\}.\]
If $\mathcal{M}_1\not= \emptyset$ then there is a unique element $n_0$ of $\bigcup_i D_i$ such that $\sigma_{n_0}=1$, hence it follows from (38) and (39) that there is a unique element $i_0$ of $[1, m]$ such that
\begin{equation*}
\mathcal{M}_1=\{I \in  \Lambda_{i_0}(\mathcal{K}): n_0 \in I\}.\]
It can then be shown that if
\[
\sigma=\sum_i|D_i|\ \]
and
\[
m=\textrm{ the number of blocks in the quotient diagram of  $(\textbf{a}, \textbf{b})$},\]
then the density of $\Pi_+(\textbf{a}, \textbf{b})$ is
\begin{equation*}
2^{m- \sigma},\  \textrm{if $\mathcal{M}_1= \emptyset$ or  $\mathcal{M}_1=\Lambda_{i_0}(\mathcal{K}) $,}\ \textrm{or} \tag{40}
\]
\begin{equation*}
2^{1- \sigma}(2^m-1),\   \textrm{if  $\emptyset \not= \mathcal{M}_1\not=\Lambda_{i_0}(\mathcal{K})$.} \tag{41}
\]

It follows that whenever $(\textbf{a}, \textbf{b})$ is an admissible $2k$-tuple for which the square-free parts of the coordinates of $\textbf{b}$ are distinct and satisfy condition (36), the cardinality of $\bigcup_i D_i$, the number of blocks $m$ in the quotient diagram, and the set $\mathcal{M}_1$ completely determine the density of $\Pi_+(\textbf{a}, \textbf{b})$ by means of formulae (40) and (41). Those formulae show that each element of $\bigcup_i D_i$ contributes a factor of 1/2 to the density of $\Pi_+(\textbf{a}, \textbf{b})$ and each block of the quotient diagram of  $(\textbf{a}, \textbf{b})$ contributes essentially a factor of 2 to the density. Because $|V_i| \geq 2$ for all $i$, it follows that $|D_i| \geq 2$ for all $i$ and so $\sigma \geq 2m$; in particular, the density of $\Pi_+(\textbf{a}, \textbf{b})$ is at most $2^{-m}$ whenever  $\mathcal{M}_1= \emptyset$ or $\mathcal{M}_1=\Lambda_{i_0}(\mathcal{K}) $ and is at most $(2^m-1)/2^{2m-1}$, otherwise. This gives an interesting number-theoretic interpretation to the number of blocks in the quotient diagram. In fact, if for each $k \in [2, \infty)$, we let $\mathcal{A}_k$ denote the set of all admissible $2k$-tuples which satisfy condition (36), set $\mathcal{A}=\bigcup_{k \in [2, \infty)} \mathcal{A}_k$, and take  $m \in [1, \infty)$, then Lemma 9.10 can be used to show that there exists infinitely many elements $(\textbf{a}, \textbf{b})$ of $\mathcal{A}$ such that the quotient diagram of  $(\textbf{a}, \textbf{b})$ has $m$ blocks and the density of $\Pi_+(\textbf{a}, \textbf{b})$ is $2^{-m}$ (respectively, $(2^m-1)/2^{2m-1})$. One can also show that if $\{l, n\} \subseteq [1, \infty),$ with $l \geq 2n,$ then there are infinitely many elements $(\textbf{a}, \textbf{b})$ of $\mathcal{A}$ such that  the density of $\Pi_+(\textbf{a}, \textbf{b})$ is $2^{1-l}(2^n-1)$.

For more details in this situation and for what transpires for arbitrary standard $2m$-tuples, we refer the interested reader to Wright [63].
\newpage
\chapter{Are quadratic residues randomly distributed?}
The purpose of this chapter is to provide evidence that the answer to the question in the title is yes. By examining tables of residues and non-residues of certain primes in section 1, we observe that residues can occur in very  irregular patterns. In section 2, we will show how to view sums of the values of Legendre symbols $\chi_p$ as random variables and then we will employ the Central Limit Theorem from probability theory to determine a condition under which, at least when $p$ is sufficiently large
, the values of $\chi_p$ can be interpreted to behave randomly and independently. In section 3, a very interesting result of Davenport and Erd$\ddot{\textrm{o}}$s on the the distribution of residues will then be employed to verify that the condition from section 2 that detects random behavior of residues and non-residues does indeed hold. Interestingly enough, the Weil-sum estimates from Theorem 9.1, which were so useful in our work in Chapter 9, will also be very useful in our proof of Davenport and Erd$\ddot{\textrm{o}}$s' result.

\section{Irregularity of the Distribution of Quadratic Residues}
Extensive numerical calculations performed over the years indicate that, at least in certain subintervals of $[1, p-1]$, residues and non-residues of $p$ occur in very irregular patterns. For example, we present below four tables which exhibit the residues and non-residues of the primes 41, 79, 101, and 139. A 0 indicates that the corresponding entry is a residue of the indicated prime and 1 indicates that the corresponding entry is a non-residue. The tables for 41 and 101 are palindromic, i.e., they read the same from left to right, starting from the first entry in the table, as from right to left, starting from the last entry (note that 41 and 101 are congruent to 1 mod 4), and the tables for 79 and 139 become palindromic if the 0 and 1 entries in the last half of the tables are switched to 1 and 0, respectively (note that 79 and 139 are congruent to 3 mod 4). However, the entries in various subintervals of consecutive integers in the first half of the tables are fairly irregular and do not appear to exhibit any predictable pattern.
 
\newpage
\vspace{0.5cm}

\hspace{3.7cm} 1-20: 0  0  1  0  0  1  1  0  0  0  1  1  1  1  1  0  1   0  1  0

\hspace{3.5cm} 21-40: 0  1  0  1  0  1  1  1  1  1  0  0  0  1  1  0  0  1  0  0

 \vspace{0.3cm}
\begin{center}
Table 1: Residues and Non-residues of 41
\end{center}

 \vspace{0.5cm}

\hspace{3.7cm} 1-20:  0  0  1  0  0  1  1  0  0  0  0  1  0  1  1  0  1  0  0  0

 \hspace{3.5cm} 21-40: 0  0  0  1  0  0  1  1  1  1  0  0  1  1  1  0  1  0  1  0
 
 \hspace{3.5cm} 41-60: 1  0  1  0  0  0  1  1  0  0  0  0  1  1  0  1  1  1  1  1

 \hspace{3.5cm} 61-78: 1  0  1  0  0  1  0  1  1  1  1  0  0  1  1  0  1  1

 \vspace{0.3cm}
\begin{center}
Table 2: Residues and Non-residues of 79
 \end{center}

 \vspace{0.5cm}
\hspace{3.7cm} 1-20: 0  1  1  0  0  0  1  1  0  1  1  1  0  0  1  0  0  1  0  0

 \hspace{3.5cm} 21-40: 0  0  0  0  0  1  1  1  1  0  0  1  0  1  1  0  0  1  1  1  
 
 \hspace{3.5cm} 41-60: 1  1  0  1  0  1  0  1  0  1  1  0  1  0  1  0  1  0  1  1   

 \hspace{3.5cm} 61-80: 1  1  1  0  0  1  1  0  1  0  0  1  1  1  1  0  0  0  0  0

 \hspace{3.28cm} 81-100: 0  0  1  0  0  1  0  0  1  1  1  0  1  1  0  0  0  1  1  0 

 \vspace{0.3cm}
\begin{center}
Table 3: Residues and Non-residues of 101
 \end{center}

 \vspace{0.5cm}

\hspace{3.7cm} 1-20: 0  1  1  0  0  0  0  1  0  1  0  1  0  1  1  0  1  1  1  0

 \hspace{3.5cm} 21-40: 1  1  1  0  0  1  1  0  0  0  0  1  1  0  0  0  0  0  1  1
 
 \hspace{3.5cm} 41-60: 0  0  1  0  0  0  0  1  0  1  0  0  1  0  0  1  0  1  1  1

 \hspace{3.5cm} 61-80: 1  1  0  0  0  0  0  1  0  1  0  1  1  1  1  1  0  0  0  0

 \hspace{3.28cm} 81-100: 0  1  0  1  1  0  1  1  0  1  0  1  1  1  1  0  1  1  0  0

 \hspace{3.1cm} 101-120: 1  1  1  1  1  0  0  1  1  1  1  0  0  1  1  0  0  0  1  0

 \hspace{3.1cm} 121-138: 0 0  1  0  0  1  0  1  0  1  0  1  1  1  1  0  0  1

 \vspace{0.3cm}
\begin{center}
Table 4: Residues and Non-residues of 139
 \end{center}

\vspace{0,5cm}
\noindent This has led to speculation about whether residues occur more or less randomly in certain intervals of consecutive integers. In the following section, we set up a procedure which can be used to provide positive evidence for the contention that residues are in fact distributed in this manner.

\section{Detecting Random Behavior Using the Central Limit Theorem}

The method which we will use to detect random behavior in the distribution of residues employs the Central Limit Theorem from the mathematical theory of probability. In order to set the stage for that result, we will briefly review some basic facts and terminology from probability theory. 

One starts with a \emph{probability space}, i.e., a triple $(\Omega, \mathcal{M}, \mu)$ consisting of a set $\Omega$ (the sample space), a distinguished $\sigma$-algebra $\mathcal{M}$ of subsets of $\Omega$ (the events), and a non-negative, countably additive measure $\mu$ defined on $\mathcal{M}$ such that $\mu(\Omega)=1$ (the probability measure). A \emph{random variable X on} $\Omega$ is an extended real-valued function defined on $\Omega$ such that for each real number $r$, the set $\{\omega\in \Omega: X(\omega)<r\}$ is in $\mathcal{M}$, i.e., $X$ is \emph{measurable with respect to} $\mathcal{M}$. The \emph{mean} and \emph{variance} of a random variable $X$ is the value of the integral $\overline{X}=\displaystyle{\int_{\Omega} X\ d\mu}$ and $\displaystyle{\int_{\Omega} (X-\overline{X})^2\ d\mu}$, respectively. The \emph{distribution function of X} is the function defined on the real line $\mathbb{R}$ by
\[
\lambda \rightarrow \mu\big(\{\omega\in \Omega: X(\omega)\leq \lambda\}\big),\ \lambda \in \mathbb{R}.\]
It can be shown that a non-negative function $F$ defined on $\mathbb{R}$ is the distribution function of a random variable if and only if $F$ is non-decreasing, right continuous, $\lim_{\lambda \rightarrow -\infty} F(\lambda)=0$ and $\lim_{\lambda \rightarrow +\infty} F(\lambda)=1$ (Chung [4], Theorem 2.2.4). 
 
A set $\{X_1,\dots,X_n\}$ of random variables on the probability space $(\Omega, \mathcal{M}, \mu)$ is (stochastically) \emph{independent} if for any $n$-tuple $(B_1,\dots,B_n)$ of Borel subsets of the real line, we have that
\[
\mu\Big(\bigcap_{i=1}^n\{\omega\in \Omega: X_i(\omega)\in B_i\}\Big)=\prod_{i=1}^n \mu\big(\{\omega\in \Omega: X_i(\omega)\in B_i\}\big).\]
An infinite sequence $(X_n)$ of random variables is independent if every finite subset of the $X_i$'s is independent.

Stochastic independence is a way of making mathematically precise the intuitive notion of describing how events are determined by the \emph{outcomes} of random \emph{trials}. An unbiased coin is tossed, the two possible outcomes are recorded as 0 or 1, with roughly probabilities of 1/2 each, and repeated tossings generate a sequence of outcomes. In a similar way, we may repeatedly cast a fair die, or draw colored beads (with replacement) from an urn, or take repeated measurements of a certain quantity from a sample population, each process generating values for a sequence of random variables. As Chung [4] states, ``it is very easy to conceive of undertaking these various trials under conditions such that their respective outcomes do not appreciably effect each other.... In this circumstance, idealized trials are carried out ``independently of one another" and the corresponding [random variables] are ``independent" according to definition". In fact, according to Chung, ``it may be said that no one could have learned the subject of probability properly without acquiring some feeling for the intuitive content of the concept of stochastic independence."

Sequences of independent random variables exist in abundance. One can in fact prove that if, for each positive integer $n$, $Y_n$ is a fixed but arbitrary random variable defined on a fixed but arbitrary probability space $(\Omega_n, \mathcal{M}_n, \mu_n)$, then there exists a probability space $(\Omega, \mathcal{M}, \mu)$ and a sequence of independent random variables $(X_n)$ on  $(\Omega, \mathcal{M}, \mu)$ such that, for all $n$, $X_n$ and $Y_n$ have the same mean, variance, and distribution function. The construction of $(\Omega, \mathcal{M}, \mu)$ and the sequence of independent random variable $(X_n)$ uses the measure-theoretic infinite product of the probability spaces  $(\Omega_n, \mathcal{M}_n, \mu_n)$ and is hence too involved to go into further here; for complete details of this construction, we refer the reader interested in them to Chung [4], Theorem 3.3.4 and its proof.

Now, suppose that $X_1, X_2,\dots$ is a sequence of random variables defined on a probability space ($\Omega, \mathcal{M}, \mu)$ which is independent, identically distributed, i.e., all $X_n$'s have the same distribution function, and each random variable has mean 0 and variance 1. If we set
\[
S_n=\sum_{k=1}^n X_k,\ n \in [1, \infty),\]
then the Central Limit Theorem (Chung [4], Theorem 6.4.4) asserts that for each real number $\lambda$,
\begin{equation*}
\lim_{n \rightarrow +\infty}\mu\Big(\Big\{\omega \in \Omega: \frac{S_n(\omega)}{\sqrt n} \leq \lambda\Big\}\Big)=\frac{1}{\sqrt{2\pi}} \int_{-\infty}^{\lambda} e^{-t^2/2} dt, \tag{1} 
\]
i.e., as $n \rightarrow +\infty, S_n/\sqrt n$ tends to becomes normally distributed with mean 0 and variance 1.

Now let $p$ be a prime. We convert the set $[0, p-1]$ into a (discrete and finite) probability space by assigning probability $1/p$ to each element of $[0, p-1]$. This induces the probability measure $\mu_p$ on $[0, p-1]$ defined by
\begin{equation*}
\mu_p(S)=\frac{|S|}{p},\ S \subseteq [0, p-1]. \tag{2}\]
For each positive integer $h<p$, consider the sums
\[
S_h(x)=\sum_{n=x+1}^{x+h} \chi_p(n),\ x=0,\dots,p-1,\]
which is just the quadratic excess of the interval $(x, x+h+1)$ that we studied in Chapter 7. The function $S_h$ is a random variable on $([0, p-1], \mu_p)$, and so by way of analogy with (1), we consider the distribution function
\begin{equation*}
\lambda \rightarrow \mu_p \Big(\Big\{x \in [0, p-1]: \frac{S_h(x)}{\sqrt h} \leq \lambda\Big\}\Big),\ \lambda \in \mathbb{R}, \tag{3}\]
of $S_h/\sqrt h$.

We next let $h=h(p)$ be a function of $p$ and look for conditions on the growth of $h(p)$ which guarantee that for each real number $\lambda$,
\begin{equation*}
\lim_{p \rightarrow +\infty} \frac{1}{p}\ \Big|\Big\{x \in [0,p-1]: \frac{S_{h(p)}(x)}{\sqrt {h(p)}} \leq \lambda\Big\}\Big|=\frac{1}{\sqrt{2\pi}} \int_{-\infty}^{\lambda} e^{-t^2/2} dt, \tag{4} \]
It is easy to see that a necessary condition for (4) to occur is that $\lim_{p \rightarrow +\infty} h(p)=+\infty$. If (4) is valid then, as we see from (2) and (3), when $p \rightarrow +\infty$ the sums $S_{h(p)}$ satisfy a ``central limit theorem" relative to the probability spaces $([0, p-1], \mu_p)$. If (4) can be verified, then upon comparing it to (1), we conclude that for $p$ sufficiently large, at least with respect to sampling using $\chi_p$ in the intervals $[x+1, x+h(p)],\ x=0,1,\dots,p-1$, residues and non-residues of $p$ appear to behave as if they are distributed randomly and independently!

\section{Verifying Random Behavior via a Result of  Davenport and Erd$\ddot{\textrm{o}}$s}
The following theorem of Davenport and Erd$\ddot{\textrm{o}}$s  ([7], Theorem 5) provides conditions on $h(p)$ which imply that (4) is true:
\begin{thm}
If $h:P \rightarrow [1, \infty)$ is any function such that
\[
\lim_{q \rightarrow +\infty} h(q)=+\infty,\ \lim_{q \rightarrow +\infty} \frac{h(q)^r}{\sqrt q}=0,\ \textrm{for all}\ r \in [1, \infty)\]
$($e.g., $h(q)=[\log^N q]$, where $N$ is any fixed positive integer$)$, then for each real number $\lambda$,
\[
\lim_{p \rightarrow +\infty} \frac{1}{p}\ \Big|\Big\{x \in [0,p-1]: \frac{S_{h(p)}(x)}{\sqrt {h(p)}} \leq \lambda\Big\}\Big|=\frac{1}{\sqrt{2\pi}} \int_{-\infty}^{\lambda} e^{-t^2/2} dt.\]
\end{thm}
\noindent As a consequence of this theorem and our discussion in section 2, we conclude that, at least for sufficiently large primes $p$, residues and non-residues do appear to be distributed randomly within intervals whose length does not increase too fast as $p\rightarrow +\infty$.

The proof of Theorem 10.1 relies on the following lemma: we will first state the lemma, use it to prove Theorem 10.1, and then prove the lemma.
\begin{lem}
Let $r$ be a fixed positive integer, and let $h$ be an integer and $p$ a prime such that $r<h<p$. Then there exists numbers $0 \leq \theta \leq 1, 0 \leq \theta^{\prime} \leq 1$ such that 
\begin{equation*}
\Big|\sum_{x=0}^{p-1} S_h(x)^{2r}\ -(p-\theta r)(h-\theta^{\prime}r)^r \prod_{i=1}^r (2i-1) \Big| \leq 2rh^{2r} \sqrt p, \tag{5}\]
\begin{equation*}
\Big|\sum_{x=0}^{p-1} S_h(x)^{2r-1} \Big| \leq 2rh^{2r} \sqrt p. \tag{6}\]
\end{lem}

\emph{Proof of Theorem} 10.1. Let $r$ be a fixed positive integer. Then by the hypotheses satisfied by $h(p)$, we have that $r<h(p)<p$ for all $p$ sufficiently large, hence Lemma 10.2 implies that for all such $p$,
\[
\Big|\frac{1}{p}\ \sum_{x=0}^{p-1} (h(p)^{-1/2}S_{h(p)}(x))^{2r}\ -\Big(1-\frac{\theta r}{p}\Big) \Big(1-\frac{\theta^{\prime}r}{h(p)} \Big)^r  \prod_{i=1}^r (2i-1) \Big| \leq 2r\frac{h(p)^r}{\sqrt p},\]
\[
\Big|\frac{1}{p}\ \sum_{x=0}^{p-1} (h(p)^{-1/2}S_{h(p)}(x))^{2r-1} \Big| \leq 2r\frac{h(p)^r}{\sqrt p}.\]
Letting $p \rightarrow +\infty$ in these inequalities, we deduce from the growth conditions on $h(p)$ that if $r$ is any positive integer and
\[
\mu_r=\left\{\begin{array}{rl}\displaystyle  \prod_{i=1}^{r/2} (2i-1),& \textrm{if $r$ is even,}\\
0,& \textrm{if $r$ is odd,} 
\\\end{array}\right. 
\]
then
\begin{equation*}
\lim_{p \rightarrow +\infty} \frac{1}{p}\ \sum_{x=0}^{p-1} (h(p)^{-1/2}S_{h(p)}(x))^r=\mu_r. \tag{7} \]

Now for each real number $s$, let
\[
N_p(s)=\frac{1}{p}\ \big|\big\{x \in [0, p-1]: S_{h(p)}(x) \leq s \big\}\big|.\]
The function $N_p$ is nondecreasing in $s$, constant except for possible discontinuities at certain integral values of $s$, and is right-continuous at every value of $s$. Because
\[
\big|S_{h(p)}(x)\big| \leq h(p),\ \textrm{for all}\ x,\]
it follows that
\[
N_p(s)=\left\{\begin{array}{rl} 0,& \textrm{if $s<-h(p)$ ,}\\
1,& \textrm{if $s \geq h(p)$.} \\\end{array}\right. \]
We also have that
\begin{eqnarray*}
(8)\hspace{2cm} \frac{1}{p} \sum_x \big(h(p)^{-1/2}S_{h(p)}(x)\big)^r&=&\frac{1}{p} \sum_{s=-h(p)}^{h(p)} \Big(\sum_{x:S_{h(p)}(x)=s} (h(p)^{-1/2}s)^r\Big)\\
&=&\frac{1}{p} \sum_{s=-h(p)}^{h(p)} (h(p)^{-1/2}s)^r|\{x: S_{h(p)}(x)=s\}|\\
&=& \sum_{s=-h(p)}^{h(p)} (h(p)^{-1/2}s)^r(N_p(s)-N_p(s-1)),
\end{eqnarray*}
and so if we let \[
\Phi_p(t)=N_p(th(p)^{-1/2}),\]
then the last sum in (8) can be written as the Stieltjes integral
\[
\int_{-\infty}^{\infty} t^rd\Phi_p(t).\]
Putting
\[
\Phi(t)=\frac{1}{\sqrt{2\pi}} \int_{-\infty}^t e^{-u^2/2} du,\]
we have
\[ 
\int_{-\infty}^{\infty} t^rd\Phi(t)=\frac{1}{\sqrt{2\pi}}\int_{-\infty}^{\infty} t^re^{-t^2/2} dt=\mu_r,\]
hence (7), (8) imply that
\begin{equation*}
\lim_{p \rightarrow +\infty} \int_{-\infty}^{\infty} t^rd\Phi_p(t)=\int_{-\infty}^{\infty} t^rd\Phi(t),\ \textrm{for all}\ r \in [0, \infty). \tag{9}\]

By virtue of the definition of $\Phi_p$, the conclusion of Theorem 10.1 can be stated as
\begin{equation*}
\lim_{p \rightarrow +\infty} \Phi_p(\lambda)=\Phi(\lambda),\ \textrm{for all real numbers}\ \lambda.\tag{10}\]
We will deduce (10) from (9) by an appeal to the classical theory of moments.

Suppose by way of contradiction that (10) is false for some $\lambda$; then there exists $\delta>0$ such that
\begin{equation*}
|\Phi_p(\lambda)-\Phi(\lambda)| \geq \delta\ \textrm{for infinitely many}\ p. \tag{11}\]
Using the first and second Helly selection theorems (Shohat and Tamarkin, [54], Introduction, section 3),we find a subsequence of these $p$, say $p^{\prime}$, and a nondeceasing real-valued function $\Phi^{*}$ defined on $\mathbb{R}$ such that
\begin{equation*}
\lim_{t \rightarrow -\infty} \Phi^{*}(t)=0,\ \lim_{t \rightarrow +\infty} \Phi^{*}(t)=1, \tag{12}\]
\begin{equation*}
\Phi^{*}\ \textrm{ is right-continuous at all points of}\ \mathbb{R}, \tag{13}\]
\begin{equation*}
\lim_{p^{\prime} \rightarrow +\infty} \Phi_{p^{\prime}}(t)= \Phi^{*}(t),\ \textrm{for all points $t$ at which $\Phi^{*}$ is continuous}, \tag{14}\]
and
\begin{equation*}
\lim_{p^{\prime} \rightarrow +\infty} \int_{-\infty}^{\infty} t^rd\Phi_{p^{\prime}}(t)=\int_{-\infty}^{\infty} t^rd\Phi^{*}(t),\ \textrm{for all}\ r \in [0, \infty). \tag{15}\]
By way of (9) and (15),
\begin{equation*}
\int_{-\infty}^{\infty} t^rd\Phi^{*}(t)=\int_{-\infty}^{\infty} t^rd\Phi(t),\ \textrm{for all}\ r \in [0, \infty). \tag{16}\]
The Weierstrass approximation theorem, which asserts that each function continuous on a closed and bounded interval of the real line is the uniform limit on that interval of a sequence of polynomials, and (16) imply
\begin{equation*}
\int_{-\infty}^{\infty} fd\Phi^{*}(t)=\int_{-\infty}^{\infty} fd\Phi(t), \tag{17}\]
for all real-valued functions $f$ continuous on  $\mathbb{R}$ of compact support. Equations (12), (13), and (17) imply that 
\begin{equation*}
\Phi^{*}(t)=\Phi(t),\ \textrm{for all}\ t \in \mathbb{R}. \tag{18}\]
Hence $\Phi^{*}$ is continuous everywhere in $\mathbb{R}$, and so by (14) and (18),
\[
\lim_{p^{\prime} \rightarrow +\infty} \Phi_{p^{\prime}}(\lambda)= \Phi(\lambda),\]
and this contradicts (11).

It remains to prove Lemma 10.2. The argument here makes use of another interesting application of the Weil-sum estimates available from Theorem 9.1.

Consider first the case with $2r$ as the exponent. We have that
\begin{equation*}
\sum_{x=0}^{p-1} (S_h(x))^{2r}= \sum_{(n_1,\dots,n_r) \in [1, h]^{2r}} \sum_{x=0}^{p-1} \chi_p\Big(\prod_{i=1}^{2r} (x+n_i)\Big). \tag{19}\]
In order to estimate the absolute value of this sum, we divide the elements $(n_1,\dots,n_{2r})$ of $[1, h]^{2r}$ into two types: $(n_1,\dots,n_{2r})$ is of \emph{type} 1 if it has at most $r$ distinct coordinates, each of which occurs an even number of times; all other elements of  $[1, h]^{2r}$ are of \emph{type} 2.

If $(n_1,\dots,n_{2r})$ is of type 1 then the polynomial $\prod_i (x+n_i)$ is a perfect square in $(\mathbb{Z}/p\mathbb{Z})[x]$. If $s$ is the number of distinct coordinates of $(n_1,\dots,n_{2r})$, then $\chi_p\Big(\prod_i (x+n_i)\Big)=0$ whenever there is a distinct coordinate $n_j$ of  $(n_1,\dots,n_{2r})$ such that $x \equiv -n_j$ mod $p$, and $\chi_p\big(\prod_i (x+n_i)\big)=1$ otherwise.
It follows that the value of the sum
\[
\sum_{x=0}^{p-1} \chi_p\Big(\prod_{i=1}^{2r} (x+n_i)\Big)\]
is at least $p-r$, and this value is clearly at most $p$. Hence there exists a number $0 \leq \theta \leq 1$ such that the sum (19) is
\[
F(h, r)(p-\theta r),\]
where $F(h, r)$ denotes the cardinality of the set of all elements of $[1, h]^{2r}$ of type 1. 

On the other hand, if  $(n_1,\dots,n_{2r})$ is of type 2 then the polynomial  $\prod_i (x+n_i)$ reduces modulo $p$ to a product of at least one and at most $2r$ distinct linear factors over $\mathbb{Z}/p\mathbb{Z}$, hence Theorem 9.1 implies that
\[
\Big|\sum_{x=0}^{p-1} \chi_p\Big(\prod_{i=1}^{2r} (x+n_i)\Big) \Big| \leq 2r\sqrt p.\]
Hence the contribution of the elements of type 2 to the sum (19) has an absolute value that does not exceed $2rh^{2r} \sqrt p.$

An appropriate estimate of the size of $F(h, r)$ is now required. Following Davenport and Erd$\ddot{\textrm{o}}$s, we  note first that the number of ways of choosing \emph{exactly} $r$ distinct integers from $[1, h]$ is $h(h-1)\cdots(h-r+1)$, and the number of ways of arranging these as $r$ pairs is $\prod_{i=1}^r (2i-1)$. Hence
\begin{eqnarray*}
F(h, r) &\geq& h(h-1)\dots(h-r+1) \prod_{i=1}^r (2i-1)\\
&>&(h-r)^r \prod_{i=1}^r (2i-1).
\end{eqnarray*}
On the other hand, the number of ways of choosing at most $r$ distinct elements from $[1, h]$ is at most $h^r$, and when these have been chosen, the number of different ways of arranging them in $2r$ places is at most $\prod_{i=1}^r (2i-1)$. Hence
\[
F(r, h) \leq h^r \prod_{i=1}^r (2i-1).\]
Hence there is a number $0 \leq \theta^{\prime} \leq 1$ such that
\[
F(r, h)=(h-\theta^{\prime} r)^r  \prod_{i=1}^r (2i-1).
\]
The conclusion of Lemma 10.2 for odd exponents follows from these estimates, and when the sum has an even exponent, the desired conclusion is now obvious, because in this case there are no elements of type 1. $\hspace{10.7cm} \textrm{QED}$

\emph{Remark}. More recently, Kurlberg and Rudnick [32] and Kurlberg [31] have provided further evidence of the random behavior of quadratic residues by computing the limiting distribution of normalized consecutive spacings between representatives of the squares in $\mathbb{Z}/n\mathbb{Z}$ as $|\pi(n)| \rightarrow +\infty$. In order to describe their work there, let $S_n \subseteq [0, n-1]$ denote the set of representatives of the squares in $\mathbb{Z}/n\mathbb{Z}$, i.e., the set of quadratic residues modulo $n$ inside $[0, n-1]$ (N.B. It is \emph{not} assumed here that a quadratic residue mod $n$ is relatively prime to $n$). Order the elements of $S_n$ as $r_1<\dots<r_N$ and then let $x_i=(r_{i+1}-r_i)/s$, where $s=(r_N-r_1)/N$ is the mean spacing; $x_i, i=1,\dots,N-1$, are the distances between consecutive elements of $S_n$ normalized to have mean distance 1. If $t$ is any fixed positive real number then it is shown in [31] and [32] that
\[
\lim_{|\pi(n)| \rightarrow +\infty}\frac{|\{x_i: x_i \leq t\}|}{|S_n|-1}=1-e^{-t},\]
i.e., for all $n$ with $|\pi(n)|$ large enough, the normalized spacings between quadratic residues of $n$ follow (approximately) a Poisson distribution. Among many other things, the Poisson distribution governs the number of customers and their arrival times in queueing theory, and so the results of Kurlburg and Rudnick can be interpreted to say that if the number of prime factors of $n$ is sufficiently large then quadratic residues of $n$ appear consecutively in the set $[0, n-1]$ in the same way as customers arriving randomly to join a queue. 
\clearpage
\fancyhead[LO]{ \fontsize{13}{12} \selectfont \nouppercase  Bibliography}
\fancyhead[LE]{ \fontsize{13}{12} \selectfont \nouppercase  Bibliography}

\newpage
\afterpage{\null\newpage}
\thispagestyle{empty}

\clearpage
\fancyhead[LO]{ \fontsize{13}{12} \selectfont \nouppercase  Index}
\fancyhead[LE]{ \fontsize{13}{12} \selectfont \nouppercase  Index}

\begin{theindex}
\item abelian extension, 66
\item admissible $2k$-tuple, 234
\item algebraic curve, 211
\subitem estimate of the number of rational points on a non-singular, 212
\subitem non-singular, 212
\subitem rational point of an, 211
\item algebraic integer, 43
\item algebraic number, 37
   \subitem degree of an, 38
\item algebraic number field, 50
\subitem embedding of an, 115
\subitem zeta function of an, 120
 \subsubitem Euler-Dedekind product expansion of the, 121
 \subsubitem elementary factors of the, 124
\item algebraic number theory, 37-44, 50-57, 107-116
\item al-Hasan ibn al-Haythem, A. A., 12
\item allowable prime, 226
\item analytic function, 157
\subitem Taylor-series expansion of an, 157
\item analytic number theory, 82-88, 116-131, 147-177
\item arithmetic algebraic geometry, 213
\item arithmetic progression, 81
\item asymptotic density, 88
\item asymptotic functions, 88
\item basic modulus, 185
\item Basic Problem, 13
\subitem solution of the, 75-79
\item Berndt, B., 148, 157, 174, 257
\item Bessel's inequality, 165

\item $(B, \textbf{S})$-signature, 226 

\item Cauchy's Integral Theorem, 158
\item Central Limit Theorem, 250 
\item character, 11
  \subitem additive, 214
  \subsubitem orthogonality relations for an, 214
  \subitem basic, 185
  \subitem Dirichlet, 84
  \subsubitem conductor of a, 184
  \subsubitem induced modulus of a, 184
  \subsubitem orthogonality relations for a, 85
  \subsubitem primitive, 184
  \subsubitem principal, 84
  \subsubitem real, 84
\item Chinese remainder theorem, 7
\item Chung, K.-L., 249, 250, 257 
\item circle group, 11, 84
\item class field theory, 23, 70
\item class number 
\subitem of a fundamental discriminant, 190
\subitem of a number field, 53
\item Clay Mathematics Institute, 87
\item combinatorial number theory, 138-144, 217-246
\item complete Weil sum, 211
\subitem estimate of a, 213
\item computational complexity, 105
\item contour, 158
\subitem closed, 158
\subitem Jordan, 158
\subsubitem exterior of a, 159
\subsubitem interior of a, 159
\subsubitem positively oriented, 159
\item contour integral, 158
\item Cohen, H., 187, 203, 257
\item Conway, J. B., 158, 160, 257

\item cyclotomic number field, 67
\item cyclotomy, 5, 6 
\item Davenport, H., 84, 174, 183, 185, 191, 207, 208, 257
\item Davenport, H. and P. Erd$\ddot{\textrm{o}}$s, 247, 251, 257
\item Dedekind, R., 37, 52, 108, 135, 257
\item Dedekind's Ideal-Distribution Theorem, 113, 114
\item Dirichlet, P. G. L., 24, 77, 82-83, 116, 131, 147, 162, 165, 257 
\item Dirichlet-Hilbert trick, 207, 208
\item Dirichlet kernel, 163
\item Dirichlet $L$-function, 85, 150
\item Dirichlet's class-number formula, 192
\item Dirichlet series, 116, 117
\subitem convergence theorem for, 117

 \item Dirichlet's theorem on primes in arithmetic progression, 81 
\subitem elementary proof of, 137
\subitem proof of, 83-86
\item \emph{Disquisitiones Arithmeticae}, 1, 4-6, 30, 46, 50, 65, 83, 257 
\item Dugundji, J., 159, 257 
\item Eisenstein, M., 32, 50
\item Eisenstein's criterion, 39
\item elementary number theory, 6-8, 137-145, 174

\item elementary symmetric polynomial, 41
\item entire function, 157
\item Erd$\ddot{\textrm{o}}$s, P., 137, 257 
\item Euclidean algorithm, 7
\item Euler-Dirichlet product formula, 85-86, 152
\item Euler, L., 11, 26, 83, 81, 137, 257
\item Euler's constant, 216
\item Euler's criterion, 11
\item Euler's totient function, 84 
\item Fermat, P., 20, 26
\item Fermat's Last Theorem, 213
\item field of complex numbers, 37
\subitem degree of a, 50
\item Filaseta, M. and D. Richman, 91, 131, 257
\item finite extension, 66
\item Fourier series, 160-161
\subitem complex form of a, 178
\subitem convergence theorem for, 162
\subitem cosine coefficient of a, 161
\subitem finite, 177
\subsubitem Fourier coefficients of a, 179
\subitem sine coefficient of a, 161
\item function of bounded variation, 165
\item fundamental discriminant, 62
\item Fundamental Problem, 15 
\subitem solution of the Fundamental Problem for the prime 2, 15-18
\subitem solution of the Fundamental Problem for odd primes, 72-74
\item Fundamental Theorem of Ideal Theory, 51
\subitem proof of the, 132-135
\item fundamental unit, 115
\item Galois automorphism, 65
\item Galois field $GF(2)$ of order 2, 89
\item Galois group, 65
\item Gauss, C. F., 3, 11, 15, 18, 19, 29, 30, 36, 37, 46, 50, 65, 94, 154, 257-258 
\item Gauss' lemma, 16, 40 
\item Gauss sum, 46, 47, 50, 203
\subitem theorem on the value of a, 154
\item Generalized Riemann Hypothesis, 87, 214
\item  Gr$\ddot{\textrm{o}}$ger, D., 29, 258
\item group of units, 84
\item Hecke, E., 30, 52-55, 58, 63-65, 111, 113, 114, 258
\item higher reciprocity laws, 22, 30 
\item Hilbert, D., 87, 107, 131, 258
\item Hungerford, T., 23, 51, 66, 67, 111, 258
\item hybrid or mixed Weil sum, 216 

\item ideal(s), 50
\subitem equivalent, 53 
\subitem genus of, 64
\subitem maximal, 51
\subitem narrow equivalence of, 31, 57, 58

\subitem norm of an, 54, 111
\subitem prime, 51
\subsubitem degree of a, 51, 123
\subitem product of, 51
\item ideal class, 53
\item ideal-class group, 53
\item incomplete Weil sum, 213
\subitem estimate of an, 214
\item infinite product, 121
\subitem absolute convergence of an, 122
\subitem convergence of an, 122
\item integral basis, 54, 108
\item inverse modulo $m$, 7
\subitem existence and uniqueness theorem for an, 7
\item Ireland, K. and M. Rosen, 7, 36, 70, 154, 258
\item isolated singularity, 159
\item Jacobi, C. G. J., 50, 99
\item Jacobi symbol, 99
\subitem Reciprocity Law for the, 100

\item Jordan curve theorem, 159
\item Kronecker, L., 154
\item Kronecker symbol, 187
\item Kurlberg, P., 256, 258
\item Kurlberg, P. and Z. Rudnick, 256, 258
\item Lagrange, J. L., 26, 29, 63, 258
\item Landau, E., 62, 63, 186, 187, 189-191, 258
\item Law of Quadratic Reciprocity (LQR), 23

\subitem Gauss' first proof of the, 30
\subitem Gauss' second proof of the, 30, 57-65
\subitem Gauss' third proof of the, 15, 31-34
\subitem Gauss' fourth proof of the, 31
\subitem Gauss' fifth proof of the, 31, 35
\subitem Gauss' sixth proof of the, 31, 35, 44-49
\subitem Gauss' seventh proof of the, 29
\subitem Gauss' eighth proof of the, 29

\item Legendre, A. M., 28, 29, 258
\item Legendre symbol, 10
\item Lemmermeyer, F., 26, 30, 35, 258
\item LeVeque, W., 88, 92, 258
\item linear Diophantine equation, 7
\subitem solution of a, 7
\item logarithmic integral, 87
\item Marcus, D., 55, 67, 258
\item method of successive substitution, 75, 76
\item Millennium Prize Problems, 87
\item minimal polynomial, 37
\item Minkowski's constant, 55
\item modular substitution, 63
\item Montgomery, M. and R. Vaughan, 88, 92, 214, 258
\item narrow class number, 58
\item narrow ideal class, 58
\item narrow ideal-class group, 58
\item Nevenlinna, R. and V. Paatero, 122, 258
\item norm of a field element, 54

\item normal distribution, 6, 250
\item notation
\subitem $P$, 6
\subitem $\mathbb{Z}$, 6
\subitem $\mathbb{Q}$, 6
\subitem $\mathbb{R}$, 6
\subitem $[m, n], m$ and $n$ integers, $m\leq n$, 6
\subitem $[m, \infty), m$ an integer, 6
\subitem $|A|, A$ a set, 6
\subitem $2^A, A$ a set, 6
\subitem gcd$(m, n)$, 6
\subitem $\pi(z), z$ an integer, 6
\subitem $U(n), n$ a positive integer, 6
\subitem for the Legendre symbol, 10
\subitem $X_{\pm}(z), z$ an integer, 13
\subitem $\pi_{\textnormal{even}}(z), z$ an integer, 13
\subitem $\pi_{\textnormal{odd}}(z), z$ an integer, 13
\subitem $\mathbf{C}$, 37
\subitem $A[x], A$ a commutative ring, 37
\subitem $\mathcal{A}(F), F$ a field, 43
\subitem $\mathcal{R}$, 44
\subitem $[a, b, c]$, 62
\subitem $\mathcal{Q}(d)$, 62
\subitem $\mathcal{I}(d)$, 63
\subitem for the Jacobi symbol, 99
\subitem $\mathcal{I}$, 111
\subitem $Z(n)$, 111
\subitem $\mathcal{Q}$, 121
\subitem $q(I), I$ a subinterval of the real line, 148
\subitem $G(n, p)$, 154

\subitem $R(n)$, 190, 191
\subitem $AP(\textbf{b}; s)$, 210
\subitem $AP(\textbf{a}, \textbf{b}; s)$, 210
\subitem $q_{\varepsilon}(p), p$ an odd prime, 220, 221
\subitem $AP(B, \textbf{S})$, 222
\subitem $\mathcal{E}(A), A$ a finite set, 225
\subitem $\mathcal{K}_{\max}$, 225
\subitem $\Lambda(\mathcal{K})$, 226
\subitem $\Pi_{\pm}(B, \textbf{S})$, 226
\subitem $\Pi_{\pm}(\mathbf{a}, \mathbf{b})$, 237
\item ordinary residue, 6
\subitem minimal non-negative, 6
\item Ore, O., 132, 258

\item overlap diagram, 235
\subitem block of an, 235
\subitem column of an, 235
\subitem gap sequence of an, 235
\item Paley, R. E. A. C., 214
\item Perel'muter, G., 4, 214, 258
\item piecewise differentiable function, 161
\item Poisson distribution, 256
\item pole, 159
\subitem order of a, 159
\subitem simple, 159
\item Polya, G., 213
\item Poussin, C. de la Vall$\acute{\textrm{e}}$e, 86, 174, 258 
\item primary discriminant, 187
\item primary factorization, 187
\item Prime Number Theorem, 88 
\subitem elementary proof of the, 137
\subitem optimal error estimate for the, 87
\item Prime Number Theorem on primes in arithmetic progression, 91
\subitem elementary proof of the, 137
\item primes in arithmetic progression, 77
\item principal genus, 64
\item principal ideal class, 53
\item probability space, 249
\item quadratic congruence, 1
\item quadratic excess, 148
\item quadratic form(s), 5,
\subitem automorph of a, 188
\subitem composition of, 5
\subitem discriminant of a, 62
\subitem equivalent, 63
\subitem genus determined by a, 64
\subitem indefinite, 63
\subitem irreducible, 62
\subitem negative definite, 63
\subitem positive definite, 63
\subitem primary representation by a, 191
\subitem primitive, 188
\subitem representation by a, 62
\subitem representative system of, 190

\item quadratic non-residue, 3
\item quadratic number field, 52
\subitem discriminant of a, 55
\subitem algebraic integers in a, 52
\subitem decomposition law in a, 52-53
\subitem regulator of a, 115
\subitem zeta function of a, 125

\item quadratic residue, 3
\item quadratic residuosity problem, 98
\item quotient diagram, 236
\subitem block of a, 236
\item Rademacher, H., 177, 258
\item random variable(s), 249
\subitem distribution function of a, 249
\subitem independent, 249
\subitem mean of a, 249
\subitem variance of a, 249
\item reciprocity law for polynomials, 22, 23
\item residue (at a pole of an analytic function), 159 
\item residue pattern, 217  
\item  residue (non-residue) support property, 218
\item residue (non-residue) support set, 217-218 
\item residue theorem, 160
\item Riemann, G. F. B., 87, 259
\item Riemann Hypothesis, 87
\item Riemann-Lebesgue lemma, 165
\item Riemann zeta function, 87
\subitem Euler-product expansion of the, 123
\item Rosen, K., 7, 97, 259
\item Rosenberg, J., 30, 259
\item Schmidt, W., 213, 216, 259
\item Selberg, A., 137, 259
\item Shamir, A. 97, 259
\item Shamir's algorithm, 97, 98
\item \item Shohat, J. and J. D. Tamarkin, 253, 259
\item splitting field, 66
\item splitting modulus, 21
\item square-free integer, 39
\item square-free part, 138
\item standard $2m$-tuple, 220
\item Supplement X, 52
\item supports all patterns, set which, 130

\item symmetric difference, 139
\item Taylor, R., 213
\item Taylor, R. and A. Wiles, 213, 259
\item \emph{theorema aureum}, 18, 19
\item unit (in a ring), 84
\item universal pattern property, 217
\item Urbanowicz, J. and K. S. Williams, 203, 259
\item Vinogradov, I. M., 214
\item \emph{Vorlesungen $\ddot{\textrm{u}}$ber Zahlentheorie}, 24, 30, 52, 83, 259
\item Weierstrass approximation theorem, 254
\item Weil, A., 211, 213, 259
\item Weisner, L., 42, 43, 259
\item Wiles, A., 213, 259
\item Wilson, J., 12
\item Wilson's theorem, 12
\item Wright, S., 145, 219, 231, 231, 237, 246, 259
\item Wyman, B. F., 20, 22, 259
\item zero-knowledge proof, 97
\item Zygmund, A., 166, 259

\end{theindex}

\end{document}